# The Distribution of Prime Numbers on the Square Root Spiral

- Harry K. Hahn -

Germany

contribution from " The Number Spiral " by

- Robert Sachs -

**30. June 2007**

## Abstract :

Prime Numbers accumulate on defined spiral graphs, which run through the Square Root Spiral.
These spiral graphs can be assigned to different spiral-systems, in which all spiral-graphs have the same direction of rotation and the same "second difference" between the numbers, which lie on these spiral-graphs. A mathematical analysis shows, that these spiral-graphs are caused exclusively by quadratic polynomials. For example the well known Euler Polynomial $x^2+x+41$ appears on the Square Root Spiral in the form of three spiral-graphs, which are defined by three different quadratic polynomials.

All natural numbers, divisible by a certain prime factor, also lie on defined spiral graphs on the Square Root Spiral ( or "Spiral of Theodorus", or "Wurzelspirale" ). And the Square Numbers 4, 9, 16, 25, 36 … even form a highly three-symmetrical system of three spiral graphs, which divides the square root spiral into three equal areas. Fibonacci number sequences also play a part in the structure of the Square Root Spiral. To learn more about these amazing facts, see my detailed introduction to the Square Root Spiral :

→ " **The ordered distribution of natural numbers on the Square Root Spiral** "

With the help of the "Number-Spiral" , described by Mr. Robert Sachs, a comparison can be drawn between the Square Root Spiral and the Ulam Spiral.

With the kind permission of Mr Robert Sachs, I show some sections of his webside : www.numberspiral.com in this study. These sections contain interesting diagrams, which are related to my analysis results, especially in regards to the distribution of prime numbers.

## Contents Page

# 1 Introduction to the Square Root Spiral :

The Square Root Spiral ( or "Spiral of Theodorus" or "Einstein Spiral" ) is a very interesting geometrical structure in which the square roots of all natural numbers have a clear defined orientation to each other. This enables the attentive viewer to find many spatial interdependencies between natural numbers, by applying simple graphical analysis techniques. Therefore the Square Root Spiral should be an important research object for professionals, who work in the field of number theory !

Here a first impressive image of the Square Root Spiral :

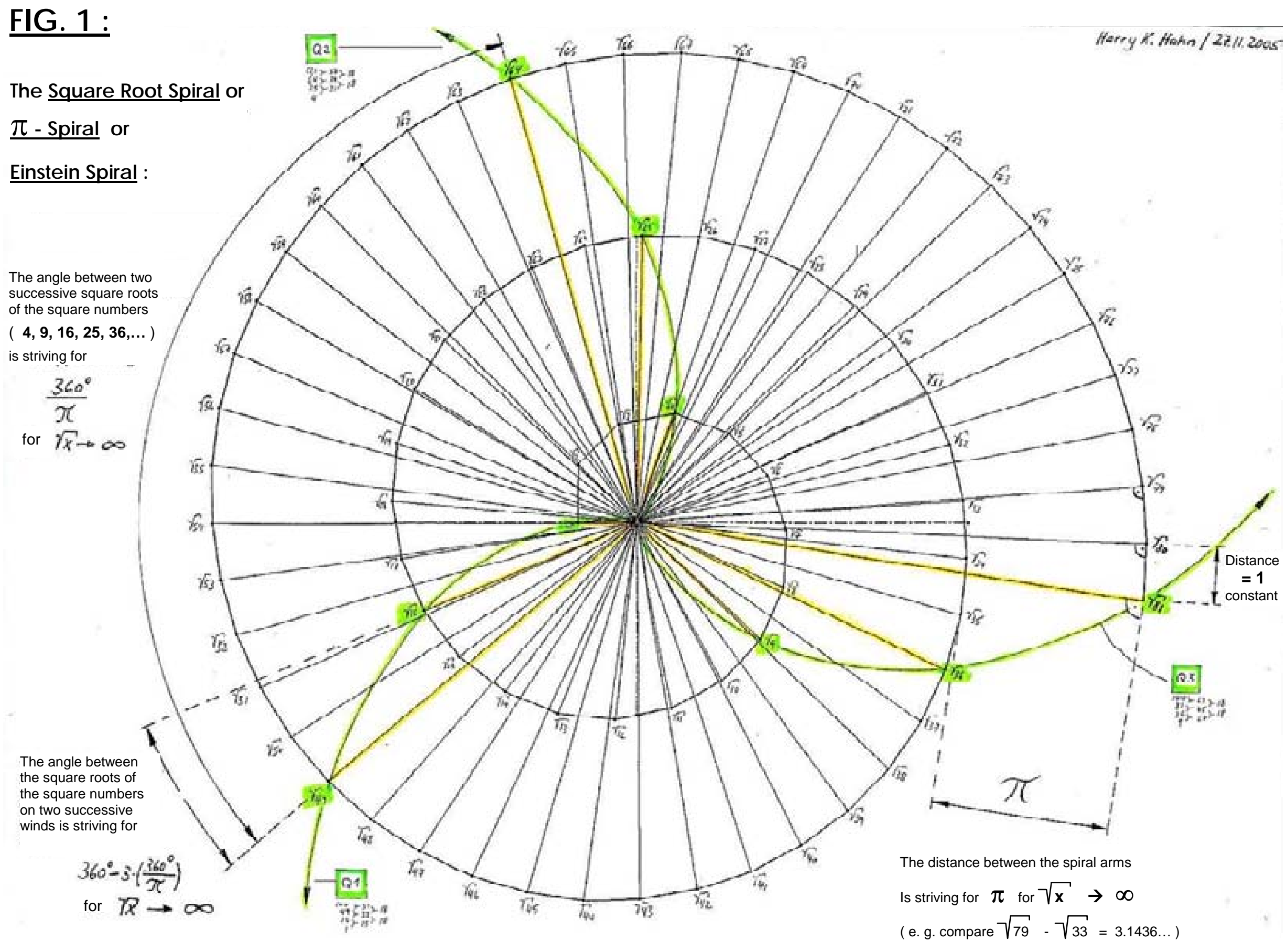

The most amazing property of the square root spiral is surely the fact, that the distance between two successive winds of the Square Root Spiral quickly strives for the well known geometrical constant $\pi$ !!

Mathematical proof that this statement is correct is shown in **Chapter 1 " The correlation to $\pi$ "** in the mathematical section of my detailed introduction to the Square Root Spiral ( → previous study ! ) :

**→ Title : " The ordered distribution of the natural numbers on the Square Root Spiral " → see ArXiv-achive**

Another striking property of the Square Root Spiral is the fact, that the square roots of all square numbers ( 4, 9, 16, 25, 36… ) lie on three highly symmetrical spiral graphs which divide the square root spiral into three equal areas. ( → see FIG.1 : graphs **Q1, Q2** and **Q3** drawn in green ). For these three graphs the following rules apply :

1.) The angle between successive Square Numbers ( on the "Einstein-Spiral" ) is striving for $360°/\pi$ for sqrt( X ) → ∞

2.) The angle between the Square Numbers on two successive winds of the "Einstein-Spiral" is striving for $360° - 3x(360°/\pi)$ for sqrt( X ) → ∞

Proof that these propositions are correct, shows **Chapter 2 " The Spiral Arms"** in the mathematical section of the above mentioned introduction study to the Square-Root Spiral.

**The Square Root Spiral** develops from a right angled base triangle ( **P1** ) with the two legs ( cathets ) having the length 1, and with the long side ( hypotenuse ) having a length which is equal to the square root of 2.

**→ see FIG. 2 and 3**

The square root spiral is formed by further adding right angled triangles to the base triangle **P1** ( see FIG 3)
In this process the longer legs of the next triangles always attach to the hypotenuses of the previous triangles. And the longer leg of the next triangle always has the same length as the hypotenuse of the previous triangle, and the shorter leg always has the length 1.
In this way a spiral structure is developing in which the spiral is created by the shorter legs of the triangles which have the constant length of 1 and where the lengths of the radial rays ( or spokes ) coming from the centre of this spiral are the square roots of the natural numbers ( sqrt 2 , sqrt 3, sqrt 4, sqrt 5 …. ).

**→ see FIG. 3**

The special property of this infinite chain of triangles is the fact that all triangles are also linked through the Pythagorean Theorem of the right angled triangle. This means that there is also a logical relationship between the imaginary square areas which can be linked up with the cathets and hypotenuses of this infinite chain of triangles ( → all square areas are multiples of the base area 1 , and these square areas represent the natural numbers N = 1, 2, 3, 4,…..) → **see FIG. 2 and 3**. This is an important property of the Square Root Spiral, which might turn out one day to be a "golden key" to number theory !

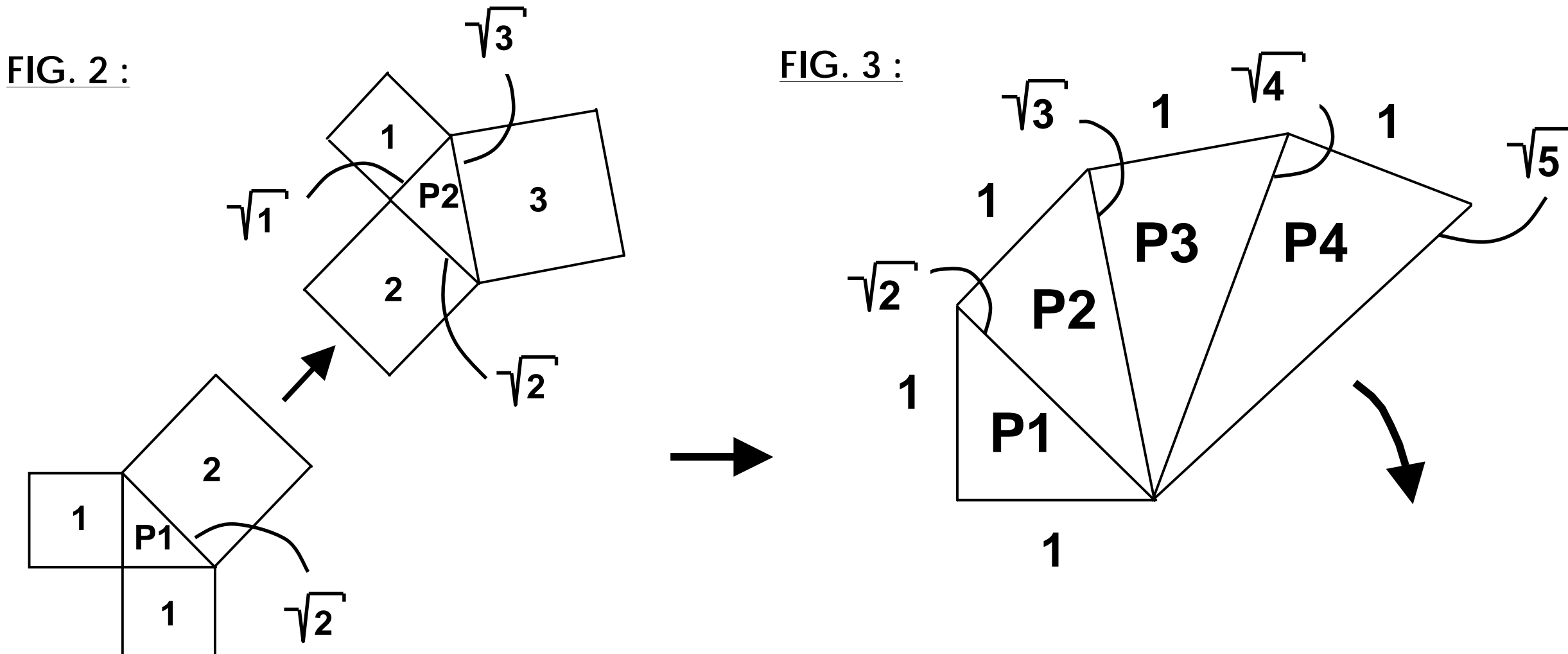

By the way, the first two triangles **P1** and **P2** , which essentially define the structure of the complete Square Root Spiral ad infinitum, are also responsible for the definition of the cube structure. → **see FIG. 4**
Here the triangle **P1** defines the geometry of the area diagonal of the cube, whereas triangle **P2** defines the geometry of the space diagonal of the cube.

→ A cube with the edge length of **1** can be considered as base unit of space itself.

**FIG. 4 :**

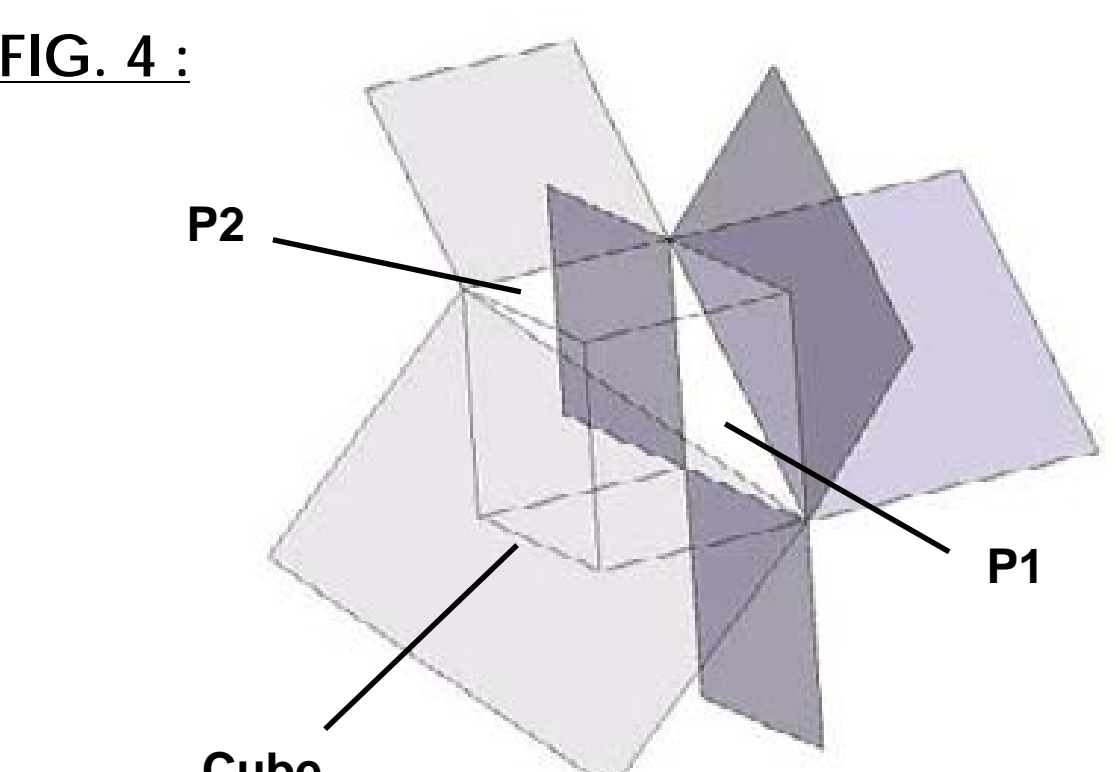

**FIG. 1** shows the further development of the Square Root Spiral ( or Einstein-Spiral ) if one rectangular triangle after the other is added to the growing chain of triangles as described in FIG. 3.

For my further analysis I have created a square root spiral consisting of nearly 300 precise constructed triangles. For this I used the CAD Software SolidWorks. The length of the hypotenuses of these triangles which represent the square roots from the natural numbers 1 to nearly 300, has an accuracy of 8 places after the decimal point. Therefore, the precision of the square root spiral used for the further analysis can be considered to be very high.

The lengths of the radial rays ( or spokes ) coming from the centre of the square root spiral represent the square roots of the natural numbers ( n = { 1, 2, 3, 4,...} ) in reference to the length 1 of the cathets of the base triangle P1 ( see FIG. 3 ). And the natural numbers themselves are imaginable by the areas of "imaginary squares", which stay vertically on these "square root rays". → **see FIG. 5** (compare with FIG.2 )

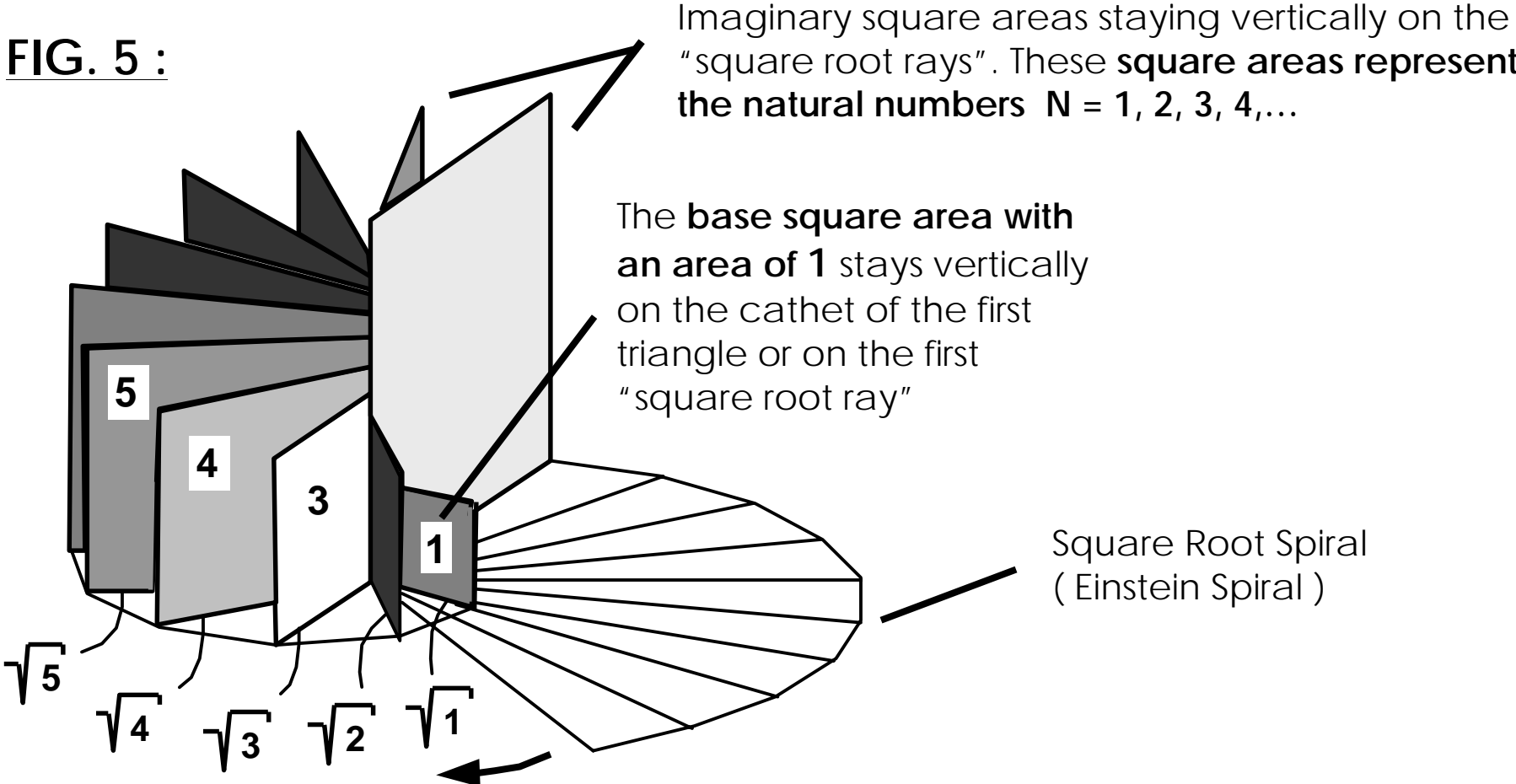

→ The "square root rays" of the Einstein Spiral can simply be seen as a projection of these spatial arranged "imaginary square areas", shown in FIG. 5, onto a 2-dimensional plane.

## 2 Mathematical description of the Square Root Spiral

Comparing the Square Root Spiral with different types of spirals ( e.g. logarithmic-, hyperbolic-, parabolic- and Archimedes- Spirals ), then the Square Root Spiral obviously seems to belong to the Archimedes Spirals.

An Archimedes Spiral is the curve ( or graph ) of a point which moves with a constant angle velocity around the centre of the coordinate system and at the same time with a constant radial velocity away from the centre. Or in other words, the radius of this spiral grows proportional to its rotary angle.

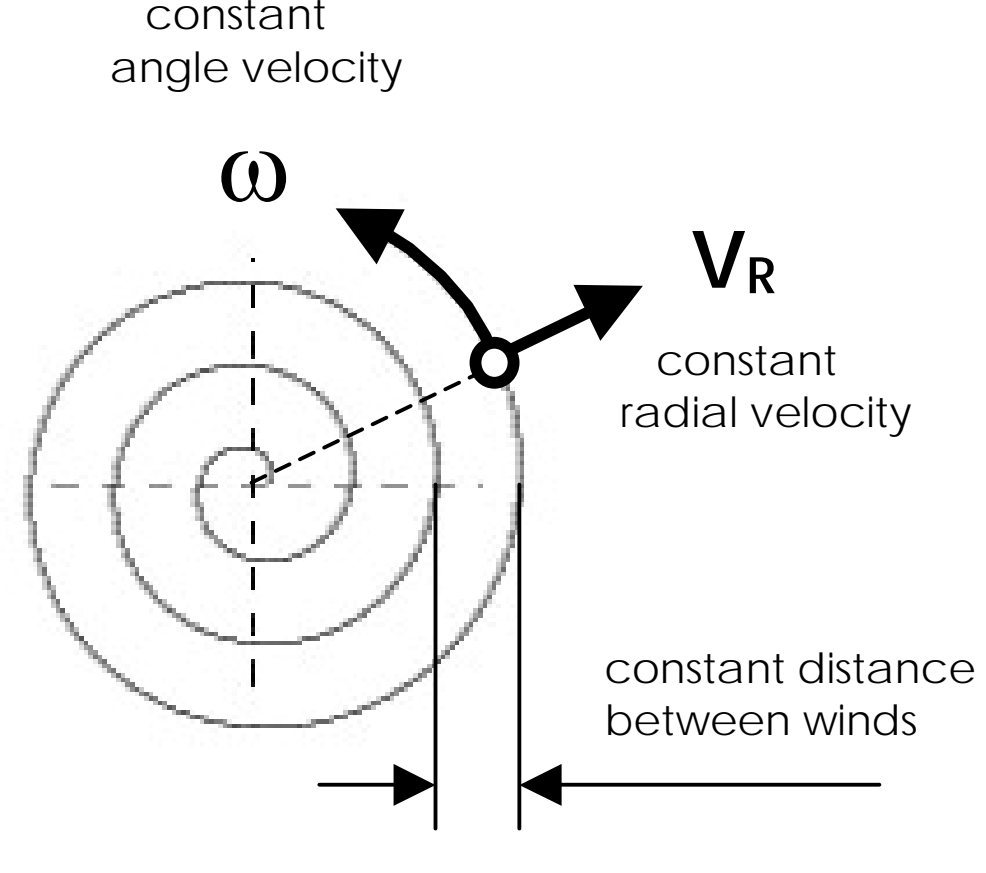

**Archimedes Spiral**

In polar coordinate style the definition of an Archimedes Spiral reads as follows :

$$r(\varphi) = a\varphi \qquad \text{with} \qquad a = \text{const.} = \frac{V_R}{\omega} > 0$$

for $r \to \infty$ the Square Root Spiral is an Archimedes Spiral with the following definition :

$$r(\varphi) = a\varphi + b$$

with $a$ = const. and $b$ = const.

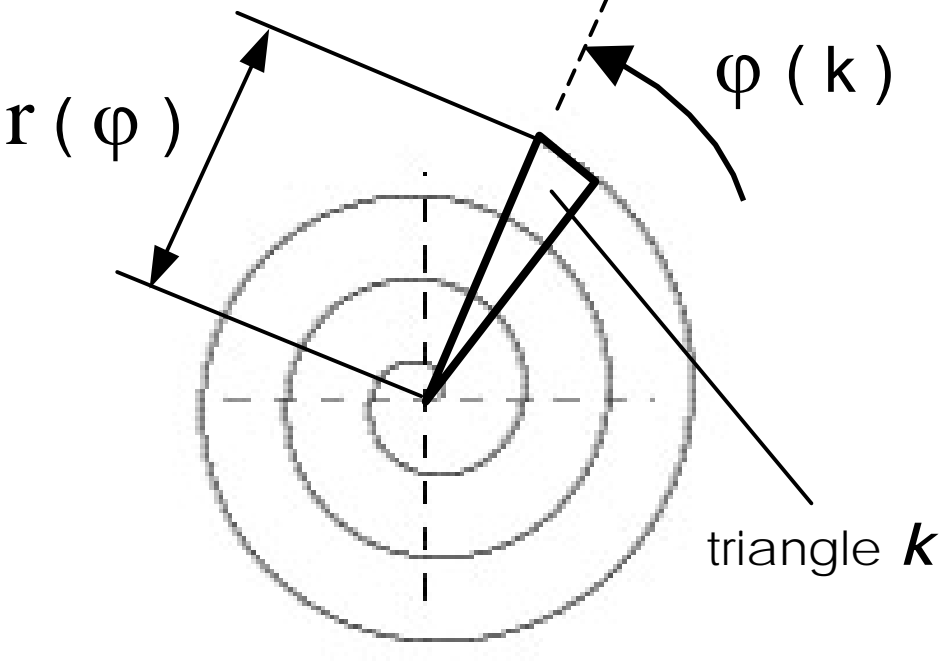

The values of the parameters $a$ and $b$ are

$$a = \frac{1}{2} \quad \text{and} \quad b = -\frac{c_2}{2} \quad ; \text{ with } \quad c_2 = \text{Square Root Spiral Constant}$$

$$c_2 = -2.157782996659....$$

Hence the following formula applies for the Square Root Spiral :

$$r(\varphi) = \frac{1}{2}\varphi + 1.078891498..... \qquad \text{for } r \to \infty$$

for $r \to \infty$ therefore the growth of the radius of the Square Root Spiral after a full rotation is striving for $\pi$ ( corresponding to the angle of a full rotation which is $2\pi$ )

<u>Note</u> : The mathematical definitions shown on this page and on the following page can also be found either in the mathematical section of my introduction study to the Square Root Spiral
→ " The ordered distribution of natural numbers on the Square Root Spiral "
or in other studies referring to the Square Root Spiral.
→ e.g. a mathematical analysis of the Square Root Spiral is available on the following website : → **http://kociemba.org/themen/spirale/spirale.htm**

**Further dependencies in the Square Root Spiral :**

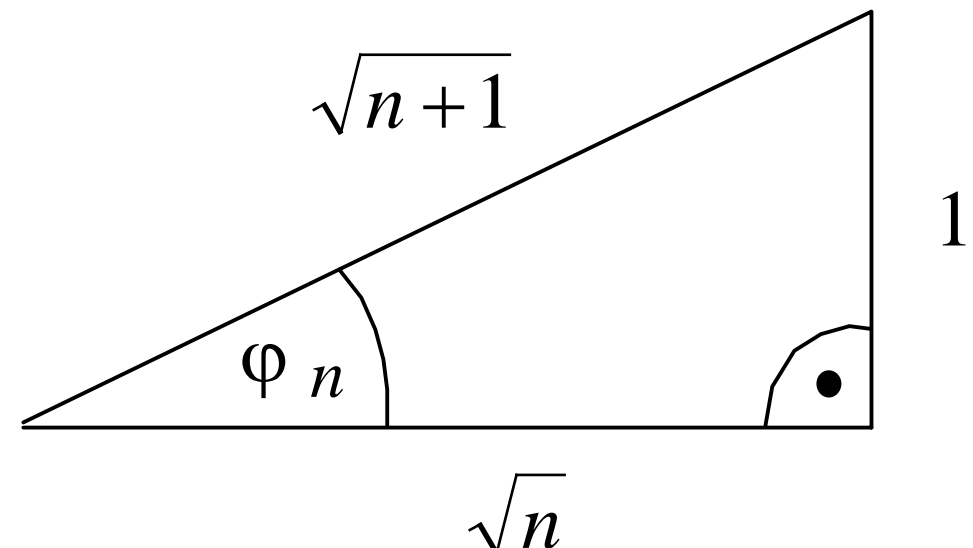

If $\varphi_n$ is the angle of the $n$th spiral segment ( or triangle ) of the Square Root Spiral, then

$\tan(\varphi_n) = \dfrac{1}{\sqrt{n}}$ ; ( ratio $\dfrac{\text{counter cathet}}{\text{cathet}}$ )

If the $n$th triangle is added to the Square Root Spiral the growth of the angle is

$$\varphi_n = \arctan\left(\frac{1}{\sqrt{n}}\right)$$ ; Note : angle in radian

The total angle $\varphi(k)$ of a number of $k$ triangles is

$\varphi(k) = \sum_{n=1}^{k} \varphi_n$ or described by an integral $\int_0^k \arctan\left(\frac{1}{\sqrt{n}}\right) dn + c_1(k)$

$\Rightarrow \quad \varphi(k) = 2\sqrt{k} + c_2(k)$ with $\lim_{k \to \infty} c_2(k) = \text{const.} = -2.157782996659.....$

$c_2$ = Square Root Spiral Constant

The growth of the radius of the Square Root Spiral at a certain triangle $n$ is

$$\Delta r = \sqrt{n+1} - \sqrt{n}$$

The radius $r$ of the Square Root Spiral ( i.e. the big cathet of triangle $k$ ) is

$r(k) = \sqrt{k}$ and by converting the above shown equation for $\varphi(k)$ it applies that

$$r(k(\varphi)) = r(\varphi) = \sqrt{\frac{1}{4}(\varphi - c_2(\varphi))^2} = \frac{1}{2}\varphi - \frac{c_2}{2}$$

For large $n$ it also applies that $\varphi_n$ is approximately $\dfrac{1}{\sqrt{n}}$ and $\Delta r$ has pretty well half of this value, that is $\dfrac{1}{2\sqrt{n}}$ , what can be proven with the help of a Taylor Sequence.

## 3 The distribution of Prime Numbers on defined spiral-graphs :

In a similar way as the Square Numbers shown in **FIG.1** , Prime Numbers also accumulate on certain spiral-graphs, which run through the Square Root Spiral ( Einstein Spiral ).

After marking all Prime Numbers with yellow color on the Square Root Spiral , it is easy to see, that the Prime Numbers accumulate on certain spiral-shaped graphs. I want to call these spiral-graphs with a high share in prime numbers "Prime Number Spiral-Graphs".

To identify these "Prime Number Spiral-Graphs", I marked the most conspicuous spiral graphs with a high share of Prime Numbers, then I tried to make sense out of their arrangement on the Square Root Spiral .

→ **see FIG. 6-A** to **6-B** on the following pages.

A closer look to these " Prime Number Spiral-Graphs" revealed the following properties :

- It is obvious that certain Prime Number Spiral-Graphs are related to each other and that they belong to the same "spiral-graph-system".
- By calculating the differences ( first difference ) of the consecutive numbers lying on one of the found "Prime Number Spiral-Graphs", and by further calculating the differences of these differences ( second difference ) we always obtain one of the following three numbers :

  **18, 20** or **22** → I called these numbers the " **2. Differential** " of the spiral graphs.

  → see for example the difference calculation in FIG. 6-A for the exemplary spiralarm **A3** : ( see PNS-P18-A )

  The calculation of the differences between the numbers 11, 41, 89, 155, 239,.... , which lie on the spiralarm A3 results in the following numbers : 30, 48, 66, 84,..... And the calculation of the differences between these numbers results in the constant value **18**.
  And this number represents the " **2. Differential** " of this spiralarm **A3**.

It is notable that the 2. Differential of the Prime Number Spiral-Graphs is always an even number .
And it seems that the 2. Differential only takes on one of these three values : **18, 20** or **22**

That is the reason why I used these three different possible values of the **2. Differential** as distinguishing property for the graphical representation of the Prime Number Spiral-Graphs shown in FIG. 6-A to 6-C.

→ Therefore in **FIG. 6-A to 6-C** the following assignment applies :

**FIG. 6-A** : shows only Prime Number Spiral-Graphs which have a **2. Differential** of **18**

**FIG. 6-B** : shows only Prime Number Spiral-Graphs which have a **2. Differential** of **20**

**FIG. 6-C** : shows only Prime Number Spiral-Graphs which have a **2. Differential** of **22**

## 3.1 The found Spiral-Graphs can be assigned to different Spiral-Graph-Systems

In my attempt to make sense out of the distribution of the Prime-Numbers on the Square Root Spiral, I first tried to establish order under the found Prime Number Spiral Graphs.

By doing this, I realized that the Prime Number Spiral Graphs are arranged in different " systems" .

As best example I want to refer to **FIG. 6-A → see following pages !**

→ **FIG. 6-A** : The **3** Prime Number Spiral Systems shown in this diagram all have a **2. Differential** of **18** !!

→ see difference calculation for the three exemplary spiralarms **A3, B5** and **C12**

The diagram shows how the Prime Numbers are clearly distributed on **3** defined spiral graph systems, which are arranged in a highly symmetrical manner ( in an angle of around 120° to each other ) around the centre of the Square Root Spiral.

On the shown **3 Prime Number-(Spiral)-Systems** ( **PNS** ) : **P18-A, P18-C** and **P18-C** , the Prime Numbers are located on pairs of spiral arms, which are separated by three numbers in between. And two spiral arms of one such pair of spiral arms, are separated by one number in between.

All spiral-graphs of the shown **3** Prime Number-( Spiral )-Systems ( **PNS** ) have a **positive rotation direction** ( **P** ) and, as already mentioned before, the **2. Differential** of all spiral-graphs has the constant value of **18**. That's why the first part of the naming of the **3** Prime Number-Spiral-Systems ( **PNS** ) is **P18**.
The **3** spiral-graph systems **A** ( drawn in orange ), **B** ( drawn in pink ) and **C** ( drawn in blue ) have further spiralarms. But for clearness there are only around 10 spiralarms drawn per system.

One important property of all Prime Number spiral-graphs ( shown in FIG 6-A ) is the obvious missing of numbers which are divisible by **2** or **3** in these graphs ! That means, that the smallest possible prime factor of the " Non-Prime Numbers" , which lie on these spiralarms, is **5**.

→ **FIG. 6-B** , on the following pages, shows a diagram with another set of **12** Prime Number Spiral Systems In this diagram all Spiral-Graphs have a **2. Differential** of **20** !!

→ see difference calculation for the four exemplary spiralarms **D8, F2, G5** and **I5**

On the shown **12 Prime Number Spiral Systems** ( **PNS** ) : **N20-D** to **N20-I** , and **P20-D** to **P20-I** , the Prime Numbers are again located on pairs of spiral arms, which are separated by three number in between. And two spiral arms of one such pair of spiral arms, are again separated by one number in between.

**6** of the shown Prime-Number-Spiral-Systems ( **PNS** ) have a **positive rotation direction** ( **P** ) and the other **6** Prime-Number-Spiral-Systems ( **PNS** ) have a **negative rotation direction** ( **N** ). The spiral-graph systems of these two groups are arranged in a symmetrical manner around the centre of the Square Root Spiral, in an angle of approx. 60° to each other. And two systems at a time are approx. point-symmetrical to each other ( in reference to the centre of the Square Root Spiral ). For example the two systems **N20-I** & **N20-F**

For clearness only **4** Prime-Number-Spiral-Systems ( **PNS** ) with a **negative rotation direction** ( **N** ) are drawn in color !! These are the **4** systems : **N20-D** (drawn in orange) , **N20-F** (drawn in red) , **N20-G** (drawn in blue ), **N20-I** (drawn in pink). There are only around 10 spiralarms drawn ofr each of these 4 spiral systems.

Note, that the **other 8** Prime Number Spiral Systems are all drawn in light grey color, and that there are only 2 to 6 spiralarms drawn of each of these systems, for clearness ! Please also note, that only the Spiral-Graphs with a **negative rotation direction** ( **N** ) are named in FIG 6-B. The naming of the Spiral-Graphs with a **positive rotation direction** ( **P** ) P-20… would in principle just mirror the naming of the "N20…-spiralarms". → For example the spiralarm with the " P20-G" -mark attached to it ( = name of this system ), would be named G1 and the next spiralarm on the left G2 and so on.

→ **FIG. 6-C**, shows the third group of Spiral-Systems , which contains altogether **11** Prime Number Spiral Systems , which all have a **2. Differential** of **22** !!

→ see difference calculation for the 4 exemplary spiralarms **J10, L6, N11** and **Q7**

On the shown **11 Prime Number Spiral Systems** ( **PNS** ) : **N22-J** to **N22-T** , the Prime Numbers are again located on pairs of spiralarms, which are separated by three numbers in between. And two spiralarms of one such pair of spiralarms, are again separated by one number in between.

All **11** shown Prime Number Spiral Systems ( **PNS** ) have a **negative rotation direction** ( **N** ). And it seems that the spiral graph systems are arranged in a symmetrical manner around the centre of the Square Root Spiral ( in an angle of ~ 360°/11 to each other, in reference to the centre of the Square Root Spiral ).

For clearness only **4** Prime Number Spiral Systems ( **PNS** ) are drawn in color !! These are the following **4** systems : **N22-J** ( drawn in orange ) , **N22-L** ( drawn in pink ), **N22-N** ( drawn in blue ) , **N20-Q** ( drawn in red ). Note, that there are only around 10 spiralarms drawn for each of these four spiral-systems.

Note, that the **other 7** Prime Number Spiral Systems are all drawn in light grey color, and that there are only around 4 to 6 spiralarms drawn of each of these systems !
For clearness, the spiralarms of these 7 systems are not named on FIG. 6-C. But the naming of these spiralarms would start on the spiralarm which has the naming of the system attached to it. For example the spiralarm with the " N22-M" -mark attached to it would be named M1, and the next spiralarm below would be named M2 and so on.

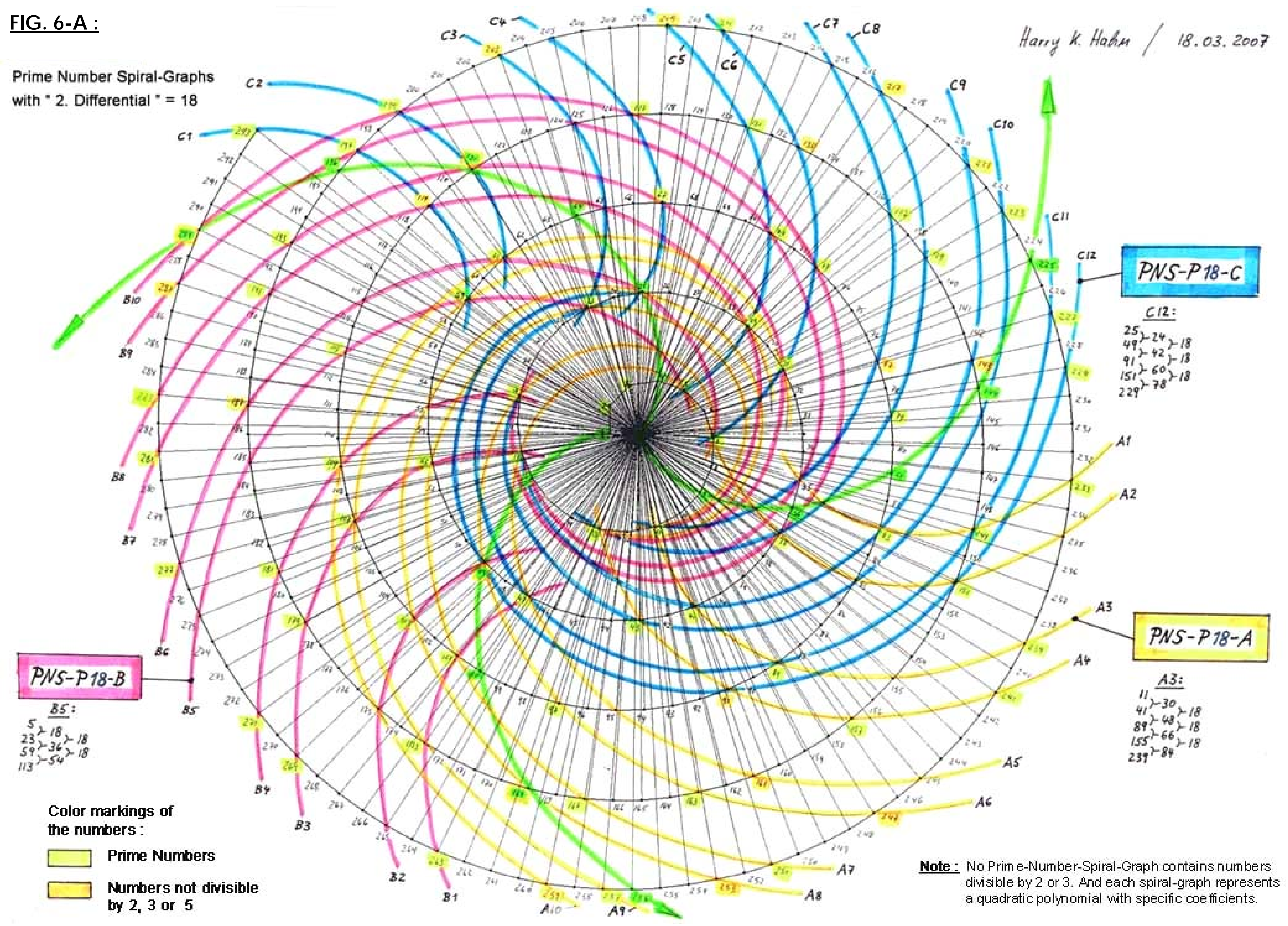

FIG. 6-A :
Prime Number Spiral-Graphs
with " 2. Differential " = 18
Harry K. Hahn / 18.03.2007
PNS-P18-C
C12:
25 > 24 > 18
49 > 42 > 18
91 > 60 > 18
151 > 78 > 18
229
PNS-P18-A
A3:
11 > 30 > 18
41 > 48 > 18
89 > 66 > 18
155 > 84 > 18
239
PNS-P18-B
B5:
5 > 18 > 18
23 > 36 > 18
59 > 54
113
Color markings of the numbers :
Prime Numbers
Numbers not divisible by 2, 3 or 5
Note : No Prime-Number-Spiral-Graph contains numbers divisible by 2 or 3. And each spiral-graph represents a quadratic polynomial with specific coefficients.

**FIG. 6-B :**

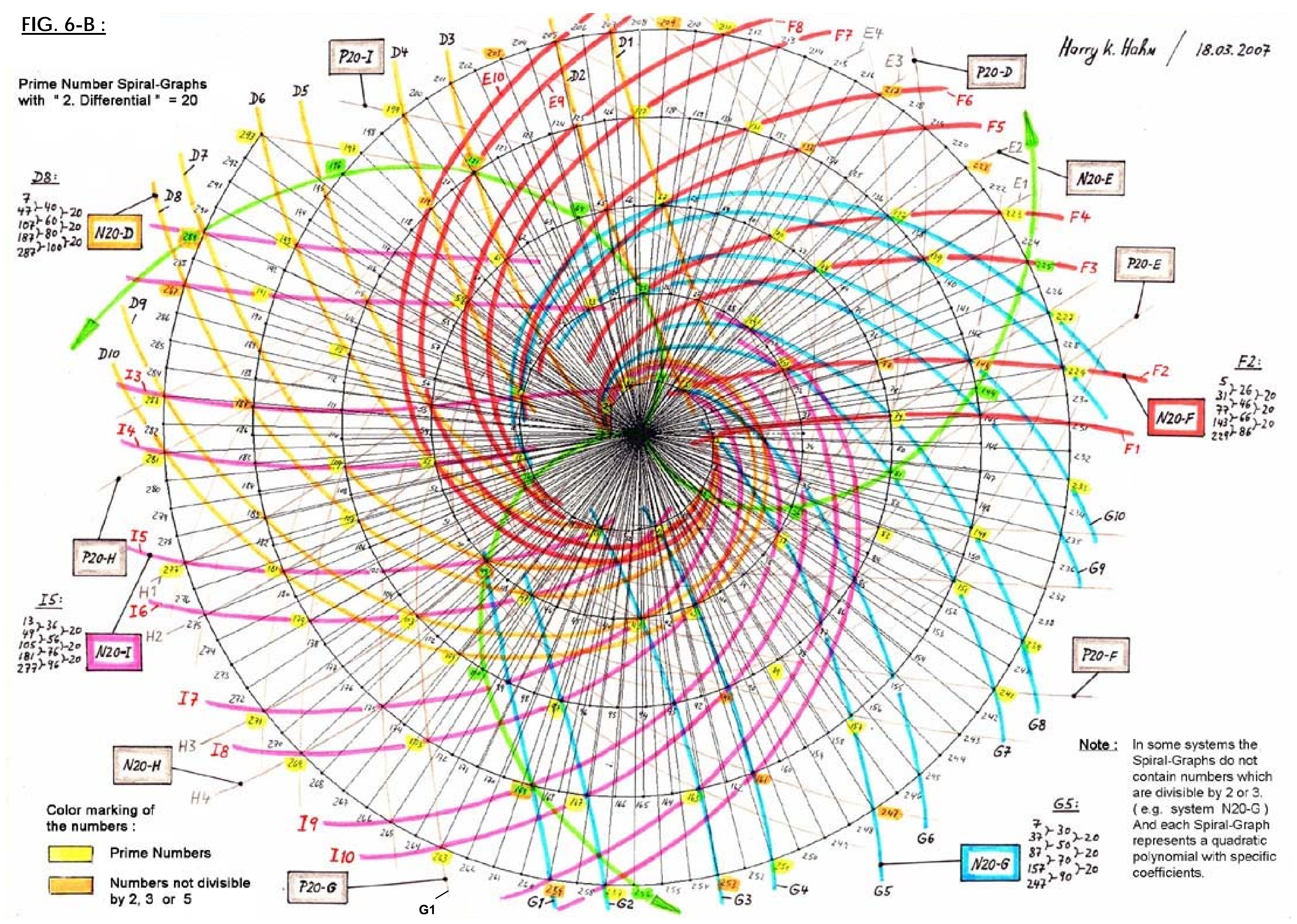

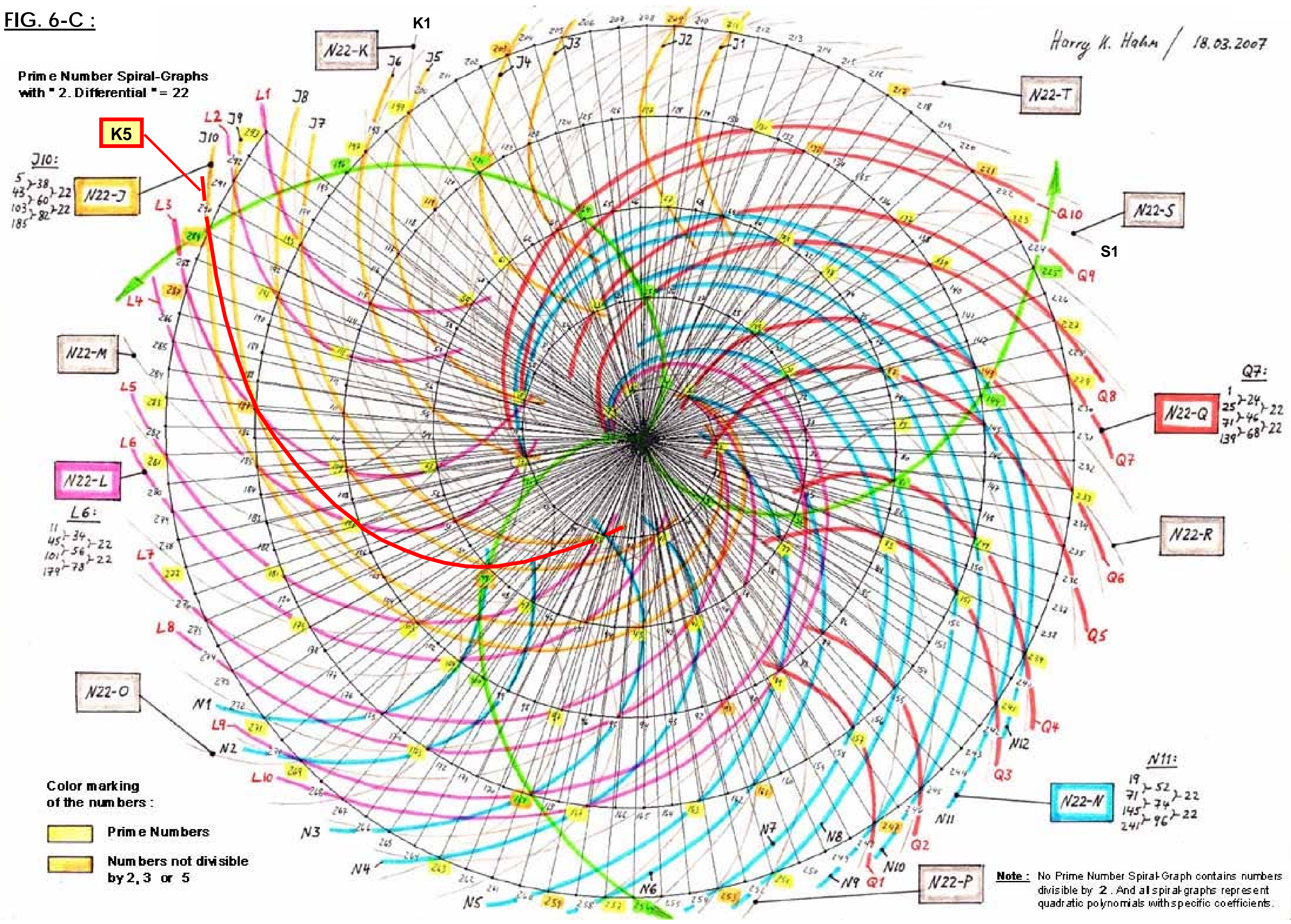
FIG. 6-C :
Prime Number Spiral-Graphs with " 2. Differential " = 22
Harry K. Hahn / 18.03.2007
K1
K5
S1
N22-K
N22-T
N22-S
N22-Q
N22-R
N22-N
N22-P
N22-O
N22-L
N22-M
N22-J
J1 J2 J3 J4 J5 J6 J7 J8 J9 J10
L1 L2 L3 L4 L5 L6 L7 L8 L9 L10
Q1 Q2 Q3 Q4 Q5 Q6 Q7 Q8 Q9 Q10
N1 N2 N3 N4 N5 N6 N7 N8 N9 N10 N11 N12
J10:
5, 43, 103, 185 ; 38, 60, 82 ; 22, 22
L6:
11, 45, 101, 179 ; 34, 56, 78 ; 22, 22
Q7:
1, 25, 71, 139 ; 24, 46, 68 ; 22, 22
N11:
19, 71, 145, 241 ; 52, 74, 96 ; 22, 22
Color marking of the numbers :
Prime Numbers
Numbers not divisible by 2, 3 or 5
Note : No Prime Number Spiral-Graph contains numbers divisible by 2. And all spiral-graphs represent quadratic polynomials with specific coefficients.

## 4 Number Sequences derived from the Spiral-Graphs shown in FIG. 6-A to 6-C

In the following section I want to show an analysis of the Number Sequences which belong to the " Prime Number Spiral Graphs" shown in FIG 6-A to 6-C.

In the Tables **6-A1** to **6-C1** which refer to the diagrams 6-A to 6-C, I have set-up three columns for each analysed "Prime Number Spiral Graph".
( → Tables **6-B1**and **6-C1** can be found **in the Appendix** ! )

The three colums in Table **6-A1**, which refer to the exemplary Spiral-Graph **A3** ( → see **FIG 6-A** ) are named A3, A3′ and A3′′.

In the left column A3 , I listed all " known" numbers which lie on this spiral-graph ( see FIG 6-A ). These are the numbers 11, 41, 89, 155 and 239.

For each analysed Spiral-Graph, there are 4 to 6 numbes available, which can be extracted from the graphs drawn in FIG 6-A to FIG 6-C.

With the help of at least three numbers of each Prime Number Spiral Graph the "first differences" and the "second differences" ( 2. Differential ) of these numbers can be calculated, which then can be used for the further development of the number sequence belonging to the analysed spiral-graph.

For the exemplary spiral-graph A3, the "first differences" between the known five numbers 11, 41, 89, 155, 239 are 30, 48, 66, 84.
And the "second differences" ( 2. Differential ) for this number sequence result in the constant value 18.

Now the "first differences" in column A3′ can be extended, with the help of the calculated 2. Differential = A3′′ = 18. And with the help of the extended sequence A3′ the number sequence of the Prime Number Spiral Graph A3 can be further extended too.

In this way, I have extended the number-sequences derived from the spiral-graphs shown in FIG 6-A to FIG 6-C up to numbers > 100 000.

### 4.1 Quadratic Polynomials are the foundation of the Spiral-Graphs :

The calculated " first differences" and " second differences" can also be used to determine the quadratic polynomials which define these Prime Number Spiral Graphs. An explanation of the calculation procedure ( → Newton Interpolation Polynomial ), to determine the quadratic polynomials of the Prime Number Spiral Graphs, is available in my introduction to the Square Root Spiral :

- → **"The ordered distribution of the natural numbers on the Square Root Spiral"**
- → This study can be found under my author-name in the arXiv - archive

Please also have a read through this study ! There is a mathematical section included in this study, which describes the shown Spiral-Graphs from the mathematical point of view.

The calculated **quadratic polynomials**, which define the Prime Number Spiral Graphs, and the number sequences belonging to it, can be found in the **Tables 6-A2** to **6-C2** → **see next pages ( Table 6-A2 )**
**and Appendix ( Table 6-B2 and Table 6-C2 ) !**

→ Referring to the general quadratic polynomial : $f(x) = \mathbf{a}x^2 + \mathbf{b}x + \mathbf{c}$

the following rules apply for the quadratic polynomials, belonging to the analysed spiral-graphs ( → shown in Tables 6-A2 to 6-C2 ) :

Rules for coefficients **a, b** and **c** :
( or meaning of coefficients )

- **a** → equivalent to the **2.Differential** of the Spiral-Graph divided by **2**
- **b** → this coefficient ( or sequence of coefficients ), refers to the system the Spiral-Graph is belonging to.
- **c** → consecutive parallel distance between the Spiral-Graphs , belonging to the same system.

The next step in my analysis of the number sequences shown in the tables 6-A1 to 6-C1, was the marking of the Prime Numbers in the first 25 numbers of each sequence. I used the color yellow to mark the Prime Numbers. I also marked the numbers which are not divisible by 2, 3 or 5 in the first 25 numbers of these sequences with red color.

Additional to the described three columns per number-sequence , I added another column to each analysed sequence ( → see Tables 6-A1to 6-C1 ).
This additional column has the naming " SD " , which means " sum of the digits" and it shows the sums of the digits of the first 25 numbers of each number sequence.
For example the " sums of the digits" of the first 3 numbers of the A3 - number sequence are 2 ( = 1 + 1 → 11 ) ; 5 ( = 4 + 1 → 41 ) ; 17 ( = 8 + 9 → 89 ) etc.

By looking over tables 6-A1to 6-C1 , it is easy to see, that there are differences in the distribution of the color marked cells, and that in many of the shown number sequences, the distribution of the color-marked cells has a clear periodic character !

Now we want to start with a more detailed analysis of the number sequences shown in the Tables 6-A1, 6-B1 and 6-C1 :

### 4.2 Differences in the Number-Sequences shown in Table 6-A1 to 6-C1

→ **Here is a list of some differences in these number-sequences, which are easy to notice :**

- Approximately **30 %** of all analysed number sequences only contain numbers which are marked in yellow and red. This means they do not contain numbers which are divisible by 2, 3 or 5.
Therefore all non-prime numbers ( marked in red ) in these number sequences consist of prime factors ≥ 7
- It is easy to see, that the white cells and the colored cells (yellow or red) periodically alternate in various ways, in the most number-sequences.
- A closer look shows, that there are sequences, where numbers divisible by 3 or 5 occur in a periodic manner.
For example in the number-sequences A3, A10 and B4 ( in Table 6-A1 ) a number divisible by 5 periodically alternates with 4 numbers which are not divisible by 2, 3 or 5. Or in the number-sequences K2, L2 and M2 ( see Table 6-C1 in the Appendix ) a number divisible by 3 periodically alternates with two numbers not divisible by 2, 3, or 5 .
- All Spiral-Graphs of the spiral-graph systems N20-D, N20-G, P20-D and P20-G ( → see FIG 6-B ) only contain numbers with the same ending. For example the Spiral-Graph **D8** in FIG 6-B only contains numbers with the ending **7**.
- The number sequences of all other spiral-graphs in FIG 6-A to FIG 6-C contain numbers, whose endings always follow a sequence of exactly 5 numbers, which is periodically recurring ad infinitum.

  For example Spiral-Graph **Q3** ( → see FIG 6-C and Table 6-C1 ) only contains numbers with the ending 1, 3 or 7. And these number-endings always occur in the following sequence : **…3, 1, 1, 3, 7,…..**, which is recurring in a periodic manner.

In principle these were the most noticeable differences of the number-sequences listed in Table 6-A1 to 6-C1.

### 4.3 Analysis of the number-endings :

A closer look shows, that the number-endings, which occur in the number-sequences depend on the "2.Differential" of the number-sequence ( or spiral-graph ). This connection is shown in **Table 3** ( → see next pages ).

In this table the prime number sequences ( shown in table 6-A1 to 6-C1 ) are ordered according to the number-endings, which occur in these sequences.

It is easy noticeable that in number-sequences, which have the 2. Differential 18 or 22 ( see first column of table ), the number-endings are defined by the same groups of 3 numbers.

For example the number-endings 1, 5, 7 or 1, 3, 7 occur in number-sequences which have either the 2. Differential 18 or 22.

The only differences are the "Number-Endings-Sequences" in which these groups of three numbers occur !

For example in number-sequences which have the 2. Differential 18 , the number-endings 1, 5, 7 or 1, 3, 7 occur as "Number-Endings-Sequences" …7, 7, 5, 1, 5,… and …3, 3, 1, 7, 1,…., whereas the same number-endings occur as quite different number-endings-sequences …7, 5, 5, 7, 1,… and …3, 1, 1, 3, 7,….. in the number-sequences with the 2. Differential 22.

In number-sequences which have the 2. Differential 20 all numbers either have the same ending, which is one of the odd numbers 1, 3, 5, 7 or 9 , or the numbers have as endings the mentioned 5 odd numbers all in the same sequence. In number-sequences, where all 5 odd numbers occur, one of the following two number-endings-sequences appears : …1, 5, 9, 3, 7,…. or …7, 3, 9, 5, 1,…. → Noticeable is the **opposite** direction of these two sequences !

The most remarkable property of all number-endings-sequences listed in Table 3, is certainly the fact, that all number-endings-sequences cover exactly **5** consecutive numbers !

This indicates, that there is a kind of "higher basic-oscillation" acting in the Square Root Spiral, which seems to cover exactly 5 windings of the Square Root Spiral per oscillation, and which interacts with all Prime Number Spiral Graphs shown in FIG 6-A to 6-C ! This amazing fact is worth an own analysis !!

From the mathematical point of view the following explanation can be given for the periodic occurrence of the number-endings described in Table 3 :

"…… Because of the recursive representation

$$p(t+1) = p(t) + a(2t+1) + b = f(p(t), t)$$

the sequence ( $p(t)$ mod $k$ ) recurs for each natural $k$, since only a maximum of $k^2$ pairs ( $p(t), t$ ) are available. This explains on the one hand the recurrence of the last sequence of figures ( $k$ = 10 ) and on the other hand the regular occurrence of certain $k$ factors ( as shown e.g. in Table 4, 5-A & 5-B ).
This applies to all $k$, not only to the prime numbers.
The sequence of numbers not divisible by 2, 3 or 5 is recurringly modulo 30 and hence a repetend of the corresponding figures also exists in the series ( number-sequences ) shown in Table 6-A1to 6-C1.
The length of this repetend can theoretically be $30^2$ = 900 ……."

This mathematical explanation for the periodic occurrence of the number-endings shown in Table 3 is from Mr. Kay Schoenberger, who also contributed the mathematical section of my introduction to the Square Root Spiral.

→ **"The ordered distribution of the natural numbers on the Square Root Spiral"**

( → This study can be found under my author-name in the arXiv – data bank )

Please have a read through the mathematical section included in this study, which describes the shown spiral-graphs from the mathematical point of view.

Table 3 gives a general information about the periodic occurrence of numbers divisible by 3 or 5 in the number-sequences derived from the Prime Number Spiral Graphs → see last two columns on the righthand side of the table .

### 4.4 Analysis of the "sums of the digits" :

Another remarkable property of the number-sequences listed in Table 6-A1 to 6-C1 is the fact, that only defined "sums of the digits" occur in every number-sequence. In the following we want to carry out a simple analysis of the "sums of the digits" which occur in these number sequences.

As described before, I added a column with the naming " SD " ( which means " sum of the digits" ) to each number-sequence listed in Table 6-A1 to 6-C1. And the numbers in this column represent the sums of the digits of the first 25 numbers of each number sequence. → For example the "sums of the digits" of the first three numbers of the **A3** – number-sequence in Table 6-A1 are **2** ( = 1 + 1 → 11 ) ; **5** ( = 4 + 1 → 41 ) ; **17** ( = 8 + 9 → 89 ) etc.

By looking over the numbers listed in the columns "SD" it is easy to see, that there are only defined "SD" - numbers occurring in each number-sequence.

And by putting these sums of the digits **in order** ( according to their value ! ) , a certain "sums of the digits – sequence" for every number-sequence in Table 6-A1 to 6-C1 can be found.

**Table 2** on the next pages shows the results of this analysis !

In this table the number-sequences ( shown in table 6-A1 to 6-C1 ) are ordered according to the "sums of the digits – sequences" , which occur in these sequences.

Here the following classification can be made :

Number-sequences with the 2. Differential **18** produce sums of the digits – sequences with either the distance of **3** or **9** between two consecutive numbers of the sums of the digits – sequence.

For example the number-sequences A2, A4, A6,…etc. produce the sums of the digits – sequence 4, 7, 10, 13, 16,… in which the distance between two consecutive numbers is **3**. Or the number-sequences B1, B7, B13,…, as another example, produce the sums of the digits – sequence 2, 11, 20, 29, 38,… in which the distance between two consecutive numbers is **9**.

Number-sequences with the 2. Differential 20 or 22 produce "sums of the digits" – sequences which show a kind of periodic behavior similar to the "number-endings-sequences" in Table 3.

Here the differences between **5** consecutive numbers of the sums of the digits – sequences form a periodic sequence of four numbers which recur ad infinitum.

For example the number-sequences N20-D1 and N20-F3 with the 2. Differential **20**, produce the following **ordered** sums of the digits – sequence : 1, 4, 7, 9, 10, 13, 16, 18, 19,… in which the periodic recurring differences **…..3, 3, 2, 1,….** occur.

And as another example, the number-sequences K3 and L1 with the 2. Differential **22**, produce the following **ordered** sums of the digits – sequence : 4, 5, 7, 10, 13, 14, 16, 19, 22,….. in which the periodic recurring differences **…..1, 2, 3, 3,….** occur.

Remarkable is here the fact, that the direction of the periodic recurring differences …..3, 3, 2, 1,…. and …..1, 2, 3, 3,…. in the sums of the digits – sequences, is exactly opposite in the number-sequences with the different "2. Differentials" 20 or 22.

But even more remarkable is the fact, that these mentioned "periodic-recurring-difference-sequences" only occur in the **ordered** sums of the digits – sequences shown in Table 2, but no periodic behavior at all, of the **unordered** " sums of the digits" can be noticed in the "SD" - colums in Table 6-A1 to 6-C1 !!

I haven't found an explanation for this strange characteristic yet !

**Table 6-A1 :** “Prime Number Sequences” derived from the graphs shown in the “Prime Number Spiral Systems” **P18-A ; P18-B** and **P18-C** → see **FIG. 6-A**

| | A1 | | | | A2 | | | | A3 | | | | A4 | | |
|---|---|---|---|---|---|---|---|---|---|---|---|---|---|---|---|
| SD | A1 | A1' | A1" | SD | A2 | A2' | A2" | SD | A3 | A3' | A3" | SD | A4 | A4' | A4" |
| 5 | 5 | | | 7 | 7 | | | 2 | 11 | | | 4 | 13 | | |
| 8 | 35 | 30 | | 10 | 37 | 30 | | 5 | 41 | 30 | | 7 | 43 | 30 | |
| 11 | 83 | 48 | 18 | 13 | 85 | 48 | 18 | 17 | 89 | 48 | 18 | 10 | 91 | 48 | 18 |
| 14 | 149 | 66 | 18 | 7 | 151 | 66 | 18 | 11 | 155 | 66 | 18 | 13 | 157 | 66 | 18 |
| 8 | 233 | 84 | 18 | 10 | 235 | 84 | 18 | 14 | 239 | 84 | 18 | 7 | 241 | 84 | 18 |
| 11 | 335 | 102 | 18 | 13 | 337 | 102 | 18 | 8 | 341 | 102 | 18 | 10 | 343 | 102 | 18 |
| 14 | 455 | 120 | 18 | 16 | 457 | 120 | 18 | 11 | 461 | 120 | 18 | 13 | 463 | 120 | 18 |
| 17 | 593 | 138 | 18 | 19 | 595 | 138 | 18 | 23 | 599 | 138 | 18 | 7 | 601 | 138 | 18 |
| 20 | 749 | 156 | 18 | 13 | 751 | 156 | 18 | 17 | 755 | 156 | 18 | 19 | 757 | 156 | 18 |
| 14 | 923 | 174 | 18 | 16 | 925 | 174 | 18 | 20 | 929 | 174 | 18 | 13 | 931 | 174 | 18 |
| 8 | 1115 | 192 | 18 | 10 | 1117 | 192 | 18 | 5 | 1121 | 192 | 18 | 7 | 1123 | 192 | 18 |
| 11 | 1325 | 210 | 18 | 13 | 1327 | 210 | 18 | 8 | 1331 | 210 | 18 | 10 | 1333 | 210 | 18 |
| 14 | 1553 | 228 | 18 | 16 | 1555 | 228 | 18 | 20 | 1559 | 228 | 18 | 13 | 1561 | 228 | 18 |
| 26 | 1799 | 246 | 18 | 10 | 1801 | 246 | 18 | 14 | 1805 | 246 | 18 | 16 | 1807 | 246 | 18 |
| 11 | 2063 | 264 | 18 | 13 | 2065 | 264 | 18 | 17 | 2069 | 264 | 18 | 10 | 2071 | 264 | 18 |
| 14 | 2345 | 282 | 18 | 16 | 2347 | 282 | 18 | 11 | 2351 | 282 | 18 | 13 | 2353 | 282 | 18 |
| 17 | 2645 | 300 | 18 | 19 | 2647 | 300 | 18 | 14 | 2651 | 300 | 18 | 16 | 2653 | 300 | 18 |
| 20 | 2963 | 318 | 18 | 22 | 2965 | 318 | 18 | 26 | 2969 | 318 | 18 | 19 | 2971 | 318 | 18 |
| 23 | 3299 | 336 | 18 | 7 | 3301 | 336 | 18 | 11 | 3305 | 336 | 18 | 13 | 3307 | 336 | 18 |
| 17 | 3653 | 354 | 18 | 19 | 3655 | 354 | 18 | 23 | 3659 | 354 | 18 | 16 | 3661 | 354 | 18 |
| 11 | 4025 | 372 | 18 | 13 | 4027 | 372 | 18 | 8 | 4031 | 372 | 18 | 10 | 4033 | 372 | 18 |
| 14 | 4415 | 390 | 18 | 16 | 4417 | 390 | 18 | 11 | 4421 | 390 | 18 | 13 | 4423 | 390 | 18 |
| 17 | 4823 | 408 | 18 | 19 | 4825 | 408 | 18 | 23 | 4829 | 408 | 18 | 16 | 4831 | 408 | 18 |
| 20 | 5249 | 426 | 18 | 13 | 5251 | 426 | 18 | 17 | 5255 | 426 | 18 | 19 | 5257 | 426 | 18 |
| 23 | 5693 | 444 | 18 | 25 | 5695 | 444 | 18 | 29 | 5699 | 444 | 18 | 13 | 5701 | 444 | 18 |
| 17 | 6155 | 462 | 18 | 19 | 6157 | 462 | 18 | 14 | 6161 | 462 | 18 | 16 | 6163 | 462 | 18 |
| 20 | 6635 | 480 | 18 | 22 | 6637 | 480 | 18 | 17 | 6641 | 480 | 18 | 19 | 6643 | 480 | 18 |

| | A5 | | | | A6 | | | | A7 | | | | A8 | | |
|---|---|---|---|---|---|---|---|---|---|---|---|---|---|---|---|
| SD | A5 | A5' | A5" | SD | A6 | A6' | A6" | SD | A7 | A7' | A7" | SD | A8 | A8' | A8" |
| 5 | 5 | | | 7 | 7 | | | 2 | 11 | | | 4 | 13 | | |
| 8 | 17 | 12 | | 10 | 19 | 12 | | 5 | 23 | 12 | | 7 | 25 | 12 | |
| 11 | 47 | 30 | 18 | 13 | 49 | 30 | 18 | 8 | 53 | 30 | 18 | 10 | 55 | 30 | 18 |
| 14 | 95 | 48 | 18 | 16 | 97 | 48 | 18 | 2 | 101 | 48 | 18 | 4 | 103 | 48 | 18 |
| 8 | 161 | 66 | 18 | 10 | 163 | 66 | 18 | 14 | 167 | 66 | 18 | 16 | 169 | 66 | 18 |
| 11 | 245 | 84 | 18 | 13 | 247 | 84 | 18 | 8 | 251 | 84 | 18 | 10 | 253 | 84 | 18 |
| 14 | 347 | 102 | 18 | 16 | 349 | 102 | 18 | 11 | 353 | 102 | 18 | 13 | 355 | 102 | 18 |
| 17 | 467 | 120 | 18 | 19 | 469 | 120 | 18 | 14 | 473 | 120 | 18 | 16 | 475 | 120 | 18 |
| 11 | 605 | 138 | 18 | 13 | 607 | 138 | 18 | 8 | 611 | 138 | 18 | 10 | 613 | 138 | 18 |
| 14 | 761 | 156 | 18 | 16 | 763 | 156 | 18 | 20 | 767 | 156 | 18 | 22 | 769 | 156 | 18 |
| 17 | 935 | 174 | 18 | 19 | 937 | 174 | 18 | 14 | 941 | 174 | 18 | 16 | 943 | 174 | 18 |
| 11 | 1127 | 192 | 18 | 13 | 1129 | 192 | 18 | 8 | 1133 | 192 | 18 | 10 | 1135 | 192 | 18 |
| 14 | 1337 | 210 | 18 | 16 | 1339 | 210 | 18 | 11 | 1343 | 210 | 18 | 13 | 1345 | 210 | 18 |
| 17 | 1565 | 228 | 18 | 19 | 1567 | 228 | 18 | 14 | 1571 | 228 | 18 | 16 | 1573 | 228 | 18 |
| 11 | 1811 | 246 | 18 | 13 | 1813 | 246 | 18 | 17 | 1817 | 246 | 18 | 19 | 1819 | 246 | 18 |
| 14 | 2075 | 264 | 18 | 16 | 2077 | 264 | 18 | 11 | 2081 | 264 | 18 | 13 | 2083 | 264 | 18 |
| 17 | 2357 | 282 | 18 | 19 | 2359 | 282 | 18 | 14 | 2363 | 282 | 18 | 16 | 2365 | 282 | 18 |
| 20 | 2657 | 300 | 18 | 22 | 2659 | 300 | 18 | 17 | 2663 | 300 | 18 | 19 | 2665 | 300 | 18 |
| 23 | 2975 | 318 | 18 | 25 | 2977 | 318 | 18 | 20 | 2981 | 318 | 18 | 22 | 2983 | 318 | 18 |
| 8 | 3311 | 336 | 18 | 10 | 3313 | 336 | 18 | 14 | 3317 | 336 | 18 | 16 | 3319 | 336 | 18 |
| 20 | 3665 | 354 | 18 | 22 | 3667 | 354 | 18 | 19 | 3671 | 354 | 18 | 19 | 3673 | 354 | 18 |
| 14 | 4037 | 372 | 18 | 16 | 4039 | 372 | 18 | 11 | 4043 | 372 | 18 | 13 | 4045 | 372 | 18 |
| 17 | 4427 | 390 | 18 | 19 | 4429 | 390 | 18 | 14 | 4433 | 390 | 18 | 16 | 4435 | 390 | 18 |
| 20 | 4835 | 408 | 18 | 22 | 4837 | 408 | 18 | 17 | 4841 | 408 | 18 | 19 | 4843 | 408 | 18 |
| 14 | 5261 | 426 | 18 | 16 | 5263 | 426 | 18 | 20 | 5267 | 426 | 18 | 22 | 5269 | 426 | 18 |
| 17 | 5705 | 444 | 18 | 19 | 5707 | 444 | 18 | 14 | 5711 | 444 | 18 | 16 | 5713 | 444 | 18 |
| 20 | 6167 | 462 | 18 | 22 | 6169 | 462 | 18 | 17 | 6173 | 462 | 18 | 19 | 6175 | 462 | 18 |

| | A9 | | | | A10 | | | | A11 | | | | A12 | | |
|---|---|---|---|---|---|---|---|---|---|---|---|---|---|---|---|
| SD | A9 | A9' | A9" | SD | A10 | A10' | A10" | SD | A11 | A11' | A11" | SD | A12 | A12' | A12" |
| 8 | 17 | | | 10 | 19 | | | 5 | 23 | | | 7 | 25 | | |
| 11 | 29 | 12 | | 4 | 31 | 12 | | 8 | 35 | 12 | | 10 | 37 | 12 | |
| 14 | 59 | 30 | 18 | 7 | 61 | 30 | 18 | 11 | 65 | 30 | 18 | 13 | 67 | 30 | 18 |
| 8 | 107 | 48 | 18 | 10 | 109 | 48 | 18 | 5 | 113 | 48 | 18 | 7 | 115 | 48 | 18 |
| 11 | 173 | 66 | 18 | 13 | 175 | 66 | 18 | 17 | 179 | 66 | 18 | 10 | 181 | 66 | 18 |
| 14 | 257 | 84 | 18 | 16 | 259 | 84 | 18 | 11 | 263 | 84 | 18 | 13 | 265 | 84 | 18 |
| 17 | 359 | 102 | 18 | 10 | 361 | 102 | 18 | 14 | 365 | 102 | 18 | 16 | 367 | 102 | 18 |
| 20 | 479 | 120 | 18 | 13 | 481 | 120 | 18 | 17 | 485 | 120 | 18 | 19 | 487 | 120 | 18 |
| 14 | 617 | 138 | 18 | 16 | 619 | 138 | 18 | 11 | 623 | 138 | 18 | 13 | 625 | 138 | 18 |
| 17 | 773 | 156 | 18 | 19 | 775 | 156 | 18 | 23 | 779 | 156 | 18 | 16 | 781 | 156 | 18 |
| 20 | 947 | 174 | 18 | 22 | 949 | 174 | 18 | 17 | 953 | 174 | 18 | 19 | 955 | 174 | 18 |
| 14 | 1139 | 192 | 18 | 7 | 1141 | 192 | 18 | 11 | 1145 | 192 | 18 | 13 | 1147 | 192 | 18 |
| 17 | 1349 | 210 | 18 | 10 | 1351 | 210 | 18 | 14 | 1355 | 210 | 18 | 16 | 1357 | 210 | 18 |
| 20 | 1577 | 228 | 18 | 22 | 1579 | 228 | 18 | 17 | 1583 | 228 | 18 | 19 | 1585 | 228 | 18 |
| 14 | 1823 | 246 | 18 | 16 | 1825 | 246 | 18 | 20 | 1829 | 246 | 18 | 13 | 1831 | 246 | 18 |
| 17 | 2087 | 264 | 18 | 19 | 2089 | 264 | 18 | 14 | 2093 | 264 | 18 | 16 | 2095 | 264 | 18 |
| 20 | 2369 | 282 | 18 | 13 | 2371 | 282 | 18 | 17 | 2375 | 282 | 18 | 19 | 2377 | 282 | 18 |
| 23 | 2669 | 300 | 18 | 16 | 2671 | 300 | 18 | 20 | 2675 | 300 | 18 | 22 | 2677 | 300 | 18 |
| 26 | 2987 | 318 | 18 | 28 | 2989 | 318 | 18 | 23 | 2993 | 318 | 18 | 25 | 2995 | 318 | 18 |
| 11 | 3323 | 336 | 18 | 13 | 3325 | 336 | 18 | 17 | 3329 | 336 | 18 | 10 | 3331 | 336 | 18 |
| 23 | 3677 | 354 | 18 | 25 | 3679 | 354 | 18 | 20 | 3683 | 354 | 18 | 22 | 3685 | 354 | 18 |
| 17 | 4049 | 372 | 18 | 10 | 4051 | 372 | 18 | 14 | 4055 | 372 | 18 | 16 | 4057 | 372 | 18 |
| 20 | 4439 | 390 | 18 | 13 | 4441 | 390 | 18 | 17 | 4445 | 390 | 18 | 19 | 4447 | 390 | 18 |
| 23 | 4847 | 408 | 18 | 25 | 4849 | 408 | 18 | 20 | 4853 | 408 | 18 | 22 | 4855 | 408 | 18 |
| 17 | 5273 | 426 | 18 | 19 | 5275 | 426 | 18 | 23 | 5279 | 426 | 18 | 16 | 5281 | 426 | 18 |
| 20 | 5717 | 444 | 18 | 22 | 5719 | 444 | 18 | 17 | 5723 | 444 | 18 | 19 | 5725 | 444 | 18 |
| 23 | 6179 | 462 | 18 | 16 | 6181 | 462 | 18 | 20 | 6185 | 462 | 18 | 22 | 6187 | 462 | 18 |

+ 4 →   + 2 →   + 4 → ....

+ 4 →   + 2 →   + 4 → ....

| | B1 | | | | B2 | | | | B3 | | | | B4 | | |
|---|---|---|---|---|---|---|---|---|---|---|---|---|---|---|---|
| SD | B1 | B1' | B1" | SD | B2 | B2' | B2" | SD | B3 | B3' | B3" | SD | B4 | B4' | B4" |
| 2 | 11 | | | 4 | 13 | | | 8 | 17 | | | 10 | 19 | | |
| 11 | 47 | 36 | | 13 | 49 | 36 | | 8 | 53 | 36 | | 10 | 55 | 36 | |
| 2 | 101 | 54 | 18 | 4 | 103 | 54 | 18 | 8 | 107 | 54 | 18 | 10 | 109 | 54 | 18 |
| 11 | 173 | 72 | 18 | 13 | 175 | 72 | 18 | 17 | 179 | 72 | 18 | 10 | 181 | 72 | 18 |
| 11 | 263 | 90 | 18 | 13 | 265 | 90 | 18 | 17 | 269 | 90 | 18 | 10 | 271 | 90 | 18 |
| 11 | 371 | 108 | 18 | 13 | 373 | 108 | 18 | 17 | 377 | 108 | 18 | 19 | 379 | 108 | 18 |
| 20 | 497 | 126 | 18 | 22 | 499 | 126 | 18 | 8 | 503 | 126 | 18 | 10 | 505 | 126 | 18 |
| 11 | 641 | 144 | 18 | 13 | 643 | 144 | 18 | 17 | 647 | 144 | 18 | 19 | 649 | 144 | 18 |
| 11 | 803 | 162 | 18 | 13 | 805 | 162 | 18 | 17 | 809 | 162 | 18 | 10 | 811 | 162 | 18 |
| 20 | 983 | 180 | 18 | 22 | 985 | 180 | 18 | 26 | 989 | 180 | 18 | 19 | 991 | 180 | 18 |
| 11 | 1181 | 198 | 18 | 13 | 1183 | 198 | 18 | 17 | 1187 | 198 | 18 | 19 | 1189 | 198 | 18 |
| 20 | 1397 | 216 | 18 | 22 | 1399 | 216 | 18 | 8 | 1403 | 216 | 18 | 10 | 1405 | 216 | 18 |
| 11 | 1631 | 234 | 18 | 13 | 1633 | 234 | 18 | 17 | 1637 | 234 | 18 | 19 | 1639 | 234 | 18 |
| 20 | 1883 | 252 | 18 | 22 | 1885 | 252 | 18 | 26 | 1889 | 252 | 18 | 19 | 1891 | 252 | 18 |
| 11 | 2153 | 270 | 18 | 13 | 2155 | 270 | 18 | 17 | 2159 | 270 | 18 | 10 | 2161 | 270 | 18 |
| 11 | 2441 | 288 | 18 | 13 | 2443 | 288 | 18 | 17 | 2447 | 288 | 18 | 19 | 2449 | 288 | 18 |
| 20 | 2747 | 306 | 18 | 22 | 2749 | 306 | 18 | 17 | 2753 | 306 | 18 | 19 | 2755 | 306 | 18 |
| 11 | 3071 | 324 | 18 | 13 | 3073 | 324 | 18 | 17 | 3077 | 324 | 18 | 19 | 3079 | 324 | 18 |
| 11 | 3413 | 342 | 18 | 13 | 3415 | 342 | 18 | 17 | 3419 | 342 | 18 | 10 | 3421 | 342 | 18 |
| 20 | 3773 | 360 | 18 | 22 | 3775 | 360 | 18 | 26 | 3779 | 360 | 18 | 19 | 3781 | 360 | 18 |
| 11 | 4151 | 378 | 18 | 13 | 4153 | 378 | 18 | 17 | 4157 | 378 | 18 | 19 | 4159 | 378 | 18 |
| 20 | 4547 | 396 | 18 | 22 | 4549 | 396 | 18 | 17 | 4553 | 396 | 18 | 19 | 4555 | 396 | 18 |
| 20 | 4961 | 414 | 18 | 22 | 4963 | 414 | 18 | 26 | 4967 | 414 | 18 | 28 | 4969 | 414 | 18 |
| 20 | 5393 | 432 | 18 | 22 | 5395 | 432 | 18 | 26 | 5399 | 432 | 18 | 10 | 5401 | 432 | 18 |
| 20 | 5843 | 450 | 18 | 22 | 5845 | 450 | 18 | 26 | 5849 | 450 | 18 | 19 | 5851 | 450 | 18 |
| 11 | 6311 | 468 | 18 | 13 | 6313 | 468 | 18 | 17 | 6317 | 468 | 18 | 19 | 6319 | 468 | 18 |
| 29 | 6797 | 486 | 18 | 31 | 6799 | 486 | 18 | 17 | 6803 | 486 | 18 | 19 | 6805 | 486 | 18 |

| | B5 | | | | B6 | | | | B7 | | | | B8 | | |
|---|---|---|---|---|---|---|---|---|---|---|---|---|---|---|---|
| SD | B5 | B5' | B5" | SD | B6 | B6' | B6" | SD | B7 | B7' | B7" | SD | B8 | B8' | B8" |
| 5 | 5 | | | 7 | 7 | | | 2 | 11 | | | 4 | 13 | | |
| 5 | 23 | 18 | | 7 | 25 | 18 | | 11 | 29 | 18 | | 4 | 31 | 18 | |
| 14 | 59 | 36 | 18 | 7 | 61 | 36 | 18 | 11 | 65 | 36 | 18 | 13 | 67 | 36 | 18 |
| 5 | 113 | 54 | 18 | 7 | 115 | 54 | 18 | 11 | 119 | 54 | 18 | 4 | 121 | 54 | 18 |
| 14 | 185 | 72 | 18 | 16 | 187 | 72 | 18 | 11 | 191 | 72 | 18 | 13 | 193 | 72 | 18 |
| 14 | 275 | 90 | 18 | 16 | 277 | 90 | 18 | 11 | 281 | 90 | 18 | 13 | 283 | 90 | 18 |
| 14 | 383 | 108 | 18 | 16 | 385 | 108 | 18 | 20 | 389 | 108 | 18 | 13 | 391 | 108 | 18 |
| 14 | 509 | 126 | 18 | 7 | 511 | 126 | 18 | 11 | 515 | 126 | 18 | 13 | 517 | 126 | 18 |
| 14 | 653 | 144 | 18 | 16 | 655 | 144 | 18 | 20 | 659 | 144 | 18 | 13 | 661 | 144 | 18 |
| 14 | 815 | 162 | 18 | 16 | 817 | 162 | 18 | 11 | 821 | 162 | 18 | 13 | 823 | 162 | 18 |
| 23 | 995 | 180 | 18 | 25 | 997 | 180 | 18 | 2 | 1001 | 180 | 18 | 4 | 1003 | 180 | 18 |
| 14 | 1193 | 198 | 18 | 16 | 1195 | 198 | 18 | 20 | 1199 | 198 | 18 | 4 | 1201 | 198 | 18 |
| 14 | 1409 | 216 | 18 | 7 | 1411 | 216 | 18 | 11 | 1415 | 216 | 18 | 13 | 1417 | 216 | 18 |
| 14 | 1643 | 234 | 18 | 16 | 1645 | 234 | 18 | 20 | 1649 | 234 | 18 | 13 | 1651 | 234 | 18 |
| 23 | 1895 | 252 | 18 | 25 | 1897 | 252 | 18 | 11 | 1901 | 252 | 18 | 13 | 1903 | 252 | 18 |
| 14 | 2165 | 270 | 18 | 16 | 2167 | 270 | 18 | 11 | 2171 | 270 | 18 | 13 | 2173 | 270 | 18 |
| 14 | 2453 | 288 | 18 | 16 | 2455 | 288 | 18 | 20 | 2459 | 288 | 18 | 13 | 2461 | 288 | 18 |
| 23 | 2759 | 306 | 18 | 16 | 2761 | 306 | 18 | 20 | 2765 | 306 | 18 | 22 | 2767 | 306 | 18 |
| 14 | 3083 | 324 | 18 | 16 | 3085 | 324 | 18 | 20 | 3089 | 324 | 18 | 13 | 3091 | 324 | 18 |
| 14 | 3425 | 342 | 18 | 16 | 3427 | 342 | 18 | 11 | 3431 | 342 | 18 | 13 | 3433 | 342 | 18 |
| 23 | 3785 | 360 | 18 | 25 | 3787 | 360 | 18 | 20 | 3791 | 360 | 18 | 22 | 3793 | 360 | 18 |
| 14 | 4163 | 378 | 18 | 16 | 4165 | 378 | 18 | 20 | 4169 | 378 | 18 | 13 | 4171 | 378 | 18 |
| 23 | 4559 | 396 | 18 | 16 | 4561 | 396 | 18 | 20 | 4565 | 396 | 18 | 22 | 4567 | 396 | 18 |
| 23 | 4973 | 414 | 18 | 25 | 4975 | 414 | 18 | 29 | 4979 | 414 | 18 | 22 | 4981 | 414 | 18 |
| 14 | 5405 | 432 | 18 | 16 | 5407 | 432 | 18 | 11 | 5411 | 432 | 18 | 13 | 5413 | 432 | 18 |
| 23 | 5855 | 450 | 18 | 25 | 5857 | 450 | 18 | 20 | 5861 | 450 | 18 | 22 | 5863 | 450 | 18 |
| 14 | 6323 | 468 | 18 | 16 | 6325 | 468 | 18 | 20 | 6329 | 468 | 18 | 13 | 6331 | 468 | 18 |

| | B9 | | | | B10 | | | | B11 | | | | B12 | | |
|---|---|---|---|---|---|---|---|---|---|---|---|---|---|---|---|
| SD | B9 | B9' | B9" | SD | B10 | B10' | B10" | SD | B11 | B11' | B11" | SD | B12 | B12' | B12" |
| 8 | 17 | | | 10 | 19 | | | 5 | 23 | | | 5 | 25 | | |
| 8 | 35 | 18 | | 10 | 37 | 18 | | 5 | 41 | 18 | | 7 | 43 | 18 | |
| 8 | 71 | 36 | 18 | 10 | 73 | 36 | 18 | 14 | 77 | 36 | 18 | 16 | 79 | 36 | 18 |
| 8 | 125 | 54 | 18 | 10 | 127 | 54 | 18 | 5 | 131 | 54 | 18 | 7 | 133 | 54 | 18 |
| 17 | 197 | 72 | 18 | 19 | 199 | 72 | 18 | 5 | 203 | 72 | 18 | 7 | 205 | 72 | 18 |
| 17 | 287 | 90 | 18 | 19 | 289 | 90 | 18 | 14 | 293 | 90 | 18 | 16 | 295 | 90 | 18 |
| 17 | 395 | 108 | 18 | 19 | 397 | 108 | 18 | 5 | 401 | 108 | 18 | 7 | 403 | 108 | 18 |
| 8 | 521 | 126 | 18 | 10 | 523 | 126 | 18 | 14 | 527 | 126 | 18 | 16 | 529 | 126 | 18 |
| 17 | 665 | 144 | 18 | 19 | 667 | 144 | 18 | 14 | 671 | 144 | 18 | 16 | 673 | 144 | 18 |
| 17 | 827 | 162 | 18 | 19 | 829 | 162 | 18 | 14 | 833 | 162 | 18 | 16 | 835 | 162 | 18 |
| 8 | 1007 | 180 | 18 | 10 | 1009 | 180 | 18 | 5 | 1013 | 180 | 18 | 7 | 1015 | 180 | 18 |
| 8 | 1205 | 198 | 18 | 10 | 1207 | 198 | 18 | 5 | 1211 | 198 | 18 | 7 | 1213 | 198 | 18 |
| 8 | 1421 | 216 | 18 | 10 | 1423 | 216 | 18 | 14 | 1427 | 216 | 18 | 16 | 1429 | 216 | 18 |
| 17 | 1655 | 234 | 18 | 19 | 1657 | 234 | 18 | 14 | 1661 | 234 | 18 | 16 | 1663 | 234 | 18 |
| 17 | 1907 | 252 | 18 | 19 | 1909 | 252 | 18 | 14 | 1913 | 252 | 18 | 16 | 1915 | 252 | 18 |
| 17 | 2177 | 270 | 18 | 19 | 2179 | 270 | 18 | 14 | 2183 | 270 | 18 | 16 | 2185 | 270 | 18 |
| 17 | 2465 | 288 | 18 | 19 | 2467 | 288 | 18 | 14 | 2471 | 288 | 18 | 16 | 2473 | 288 | 18 |
| 17 | 2771 | 306 | 18 | 19 | 2773 | 306 | 18 | 23 | 2777 | 306 | 18 | 25 | 2779 | 306 | 18 |
| 17 | 3095 | 324 | 18 | 19 | 3097 | 324 | 18 | 5 | 3101 | 324 | 18 | 7 | 3103 | 324 | 18 |
| 17 | 3437 | 342 | 18 | 19 | 3439 | 342 | 18 | 14 | 3443 | 342 | 18 | 16 | 3445 | 342 | 18 |
| 26 | 3797 | 360 | 18 | 28 | 3799 | 360 | 18 | 14 | 3803 | 360 | 18 | 16 | 3805 | 360 | 18 |
| 17 | 4175 | 378 | 18 | 19 | 4177 | 378 | 18 | 14 | 4181 | 378 | 18 | 16 | 4183 | 378 | 18 |
| 17 | 4571 | 396 | 18 | 19 | 4573 | 396 | 18 | 23 | 4577 | 396 | 18 | 25 | 4579 | 396 | 18 |
| 26 | 4985 | 414 | 18 | 28 | 4987 | 414 | 18 | 23 | 4991 | 414 | 18 | 25 | 4993 | 414 | 18 |
| 17 | 5417 | 432 | 18 | 19 | 5419 | 432 | 18 | 14 | 5423 | 432 | 18 | 16 | 5425 | 432 | 18 |
| 26 | 5867 | 450 | 18 | 28 | 5869 | 450 | 18 | 23 | 5873 | 450 | 18 | 25 | 5875 | 450 | 18 |
| 17 | 6335 | 468 | 18 | 19 | 6337 | 468 | 18 | 14 | 6341 | 468 | 18 | 16 | 6343 | 468 | 18 |

+ 4 →   + 2 →   + 4 → ....

+ 4 →   + 2 →   + 4 → ....

| | C1 | | | | C2 | | | | C3 | | | | C4 | | |
|---|---|---|---|---|---|---|---|---|---|---|---|---|---|---|---|
| SD | C1 | C1' | C1" | SD | C2 | C2' | C2" | SD | C3 | C3' | C3" | SD | C4 | C4' | C4" |
| 8 | 17 | | | 10 | 19 | | | 5 | 23 | | | 7 | 25 | | |
| 16 | 59 | 42 | | 7 | 61 | 42 | | 13 | 65 | 42 | | 13 | 67 | 42 | |
| 11 | 119 | 60 | 18 | 4 | 121 | 60 | 18 | 8 | 125 | 60 | 18 | 10 | 127 | 60 | 18 |
| 17 | 197 | 78 | 18 | 19 | 199 | 78 | 18 | 5 | 203 | 78 | 18 | 7 | 205 | 78 | 18 |
| 14 | 293 | 96 | 18 | 16 | 295 | 96 | 18 | 20 | 299 | 96 | 18 | 4 | 301 | 96 | 18 |
| 13 | 407 | 114 | 18 | 13 | 409 | 114 | 18 | 8 | 413 | 114 | 18 | 10 | 415 | 114 | 18 |
| 17 | 539 | 132 | 18 | 10 | 541 | 132 | 18 | 14 | 545 | 132 | 18 | 16 | 547 | 132 | 18 |
| 23 | 689 | 150 | 18 | 16 | 691 | 150 | 18 | 20 | 695 | 150 | 18 | 22 | 697 | 150 | 18 |
| 20 | 857 | 168 | 18 | 22 | 859 | 168 | 18 | 17 | 863 | 168 | 18 | 19 | 865 | 168 | 18 |
| 8 | 1043 | 186 | 18 | 10 | 1045 | 186 | 18 | 14 | 1049 | 186 | 18 | 7 | 1051 | 186 | 18 |
| 14 | 1247 | 204 | 18 | 16 | 1249 | 204 | 18 | 11 | 1253 | 204 | 18 | 13 | 1255 | 204 | 18 |
| 20 | 1469 | 222 | 18 | 13 | 1471 | 222 | 18 | 17 | 1475 | 222 | 18 | 19 | 1477 | 222 | 18 |
| 17 | 1709 | 240 | 18 | 10 | 1711 | 240 | 18 | 14 | 1715 | 240 | 18 | 16 | 1717 | 240 | 18 |
| 23 | 1967 | 258 | 18 | 25 | 1969 | 258 | 18 | 20 | 1973 | 258 | 18 | 22 | 1975 | 258 | 18 |
| 11 | 2243 | 276 | 18 | 13 | 2245 | 276 | 18 | 17 | 2249 | 276 | 18 | 10 | 2251 | 276 | 18 |
| 17 | 2537 | 294 | 18 | 19 | 2539 | 294 | 18 | 14 | 2543 | 294 | 18 | 16 | 2545 | 294 | 18 |
| 23 | 2849 | 312 | 18 | 16 | 2851 | 312 | 18 | 20 | 2855 | 312 | 18 | 22 | 2857 | 312 | 18 |
| 20 | 3179 | 330 | 18 | 13 | 3181 | 330 | 18 | 17 | 3185 | 330 | 18 | 19 | 3187 | 330 | 18 |
| 17 | 3527 | 348 | 18 | 19 | 3529 | 348 | 18 | 14 | 3533 | 348 | 18 | 16 | 3535 | 348 | 18 |
| 23 | 3893 | 366 | 18 | 25 | 3895 | 366 | 18 | 29 | 3899 | 366 | 18 | 13 | 3901 | 366 | 18 |
| 20 | 4277 | 384 | 18 | 22 | 4279 | 384 | 18 | 17 | 4283 | 384 | 18 | 19 | 4285 | 384 | 18 |
| 26 | 4679 | 402 | 18 | 19 | 4681 | 402 | 18 | 23 | 4685 | 402 | 18 | 25 | 4687 | 402 | 18 |
| 23 | 5099 | 420 | 18 | 7 | 5101 | 420 | 18 | 11 | 5105 | 420 | 18 | 13 | 5107 | 420 | 18 |
| 20 | 5537 | 438 | 18 | 22 | 5539 | 438 | 18 | 17 | 5543 | 438 | 18 | 19 | 5545 | 438 | 18 |
| 26 | 5993 | 456 | 18 | 28 | 5995 | 456 | 18 | 32 | 5999 | 456 | 18 | 7 | 6001 | 456 | 18 |
| 23 | 6467 | 474 | 18 | 25 | 6469 | 474 | 18 | 20 | 6473 | 474 | 18 | 22 | 6475 | 474 | 18 |
| 29 | 6959 | 492 | 18 | 22 | 6961 | 492 | 18 | 26 | 6965 | 492 | 18 | 28 | 6967 | 492 | 18 |

| | C5 | | | | C6 | | | | C7 | | | | C8 | | |
|---|---|---|---|---|---|---|---|---|---|---|---|---|---|---|---|
| SD | C5 | C5' | C5" | SD | C6 | C6' | C6" | SD | C7 | C7' | C7" | SD | C8 | C8' | C8" |
| 5 | 5 | | | 7 | 7 | | | 2 | 11 | | | 4 | 13 | | |
| 11 | 29 | 24 | | 4 | 31 | 24 | | 8 | 35 | 24 | | 10 | 37 | 24 | |
| 8 | 71 | 42 | 18 | 10 | 73 | 42 | 18 | 14 | 77 | 42 | 18 | 16 | 79 | 42 | 18 |
| 5 | 131 | 60 | 18 | 7 | 133 | 60 | 18 | 11 | 137 | 60 | 18 | 13 | 139 | 60 | 18 |
| 11 | 209 | 78 | 18 | 4 | 211 | 78 | 18 | 8 | 215 | 78 | 18 | 10 | 217 | 78 | 18 |
| 8 | 305 | 96 | 18 | 10 | 307 | 96 | 18 | 5 | 311 | 96 | 18 | 7 | 313 | 96 | 18 |
| 14 | 419 | 114 | 18 | 7 | 421 | 114 | 18 | 11 | 425 | 114 | 18 | 13 | 427 | 114 | 18 |
| 11 | 551 | 132 | 18 | 13 | 553 | 132 | 18 | 17 | 557 | 132 | 18 | 19 | 559 | 132 | 18 |
| 8 | 701 | 150 | 18 | 10 | 703 | 150 | 18 | 14 | 707 | 150 | 18 | 16 | 709 | 150 | 18 |
| 23 | 869 | 168 | 18 | 16 | 871 | 168 | 18 | 20 | 875 | 168 | 18 | 22 | 877 | 168 | 18 |
| 11 | 1055 | 186 | 18 | 13 | 1057 | 186 | 18 | 8 | 1061 | 186 | 18 | 10 | 1063 | 186 | 18 |
| 17 | 1259 | 204 | 18 | 10 | 1261 | 204 | 18 | 14 | 1265 | 204 | 18 | 16 | 1267 | 204 | 18 |
| 14 | 1481 | 222 | 18 | 16 | 1483 | 222 | 18 | 20 | 1487 | 222 | 18 | 22 | 1489 | 222 | 18 |
| 11 | 1721 | 240 | 18 | 13 | 1723 | 240 | 18 | 17 | 1727 | 240 | 18 | 19 | 1729 | 240 | 18 |
| 26 | 1979 | 258 | 18 | 19 | 1981 | 258 | 18 | 23 | 1985 | 258 | 18 | 25 | 1987 | 258 | 18 |
| 14 | 2255 | 276 | 18 | 16 | 2257 | 276 | 18 | 11 | 2261 | 276 | 18 | 13 | 2263 | 276 | 18 |
| 20 | 2549 | 294 | 18 | 13 | 2551 | 294 | 18 | 17 | 2555 | 294 | 18 | 19 | 2557 | 294 | 18 |
| 17 | 2861 | 312 | 18 | 19 | 2863 | 312 | 18 | 23 | 2867 | 312 | 18 | 25 | 2869 | 312 | 18 |
| 14 | 3191 | 330 | 18 | 16 | 3193 | 330 | 18 | 20 | 3197 | 330 | 18 | 22 | 3199 | 330 | 18 |
| 20 | 3539 | 348 | 18 | 13 | 3541 | 348 | 18 | 17 | 3545 | 348 | 18 | 19 | 3547 | 348 | 18 |
| 17 | 3905 | 366 | 18 | 19 | 3907 | 366 | 18 | 14 | 3911 | 366 | 18 | 16 | 3913 | 366 | 18 |
| 23 | 4289 | 384 | 18 | 16 | 4291 | 384 | 18 | 20 | 4295 | 384 | 18 | 22 | 4297 | 384 | 18 |
| 20 | 4691 | 402 | 18 | 22 | 4693 | 402 | 18 | 26 | 4697 | 402 | 18 | 28 | 4699 | 402 | 18 |
| 8 | 5111 | 420 | 18 | 10 | 5113 | 420 | 18 | 14 | 5117 | 420 | 18 | 16 | 5119 | 420 | 18 |
| 23 | 5549 | 438 | 18 | 16 | 5551 | 438 | 18 | 20 | 5555 | 438 | 18 | 22 | 5557 | 438 | 18 |
| 11 | 6005 | 456 | 18 | 13 | 6007 | 456 | 18 | 8 | 6011 | 456 | 18 | 10 | 6013 | 456 | 18 |
| 26 | 6479 | 474 | 18 | 19 | 6481 | 474 | 18 | 23 | 6485 | 474 | 18 | 25 | 6487 | 474 | 18 |

| | C9 | | | | C10 | | | | C11 | | | | C12 | | |
|---|---|---|---|---|---|---|---|---|---|---|---|---|---|---|---|
| SD | C9 | C9' | C9" | SD | C10 | C10' | C10" | SD | C11 | C11' | C11" | SD | C12 | C12' | C12" |
| 8 | 17 | | | 10 | 19 | | | 5 | 23 | | | 7 | 25 | | |
| 5 | 41 | 24 | | 7 | 43 | 24 | | 11 | 47 | 24 | | 13 | 49 | 24 | |
| 11 | 83 | 42 | 18 | 13 | 85 | 42 | 18 | 17 | 89 | 42 | 18 | 10 | 91 | 42 | 18 |
| 8 | 143 | 60 | 18 | 10 | 145 | 60 | 18 | 14 | 149 | 60 | 18 | 7 | 151 | 60 | 18 |
| 5 | 221 | 78 | 18 | 7 | 223 | 78 | 18 | 11 | 227 | 78 | 18 | 13 | 229 | 78 | 18 |
| 11 | 317 | 96 | 18 | 13 | 319 | 96 | 18 | 8 | 323 | 96 | 18 | 10 | 325 | 96 | 18 |
| 8 | 431 | 114 | 18 | 10 | 433 | 114 | 18 | 14 | 437 | 114 | 18 | 16 | 439 | 114 | 18 |
| 14 | 563 | 132 | 18 | 16 | 565 | 132 | 18 | 20 | 569 | 132 | 18 | 13 | 571 | 132 | 18 |
| 11 | 713 | 150 | 18 | 13 | 715 | 150 | 18 | 17 | 719 | 150 | 18 | 10 | 721 | 150 | 18 |
| 17 | 881 | 168 | 18 | 19 | 883 | 168 | 18 | 23 | 887 | 168 | 18 | 25 | 889 | 168 | 18 |
| 14 | 1067 | 186 | 18 | 16 | 1069 | 186 | 18 | 11 | 1073 | 186 | 18 | 13 | 1075 | 186 | 18 |
| 11 | 1271 | 204 | 18 | 13 | 1273 | 204 | 18 | 17 | 1277 | 204 | 18 | 19 | 1279 | 204 | 18 |
| 17 | 1493 | 222 | 18 | 19 | 1495 | 222 | 18 | 23 | 1499 | 222 | 18 | 7 | 1501 | 222 | 18 |
| 14 | 1733 | 240 | 18 | 16 | 1735 | 240 | 18 | 20 | 1739 | 240 | 18 | 13 | 1741 | 240 | 18 |
| 20 | 1991 | 258 | 18 | 22 | 1993 | 258 | 18 | 26 | 1997 | 258 | 18 | 28 | 1999 | 258 | 18 |
| 17 | 2267 | 276 | 18 | 19 | 2269 | 276 | 18 | 14 | 2273 | 276 | 18 | 16 | 2275 | 276 | 18 |
| 14 | 2561 | 294 | 18 | 16 | 2563 | 294 | 18 | 20 | 2567 | 294 | 18 | 22 | 2569 | 294 | 18 |
| 20 | 2873 | 312 | 18 | 22 | 2875 | 312 | 18 | 26 | 2879 | 312 | 18 | 19 | 2881 | 312 | 18 |
| 8 | 3203 | 330 | 18 | 10 | 3205 | 330 | 18 | 14 | 3209 | 330 | 18 | 7 | 3211 | 330 | 18 |
| 14 | 3551 | 348 | 18 | 16 | 3553 | 348 | 18 | 20 | 3557 | 348 | 18 | 22 | 3559 | 348 | 18 |
| 20 | 3917 | 366 | 18 | 22 | 3919 | 366 | 18 | 17 | 3923 | 366 | 18 | 19 | 3925 | 366 | 18 |
| 8 | 4301 | 384 | 18 | 10 | 4303 | 384 | 18 | 14 | 4307 | 384 | 18 | 16 | 4309 | 384 | 18 |
| 14 | 4703 | 402 | 18 | 16 | 4705 | 402 | 18 | 20 | 4709 | 402 | 18 | 13 | 4711 | 402 | 18 |
| 11 | 5123 | 420 | 18 | 13 | 5125 | 420 | 18 | 17 | 5129 | 420 | 18 | 10 | 5131 | 420 | 18 |
| 17 | 5561 | 438 | 18 | 19 | 5563 | 438 | 18 | 23 | 5567 | 438 | 18 | 25 | 5569 | 438 | 18 |
| 14 | 6017 | 456 | 18 | 16 | 6019 | 456 | 18 | 11 | 6023 | 456 | 18 | 13 | 6025 | 456 | 18 |
| 20 | 6491 | 474 | 18 | 22 | 6493 | 474 | 18 | 26 | 6497 | 474 | 18 | 28 | 6499 | 474 | 18 |

+ 4 →   + 2 →   + 4 → ....

+ 4 →   + 2 →   + 4 → ....

**Table 6-A2** : Quadratic Polynomials of the Spiral-Graphs belonging to the "Prime Number Spiral Systems" **P18-A, P18-B** and **P18-C** ( with the **2. Differential = 18** )

| Spiral Graph System | Spiral Graph | Number Sequence of Spiral Graph | Quadratic Polynomial 1 ( calculated with the first 3 numbers of the given sequence ) | Quadratic Polynomial 2 ( calculated with 3 numbers starting with the **2.** Number of the sequence ) | Quadratic Polynomial 3 ( calculated with 3 numbers starting with the **3.** Number of the sequence ) | Quadratic Polynomial 4 ( calculated with 3 numbers starting with the **4.** Number of the sequence ) |
|---|---|---|---|---|---|---|
| **P18-A** | A1 | 5 , 35 , 83 , 149 , 233 , 335 ,...... | $f_1(x) = 9x^2 + 3x - 7$ | $f_2(x) = 9x^2 + 21x + 5$ | $f_3(x) = 9x^2 + 39x + 35$ | $f_4(x) = 9x^2 + 57x + 83$ |
| | A2 | 7 , 37 , 85 , 151 , 235 , 337 ,...... | $f_1(x) = 9x^2 + 3x - 5$ | $f_2(x) = 9x^2 + 21x + 7$ | $f_3(x) = 9x^2 + 39x + 37$ | $f_4(x) = 9x^2 + 57x + 85$ |
| | A3 | 11 , 41 , 89 , 155 , 239 , 341 ,...... | $f_1(x) = 9x^2 + 3x - 1$ | $f_2(x) = 9x^2 + 21x + 11$ | $f_3(x) = 9x^2 + 39x + 41$ | $f_4(x) = 9x^2 + 57x + 89$ |
| | A4 | 1 , 13 , 43 , 91 , 157 , 241 ,...... | $f_1(x) = 9x^2 - 15x + 7$ | $f_2(x) = 9x^2 + 3x + 1$ | $f_3(x) = 9x^2 + 21x + 13$ | $f_4(x) = 9x^2 + 39x + 43$ |
| | A5 | 5 , 17 , 47 , 95 , 161 , 245 ,...... | $f_1(x) = 9x^2 - 15x + 11$ | $f_2(x) = 9x^2 + 3x + 5$ | $f_3(x) = 9x^2 + 21x + 17$ | $f_4(x) = 9x^2 + 39x + 47$ |
| | A6 | 7 , 19 , 49 , 97 , 163 , 247 ,...... | $f_1(x) = 9x^2 - 15x + 13$ | $f_2(x) = 9x^2 + 3x + 7$ | $f_3(x) = 9x^2 + 21x + 19$ | $f_4(x) = 9x^2 + 39x + 49$ |
| | A7 | 11 , 23 , 53 , 101 , 167 , 251 ,...... | $f_1(x) = 9x^2 - 15x + 17$ | $f_2(x) = 9x^2 + 3x + 11$ | $f_3(x) = 9x^2 + 21x + 23$ | $f_4(x) = 9x^2 + 39x + 53$ |
| | A8 | 13 , 25 , 55 , 103 , 169 , 253 ,...... | $f_1(x) = 9x^2 - 15x + 19$ | $f_2(x) = 9x^2 + 3x + 13$ | $f_3(x) = 9x^2 + 21x + 25$ | $f_4(x) = 9x^2 + 39x + 55$ |
| | A9 | 17 , 29 , 59 , 107 , 173 , 257 ,...... | $f_1(x) = 9x^2 - 15x + 23$ | $f_2(x) = 9x^2 + 3x + 17$ | $f_3(x) = 9x^2 + 21x + 29$ | $f_4(x) = 9x^2 + 39x + 59$ |
| | A10 | 19 , 31 , 61 , 109 , 175 , 259 ,...... | $f_1(x) = 9x^2 - 15x + 25$ | $f_2(x) = 9x^2 + 3x + 19$ | $f_3(x) = 9x^2 + 21x + 31$ | $f_4(x) = 9x^2 + 39x + 61$ |
| | A11 | 23 , 35 , 65 , 113 , 179 , 263 ,...... | $f_1(x) = 9x^2 - 15x + 29$ | $f_2(x) = 9x^2 + 3x + 23$ | $f_3(x) = 9x^2 + 21x + 35$ | $f_4(x) = 9x^2 + 39x + 65$ |
| | A12 | 25 , 37 , 67 , 115 , 181 , 265 ,...... | $f_1(x) = 9x^2 - 15x + 31$ | $f_2(x) = 9x^2 + 3x + 25$ | $f_3(x) = 9x^2 + 21x + 37$ | $f_4(x) = 9x^2 + 39x + 67$ |
| **P18-B** | B1 | 11 , 47 , 101 , 173 , 263 , 371 ,...... | $f_1(x) = 9x^2 + 9x - 7$ | $f_2(x) = 9x^2 + 27x + 11$ | $f_3(x) = 9x^2 + 45x + 47$ | $f_4(x) = 9x^2 + 63x + 101$ |
| | B2 | 13 , 49 , 103 , 175 , 265 , 373 ,...... | $f_1(x) = 9x^2 + 9x - 5$ | $f_2(x) = 9x^2 + 27x + 13$ | $f_3(x) = 9x^2 + 45x + 49$ | $f_4(x) = 9x^2 + 63x + 103$ |
| | B3 | 17 , 53 , 107 , 179 , 269 , 377 ,...... | $f_1(x) = 9x^2 + 9x - 1$ | $f_2(x) = 9x^2 + 27x + 17$ | $f_3(x) = 9x^2 + 45x + 53$ | $f_4(x) = 9x^2 + 63x + 107$ |
| | B4 | 1 , 19 , 55 , 109 , 181 , 271 ,...... | $f_1(x) = 9x^2 - 9x + 1$ | $f_2(x) = 9x^2 + 9x + 1$ | $f_3(x) = 9x^2 + 27x + 19$ | $f_4(x) = 9x^2 + 45x + 55$ |
| | B5 | 5 , 23 , 59 , 113 , 185 , 275 ,...... | $f_1(x) = 9x^2 - 9x + 5$ | $f_2(x) = 9x^2 + 9x + 5$ | $f_3(x) = 9x^2 + 27x + 23$ | $f_4(x) = 9x^2 + 45x + 59$ |
| | B6 | 7 , 25 , 61 , 115 , 187 , 277 ,...... | $f_1(x) = 9x^2 - 9x + 7$ | $f_2(x) = 9x^2 + 9x + 7$ | $f_3(x) = 9x^2 + 27x + 25$ | $f_4(x) = 9x^2 + 45x + 61$ |
| | B7 | 11 , 29 , 65 , 119 , 191 , 281 ,...... | $f_1(x) = 9x^2 - 9x + 11$ | $f_2(x) = 9x^2 + 9x + 11$ | $f_3(x) = 9x^2 + 27x + 29$ | $f_4(x) = 9x^2 + 45x + 65$ |
| | B8 | 13 , 31 , 67 , 121 , 193 , 283 ,...... | $f_1(x) = 9x^2 - 9x + 13$ | $f_2(x) = 9x^2 + 9x + 13$ | $f_3(x) = 9x^2 + 27x + 31$ | $f_4(x) = 9x^2 + 45x + 67$ |
| | B9 | 17 , 35 , 71 , 125 , 197 , 287 ,...... | $f_1(x) = 9x^2 - 9x + 17$ | $f_2(x) = 9x^2 + 9x + 17$ | $f_3(x) = 9x^2 + 27x + 35$ | $f_4(x) = 9x^2 + 45x + 71$ |
| | B10 | 19 , 37 , 73 , 127 , 199 , 289 ,...... | $f_1(x) = 9x^2 - 9x + 19$ | $f_2(x) = 9x^2 + 9x + 19$ | $f_3(x) = 9x^2 + 27x + 37$ | $f_4(x) = 9x^2 + 45x + 73$ |
| | B11 | 23 , 41 , 77 , 131 , 203 , 293 ,...... | $f_1(x) = 9x^2 - 9x + 23$ | $f_2(x) = 9x^2 + 9x + 23$ | $f_3(x) = 9x^2 + 27x + 41$ | $f_4(x) = 9x^2 + 45x + 77$ |
| | B12 | 25 , 43 , 79 , 133 , 205 , 295 ,...... | $f_1(x) = 9x^2 - 9x + 25$ | $f_2(x) = 9x^2 + 9x + 25$ | $f_3(x) = 9x^2 + 27x + 43$ | $f_4(x) = 9x^2 + 45x + 79$ |
| **P18-C** | C1 | 17 , 59 , 119 , 197 , 293 , 407 ,...... | $f_1(x) = 9x^2 + 15x - 7$ | $f_2(x) = 9x^2 + 33x + 17$ | $f_3(x) = 9x^2 + 51x + 59$ | $f_4(x) = 9x^2 + 69x + 119$ |
| | C2 | 19 , 61 , 121 , 199 , 295 , 409 ,...... | $f_1(x) = 9x^2 + 15x - 5$ | $f_2(x) = 9x^2 + 33x + 19$ | $f_3(x) = 9x^2 + 51x + 61$ | $f_4(x) = 9x^2 + 69x + 121$ |
| | C3 | 23 , 65 , 125 , 203 , 299 , 413 ,...... | $f_1(x) = 9x^2 + 15x - 1$ | $f_2(x) = 9x^2 + 33x + 23$ | $f_3(x) = 9x^2 + 51x + 65$ | $f_4(x) = 9x^2 + 69x + 125$ |
| | C4 | 1 , 25 , 67 , 127 , 205 , 301 ,...... | $f_1(x) = 9x^2 - 3x - 5$ | $f_2(x) = 9x^2 + 15x + 1$ | $f_3(x) = 9x^2 + 33x + 25$ | $f_4(x) = 9x^2 + 51x + 67$ |
| | C5 | 5 , 29 , 71 , 131 , 209 , 305 ,...... | $f_1(x) = 9x^2 - 3x - 1$ | $f_2(x) = 9x^2 + 15x + 5$ | $f_3(x) = 9x^2 + 33x + 29$ | $f_4(x) = 9x^2 + 51x + 71$ |
| | C6 | 1 , 7 , 31 , 73 , 133 , 211 ,...... | $f_1(x) = 9x^2 - 21x + 13$ | $f_2(x) = 9x^2 - 3x + 1$ | $f_3(x) = 9x^2 + 15x + 7$ | $f_4(x) = 9x^2 + 33x + 31$ |
| | C7 | 5 , 11 , 35 , 77 , 137 , 215 ,...... | $f_1(x) = 9x^2 - 21x + 17$ | $f_2(x) = 9x^2 - 3x + 5$ | $f_3(x) = 9x^2 + 15x + 11$ | $f_4(x) = 9x^2 + 33x + 35$ |
| | C8 | 7 , 13 , 37 , 79 , 139 , 217 ,...... | $f_1(x) = 9x^2 - 21x + 19$ | $f_2(x) = 9x^2 - 3x + 7$ | $f_3(x) = 9x^2 + 15x + 13$ | $f_4(x) = 9x^2 + 33x + 37$ |
| | C9 | 11 , 17 , 41 , 83 , 143 , 221 ,...... | $f_1(x) = 9x^2 - 21x + 23$ | $f_2(x) = 9x^2 - 3x + 11$ | $f_3(x) = 9x^2 + 15x + 17$ | $f_4(x) = 9x^2 + 33x + 41$ |
| | C10 | 13 , 19 , 43 , 85 , 145 , 223 ,...... | $f_1(x) = 9x^2 - 21x + 25$ | $f_2(x) = 9x^2 - 3x + 13$ | $f_3(x) = 9x^2 + 15x + 19$ | $f_4(x) = 9x^2 + 33x + 43$ |
| | C11 | 17 , 23 , 47 , 89 , 149 , 227 ,...... | $f_1(x) = 9x^2 - 21x + 29$ | $f_2(x) = 9x^2 - 3x + 17$ | $f_3(x) = 9x^2 + 15x + 23$ | $f_4(x) = 9x^2 + 33x + 47$ |
| | C12 | 19 , 25 , 49 , 91 , 151 , 229 ,...... | $f_1(x) = 9x^2 - 21x + 31$ | $f_2(x) = 9x^2 - 3x + 19$ | $f_3(x) = 9x^2 + 15x + 25$ | $f_4(x) = 9x^2 + 33x + 49$ |

## 5 The share of Prime Numbers in the analysed number sequences

To find out how high the share of prime numbers really is, in the prime-number-sequences shown in Table 6-A1to 6-C1, I have carried out a random "spot check" in a few of these sequences.

For these analysis I chose the three Prime Number Spiral Graphs (-sequences ) **B3** , **Q3** and **P20-G1** ( see FIG 6-A to 6-C and Table 6-A1to 6-C1 ) , because of their high share in prime numbers at the beginning of the sequence.

With the help of an excel-table I then extended these three sequences up to numbers >2,500,000,000 → **see Table 1**

Then I picked out a longer section at the beginning and four sections out of the further run of these extended sequences, to analyse them in regards to their share in Prime Numbers.

For the other four sections I picked sections of 8 numbers out of the number areas : 2,500,000 ; 25,000,000 ; 250,000,000 and 2,500,000,000 for each of the three number-sequences.

I then marked the Prime Numbers in each sequence with yellow color. I also marked numbers which are not divisible by 2, 3 or 5 with red color, in the different sections of these sequences.

By looking over Table 1 it is obvious that the share in prime numbers is dropping after the beginning of the sequences, which has a share in prime numbers of around 70 – 75 %.

But it seems that the share in prime number is striving for around 25 – 35 % in the long run. The distribution of prime numbers in the number area 0 – 2,500,000,000 seems to be relatively evenly.

However there are also sections where only few prime numbers occur ( e.g. in the second sections in sequence Q3 and P20-G1, in the number area 2,500,000)

In any case, it would be worth to further extend the three prime-numbers sequences B3, Q3 and P20-G1 shown in Table 1, and a few other sequences which also have high shares in prime numbers ( e.g. the number-sequences A9, C2, C5, P20-G3, P20-H1, S1 etc. → see Table 6-A1 to 6-C1 ), to analyse them for their content of prime numbers. This could be done with special software which automatically extends and analyses these sequences for prime numbers.

## 6 To the periodic occurrence of prime factors in non-prime-numbers

When I set-up **Table 1**, I noticed that some of the prime factors of the non-prime-numbers occured at equal intervals in the sequences B3, Q3 and P20-G1.
It also appeared that only defined prime factors occur in every sequence.
I have marked some of these prime factors in red color in Table 1(→ see column "prime-factors" ). For example the prime factors 17 and 23 in sequence B3 and the prime factors 23 and 71 in sequence P20-G1occur at equal intervals.

These peculiarities of some of the prime factors shown in Table 1, forced me to do an analysis in regards to a possible "periodic behavior" of the prime factors which form the non-prime-numbers in the analysed prime number spiral graphs .

Besides the three prime-number-sequences B3, Q3 and P20-G1, I further chose the sequences K5, S1 and B33 for this analysis.
The two sequences S1 and B33 where chosen because they have the same number-endings-sequences and the same "sums of the digits"-sequences as the number sequences Q3 and B3 ( → see Table 2 and 3 ).
Therefore it was interesting to see what effects these similarities have on the periodic distribution of the prime-factors in these sequences.

Here the sequence B33 was developed, by a further extension of Table 6-A1on the righthand side , by adding alternately the numbers 2 or 4.

The number-sequence K5 was chosen, because it contains many small prime-factors at the beginning of the sequence. ( This property was used to explain the origin-principle of the prime-number-sequences shown in FIG 6-A to 6-C in a graphic way ( in FIG. 7 ). I will come back to this point later ! )

The results of this analysis regarding the periodic behavior of some small prime-factors are shown in Table 4 , Table 5-A and Table 5-B. ( → see next pages ! )

These tables give a good insight into the periodic behavior of the smallest prime-factors contained in the non-prime-numbers of the chosen sequences.

The following list describes the most remarkable properties of these periodic occurring prime-factors : ( → see Table 4, 5-A and 5-B )

1.) - All prime-factors, which form the non-prime-numbers in the analysed number-sequences, recur in defined periods. And this principle seems to apply to all prime-number-sequences derived from the Spiral-Graphs shown in FIG 6-A to FIG 6-C.

2.) - The following general rule applies for these periods :
The smaller the prime-factor, the smaller is the period in which the prime-factor occurs in each number sequence.

3.) - Further the following remarkable rule applies :

The period-length of every prime-factor, expressed in spacings ( lines in the table ), or the sum of the period-lengths (if there are 2 different ones ! ), is identical to the value of the prime-factor !

For example : Prime Factor 13 in number-sequence " B3 " occurs periodically in every 13th line of the Table 4.
Or Prime Factor 23 occurs with the two alternating period-lengths 2 and 21, which add up to 23, which again is exactly equal to the prime factor itself ! etc.

It is notable that the number of spacings between two numbers ( → lines in the table between two numbers ) corresponds with the same number of "winds" of the Square-Root-Spiral , which lie between these two numbers !

4.) - By comparing the prime factors, which occur in the non-prime numbers of the analysed number-sequences, it is notable that only specific smallest prime-factors occur in each number sequence ( → see columns " prime-factors of non-prime-numbers" in Table 4, 5-A and 5-B )

For example these are the following prime-factors :

| | | |
|---|---|---|
| In number-sequence B3 | : | 13, 17, 23, 29, 43, 53, 61,... |
| " " " **Q3** | **:** | **11, 13, 31, 37, 73, 89,....** |
| " " " P20-G1 | : | 13, 23, 31, 67, 71, 73,.... |
| " " " K5 | : | 7, 11, 13, 17, 29, 37,.... |
| " " " **S1** | **:** | **11, 13, 31, 37, 73, 89,....** |
| " " " B33 | : | 11, 13, 29, 43, 53, 59, 61.... |

This is a clear indication for a "blueprint" which controls the composition of the non-prime-numbers in each sequence. In each number-sequence specific smallest prime-factors are completely absent. This leads to the following conclusion :

" The none-prime-numbers in each number-sequence ( shown in Table 6-A1 to 6-C1 ) are formed by a defined number of specific prime-factors, which recur in defined periods ! "

5.) - A comparison of the prime- factor distribution in the number-sequences **Q3** and **S1** clearly shows, that the non-prime-numbers in these two sequences **consist of exactly the same prime-factors !! Even the periods** in which these prime factors occur **are exactly the same !!** The only difference between the non-prime numbers in these two number-sequences is that the non-prime numbers are made of different combinations of prime-factors and that they are distributed in a different way in these two number-sequence.

However this comparison shows that number-sequences of **different spiral-graph-systems**, which have the same number endings sequence and the same "sums of the digits sequence" ( see Table 2 and 3 ), seem to contain exactly the same prime-factors with the same periods in their non-prime-numbers !!!

But a comparison of the two number-sequences B3 and B33, which belong to the **same spiral-graph system** P18-B ( see FIG 6-A ), and which also have the same number endings sequences and the same sums of the digits sequences do not contain exactly the same prime-factors with the same periods !

As mentioned on the last page, I have developed the number-sequence B33 through a further extension of Table 6-A1on the righthand side , by adding alternately the numbers 2 or 4.
I did this to have another pair of sequences ( B3 and B33 ) which according to Table 2 and 3 also have the same number endings sequences and the same sums of the digits sequences".

But comparing the prime-factors, which occur in Table 4 and 5-B in the non-prime-numbers of the two sequences B3 and B33, doesn't show the same matching of prime-factors and periods as in the sequences Q3 and S1.
The most prime factors in the sequences B3 and B33 may be the same, but the periods in which they occur are clearly different. And it is also noticeable that both sequences also contain a few different prime-factors.

## 7 Graphic explanation of the origin of the periodic occuring prime factors

Now I want to show a "graphic explanation" for the origin of the periodic occurring prime factors described in Table 4, 5-A and 5-B.

For this please have a look at → **FIG 7**

Here the spiral-graph **K5** is drawn in **red** color. ( → see also FIG 6-C )
I have chosen the spiral-graph K5 because it contains many non-prime-numbers at the start.
Because of that it can be demonstrated, that the three non-prime-numbers **49, 187** and **289** are formed by "points of intersection" of the three "number-group-spiral-systems" which contain either numbers divisible by **7, 11** or **17**.

→ ( **Ref.:** As already mentioned in the abstract and in the introduction, all natural numbers divisible by the same prime number lie in defined **"number-group-spiral-systems"**. To get a better understanding of this property, please have a read through my introduction to the Square Root Spiral :

→ " **The ordered distribution of natural numbers on the Square Root Spiral "**

In this study the number-group-spiral-systems, which contain the numbers divisible by 7, 11 and 17 are shown, and the general rule which defines these spiral-systems is described → see **chapter 5.2** in the mentioned study ! )

For clearness I have only shown one "spiral-graph-system" of each number-group-spiral-system in FIG 7 And from each of these systems I have only shown 3 to 4 spiralarms !

Besides the spiral-graph K5, the number-group-spiral-system P1 is shown, which contains numbers divisible by 17 (blue), as well as the system N2 (orange), which contains numbers divisible by 7 , and the system N2 (pink), which contains numbers divisible by 11 ( see FIG 7 ).
Now it is easy to see in FIG 7, that the course ( curvature ) of the spiral-graph K5 is already fully defined by periodic points of intersection of the three mentioned number-group-spiral-systems ( which contain the numbers divisible by 7, 11 and 17 ) with the square root spiral.

Other defined number-group-spiral-systems also contribute with periodic points of intersection to the formation of non-prime-numbers on spiral-graph K5, but the course of spiral-graph 5 is already defined by the mentioned spiral-systems.

The **Appendix** shows a diagram of the prime-number-spiral-graph B3 with the specification of the exact polar coordinates of the points of intersection of this graph with the square-root-spiral ( → positions of the natural numbers which lie on this graph ). These coordinates might be helpful for an exact analysis of this spiral-graph !

## 8 Final Comment

Every prime-number-spiral-graph presented in FIG 6-A to 6-C, shows periodicities in the distribution of the prime factors of it's non-prime numbers !
Tables 4, 5-A and 5-B are a first proof for this proposition.

And the share of Prime Numbers as well as the distribution of Prime Numbers on a certain prime-number-spiral-graph is a result of the periodic occurring prime factors which form the non-prime numbers in this graph !

The distribution of the periodic occuring prime factors is defined by the number-group-spiral-systems which I have described in my introduction study to the Square Root Spiral. (→ see arXiv – Archiv ) The Title of this study is :

" **The ordered distribution of natural numbers on the Square Root Spiral " [ 1 ]**

The general rule, which defines the arrangement of the mentioned number-group-spiral-systems on the Square Root Spiral, can be found in **chapter 5.2** in the above mentioned study ! ( → A reference to these number-group-spiral-systems is also given on the previous page ! → see righthand side ).

Similar to the prime-number-spiral-graphs in this study , the spiralarms of the mentioned number-group-spiral-systems are also clear defined by quadratic polynomials.

The periodic occuring prime factors in the non-prime numbers of the prime-number-spiral-graphs can be explained by periodic points of intersection of certain spiralarms of the number-group-spiral-systems with the Square Root Spiral ( → see example in FIG 7 )

In this connection I also want to refer to my 3. Study which I intend to file with the arXiv – Archiv. The title of this study is :

→ " The logic of the prime number distribution" [ 3 ]

In this study the general distribution of the prime numbers is decribed with a simple "Wave Model" in a visual way.
The base of this "Wave Model" is the fact, that two Prime Number Sequences ( SQ1 & SQ2 ), which seem to contain all prime numbers, recur in themselves over and over again with increasing wave-lengths, in a very similar way as "Undertones" derive from a defined fundamental frequency **f**.

Undertones are the inversion of Overtones, which are known by every musician. Overtones ( harmonics ) are integer multiples of a fundamental frequency **f**.

The continuous recurrence of these number-sequences ( SQ1 + SQ2 ) in themselves can be considered as the principle of creation of the non-prime numbers in these two number sequences. Non-prime-numbers are created on places in these number sequences where there is interference caused by the recurrences of these number-sequences. On the other hand prime numbers represent places in these number-sequences ( SQ1 + SQ2 ) where there is no interference caused through the recurrences of this number-sequences.

The logic of this "Wave Model" is really easy to understand !
→ Please have a look at Table 2 in the mentioned paper !

The two Prime Number Sequences SQ1 & SQ2 which I mentioned, are actually easy to see in FIG 6-A. Following the winds of the square root spiral it is easy to see, that all prime numbers lie on two sequences of numbers, which are shifted to each other by two numbers, and where the difference between two successive numbers in each sequence is always 6.
By the way,… the "Prime Number Spiral Graphs" shown in FIG 6-A contain the same numbers as the mentioned Number Sequences SQ1 & SQ2 ! These are the natural numbers which are not divisible by 2 and 3.

The "periodic phenomenons" described in this work and in the other mentioned study, which are responsible for the distribution of the non-prime-numbers and prime numbers, should definitely be further analysed in more detail !

From the point of view of the Square Root Spiral, the distribution of the prime numbers seems to be clearly ordered. However this ordering principle is hard to grasp, because it is defined by the spatial complex interference of the mentioned number-group-spiral-systems , which is rapidly increasing from the centre of the square root spiral towards infinity.

In this connection I want to refer to another study, which I intend to file with the ArXiv-archive, which shows surprising similarities between the periodic behavior of the prime factors in the non-prime-numbers of the spiral-graphs shown in FIG 6-A to 6-C and the periodic behavior of prime factors which occur in Fibonacci-Numbers !

The title of this study is :

→ **" The mathematical origin of natural Fibonacci-Sequences,
and the periodic distribution of prime factors in these sequences." [ 4 ]**

Because my graphical analysis of the Square-Root-Spiral seems to open up new territory, I can't really give many references for my work.

However I found an interesting webside, which was set up by Mr. Robert Sachs.

And this webside deals with a special "Number Spiral", which is closely related to the well known Ulam-Spiral and to the Square Root Spiral as well. And on this "Number-Spiral" prime numbers also accumulate on defined graphs in a very similar fashion than shown in my study in FIG 6-A to 6-C.

The Number Spiral, analysed by Mr. Sachs, is approximately winded three times tighter in comparison with the Square Root Spiral, in a way that the three square-number spiralarms Q1, Q2 and Q3 of the Square Root Spiral ( see FIG. 1 ) are congruent on one straight line !
Therefore a comparison of my analysis results with the results of the study from Mr. Robert Sachs could be very helpful for further discoveries !

Mr. Sachs gave me permission to show some sections of his analysis in my paper.

From the 10 chapters shown on the webside of Mr. Sachs I want to show the chapters which have a close connection to my findings. These are as follows :

1.) – Introduction ; 2.) – Product Curves ; 3.) – Offset Curves ; 5.) – Quadratic Polynomials ; 7.) – Prime Numbers ; 9.) - Formulas

Please have a look to **chapter 9** with the title **" The Number Spiral"** on **page 28** which shows the above mentioned chapters.

Images from the analysis of Mr. Sachs are named as follows: FIG. NS-1 to NS-18

The webside of Mr. Sachs, which deals with the mentioned Number Spiral can be found under the following weblink : → www.numberspiral.com

Following **chapter 9,** I then compared the Square Root Spiral with the Number Spiral and the Ulam Spiral, in regards to the arrangement of some selected "Reference Graphs". This comparison is shown in **chapter 10** → see **page 36**

I consider this comparison only as a first little step of a much more extensive analysis, which should be carried out here, in regards to the arrangement of prime numbers and non-prime-numbers on spirals with different spiral intensities, where the square numbers are located on a defined number of graphs.

In the following I want to show a priority list of some discoveries shown in my studies 1 to 4, where further mathematical analyses should be done, to explain this findings !

Special attention should here be paid to the distribution of prime factors in the non-prime numbers of the analysed spiral graphs and number sequences ! :

**Priority List of discoveries, suggested for further mathematical analysis :**

|  |  | Study |
|---|---|---|
| **1** | - General rule which defines the Number-Group-Spiral-Systems → chapter 5.2 | **- [ 1 ]** |
| **2** | - Prime Number Spiral Systems shown in FIG 6-A to 6-C → chapter 3 – in this study | **- [ 2 ]** |
| **3** | - Periodic occurrence of the prime factors in the non-prime-numbers of the Prime Number Spiral Graphs & their period lengths (→ Chapt. 6 /Tab. 4, 5-A, 5-B) | **- [ 2 ]** |
| **4** | - The meaning of the "Sums of the digits – Sequences" and "Number-endings-Sequences" in the whole context (→ described in Chapter 4.3 and 4.4 ) | **- [ 2 ]** |
| **5** | - Why do the number-endings-sequences and "sums of the digits sequences" usually cover **5** successive numbers of the analysed sequences ? (→ 4.3/4.4 ) | **- [ 2 ]** |
| **6** | - The interlaced occurrence of the special Prime Number Sequences SQ1 + SQ2 on the Square Root Spiral → ( see FIG 8 ) **- [ 2 ]** ; → also see my 3. Study ! | **- [ 3 ]** |
| **7** | - Pronics-Graphs and Product Curves on the Square Root Spiral & Number Spiral | **- [ 2 ]** |
| **8** | - The periodic occurrence of prime factors in the natural Fibonacci-Sequences | **- [ 4 ]** |
| **9** | - The distribution of the Square Numbers on 3 highly symmetrical spiral graphs | **- [ 1 ]** |
| **10** | - Interval $\pi$ between winds of the Square Root Spiral for sqrt → ∞ ( Chapt. 1 & 1' ) | **- [ 2 ]** |
| **11** | - Difference-Graph $f(x) = 2 ( 5x^2 - 7x + 3 )$ shown in FIG. 16 / in the Appendix | **- [ 1 ]** |

In December 2005 and June 2006 I sent the most findings shown in this study here to some universities in Germany for an assessment. But unfortunately there wasn't much response ! That's why I decided to publish my discoveries here !

Prof. S.J. Patterson from the University of Goettingen found some of my discoveries very interesting. He was especially interested in the spiral graphs which contain the Prime Numbers ( shown in FIG 6-A to 6-C ).
These spiral graphs are special quadratic polynomials, which are of great interest to Prime Number Theory. For example the quadratic polynomial B3 in FIG 6-A → $B3 = F(x) = 9x^2 + 27x + 17$ ( or $9x^2 + 9x - 1$ )

Prof. Ernst Wilhelm Zink from the Humboldt-University in Berlin also found my study very interesting and he organized a mathematical analysis of the spiral-graphs shown in FIG 6-A to FIG 6-C.

This mathematical analysis was carried out by Mr. Kay Schoenberger, a student of mathematics on the Humboldt-University of Berlin, who was working on his dissertation. The results of this mathematical analysis is shown in my first study :

→ " **The ordered distribution of natural numbers on the Square Root Spiral " [ 1 ]**

The Square Root Spiral and the Number Spiral are interesting tools to find ( spatial ) interdependencies between natural numbers and to understand the distribution of certain number groups in a visual way.

Therefore I want to ask mathematicians who read my study, to continue my work and to do more extensive graphical analyses of similar kind by using more advanced analysis-techniques and specialized computer software.

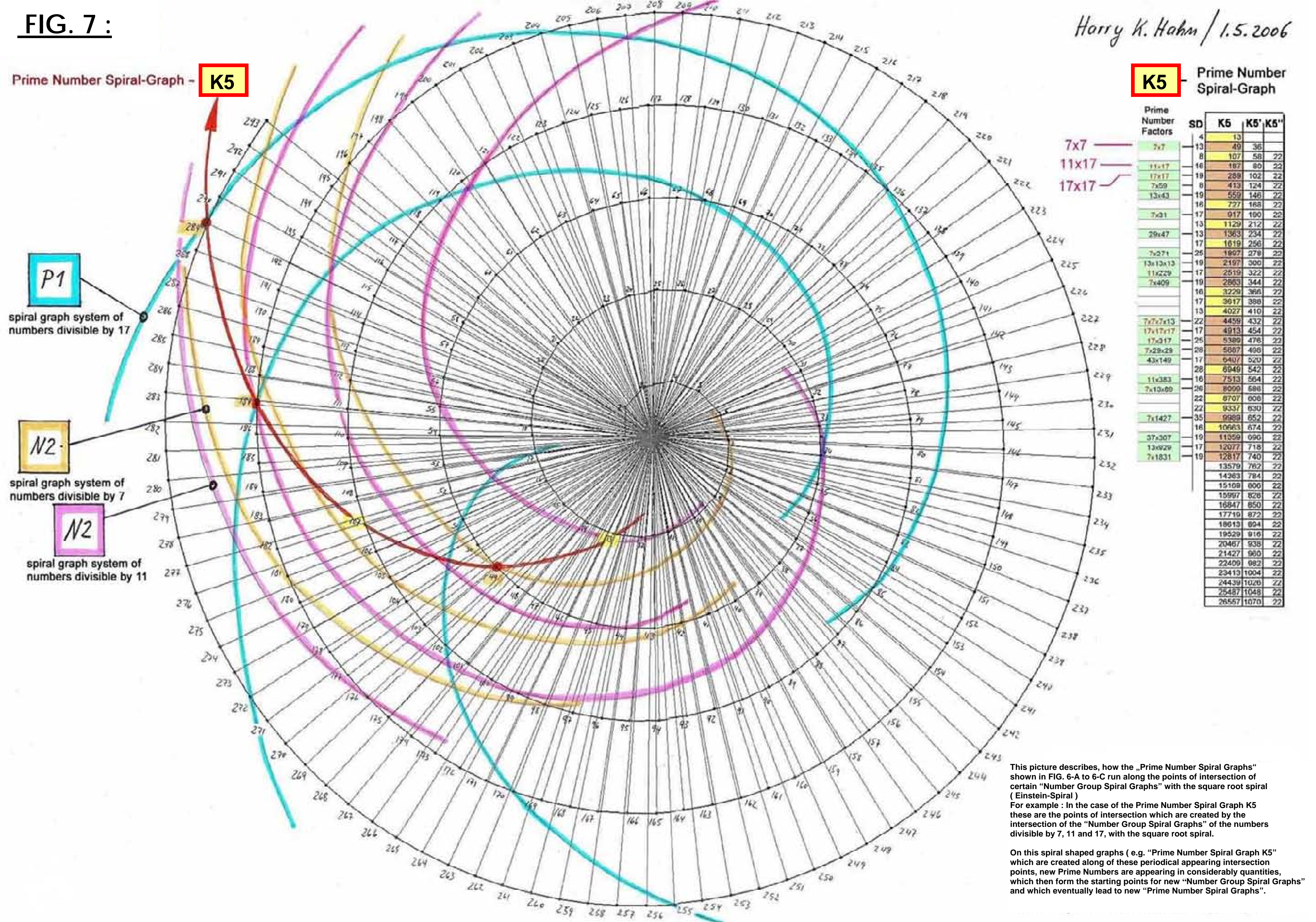

| Prime Number Factors | SD | K5 | K5' | K5'' |
|---|---|---|---|---|
| | 4 | 13 | | |
| 7x7 | 13 | 49 | 36 | |
| | 8 | 107 | 58 | 22 |
| 11x17 | 16 | 187 | 80 | 22 |
| 17x17 | 19 | 289 | 102 | 22 |
| 7x59 | 8 | 413 | 124 | 22 |
| 13x43 | 19 | 559 | 146 | 22 |
| | 16 | 727 | 168 | 22 |
| 7x31 | 17 | 917 | 190 | 22 |
| | 13 | 1129 | 212 | 22 |
| 29x47 | 13 | 1363 | 234 | 22 |
| | 17 | 1619 | 256 | 22 |
| 7x271 | 25 | 1897 | 278 | 22 |
| 13x13x13 | 19 | 2197 | 300 | 22 |
| 11x229 | 17 | 2519 | 322 | 22 |
| 7x409 | 19 | 2863 | 344 | 22 |
| | 16 | 3229 | 366 | 22 |
| | 17 | 3617 | 388 | 22 |
| | 13 | 4027 | 410 | 22 |
| 7x7x7x13 | 22 | 4459 | 432 | 22 |
| 17x17x17 | 17 | 4913 | 454 | 22 |
| 17x317 | 25 | 5389 | 476 | 22 |
| 7x29x29 | 28 | 5887 | 498 | 22 |
| 43x149 | 17 | 6407 | 520 | 22 |
| | 28 | 6949 | 542 | 22 |
| 11x383 | 16 | 7513 | 564 | 22 |
| 7x13x89 | 26 | 8099 | 586 | 22 |
| | 22 | 8707 | 608 | 22 |
| | 22 | 9337 | 630 | 22 |
| 7x1427 | 35 | 9989 | 652 | 22 |
| | 16 | 10663 | 674 | 22 |
| 37x307 | 19 | 11359 | 696 | 22 |
| 13x929 | 17 | 12077 | 718 | 22 |
| 7x1831 | 19 | 12817 | 740 | 22 |
| | | 13579 | 762 | 22 |
| | | 14363 | 784 | 22 |
| | | 15169 | 806 | 22 |
| | | 15997 | 828 | 22 |
| | | 16847 | 850 | 22 |
| | | 17719 | 872 | 22 |
| | | 18613 | 894 | 22 |
| | | 19529 | 916 | 22 |
| | | 20467 | 938 | 22 |
| | | 21427 | 960 | 22 |
| | | 22409 | 982 | 22 |
| | | 23413 | 1004 | 22 |
| | | 24439 | 1026 | 22 |
| | | 25487 | 1048 | 22 |
| | | 26557 | 1070 | 22 |

**This picture describes, how the „Prime Number Spiral Graphs“ shown in FIG. 6-A to 6-C run along the points of intersection of certain “Number Group Spiral Graphs” with the square root spiral ( Einstein-Spiral )**
**For example : In the case of the Prime Number Spiral Graph K5 these are the points of intersection which are created by the intersection of the “Number Group Spiral Graphs” of the numbers divisible by 7, 11 and 17, with the square root spiral.**

**On this spiral shaped graphs ( e.g. “Prime Number Spiral Graph K5” which are created along of these periodical appearing intersection points, new Prime Numbers are appearing in considerably quantities, which then form the starting points for new “Number Group Spiral Graphs” and which eventually lead to new “Prime Number Spiral Graphs”.**

**Table 1 :** Random analysis of the numbers of the Prime Number Spiral-Graphs (-sequences ) **B3** , **Q3** and **P20-G1** ( see FIG. 6-A to 6-C ) in regards to their share in Prime Numbers. → **Spot checks carried out in the number area : 0 - 2 500 000 000**
Note : The selected Prime Number Spiral-Graphs contain a particular high share of Prime Numbers.

| B3 | | | | |
|---|---|---|---|---|
| SD | Prime Factors | B3 | B3’ | B3’’ |
| 8 | | 17 | | |
| 8 | | 53 | 36 | |
| 8 | | 107 | 54 | 18 |
| 17 | | 179 | 72 | 18 |
| 17 | | 269 | 90 | 18 |
| 17 | 13x29 | 377 | 108 | 18 |
| 8 | | 503 | 126 | 18 |
| 17 | | 647 | 144 | 18 |
| 17 | | 809 | 162 | 18 |
| 26 | 23x43 | 989 | 180 | 18 |
| 17 | | 1187 | 198 | 18 |
| 8 | 23x61 | 1403 | 216 | 18 |
| 17 | | 1637 | 234 | 18 |
| 26 | | 1889 | 252 | 18 |
| 17 | 17x127 | 2159 | 270 | 18 |
| 17 | | 2447 | 288 | 18 |
| 17 | | 2753 | 306 | 18 |
| 17 | 17x181 | 3077 | 324 | 18 |
| 17 | 13x263 | 3419 | 342 | 18 |
| 26 | | 3779 | 360 | 18 |
| 17 | | 4157 | 378 | 18 |
| 17 | 29x157 | 4553 | 396 | 18 |
| 26 | | 4967 | 414 | 18 |
| 26 | | 5399 | 432 | 18 |
| 26 | | 5849 | 450 | 18 |
| 17 | | 6317 | 468 | 18 |
| 17 | | 6803 | 486 | 18 |
| 17 | | 7307 | 504 | 18 |
| 26 | | 7829 | 522 | 18 |
| 26 | | 8369 | 540 | 18 |
| 26 | 79x113 | 8927 | 558 | 18 |
| 17 | 13x17x43 | 9503 | 576 | 18 |
| 17 | 23x439 | 10097 | 594 | 18 |
| 17 | | 10709 | 612 | 18 |
| 17 | 17x23x29 | 11339 | 630 | 18 |
| 26 | | 11987 | 648 | 18 |
| 17 | | 12653 | 666 | 18 |
| 17 | | 13337 | 684 | 18 |
| 17 | 53x47251 | 2504303 | 9486 | 18 |
| 26 | 17x29x5099 | 2513807 | 9504 | 18 |
| 26 | | 2523329 | 9522 | 18 |
| 35 | | 2532869 | 9540 | 18 |
| 26 | 107x23761 | 2542427 | 9558 | 18 |
| 17 | 53x179x269 | 2552003 | 9576 | 18 |
| 35 | 113x22669 | 2561597 | 9594 | 18 |
| 26 | | 2571209 | 9612 | 18 |
| 17 | 17x797x1847 | 25025003 | 30006 | 18 |
| 26 | 23x23x47363 | 25055027 | 30024 | 18 |
| 35 | | 25085069 | 30042 | 18 |
| 26 | 13x1931933 | 25115129 | 30060 | 18 |
| 26 | 2551x9857 | 25145207 | 30078 | 18 |
| 26 | | 25175303 | 30096 | 18 |
| 26 | | 25205417 | 30114 | 18 |
| 35 | | 25235549 | 30132 | 18 |
| 26 | 191x1308919 | 250003529 | 94860 | 18 |
| 35 | 13x17x29x39023 | 250098407 | 94878 | 18 |
| 26 | 233x1073791 | 250193303 | 94896 | 18 |
| 35 | | 250288217 | 94914 | 18 |
| 35 | | 250383149 | 94932 | 18 |
| 44 | | 250478099 | 94950 | 18 |
| 35 | 10253x24439 | 250573067 | 94968 | 18 |
| 35 | 23x10898611 | 250668053 | 94986 | 18 |
| 17 | | 2500250003 | 300006 | 18 |
| 26 | 29x2731x31573 | 2500550027 | 300024 | 18 |
| 35 | 3793x659333 | 2500850069 | 300042 | 18 |
| 26 | 181x13818509 | 2501150129 | 300060 | 18 |
| 26 | | 2501450207 | 300078 | 18 |
| 26 | 13x269x673x1063 | 2501750303 | 300096 | 18 |
| 26 | 43x58187219 | 2502050417 | 300114 | 18 |
| 35 | | 2502350549 | 300132 | 18 |

| Q3 | | | | |
|---|---|---|---|---|
| SD | Prime Factors | Q3 | Q3’ | Q3’’ |
| 4 | | 13 | | |
| 10 | | 37 | 24 | |
| 11 | | 83 | 46 | 22 |
| 7 | | 151 | 68 | 22 |
| 7 | | 241 | 90 | 22 |
| 11 | | 353 | 112 | 22 |
| 19 | | 487 | 134 | 22 |
| 13 | | 643 | 156 | 22 |
| 11 | | 821 | 178 | 22 |
| 4 | | 1021 | 200 | 22 |
| 10 | 11x113 | 1243 | 222 | 22 |
| 20 | | 1487 | 244 | 22 |
| 16 | | 1753 | 266 | 22 |
| 7 | 13x157 | 2041 | 288 | 22 |
| 11 | | 2351 | 310 | 22 |
| 19 | | 2683 | 332 | 22 |
| 13 | | 3037 | 354 | 22 |
| 11 | | 3413 | 376 | 22 |
| 13 | 37x103 | 3811 | 398 | 22 |
| 10 | | 4231 | 420 | 22 |
| 20 | | 4673 | 442 | 22 |
| 16 | 11x467 | 5137 | 464 | 22 |
| 16 | | 5623 | 486 | 22 |
| 11 | | 6131 | 508 | 22 |
| 19 | | 6661 | 530 | 22 |
| 13 | | 7213 | 552 | 22 |
| 29 | 13x599 | 7787 | 574 | 22 |
| 22 | 83x101 | 8383 | 596 | 22 |
| 10 | | 9001 | 618 | 22 |
| 20 | 31x311 | 9641 | 640 | 22 |
| 7 | | 10303 | 662 | 22 |
| 25 | | 10987 | 684 | 22 |
| 20 | 11x1063 | 11693 | 706 | 22 |
| 10 | | 12421 | 728 | 22 |
| 13 | | 13171 | 750 | 22 |
| 20 | 73x191 | 13943 | 772 | 22 |
| 22 | | 14737 | 794 | 22 |
| 19 | 103x151 | 15553 | 816 | 22 |
| 25 | 947x3221 | 3050287 | 11574 | 22 |
| 29 | 11x278353 | 3061883 | 11596 | 22 |
| 19 | 739x4159 | 3073501 | 11618 | 22 |
| 22 | 1319x2339 | 3085141 | 11640 | 22 |
| 29 | 233x13291 | 3096803 | 11662 | 22 |
| 31 | | 3108487 | 11684 | 22 |
| 19 | 101x30893 | 3120193 | 11706 | 22 |
| 20 | 13x103x2339 | 3131921 | 11728 | 22 |
| 37 | 37x103x8017 | 30552787 | 36654 | 22 |
| 38 | | 30589463 | 36676 | 22 |
| 25 | 503x60887 | 30626161 | 36698 | 22 |
| 34 | | 30662881 | 36720 | 22 |
| 38 | 157x195539 | 30699623 | 36742 | 22 |
| 37 | 11x2794217 | 30736387 | 36764 | 22 |
| 31 | 31x197x5039 | 30773173 | 36786 | 22 |
| 38 | 743x41467 | 30809981 | 36808 | 22 |
| 34 | 113x2703137 | 305454481 | 115920 | 22 |
| 29 | 3407x89689 | 305570423 | 115942 | 22 |
| 46 | 103x1693x1753 | 305686387 | 115964 | 22 |
| 31 | 37x1303x6343 | 305802373 | 115986 | 22 |
| 38 | | 305918381 | 116008 | 22 |
| 22 | 17419x17569 | 306034411 | 116030 | 22 |
| 28 | | 306150463 | 116052 | 22 |
| 38 | | 306266537 | 116074 | 22 |
| 49 | 13x235040599 | 3055527787 | 366654 | 22 |
| 47 | 137x499x44701 | 3055894463 | 366676 | 22 |
| 31 | 37x223x370411 | 3056261161 | 366698 | 22 |
| 46 | | 3056627881 | 366720 | 22 |
| 47 | 2633x1161031 | 3056994623 | 366742 | 22 |
| 43 | 1433x2133539 | 3057361387 | 366764 | 22 |
| 43 | | 3057728173 | 366786 | 22 |
| 47 | 3803x804127 | 3058094981 | 366808 | 22 |

| P20 - G1 | | | | |
|---|---|---|---|---|
| SD | Prime Factors | G1 | G1’ | G1’’ |
| 4 | | 103 | | |
| 11 | | 173 | 70 | |
| 11 | | 263 | 90 | 20 |
| 13 | | 373 | 110 | 20 |
| 8 | | 503 | 130 | 20 |
| 14 | | 653 | 150 | 20 |
| 13 | | 823 | 170 | 20 |
| 5 | | 1013 | 190 | 20 |
| 8 | | 1223 | 210 | 20 |
| 13 | | 1453 | 230 | 20 |
| 11 | 13x131 | 1703 | 250 | 20 |
| 20 | | 1973 | 270 | 20 |
| 13 | 31x73 | 2263 | 290 | 20 |
| 17 | 31x83 | 2573 | 310 | 20 |
| 14 | | 2903 | 330 | 20 |
| 13 | | 3253 | 350 | 20 |
| 14 | | 3623 | 370 | 20 |
| 8 | | 4013 | 390 | 20 |
| 13 | | 4423 | 410 | 20 |
| 20 | 23x211 | 4853 | 430 | 20 |
| 11 | | 5303 | 450 | 20 |
| 22 | 23x251 | 5773 | 470 | 20 |
| 17 | | 6263 | 490 | 20 |
| 23 | 13x521 | 6773 | 510 | 20 |
| 13 | 67x109 | 7303 | 530 | 20 |
| 23 | | 7853 | 550 | 20 |
| 17 | | 8423 | 570 | 20 |
| 13 | | 9013 | 590 | 20 |
| 20 | | 9623 | 610 | 20 |
| 11 | | 10253 | 630 | 20 |
| 13 | | 10903 | 650 | 20 |
| 17 | 71x163 | 11573 | 670 | 20 |
| 14 | | 12263 | 690 | 20 |
| 22 | | 12973 | 710 | 20 |
| 14 | 71x193 | 13703 | 730 | 20 |
| 17 | 97x149 | 14453 | 750 | 20 |
| 13 | 13x1171 | 15223 | 770 | 20 |
| 11 | 67x239 | 16013 | 790 | 20 |
| 35 | 709x3947 | 2798423 | 10570 | 20 |
| 23 | 23x122131 | 2809013 | 10590 | 20 |
| 31 | 71x151x263 | 2819623 | 10610 | 20 |
| 23 | 53x53401 | 2830253 | 10630 | 20 |
| 26 | 13x218531 | 2840903 | 10650 | 20 |
| 31 | 71x40163 | 2851573 | 10670 | 20 |
| 29 | 617x4639 | 2862263 | 10690 | 20 |
| 38 | | 2872973 | 10710 | 20 |
| 38 | | 27855623 | 33370 | 20 |
| 38 | 71x392803 | 27889013 | 33390 | 20 |
| 31 | | 27922423 | 33410 | 20 |
| 44 | 3803x7351 | 27955853 | 33430 | 20 |
| 41 | | 27989303 | 33450 | 20 |
| 31 | 157x178489 | 28022773 | 33470 | 20 |
| 32 | 2161x12983 | 28056263 | 33490 | 20 |
| 44 | 83x338431 | 28089773 | 33510 | 20 |
| 53 | 2711x102523 | 277939853 | 105430 | 20 |
| 32 | | 278045303 | 105450 | 20 |
| 40 | | 278150773 | 105470 | 20 |
| 41 | | 278256263 | 105490 | 20 |
| 44 | 97x643x4463 | 278361773 | 105510 | 20 |
| 40 | | 278467303 | 105530 | 20 |
| 47 | 13x21428681 | 278572853 | 105550 | 20 |
| 47 | | 278678423 | 105570 | 20 |
| 50 | | 2778555623 | 333370 | 20 |
| 53 | 271x10254203 | 2778889013 | 333390 | 20 |
| 40 | | 2779222423 | 333410 | 20 |
| 56 | 509x5460817 | 2779555853 | 333430 | 20 |
| 56 | | 2779889303 | 333450 | 20 |
| 40 | 23x2539x47609 | 2780222773 | 333470 | 20 |
| 44 | 659x4219357 | 2780556263 | 333490 | 20 |
| 59 | 23x2741x44111 | 2780889773 | 333510 | 20 |

**Note :** Remarkable is the repeated occurrence of certain Prime-Factors ( marked in red ) !!
**SD** = sum of the digits

**Table 2 :** Analysis of the "Sums of the digits" which result from the numbers of the "Prime Number Sequences" shown in tables **6-A1** to **6-C1** ( --> see columns -**SD**- )

| **2. Differential** | **"Prime Number Sequences"** ( see tables **6-A1** to **6-C1** ) --> Sequences derived from Prime Number Spiral Graphs shown in FIG. 6-**A** to 6-**C** | **"Sum of the digits"-Sequence belonging to these Prime Number Sequences** ( see columns -**SD**- in tables 6-A1 to 6-C1 ) | **Periodicity of "Sum of the digits" Sequence** |
|---|---|---|---|
| **18** | A2, A4, A6, A8, A10, A12.....<br>C2, C4, C6, C8, C10, C12,....<br>2 2 2 | 4,7,10,13,16,19,....<br>3 3 3 | 3 |
| | A1, A3, A5, A7, A9, A11.....<br>C1, C3, C5, C7, C9, C11,.... | 2,5,8,11,14,17,.... | 3 |
| | B1, B7, B13, B19......<br>6 6 6 | 2,11,20,29,38,....<br>9 9 9 | 9 |
| | B2, B8, B14, B20,.... | 4,13,22,31,40,.... | 9 |
| | B3, B9, B15, B21,.... | 8,17,26,35,44,.... | 9 |
| | B4, B10, B16, B22,.... | 10,19,28,37,46,...... | 9 |
| | B5, B11, B17, B23,.... | 5,14,23,32,41,.... | 9 |
| | B6, B12, B18, B24,.... | 7,16,25,34,43,.... | 9 |

| **2. Differential** | | **D** | **E** | **F** | **G** | **H** | **I** | "Sum of the digits"-Sequence | Periodicity |
|---|---|---|---|---|---|---|---|---|---|
| **20** | N - | D1 | | F3 | G1 | | I1 | 1, 4, 7, 9,10,13,16,18,19,22,25,27,....<br>3 3 2 1 3 3 2 1 3 3 2 | 3 , 3 , 2 , 1 |
| | P - | | E2 | | | H2 | | | |
| | N - | D3 | E1 | | G3 | H1 | I3 | 1,3,4,7,10,12,13,16,19,21,22,.... | 3 , 3 , 2 , 1 |
| | P - | | | F2 | | | | | |
| | N - | | E3 | F1 | | H3 | | 1,4,6,7,10,13,15,16,19,22,24,.... | 3 , 3 , 2 , 1 |
| | P - | D2 | | | G2 | | I2 | | |
| | N - | | E2 | | | H2 | | 2,5,8,10,11,14,17,19,20,23,.... | 3 , 3 , 2 , 1 |
| | P - | D3 | | F1 | G3 | | I3 | | |
| | N - | | | F2 | | | | 2,4,5,8,11,13,14,17,20,22,..... | 3 , 3 , 2 , 1 |
| | P - | D1 | E3 | | G1 | H3 | I1 | | |
| | N - | D2 | | | G2 | | I2 | 2,5,7,8,11,14,16,17,20,23,.... | 3 , 3 , 2 , 1 |
| | P - | | E1 | F3 | | H1 | | | |

| **2. Differential** | **J** | **K** | **L** | **M** | **N** | **O** | **P** | **Q** | **R** | **S** | **T** | "Sum of the digits"-Sequence | Periodicity |
|---|---|---|---|---|---|---|---|---|---|---|---|---|---|
| **22** | J2 | K1 | | | N1 | O3 | P3 | Q3 | R3 | S1 | T3 | 2, 4, 7,10,11,13,16,19,20,22,25,....<br>2 3 3 1 2 3 3 1 2 3 3 | 1 , 2 , 3 , 3 |
| | | K3 | L1 | M1 | N3 | | | | | S3 | | 4,5,7,10,13,14,16,19,22,.... | 1 , 2 , 3 , 3 |
| | | | L3 | M3 | | O1 | P1 | Q1 | R1 | | T1 | 4,7,8,10,13,16,17,19,22,.... | 1 , 2 , 3 , 3 |
| | | | L2 | M2 | | | | | | | | 5,8,11,12,14,17,20,21,23,.... | 1 , 2 , 3 , 3 |
| | | | | | | O2 | P2 | Q2 | R2 | | T2 | 2,5,6,8,11,14,15,17,20,.... | 1 , 2 , 3 , 3 |
| | | K2 | | | N2 | | | | | S2 | | 5,8,9,11,14,17,18,20,23,.... | 1 , 2 , 3 , 3 |
| | J3 | | | | | | | | | | | 6,7,9,12,15,16,18,21,24,... | 1 , 2 , 3 , 3 |
| | J1 | | | | | | | | | | | 4,6,9,12,13,15,18,21,22,.... | 1 , 2 , 3 , 3 |

**Note :** „Sum of the digits“- Sequences were created by ordering the „sums of the digits“ occuring in the „prime number sequences“ according to their value. → see values in the columns -**SD**- in tables 6-A1 to 6-C1

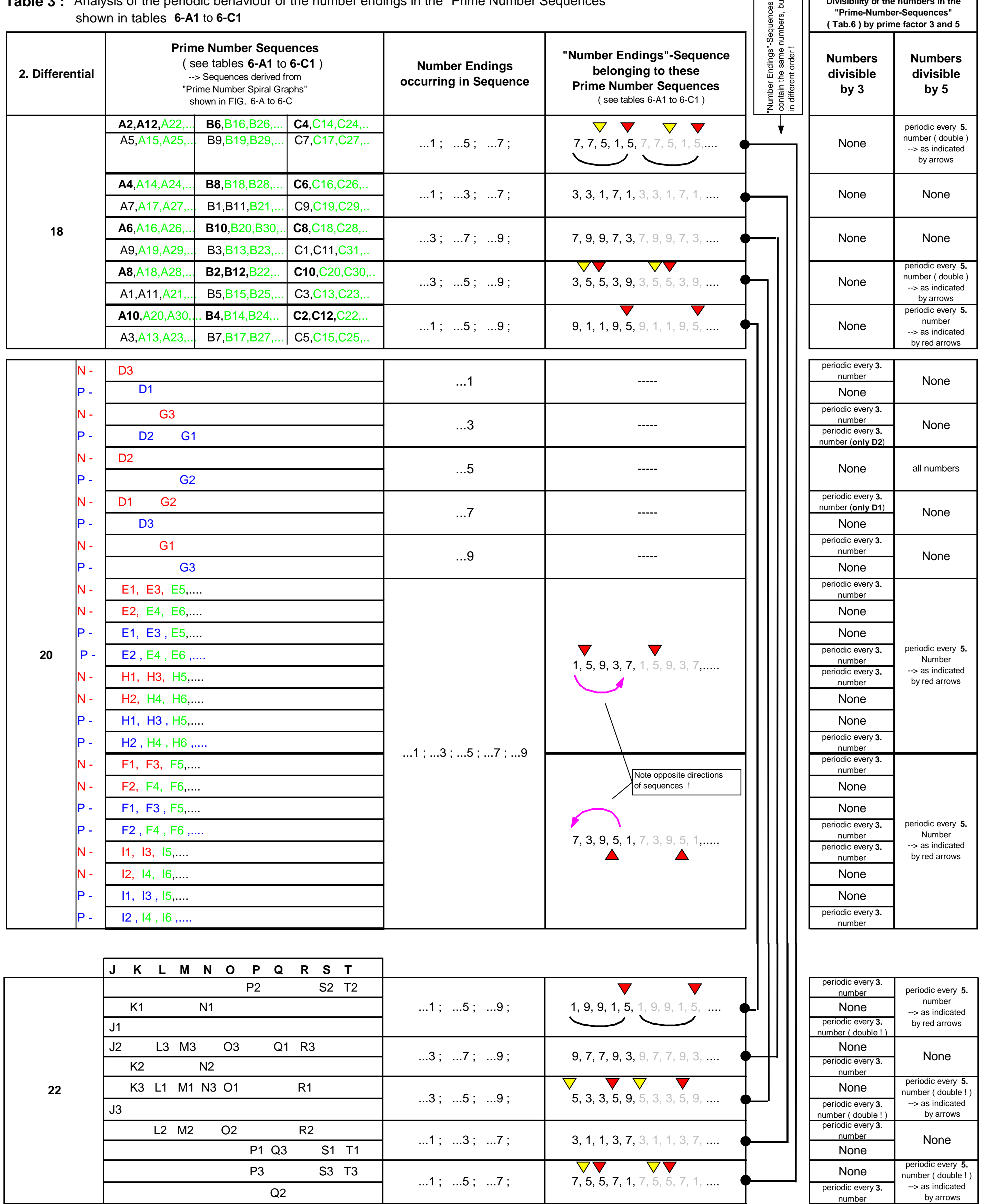

**Table 3 :** Analysis of the periodic behaviour of the number endings in the "Prime Number Sequences" shown in tables **6-A1** to **6-C1**

"Number Endings"-Sequences contain the same numbers, but in different order !

Divisibility of the numbers in the "Prime-Number-Sequences" ( Tab.6 ) by prime factor 3 and 5

| 2. Differential | Prime Number Sequences ( see tables 6-A1 to 6-C1 ) --> Sequences derived from "Prime Number Spiral Graphs" shown in FIG. 6-A to 6-C | | | Number Endings occurring in Sequence | "Number Endings"-Sequence belonging to these Prime Number Sequences ( see tables 6-A1 to 6-C1 ) | Numbers divisible by 3 | Numbers divisible by 5 |
|---|---|---|---|---|---|---|---|
| 18 | **A2,A12**,A22,... | **B6**,B16,B26,... | **C4**,C14,C24,.. | ...1 ; ...5 ; ...7 ; | 7, 7, 5, 1, 5, 7, 7, 5, 1, 5,.... | None | periodic every 5. number ( double ) --> as indicated by arrows |
| | A5,A15,A25,... | B9,B19,B29,... | C7,C17,C27,.. | | | | |
| | **A4**,A14,A24,... | **B8**,B18,B28,... | **C6**,C16,C26,.. | ...1 ; ...3 ; ...7 ; | 3, 3, 1, 7, 1, 3, 3, 1, 7, 1, .... | None | None |
| | A7,A17,A27,.. | B1,B11,B21,... | C9,C19,C29,.. | | | | |
| | **A6**,A16,A26,.. | **B10**,B20,B30,.. | **C8**,C18,C28,.. | ...3 ; ...7 ; ...9 ; | 7, 9, 9, 7, 3, 7, 9, 9, 7, 3, .... | None | None |
| | A9,A19,A29,.. | B3,B13,B23,... | C1,C11,C31,.. | | | | |
| | **A8**,A18,A28,... | **B2,B12**,B22,.. | **C10**,C20,C30,.. | ...3 ; ...5 ; ...9 ; | 3, 5, 5, 3, 9, 3, 5, 5, 3, 9, .... | None | periodic every 5. number ( double ) --> as indicated by arrows |
| | A1,A11,A21,.. | B5,B15,B25,... | C3,C13,C23,.. | | | | |
| | **A10**,A20,A30,... | **B4**,B14,B24,.. | **C2,C12**,C22,.. | ...1 ; ...5 ; ...9 ; | 9, 1, 1, 9, 5, 9, 1, 1, 9, 5, .... | None | periodic every 5. number --> as indicated by red arrows |
| | A3,A13,A23,.. | B7,B17,B27,... | C5,C15,C25,.. | | | | |

| 2. Differential | | Prime Number Sequences | Number Endings occurring in Sequence | "Number Endings"-Sequence | Numbers divisible by 3 | Numbers divisible by 5 |
|---|---|---|---|---|---|---|
| 20 | N - | D3 | ...1 | ----- | periodic every 3. number | None |
| | P - | D1 | | | None | |
| | N - | G3 | ...3 | ----- | periodic every 3. number | None |
| | P - | D2 G1 | | | periodic every 3. number (only D2) | |
| | N - | D2 | ...5 | ----- | None | all numbers |
| | P - | G2 | | | | |
| | N - | D1 G2 | ...7 | ----- | periodic every 3. number (only D1) | None |
| | P - | D3 | | | None | |
| | N - | G1 | ...9 | ----- | periodic every 3. number | None |
| | P - | G3 | | | None | |
| | N - | E1, E3, E5,.... | ...1 ; ...3 ; ...5 ; ...7 ; ...9 | 1, 5, 9, 3, 7, 1, 5, 9, 3, 7,..... | periodic every 3. number | periodic every 5. Number --> as indicated by red arrows |
| | N - | E2, E4, E6,.... | | | None | |
| | P - | E1, E3 , E5,.... | | | None | |
| | P - | E2 , E4 , E6 ,.... | | | periodic every 3. number | |
| | N - | H1, H3, H5,.... | | | periodic every 3. number | |
| | N - | H2, H4, H6,.... | | | None | |
| | P - | H1, H3 , H5,.... | | | None | |
| | P - | H2 , H4 , H6 ,.... | | | periodic every 3. number | |
| | N - | F1, F3, F5,.... | | Note opposite directions of sequences ! | periodic every 3. number | periodic every 5. Number --> as indicated by red arrows |
| | N - | F2, F4, F6,.... | | | None | |
| | P - | F1, F3 , F5,.... | | | None | |
| | P - | F2 , F4 , F6 ,.... | | | periodic every 3. number | |
| | N - | I1, I3, I5,.... | | 7, 3, 9, 5, 1, 7, 3, 9, 5, 1,..... | periodic every 3. number | |
| | N - | I2, I4, I6,.... | | | None | |
| | P - | I1, I3 , I5,.... | | | None | |
| | P - | I2 , I4 , I6 ,.... | | | periodic every 3. number | |

| 2. Differential | J | K | L | M | N | O | P | Q | R | S | T | Number Endings occurring in Sequence | "Number Endings"-Sequence | Numbers divisible by 3 | Numbers divisible by 5 |
|---|---|---|---|---|---|---|---|---|---|---|---|---|---|---|---|
| 22 | | | | | | | P2 | | | S2 | T2 | ...1 ; ...5 ; ...9 ; | 1, 9, 9, 1, 5, 1, 9, 9, 1, 5, .... | periodic every 3. number | periodic every 5. number --> as indicated by red arrows |
| | | K1 | | | N1 | | | | | | | | | None | |
| | J1 | | | | | | | | | | | | | periodic every 3. number ( double ! ) | |
| | J2 | | L3 | M3 | | O3 | | Q1 | R3 | | | ...3 ; ...7 ; ...9 ; | 9, 7, 7, 9, 3, 9, 7, 7, 9, 3, .... | None | None |
| | | K2 | | | N2 | | | | | | | | | periodic every 3. number | |
| | | K3 | L1 | M1 | N3 | O1 | | | R1 | | | ...3 ; ...5 ; ...9 ; | 5, 3, 3, 5, 9, 5, 3, 3, 5, 9, .... | None | periodic every 5. number ( double ! ) --> as indicated by arrows |
| | J3 | | | | | | | | | | | | | periodic every 3. number ( double ! ) | |
| | | | L2 | M2 | | O2 | | | R2 | | | ...1 ; ...3 ; ...7 ; | 3, 1, 1, 3, 7, 3, 1, 1, 3, 7, .... | periodic every 3. number | None |
| | | | | | | | P1 | Q3 | | S1 | T1 | | | None | |
| | | | | | | | P3 | | | S3 | T3 | ...1 ; ...5 ; ...7 ; | 7, 5, 5, 7, 1, 7, 5, 5, 7, 1, .... | None | periodic every 5. number ( double ! ) --> as indicated by arrows |
| | | | | | | | | Q2 | | | | | | periodic every 3. number | |

e.g. A15, A25 - "Prime Number Sequences" marked in green not shown in tables 6-A1 to 6-C1 ( only for reference ! )

e.g. N - D3 marking in red indicates that the "Prime Number Spiral Graphs" have a **negative** (N) rotation direction
P - D1 marking in blue indicates that the "Prime Number Spiral Graphs" have a **positive** (P) rotation direction

indicates periodic occurence of numbers divisible by **5** in "Prime Number Sequences"

**Table 4 :** Periodic occurring Prime Factors in the „Prime-Number-Spiral-Graphs" (-number sequences ) **B3** and **K5**
( see also FIG. 6-A / 6-C & 7 and Tables 6-A1 / 6-C1 )

(13) (2/21) (3/14) (16/13) (22/21) (36/25)

| some Prime Factors presented in tabular form | | | | | | <-------- periodic occurrence of Individual Prime Factors → expressed through the number of spacings = (X) | Prime Factors of the none Prime-Numbers | sum of the digits | **B3** | | |
|---|---|---|---|---|---|---|---|---|---|---|---|
| 61 | 43 | 29 | 23 | 17 | 13 | | | | B3 | B3' | B3'' |
| | | | | | | | | 8 | 17 | | |
| | | | | | | | | 8 | 53 | 36 | |
| | | | | | | | | 8 | 107 | 54 | 18 |
| | | | | | | | | 17 | 179 | 72 | 18 |
| | | | | | | | | 17 | 269 | 90 | 18 |
| | | 29 | | | 13 | | 13x29 | 17 | 377 | 108 | 18 |
| | | | | | | | | 8 | 503 | 126 | 18 |
| | | | | | | | | 17 | 647 | 144 | 18 |
| | | | | | | | | 17 | 809 | 162 | 18 |
| | 43 | | 23 | | | | 23x43 | 26 | 989 | 180 | 18 |
| | | | | | | | | 17 | 1187 | 198 | 18 |
| 61 | | | 23 | | | | 23x61 | 8 | 1403 | 216 | 18 |
| | | | | | | | | 17 | 1637 | 234 | 18 |
| | | | | | | | | 26 | 1889 | 252 | 18 |
| | | | | 17 | | | 17x127 | 17 | 2159 | 270 | 18 |
| | | | | | | | | 17 | 2447 | 288 | 18 |
| | | | | | | | | 17 | 2753 | 306 | 18 |
| | | | | 17 | | | 17x181 | 17 | 3077 | 324 | 18 |
| | | | | | 13 | | 13x263 | 17 | 3419 | 342 | 18 |
| | | | | | | | | 26 | 3779 | 360 | 18 |
| | | | | | | | | 17 | 4157 | 378 | 18 |
| | | 29 | | | | | 29x157 | 17 | 4553 | 396 | 18 |
| | | | | | | | | 26 | 4967 | 414 | 18 |
| | | | | | | | | 26 | 5399 | 432 | 18 |
| | | | | | | | | 26 | 5849 | 450 | 18 |
| | | | | | | | | 17 | 6317 | 468 | 18 |
| | | | | | | | | 17 | 6803 | 486 | 18 |
| | | | | | | | | 17 | 7307 | 504 | 18 |
| | | | | | | | | 26 | 7829 | 522 | 18 |
| | | | | | | | | 26 | 8369 | 540 | 18 |
| | | | | | | | 79x113 | 26 | 8927 | 558 | 18 |
| | 43 | | | 17 | 13 | | 13x17x43 | 17 | 9503 | 576 | 18 |
| | | | 23 | | | | 23x439 | 17 | 10097 | 594 | 18 |
| | | | | | | | | 17 | 10709 | 612 | 18 |
| | | 29 | 23 | 17 | | | 17x23x29 | 17 | 11339 | 630 | 18 |
| | | | | | | | | 26 | 11987 | 648 | 18 |
| | | | | | | | | 17 | 12653 | 666 | 18 |
| | | | | | | | | 17 | 13337 | 684 | 18 |
| | | | | | | | 101x139 | 17 | 14039 | 702 | 18 |
| | | | | | | | | 26 | 14759 | 720 | 18 |
| | | | | | | | | 26 | 15497 | 738 | 18 |
| | | | | | | | | 17 | 16253 | 756 | 18 |
| | | | | | | | | 17 | 17027 | 774 | 18 |
| | | | | | | | 103x173 | 26 | 17819 | 792 | 18 |
| | | | | | 13 | | 13x1433 | 26 | 18629 | 810 | 18 |
| | | | | | | | | 26 | 19457 | 828 | 18 |
| | | | | | | | 79x257 | 8 | 20303 | 846 | 18 |
| 61 | | | | | | | 61x347 | 17 | 21167 | 864 | 18 |
| | | | | 17 | | | 17x1297 | 17 | 22049 | 882 | 18 |
| | | | | | | | 53x433 | 26 | 22949 | 900 | 18 |
| | | 29 | | | | | 29x823 | 26 | 23867 | 918 | 18 |
| | | | | 17 | | | 17x1459 | 17 | 24803 | 936 | 18 |
| | 43 | | | | | | 43x599 | 26 | 25757 | 954 | 18 |
| | | | | | | | | 26 | 26729 | 972 | 18 |
| | | | | | | | 53x523 | 26 | 27719 | 990 | 18 |
| | | | 23 | | | | 23x1249 | 26 | 28727 | 1008 | 18 |
| | | | | | | | | 26 | 29753 | 1026 | 18 |
| | | | 23 | | 13 | | 13x23x103 | 26 | 30797 | 1044 | 18 |
| | | | | | | | | 26 | 31859 | 1062 | 18 |
| | | | | | | | | 26 | 32939 | 1080 | 18 |
| | | | | | | | 101x337 | 17 | 34037 | 1098 | 18 |
| | | | | | | | | 17 | 35153 | 1116 | 18 |
| | | | | | | | 131x277 | 26 | 36287 | 1134 | 18 |
| | | 29 | | | | | 29x1291 | 26 | 37439 | 1152 | 18 |
| | | | | | | | | 26 | 38609 | 1170 | 18 |
| | | | | 17 | | | 17x2341 | 35 | 39797 | 1188 | 18 |
| | | | | | | | 131x313 | 8 | 41003 | 1206 | 18 |
| | | | | | | | | 17 | 42227 | 1224 | 18 |
| | | | | 17 | | | 17x2557 | 26 | 43469 | 1242 | 18 |
| | | | | | | | | 26 | 44729 | 1260 | 18 |
| | | | | | 13 | | 13x3539 | 17 | 46007 | 1278 | 18 |
| | | | | | | | | 17 | 47303 | 1296 | 18 |
| 61 | | | | | | | 61x797 | 26 | 48617 | 1314 | 18 |
| | | | | | | | 199x251 | 35 | 49949 | 1332 | 18 |
| | 43 | | | | | | 43x1193 | 26 | 51299 | 1350 | 18 |
| | | | | | | | | 26 | 52667 | 1368 | 18 |
| | | | | | | | 191x283 | 17 | 54053 | 1386 | 18 |
| | | | | | | | | 26 | 55457 | 1404 | 18 |
| | | | 23 | | | | 23x2473 | 35 | 56879 | 1422 | 18 |
| | | 29 | | | | | 29x2011 | 26 | 58319 | 1440 | 18 |
| | | | 23 | | | | 23x23x113 | 35 | 59777 | 1458 | 18 |
| | | | | | | | | 17 | 61253 | 1476 | 18 |

(4/3) (1/16) (11) (7/6) (12/17) (3/34)

| some Prime Factors presented in tabular form | | | | | | <-------- periodic occurrence of Individual Prime Factors → expressed through the number of spacings = (X) | Prime Factors of the none Prime-Numbers | sum of the digits | **K5** | | |
|---|---|---|---|---|---|---|---|---|---|---|---|
| 37 | 29 | 17 | 13 | 11 | 7 | | | | K5 | K5' | K5'' |
| | | | | | | | | 4 | 13 | | |
| | | | | | 7 | | 7x7 | 13 | 49 | 36 | |
| | | | | | | | | 8 | 107 | 58 | 22 |
| | | 17 | | 11 | | | 11x17 | 16 | 187 | 80 | 22 |
| | | 17 | | | | | 17x17 | 19 | 289 | 102 | 22 |
| | | | | | 7 | | 7x59 | 8 | 413 | 124 | 22 |
| | | | 13 | | | | 13x43 | 19 | 559 | 146 | 22 |
| | | | | | | | | 16 | 727 | 168 | 22 |
| | | | | | 7 | | 7x31 | 17 | 917 | 190 | 22 |
| | | | | | | | | 13 | 1129 | 212 | 22 |
| | 29 | | | | | | 29x47 | 13 | 1363 | 234 | 22 |
| | | | | | | | | 17 | 1619 | 256 | 22 |
| | | | | | 7 | | 7x271 | 25 | 1897 | 278 | 22 |
| | | | 13 | | | | 13x13x13 | 19 | 2197 | 300 | 22 |
| | | | | 11 | | | 11x229 | 17 | 2519 | 322 | 22 |
| | | | | | 7 | | 7x409 | 19 | 2863 | 344 | 22 |
| | | | | | | | | 16 | 3229 | 366 | 22 |
| | | | | | | | | 17 | 3617 | 388 | 22 |
| | | | | | | | | 13 | 4027 | 410 | 22 |
| | | | 13 | | 7 | | 7x7x7x13 | 22 | 4459 | 432 | 22 |
| | | 17 | | | | | 17x17x17 | 17 | 4913 | 454 | 22 |
| | | 17 | | | | | 17x317 | 25 | 5389 | 476 | 22 |
| | 29 | | | | 7 | | 7x29x29 | 28 | 5887 | 498 | 22 |
| | | | | | | | 43x149 | 17 | 6407 | 520 | 22 |
| | | | | | | | | 28 | 6949 | 542 | 22 |
| | | | | 11 | | | 11x383 | 16 | 7513 | 564 | 22 |
| | | | 13 | | 7 | | 7x13x89 | 26 | 8099 | 586 | 22 |
| | | | | | | | | 22 | 8707 | 608 | 22 |
| | | | | | | | | 22 | 9337 | 630 | 22 |
| | | | | | 7 | | 7x1427 | 35 | 9989 | 652 | 22 |
| | | | | | | | | 16 | 10663 | 674 | 22 |
| 37 | | | | | | | 37x307 | 19 | 11359 | 696 | 22 |
| | | | 13 | | | | 13x929 | 17 | 12077 | 718 | 22 |
| | | | | | 7 | | 7x1831 | 19 | 12817 | 740 | 22 |
| 37 | | | | | | | 37x367 | 25 | 13579 | 762 | 22 |
| | | | | | | | 53x271 | 17 | 14363 | 784 | 22 |
| | | | | 11 | 7 | | 7x11x197 | 22 | 15169 | 806 | 22 |
| | | 17 | | | | | 17x941 | 31 | 15997 | 828 | 22 |
| | | 17 | | | | | 17x991 | 26 | 16847 | 850 | 22 |
| | 29 | | 13 | | | | 13x29x47 | 25 | 17719 | 872 | 22 |
| | | | | | 7 | | 7x2659 | 19 | 18613 | 894 | 22 |
| | | | | | | | 59x331 | 26 | 19529 | 916 | 22 |
| | | | | | | | 97x211 | 19 | 20467 | 938 | 22 |
| | | | | | 7 | | 7x3061 | 16 | 21427 | 960 | 22 |
| | | | | | | | | 17 | 22409 | 982 | 22 |
| | | | 13 | | | | 13x1801 | 13 | 23413 | 1004 | 22 |
| | | | | | | | | 22 | 24439 | 1026 | 22 |
| | | | | 11 | 7 | | 7x11x331 | 26 | 25487 | 1048 | 22 |
| | | | | | | | | 25 | 26557 | 1070 | 22 |
| | | | | | | | 43x643 | 28 | 27649 | 1092 | 22 |
| | | | | | 7 | | 7x7x587 | 26 | 28763 | 1114 | 22 |
| | 29 | | | | | | 29x1031 | 37 | 29899 | 1136 | 22 |
| | | | 13 | | | | 13x2389 | 16 | 31057 | 1158 | 22 |
| | | | | | | | | 17 | 32237 | 1180 | 22 |
| | | 17 | | | 7 | | 7x17x281 | 22 | 33439 | 1202 | 22 |
| | | 17 | | | | | 17x2039 | 22 | 34663 | 1224 | 22 |
| | | | | | | | 149x241 | 26 | 35909 | 1246 | 22 |
| | | | | | 7 | | 7x47x113 | 25 | 37177 | 1268 | 22 |
| | | | 13 | 11 | | | 11x13x269 | 28 | 38467 | 1290 | 22 |
| | | | | | | | | 35 | 39779 | 1312 | 22 |
| | | | | | | | | 10 | 41113 | 1334 | 22 |
| | | | | | 7 | | 7x6067 | 25 | 42469 | 1356 | 22 |
| | | | | | | | 163x269 | 26 | 43847 | 1378 | 22 |
| | | | | | | | | 22 | 45247 | 1400 | 22 |
| | | | | | 7 | | 7x59x113 | 31 | 46669 | 1422 | 22 |
| | | | 13 | | | | 13x3701 | 17 | 48113 | 1444 | 22 |
| | | | | | | | 43x1153 | 34 | 49579 | 1466 | 22 |
| | | | | | | | 223x229 | 19 | 51067 | 1488 | 22 |
| 37 | 29 | | | | 7 | | 7x7x29x37 | 26 | 52577 | 1510 | 22 |
| | | | | 11 | | | 11x4919 | 19 | 54109 | 1532 | 22 |
| | | | | | | | | 25 | 55663 | 1554 | 22 |
| 37 | | 17 | 13 | | 7 | | 7x13x17x37 | 26 | 57239 | 1576 | 22 |
| | | 17 | | | | | 17x3461 | 31 | 58837 | 1598 | 22 |
| | | | | | | | | 22 | 60457 | 1620 | 22 |
| | | | | | | | | 26 | 62099 | 1642 | 22 |
| | | | | | 7 | | 7x9109 | 25 | 63763 | 1664 | 22 |
| | | | | | | | | 28 | 65449 | 1686 | 22 |
| | | | | | | | | 26 | 67157 | 1708 | 22 |
| | | | 13 | | 7 | | 7x13x757 | 37 | 68887 | 1730 | 22 |
| | | | | | | | | 25 | 70639 | 1752 | 22 |
| | 29 | | | 11 | | | 11x29x227 | 17 | 72413 | 1774 | 22 |
| | | | | | | | | 22 | 74209 | 1796 | 22 |

**Note :** the number of "spacings" ( or lines ) which lie between two successive prime factors of the same value, corresponds to the number of successive spiral windings of the Square-Root-Spiral ( "Einstein-Spiral" ), which lie between the two numbers which contain these prime factors ( see FIG. 6-A / 6-C ) !

**Table 5-A :** Periodic occurring Prime Factors in the „Prime-Number-Spiral-Graphs" (-number sequences ) **Q3** and **P20-G1**
( see also FIG. 6-B / 6-C and Tables 6-B1 / 6-C1 )

| some Prime Factors presented in tabular form | | | | | <-------- periodic occurrence of Individual Prime Factors → expressed through the number of spacings = (X) | Prime Factors of the none Prime-Numbers | sum of the digits | Q3 | | |
|---|---|---|---|---|---|---|---|---|---|---|
| 73 | 37 | 31 | 13 | 11 | | | | Q3 | Q3' | Q3'' |
| | | | | | | | 4 | 13 | | |
| | | | | | | | 10 | 37 | 24 | |
| | | | | | | | 11 | 83 | 46 | 22 |
| | | | | | | | 7 | 151 | 68 | 22 |
| | | | | | | | 7 | 241 | 90 | 22 |
| | | | | | | | 11 | 353 | 112 | 22 |
| | | | | | | | 19 | 487 | 134 | 22 |
| | | | | | (11) | | 13 | 643 | 156 | 22 |
| | | | | | | | 11 | 821 | 178 | 22 |
| | | | | | | | 4 | 1021 | 200 | 22 |
| | | | | 11 | | 11x113 | 10 | 1243 | 222 | 22 |
| | | | | | (13) | | 20 | 1487 | 244 | 22 |
| | | | | | | | 16 | 1753 | 266 | 22 |
| | | | 13 | | | 13x157 | 7 | 2041 | 288 | 22 |
| | | | | | | | 11 | 2351 | 310 | 22 |
| | | | | | | | 19 | 2683 | 332 | 22 |
| | | | | | (20/17) | | 13 | 3037 | 354 | 22 |
| | | | | | | | 11 | 3413 | 376 | 22 |
| | 37 | | | | | 37x103 | 13 | 3811 | 398 | 22 |
| | | | | | | | 10 | 4231 | 420 | 22 |
| | | | | | | | 20 | 4673 | 442 | 22 |
| | | | | 11 | | 11x467 | 16 | 5137 | 464 | 22 |
| | | | | | | | 16 | 5623 | 486 | 22 |
| | | | | | | | 11 | 6131 | 508 | 22 |
| | | | | | | | 19 | 6661 | 530 | 22 |
| | | | | | | | 13 | 7213 | 552 | 22 |
| | | | 13 | | | 13x599 | 29 | 7787 | 574 | 22 |
| | | | | | | 83x101 | 22 | 8383 | 596 | 22 |
| | | | | | | | 10 | 9001 | 618 | 22 |
| | | 31 | | | | 31x311 | 20 | 9641 | 640 | 22 |
| | | | | | | | 7 | 10303 | 662 | 22 |
| | | | | | | | 25 | 10987 | 684 | 22 |
| | | | | 11 | | 11x1063 | 20 | 11693 | 706 | 22 |
| | | | | | | | 10 | 12421 | 728 | 22 |
| | | | | | | | 13 | 13171 | 750 | 22 |
| 73 | | | | | (35/38) | 73x191 | 20 | 13943 | 772 | 22 |
| | | | | | | | 22 | 14737 | 794 | 22 |
| | | | | | | 103x151 | 19 | 15553 | 816 | 22 |
| | 37 | | | | (31) | 37x443 | 20 | 16391 | 838 | 22 |
| | | | 13 | | | 13x1327 | 16 | 17251 | 860 | 22 |
| | | | | | | | 16 | 18133 | 882 | 22 |
| | | | | | | | 20 | 19037 | 904 | 22 |
| | | | | | | | 28 | 19963 | 926 | 22 |
| | | | | 11 | | 11x1901 | 13 | 20911 | 948 | 22 |
| | | | | | | | 20 | 21881 | 970 | 22 |
| | | | | | | 89x257 | 22 | 22873 | 992 | 22 |
| | | | | | | | 28 | 23887 | 1014 | 22 |
| | | | | | | | 20 | 24923 | 1036 | 22 |
| | | | | | | | 25 | 25981 | 1058 | 22 |
| | | | | | | | 16 | 27061 | 1080 | 22 |
| | | | | | | | 20 | 28163 | 1102 | 22 |
| | | | | | | | 28 | 29287 | 1124 | 22 |
| | | | 13 | | | 13x2341 | 13 | 30433 | 1146 | 22 |
| | | | | | | | 11 | 31601 | 1168 | 22 |
| | | | | 11 | | 11x11x271 | 22 | 32791 | 1190 | 22 |
| | 37 | | | | | 37x919 | 10 | 34003 | 1212 | 22 |
| | | | | | | 167x211 | 20 | 35237 | 1234 | 22 |
| | | | | | | | 25 | 36493 | 1256 | 22 |
| | | | | | | 107x353 | 25 | 37771 | 1278 | 22 |
| | | | | | | 89x439 | 20 | 39071 | 1300 | 22 |
| | | 31 | | | | 31x1303 | 19 | 40393 | 1322 | 22 |
| | | | | | | | 22 | 41737 | 1344 | 22 |
| | | | | | | | 11 | 43103 | 1366 | 22 |
| | | | | | | | 22 | 44491 | 1388 | 22 |
| | | | | | | 197x233 | 19 | 45901 | 1410 | 22 |
| | | | 13 | 11 | | 11x13x331 | 20 | 47333 | 1432 | 22 |
| | | | | | | | 34 | 48787 | 1454 | 22 |
| | | | | | | | 16 | 50263 | 1476 | 22 |
| | | | | | | 191x271 | 20 | 51761 | 1498 | 22 |
| | | | | | | | 19 | 53281 | 1520 | 22 |
| 73 | | | | | | 73x751 | 22 | 54823 | 1542 | 22 |
| | | | | | | 113x499 | 29 | 56387 | 1564 | 22 |
| | | | | | | | 31 | 57973 | 1586 | 22 |
| | | | | | | | 28 | 59581 | 1608 | 22 |
| | | | | | | | 11 | 61211 | 1630 | 22 |
| | 37 | | | | | 37x1699 | 25 | 62863 | 1652 | 22 |
| | | | | 11 | | 11x5867 | 25 | 64537 | 1674 | 22 |
| | | | | | | 107x619 | 20 | 66233 | 1696 | 22 |
| | | | 13 | | | 13x5227 | 28 | 67951 | 1718 | 22 |
| | | | | | | | 31 | 69691 | 1740 | 22 |
| | | | | | | | 20 | 71453 | 1762 | 22 |
| | | | | | | | 22 | 73237 | 1784 | 22 |

**Note :** the number of "spacings" ( or lines ) which lie between two successive prime factors of the same value, corresponds to the number of successive spiral windings of the Square-Root-Spiral (Einstein-Spiral) , which lie between the two numbers which contain these prime factors ( see FIG. 6-B / 6-C ) !

| some Prime Factors presented in tabular form | | | | | <-------- periodic occurrence of Individual Prime Factors → expressed through the number of spacings = (X) | Prime Factors of the none Prime-Numbers | sum of the digits | P20-G1 | | |
|---|---|---|---|---|---|---|---|---|---|---|
| 73 | 67 | 31 | 23 | 13 | | | | G1 | G1' | G1'' |
| | | | | | | | 4 | 103 | | |
| | | | | | | | 11 | 173 | 70 | |
| | | | | | | | 11 | 263 | 90 | 20 |
| | | | | | | | 13 | 373 | 110 | 20 |
| | | | | | | | 8 | 503 | 130 | 20 |
| | | | | | (13) | | 14 | 653 | 150 | 20 |
| | | | | | | | 13 | 823 | 170 | 20 |
| | | | | | (1/30) | | 5 | 1013 | 190 | 20 |
| | | | | | | | 8 | 1223 | 210 | 20 |
| | | | | | | | 13 | 1453 | 230 | 20 |
| | | | | 13 | | 13x131 | 11 | 1703 | 250 | 20 |
| | | | | | (43/30) | | 20 | 1973 | 270 | 20 |
| 73 | | 31 | | | | 31x73 | 13 | 2263 | 290 | 20 |
| | | 31 | | | | 31x83 | 17 | 2573 | 310 | 20 |
| | | | | | (2/21) | | 14 | 2903 | 330 | 20 |
| | | | | | | | 13 | 3253 | 350 | 20 |
| | | | | | | | 14 | 3623 | 370 | 20 |
| | | | | | | | 8 | 4013 | 390 | 20 |
| | | | | | (13/54) | | 13 | 4423 | 410 | 20 |
| | | | 23 | | | 23x211 | 20 | 4853 | 430 | 20 |
| | | | | | | | 11 | 5303 | 450 | 20 |
| | | | 23 | | | 23x251 | 22 | 5773 | 470 | 20 |
| | | | | | | | 17 | 6263 | 490 | 20 |
| | | | | 13 | | 13x521 | 23 | 6773 | 510 | 20 |
| | 67 | | | | | 67x109 | 13 | 7303 | 530 | 20 |
| | | | | | | | 23 | 7853 | 550 | 20 |
| | | | | | | | 17 | 8423 | 570 | 20 |
| | | | | | | | 13 | 9013 | 590 | 20 |
| | | | | | | | 20 | 9623 | 610 | 20 |
| | | | | | | | 11 | 10253 | 630 | 20 |
| | | | | | | | 13 | 10903 | 650 | 20 |
| | | | | | | 71x163 | 17 | 11573 | 670 | 20 |
| | | | | | | | 14 | 12263 | 690 | 20 |
| | | | | | | | 22 | 12973 | 710 | 20 |
| | | | | | | 71x193 | 14 | 13703 | 730 | 20 |
| | | | | | | 97x149 | 17 | 14453 | 750 | 20 |
| | | | | 13 | | 13x1171 | 13 | 15223 | 770 | 20 |
| | 67 | | | | | 67x239 | 11 | 16013 | 790 | 20 |
| | | | | | | | 20 | 16823 | 810 | 20 |
| | | | | | | 127x139 | 22 | 17653 | 830 | 20 |
| | | | | | | | 17 | 18503 | 850 | 20 |
| | | | | | | | 23 | 19373 | 870 | 20 |
| | | | 23 | | | 23x881 | 13 | 20263 | 890 | 20 |
| | | 31 | | | | 31x683 | 14 | 21173 | 910 | 20 |
| | | 31 | 23 | | | 23x31x31 | 8 | 22103 | 930 | 20 |
| | | | | | | | 13 | 23053 | 950 | 20 |
| | | | | | | | 11 | 24023 | 970 | 20 |
| | | | | | | | 11 | 25013 | 990 | 20 |
| | | | | | | 53x491 | 13 | 26023 | 1010 | 20 |
| | | | | 13 | | 13x2081 | 17 | 27053 | 1030 | 20 |
| | | | | | | 157x179 | 14 | 28103 | 1050 | 20 |
| | | | | | | | 22 | 29173 | 1070 | 20 |
| | | | | | | 53x571 | 14 | 30263 | 1090 | 20 |
| | | | | | | 137x229 | 17 | 31373 | 1110 | 20 |
| | | | | | | | 13 | 32503 | 1130 | 20 |
| 73 | | | | | | 73x461 | 20 | 33653 | 1150 | 20 |
| | | | | | | 97x359 | 20 | 34823 | 1170 | 20 |
| | | | | | | | 13 | 36013 | 1190 | 20 |
| | | | | | | | 17 | 37223 | 1210 | 20 |
| | | | | | | | 23 | 38453 | 1230 | 20 |
| | | | | | | | 22 | 39703 | 1250 | 20 |
| | | | | | | | 23 | 40973 | 1270 | 20 |
| | | | | 13 | | 13x3251 | 17 | 42263 | 1290 | 20 |
| | | | | | | | 22 | 43573 | 1310 | 20 |
| | | | | | | 83x541 | 20 | 44903 | 1330 | 20 |
| | | | 23 | | | 23x2011 | 20 | 46253 | 1350 | 20 |
| | | | | | | | 22 | 47623 | 1370 | 20 |
| | | | 23 | | | 23x2131 | 17 | 49013 | 1390 | 20 |
| | | | | | | | 14 | 50423 | 1410 | 20 |
| | | | | | | | 22 | 51853 | 1430 | 20 |
| | | | | | | 151x353 | 14 | 53303 | 1450 | 20 |
| | | | | | | | 26 | 54773 | 1470 | 20 |
| | | | | | | | 22 | 56263 | 1490 | 20 |
| | | | | | | | 29 | 57773 | 1510 | 20 |
| | | 31 | | | | 31x1913 | 20 | 59303 | 1530 | 20 |
| | | 31 | | 13 | | 13x31x151 | 22 | 60853 | 1550 | 20 |
| | | | | | | | 17 | 62423 | 1570 | 20 |
| | | | | | | | 14 | 64013 | 1590 | 20 |
| | | | | | | 137x479 | 16 | 65623 | 1610 | 20 |
| | | | | | | 109x617 | 23 | 67253 | 1630 | 20 |
| | | | | | | | 26 | 68903 | 1650 | 20 |
| | | | | | | | 22 | 70573 | 1670 | 20 |
| | | | | | | 127x569 | 20 | 72263 | 1690 | 20 |
| | | | | | | | 29 | 73973 | 1710 | 20 |
| | | | | | | | 22 | 75703 | 1730 | 20 |
| 73 | | | | | | 73x1061 | 26 | 77453 | 1750 | 20 |
| | | | | | | 227x349 | 23 | 79223 | 1770 | 20 |
| | | | | | | | 13 | 81013 | 1790 | 20 |
| | | | 23 | 13 | | 13x23x277 | 23 | 82823 | 1810 | 20 |
| | | | | | | | 26 | 84653 | 1830 | 20 |
| | | | 23 | | | 23x3761 | 22 | 86503 | 1850 | 20 |
| | 67 | | | | | 67x1319 | 29 | 88373 | 1870 | 20 |

**Table 5-B :** Periodic occurring Prime Factors in the „Prime-Number-Spiral-Graphs" (-number sequences ) **S1** and **B33**
( see also FIG. 6-A / 6-C and Tables 6-A1 / 6-C1 )

| some Prime Factors presented in tabular form: 73 | 37 | 31 | 13 | 11 | <-------- periodic occurence of individual Prime Factors --> expressed through the numbers of spacings = (X) | Prime Factors of the None-Prime-Numbers | sum of the digits | S1 | S1' | S1" |
|---|---|---|---|---|---|---|---|---|---|---|
| | | | | | | | 2 | 11 | | |
| | | | | | | | 4 | 31 | 20 | |
| | | | | | | | 10 | 73 | 42 | 22 |
| | | | | | | | 11 | 137 | 64 | 22 |
| | | | | | | | 7 | 223 | 86 | 22 |
| | | | | | (11) | | 7 | 331 | 108 | 22 |
| | | | | | | | 11 | 461 | 130 | 22 |
| | | | | | | | 10 | 613 | 152 | 22 |
| | | | | | (13) | | 22 | 787 | 174 | 22 |
| | | | | | | | 20 | 983 | 196 | 22 |
| | | | | | | | 4 | 1201 | 218 | 22 |
| | | | | 11 | | 11x131 | 10 | 1441 | 240 | 22 |
| | | | 13 | | (17/20) | 13x131 | 11 | 1703 | 262 | 22 |
| | | | | | | | 25 | 1987 | 284 | 22 |
| | | | | | | | 16 | 2293 | 306 | 22 |
| | | | | | | | 11 | 2621 | 328 | 22 |
| | | | | | | | 19 | 2971 | 350 | 22 |
| | | | | | | | 13 | 3343 | 372 | 22 |
| | 37 | | | | | 37x101 | 20 | 3737 | 394 | 22 |
| | | | | | | | 13 | 4153 | 416 | 22 |
| | | | | | | | 19 | 4591 | 438 | 22 |
| | | | | | | | 11 | 5051 | 460 | 22 |
| | | | | 11 | | 11x503 | 16 | 5533 | 482 | 22 |
| | | | | | | | 16 | 6037 | 504 | 22 |
| | | | | | | | 20 | 6563 | 526 | 22 |
| | | | 13 | | | 13x547 | 10 | 7111 | 548 | 22 |
| | | | | | | | 22 | 7681 | 570 | 22 |
| | | | | | | | 20 | 8273 | 592 | 22 |
| | | | | | | | 31 | 8887 | 614 | 22 |
| | | | | | | 89x107 | 19 | 9523 | 636 | 22 |
| | | | | | | | 11 | 10181 | 658 | 22 |
| | | | | | | | 16 | 10861 | 680 | 22 |
| | | 31 | | | (31) | 31x373 | 16 | 11563 | 702 | 22 |
| | | | | 11 | | 11x1117 | 20 | 12287 | 724 | 22 |
| | | | | | | | 10 | 13033 | 746 | 22 |
| | 37 | | | | | 37x373 | 13 | 13801 | 768 | 22 |
| | | | | | | | 20 | 14591 | 790 | 22 |
| 73 | | | | | | 73x211 | 13 | 15403 | 812 | 22 |
| | | | 13 | | | 13x1249 | 19 | 16237 | 834 | 22 |
| | | | | | (38/35) | | 20 | 17093 | 856 | 22 |
| | | | | | | | 25 | 17971 | 878 | 22 |
| | | | | | | 113x167 | 25 | 18871 | 900 | 22 |
| | | | | | | | 29 | 19793 | 922 | 22 |
| | | | | | | 89x233 | 19 | 20737 | 944 | 22 |
| | | | | 11 | | 11x1973 | 13 | 21703 | 966 | 22 |
| | | | | | | | 20 | 22691 | 988 | 22 |
| | | | | | | 137x173 | 13 | 23701 | 1010 | 22 |
| | | | | | | | 19 | 24733 | 1032 | 22 |
| | | | | | | 107x241 | 29 | 25787 | 1054 | 22 |
| | | | | | | | 25 | 26863 | 1076 | 22 |
| | | | | | | | 25 | 27961 | 1098 | 22 |
| | | | 13 | | | 13x2237 | 20 | 29081 | 1120 | 22 |
| | | | | | | | 10 | 30223 | 1142 | 22 |
| | | | | | | | 22 | 31387 | 1164 | 22 |
| | | | | | | | 20 | 32573 | 1186 | 22 |
| | 37 | | | 11 | | 11x37x83 | 22 | 33781 | 1208 | 22 |
| | | | | | | 157x223 | 10 | 35011 | 1230 | 22 |
| | | | | | | | 20 | 36263 | 1252 | 22 |
| | | | | | | | 25 | 37537 | 1274 | 22 |
| | | | | | | | 25 | 38833 | 1296 | 22 |
| | | | | | | | 11 | 40151 | 1318 | 22 |
| | | | | | | | 19 | 41491 | 1340 | 22 |
| | | | | | | | 22 | 42853 | 1362 | 22 |
| | | 31 | | | | 31x1427 | 20 | 44237 | 1384 | 22 |
| | | | 13 | | | 13x3511 | 22 | 45643 | 1406 | 22 |
| | | | | | | 103x457 | 19 | 47071 | 1428 | 22 |
| | | | | 11 | | 11x11x401 | 20 | 48521 | 1450 | 22 |
| | | | | | | | 34 | 49993 | 1472 | 22 |
| | | | | | | | 25 | 51487 | 1494 | 22 |
| | | | | | | | 11 | 53003 | 1516 | 22 |
| | | | | | | | 19 | 54541 | 1538 | 22 |
| | | | | | | | 13 | 56101 | 1560 | 22 |
| | 37 | | | | | 37x1559 | 29 | 57683 | 1582 | 22 |
| | | | | | | 101x587 | 31 | 59287 | 1604 | 22 |
| | | | | | | | 19 | 60913 | 1626 | 22 |
| 73 | | | | | | 73x857 | 20 | 62561 | 1648 | 22 |
| | | | | | | | 16 | 64231 | 1670 | 22 |
| | | | 13 | 11 | | 11x13x461 | 25 | 65923 | 1692 | 22 |
| | | | | | | 239x283 | 29 | 67637 | 1714 | 22 |
| | | | | | | 173x401 | 28 | 69373 | 1736 | 22 |
| | | | | | | 83x857 | 13 | 71131 | 1758 | 22 |
| | | | | | | | 20 | 72911 | 1780 | 22 |

| some Prime Factors presented in tabular form: 61 | 43 | 29 | 13 | 11 | <-------- periodic occurence of individual Prime Factors --> expressed through the numbers of spacings = (X) | Prime Factors of the None-Prime-Numbers | sum of the digits | B33 | B33' | B33" |
|---|---|---|---|---|---|---|---|---|---|---|
| | | | | | | | 17 | 89 | | |
| | | | | | (6/5) | | 8 | 107 | 18 | |
| | | | 13 | 11 | (8/5) | 11x13 | 8 | 143 | 36 | 18 |
| | | | | | | | 17 | 197 | 54 | 18 |
| | | | | | | | 17 | 269 | 72 | 18 |
| | | | | | | | 17 | 359 | 90 | 18 |
| | | | | | | | 17 | 467 | 108 | 18 |
| | | | | | | | 17 | 593 | 126 | 18 |
| | | | | 11 | (10/19) | 11x67 | 17 | 737 | 144 | 18 |
| | | 29 | | | | 29x31 | 26 | 899 | 162 | 18 |
| | | | 13 | | | 13x83 | 17 | 1079 | 180 | 18 |
| | | | | | | | 17 | 1277 | 198 | 18 |
| | | | | | | | 17 | 1493 | 216 | 18 |
| | | | | 11 | | 11x157 | 17 | 1727 | 234 | 18 |
| | | | | | | | 26 | 1979 | 252 | 18 |
| | | | 13 | | (10/33) | 13x173 | 17 | 2249 | 270 | 18 |
| | 43 | | | | | 43x59 | 17 | 2537 | 288 | 18 |
| | | | | | | | 17 | 2843 | 306 | 18 |
| | | | | | | | 17 | 3167 | 324 | 18 |
| | | 29 | | 11 | | 11x11x29 | 17 | 3509 | 342 | 18 |
| | | | | | | 53x73 | 26 | 3869 | 360 | 18 |
| | | | | | | 31x137 | 17 | 4247 | 378 | 18 |
| | | | | | | | 17 | 4643 | 396 | 18 |
| | | | 13 | | (6/55) | 13x389 | 17 | 5057 | 414 | 18 |
| | | | | 11 | | 11x499 | 26 | 5489 | 432 | 18 |
| | | | | | | | 26 | 5939 | 450 | 18 |
| | 43 | | | | | 43x149 | 17 | 6407 | 468 | 18 |
| 61 | | | | | | 61x113 | 26 | 6893 | 486 | 18 |
| | | | 13 | | | 13x569 | 26 | 7397 | 504 | 18 |
| | | | | | | | 26 | 7919 | 522 | 18 |
| | | | | 11 | | 11x769 | 26 | 8459 | 540 | 18 |
| | | | | | | 71x127 | 17 | 9017 | 558 | 18 |
| | | | | | | 53x181 | 26 | 9593 | 576 | 18 |
| 61 | | | | | | 61x167 | 17 | 10187 | 594 | 18 |
| | | | | | | | 26 | 10799 | 612 | 18 |
| | | | | 11 | | 11x1039 | 17 | 11429 | 630 | 18 |
| | | | 13 | | | 13x929 | 17 | 12077 | 648 | 18 |
| | | | | | | | 17 | 12743 | 666 | 18 |
| | | 29 | | | | 29x463 | 17 | 13427 | 684 | 18 |
| | | | | | | 71x199 | 17 | 14129 | 702 | 18 |
| | | | | | | 31x479 | 26 | 14849 | 720 | 18 |
| | | | 13 | 11 | | 11x13x109 | 26 | 15587 | 738 | 18 |
| | | | | | | 59x277 | 17 | 16343 | 756 | 18 |
| | | | | | | | 17 | 17117 | 774 | 18 |
| | | | | | | | 26 | 17909 | 792 | 18 |
| | | | | | | | 26 | 18719 | 810 | 18 |
| | | | | 11 | | 11x1777 | 26 | 19547 | 828 | 18 |
| | | | | | | | 17 | 20393 | 846 | 18 |
| | | 29 | | | | 29x733 | 17 | 21257 | 864 | 18 |
| | | | 13 | | | 13x13x131 | 17 | 22139 | 882 | 18 |
| | | | | | | | 17 | 23039 | 900 | 18 |
| | | | | | | | 26 | 23957 | 918 | 18 |
| | | | | 11 | | 11x31x73 | 26 | 24893 | 936 | 18 |
| | | | | | | | 26 | 25847 | 954 | 18 |
| | | | 13 | | | 13x2063 | 26 | 26819 | 972 | 18 |
| | | | | | | | 26 | 27809 | 990 | 18 |
| | | | | | | | 26 | 28817 | 1008 | 18 |
| | | | | 11 | | 11x2713 | 26 | 29843 | 1026 | 18 |
| | | | | | | 67x461 | 26 | 30887 | 1044 | 18 |
| | 43 | | | | | 43x743 | 26 | 31949 | 1062 | 18 |
| | | | | | | | 17 | 33029 | 1080 | 18 |
| | | | | | | | 17 | 34127 | 1098 | 18 |
| | | | 13 | | | 13x2711 | 17 | 35243 | 1116 | 18 |
| | | | | 11 | | 11x3307 | 26 | 36377 | 1134 | 18 |
| | | | | | | | 26 | 37529 | 1152 | 18 |
| | | | | | | | 35 | 38699 | 1170 | 18 |
| | | | | | | | 35 | 39887 | 1188 | 18 |
| | | 29 | 13 | | | 13x29x109 | 17 | 41093 | 1206 | 18 |
| | | | | 11 | | 11x3847 | 17 | 42317 | 1224 | 18 |
| | 43 | | | | | 43x1013 | 26 | 43559 | 1242 | 18 |
| | | | | | | | 26 | 44819 | 1260 | 18 |
| | | | | | | 31x1487 | 26 | 46097 | 1278 | 18 |
| | | | | | | 83x571 | 26 | 47393 | 1296 | 18 |
| | | | | | | 53x919 | 26 | 48707 | 1314 | 18 |
| | | | | 11 | | 11x4549 | 17 | 50039 | 1332 | 18 |
| | | | 13 | | | 13x59x67 | 26 | 51389 | 1350 | 18 |
| | | | | | | | 26 | 52757 | 1368 | 18 |
| | | 29 | | | | 29x1867 | 17 | 54143 | 1386 | 18 |
| | | | | | | | 26 | 55547 | 1404 | 18 |
| | | | | 11 | | 11x5179 | 35 | 56969 | 1422 | 18 |
| | | | 13 | | | 13x4493 | 26 | 58409 | 1440 | 18 |
| | | | | | | 131x457 | 35 | 59867 | 1458 | 18 |

**Note :** the number of "spacings" ( or lines ) which lie between two successive prime factors of the same value, corresponds to the number of successive spiral windings of the Square-Root-Spiral ( "Einstein-Spiral" ), which lie between the two numbers which contain these prime factors ( see FIG. 6-A / 6-C ) !

# 9 The Number Spiral - by Robert Sachs → www.numberspiral.com

The following chapter shows extracts from the analysis of the " Number Spiral ", carried out by Mr. Robert Sachs.

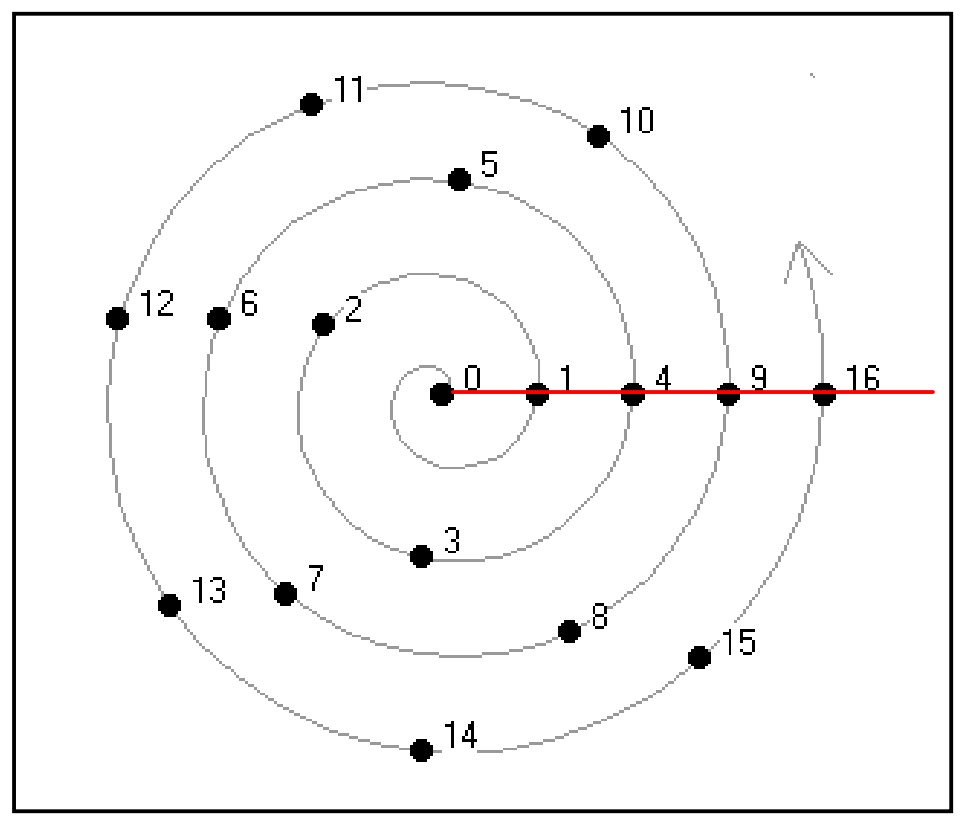

**FIG. NS-1** : Number Spiral

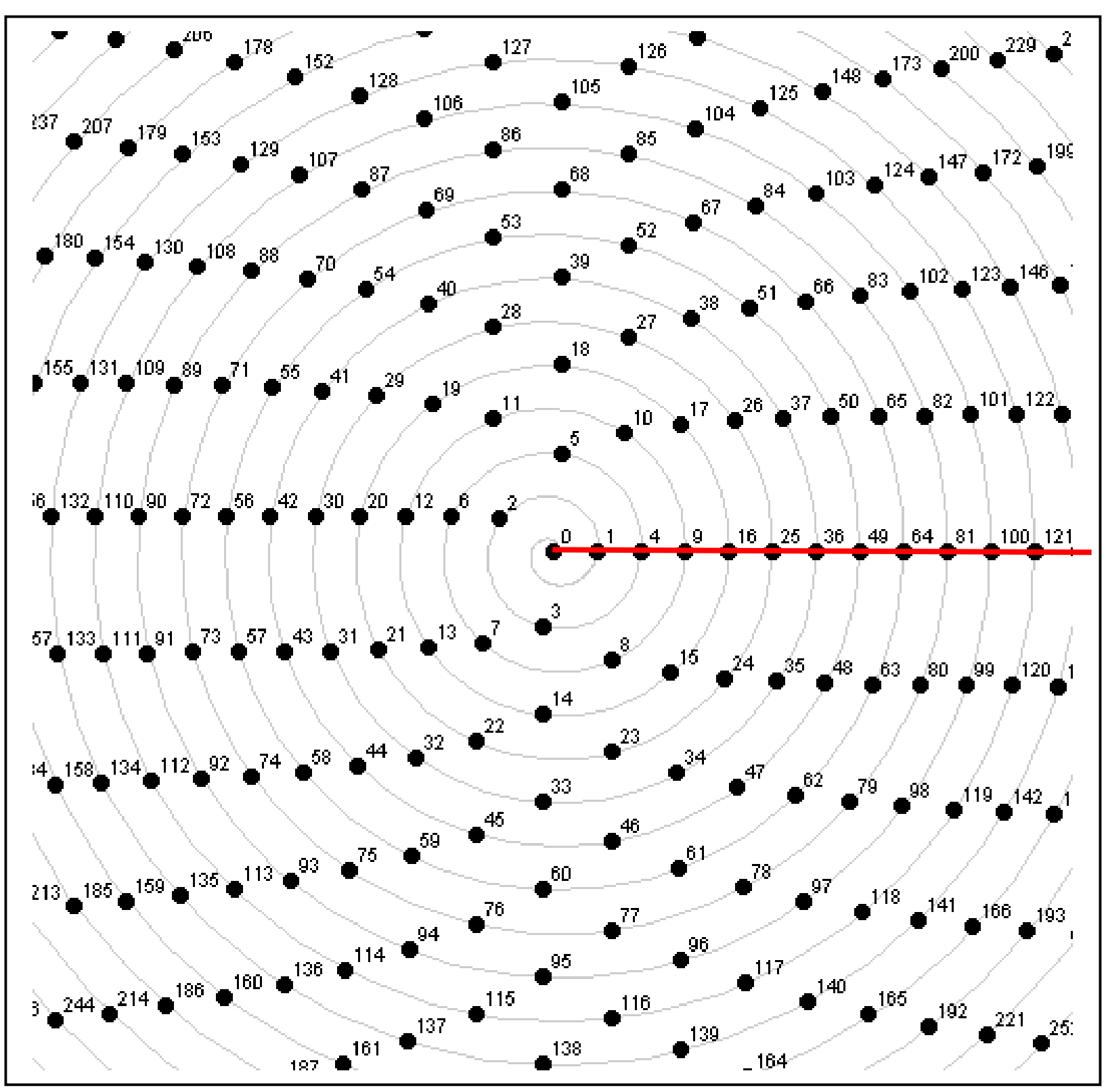

**FIG. NS-2** : Number Spiral with perfect squares lined-up on a straight line

## 9.1 Introduction

Number spirals are very simple. To make one, we just write the non-negative integers on a ribbon and roll it up with zero at the center. The trick is to arrange the spiral so all the perfect squares (1, 4, 9, 16, etc.) line up in a row on the right side: → **see Figure NS-1 & NS-2**

If we continue winding for a while and zoom out a bit, the result looks like shown on the right

If we zoom out even further and remove everything except the dots that indicate the locations of integers, we get **Figure NS-3** below. It shows 2026 dots :

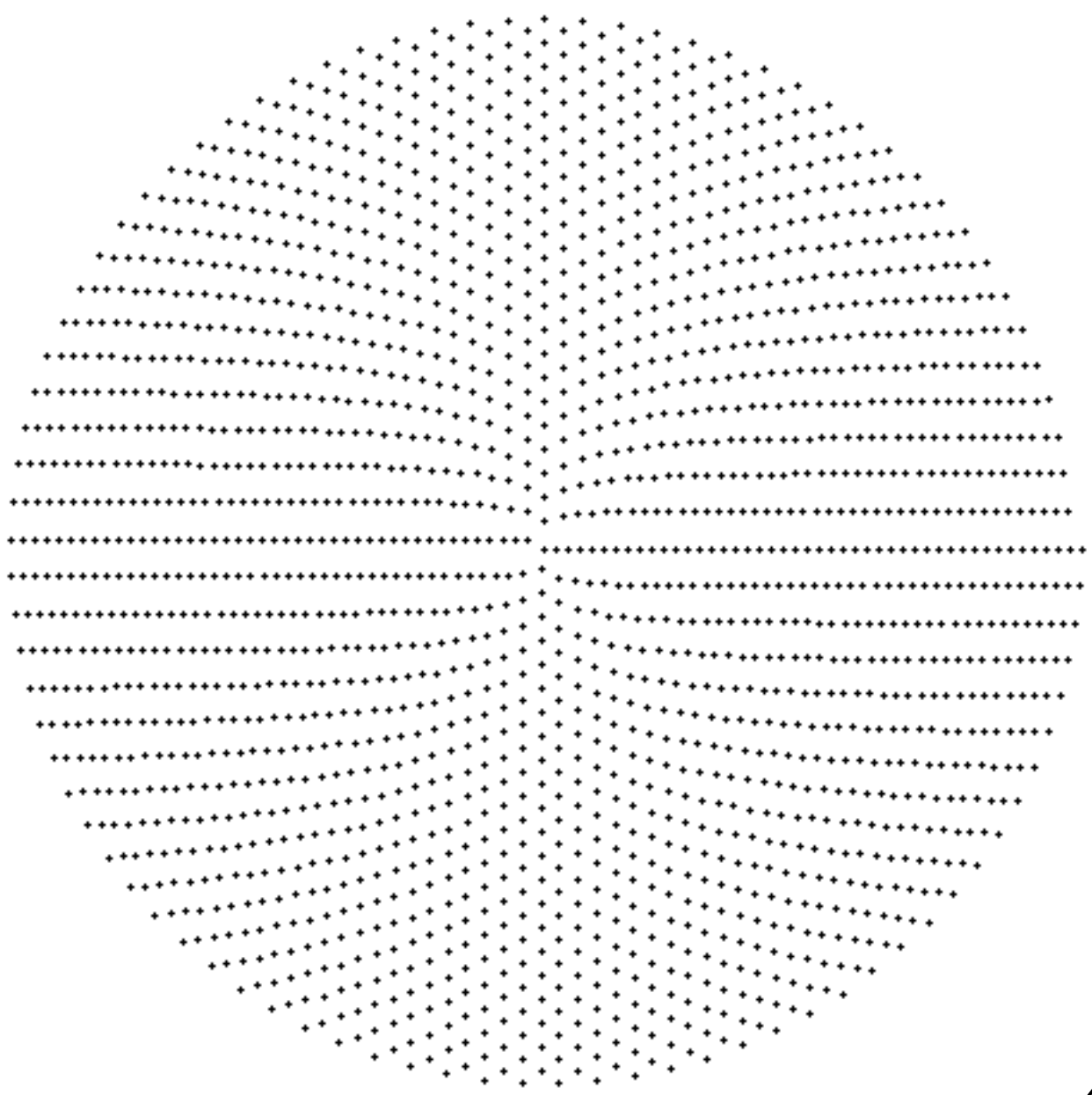

**FIG. NS-3** : Number Spiral
→ dots indicating 2026 integers

Let's try making the primes darker than the non-primes:
→ **see Figure NS-4**

The primes clearly seem to cluster along certain curves.
( see image below )

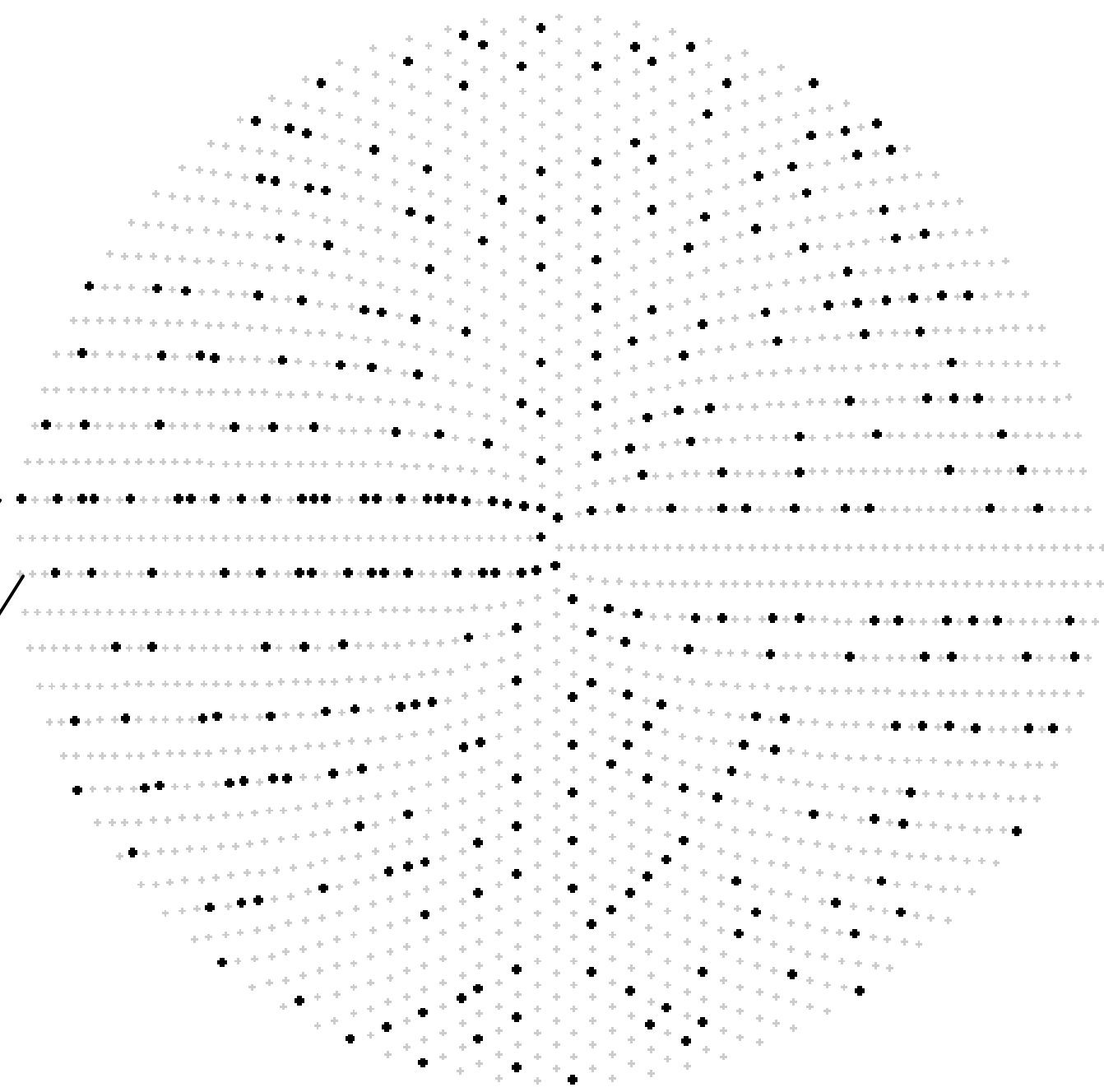

"Curve P-1" contains the numbers 5, 11, 19, 29, 41, 55, 71... which result in the quadratic polynomials :

$x^2 + 3x + 1$ or $x^2 + 5x + 5$ or $x^2 + 7x + 11$ etc.

"Curve P+1" contains the numbers 3, 7, 13, 21, 31, 43, 57... which result in the quadratic polynomials :

$x^2 + x + 1$ or $x^2 + 3x + 3$ or $x^2 + 5x + 7$ etc.

**FIG. NS-4** : Prime Numbers cluster along defined curves

Let's zoom out even further to get a better look. The following number spiral shows all the primes that occur within the first 46,656 non-negative integers. (For clarity, non-primes have been left out.)
→ **see Figure NS-5**

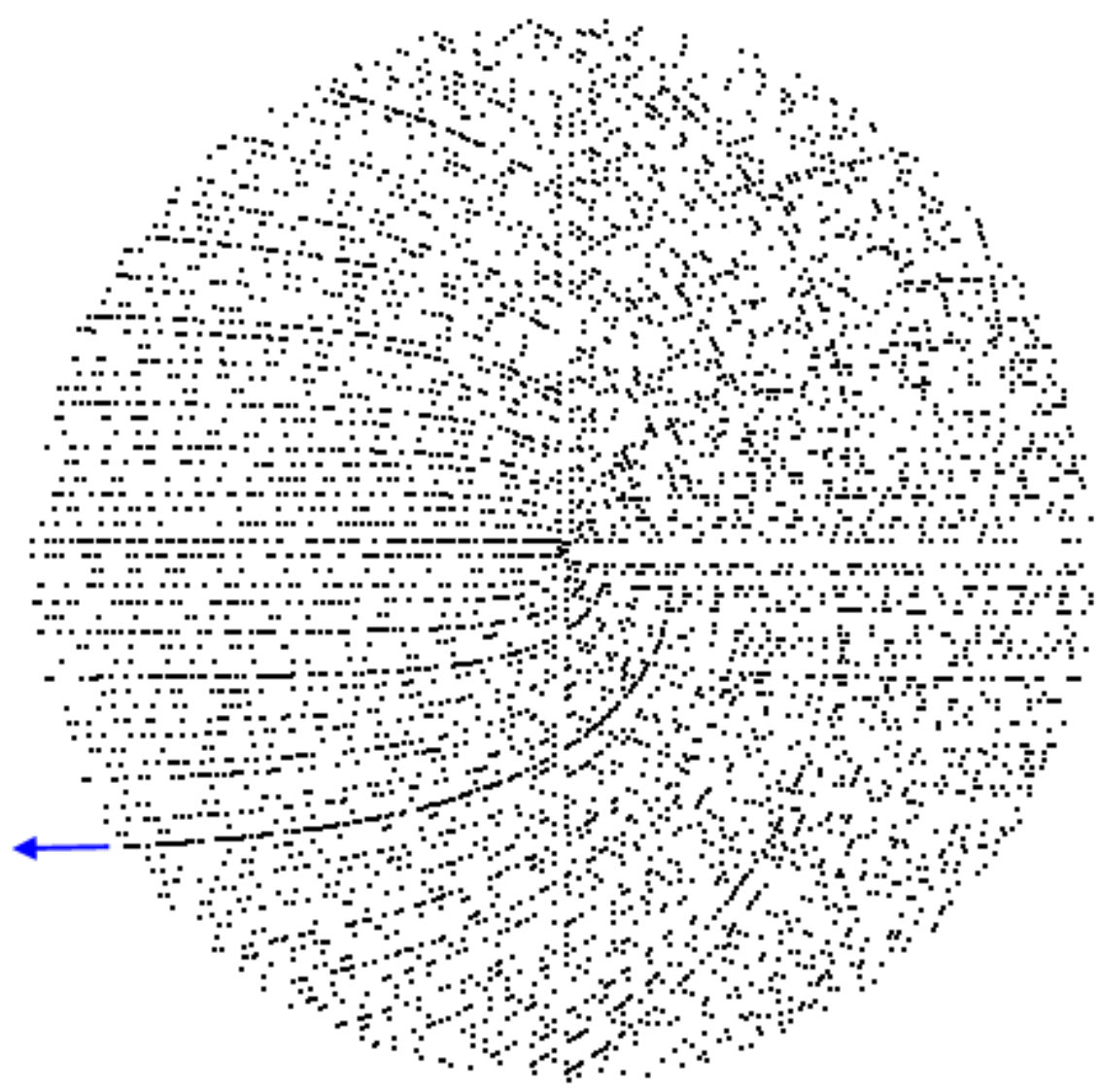

**FIG. NS-5** : Shows the Positions of the first 46656 Prime Number

It is evident that prime numbers concentrate on certain curves which run to the northwest and southwest, like the curve marked by the blue arrow.

In the following we'll investigate these patterns and try to make sense out of them.

## 9.2 Product Curves

On the previous images we saw that primes tend to line up in curves on the number spiral. In fact, the whole spiral is made of curves of this kind, and every integer belongs to an infinite number of them.

The simplest curve of this type (the one with least curvature) is the line of perfect squares, marked here in blue → **see Figure NS-6**

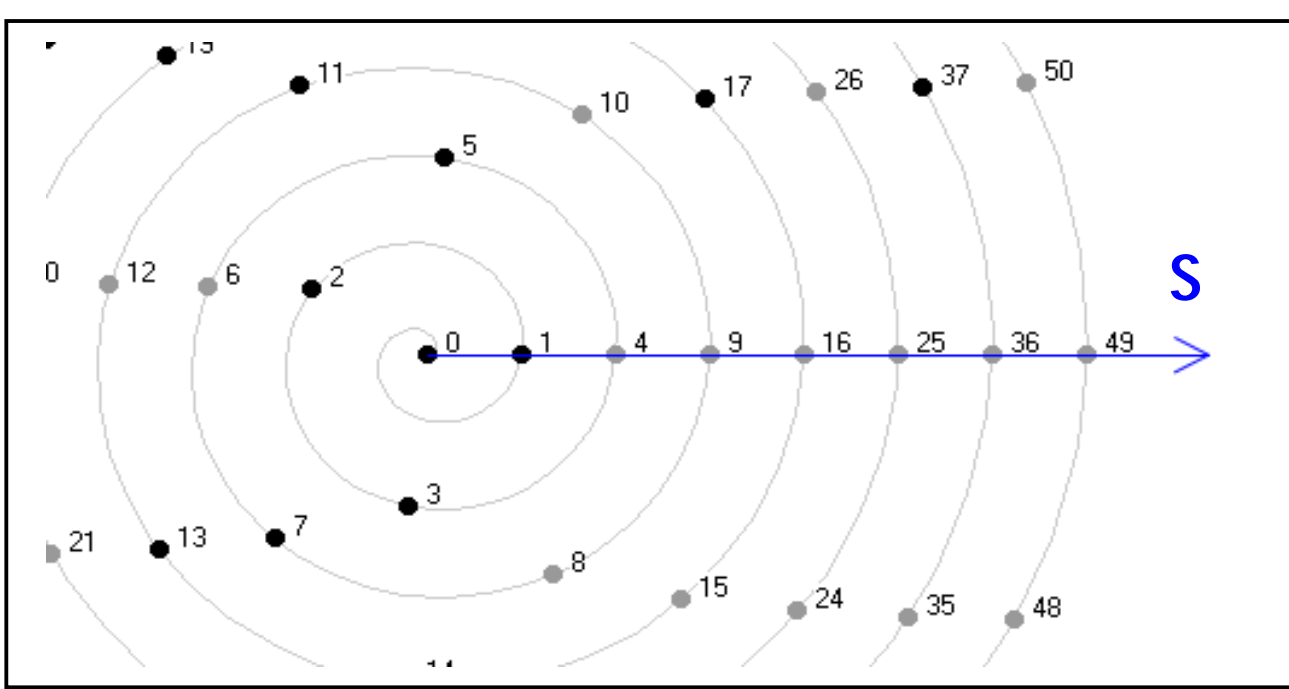

**FIG. NS-6**: "Curve S" – Line of perfect squares on the Number Spiral

For convenience, I'll call this line "curve S" for '"squares."

Here's another example → **see Figure NS-7**

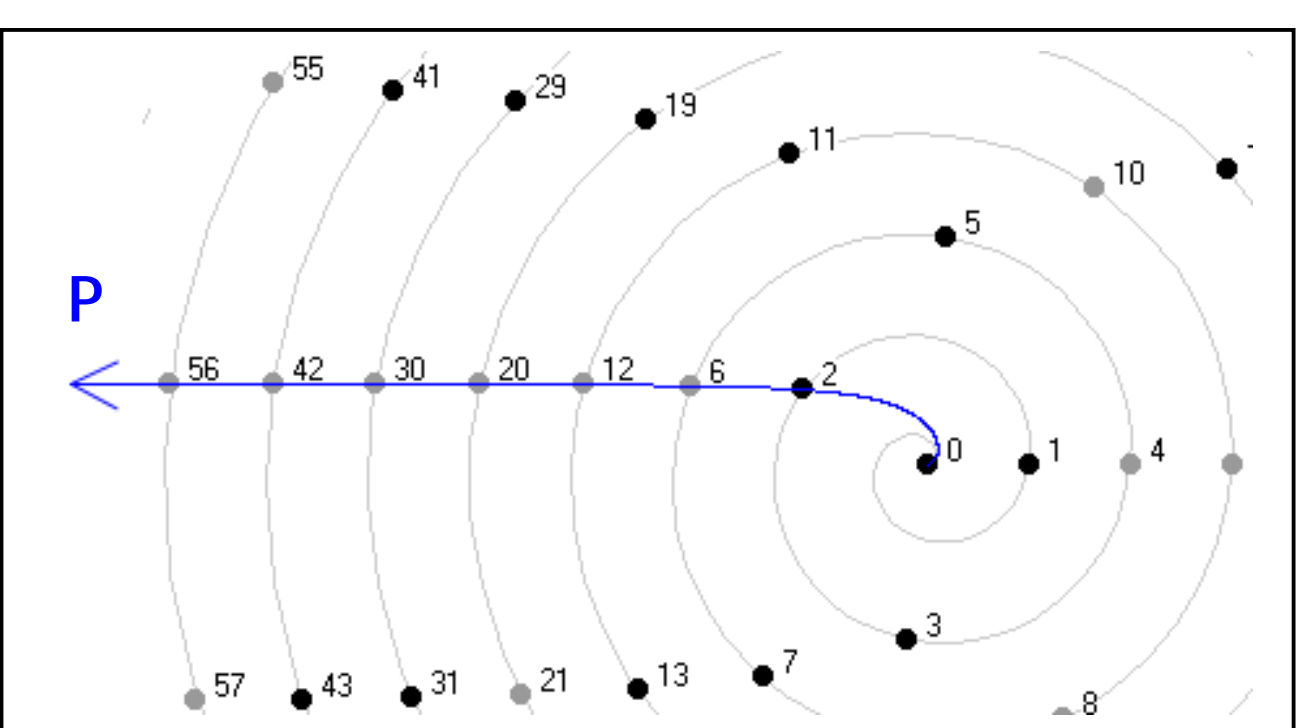

**FIG. NS-7**: "Curve P" – Line of " pronics " with difference between factors of numbers on this curve is always 1

The factors of numbers on this second curve are 1 x 2, 2 x 3, 3 x 4, etc. The difference between factors is always 1. Such numbers are called " pronics ", so I'll name this second line "Curve P".

In curve S, the difference between factors is zero. In curve P, it is one. And there are other such curves at distances of 1, 2, 4, 6 , 9, 12,..." units" to the original curves P or S .

For example curve P-6 , which is at a distance of "6 units" to P , contains the numbers 50, 66, 84, 104, 126,.... In this sequence the factors of the numbers are 5 x 10 , 6 x 11 , 7 x 12 , 8 x 13 , ...etc. , which corresponds to a difference of 5 between the factors of a number and to the difference of 1 between the factors of two successive numbers.

Or curve S-1 , which is at a distance of one unit to S , contains the numbers 15, 24, 35, 48, 63, 80,.... In this number sequence the factors of the numbers are 3 x 5 , 4 x 6 , 5 x 7 , 6 x 8 , 7 x 9 , 8 x 10 etc. , which corresponds to a difference of 2 between the factors of a number and to the difference of 1 between the factors of two successive numbers.

We can continue this way forever, increasing the distance to the curves P and S and finding new curves.

On the lefthand side in **Figure NS-8** the first ten such curves are shown:

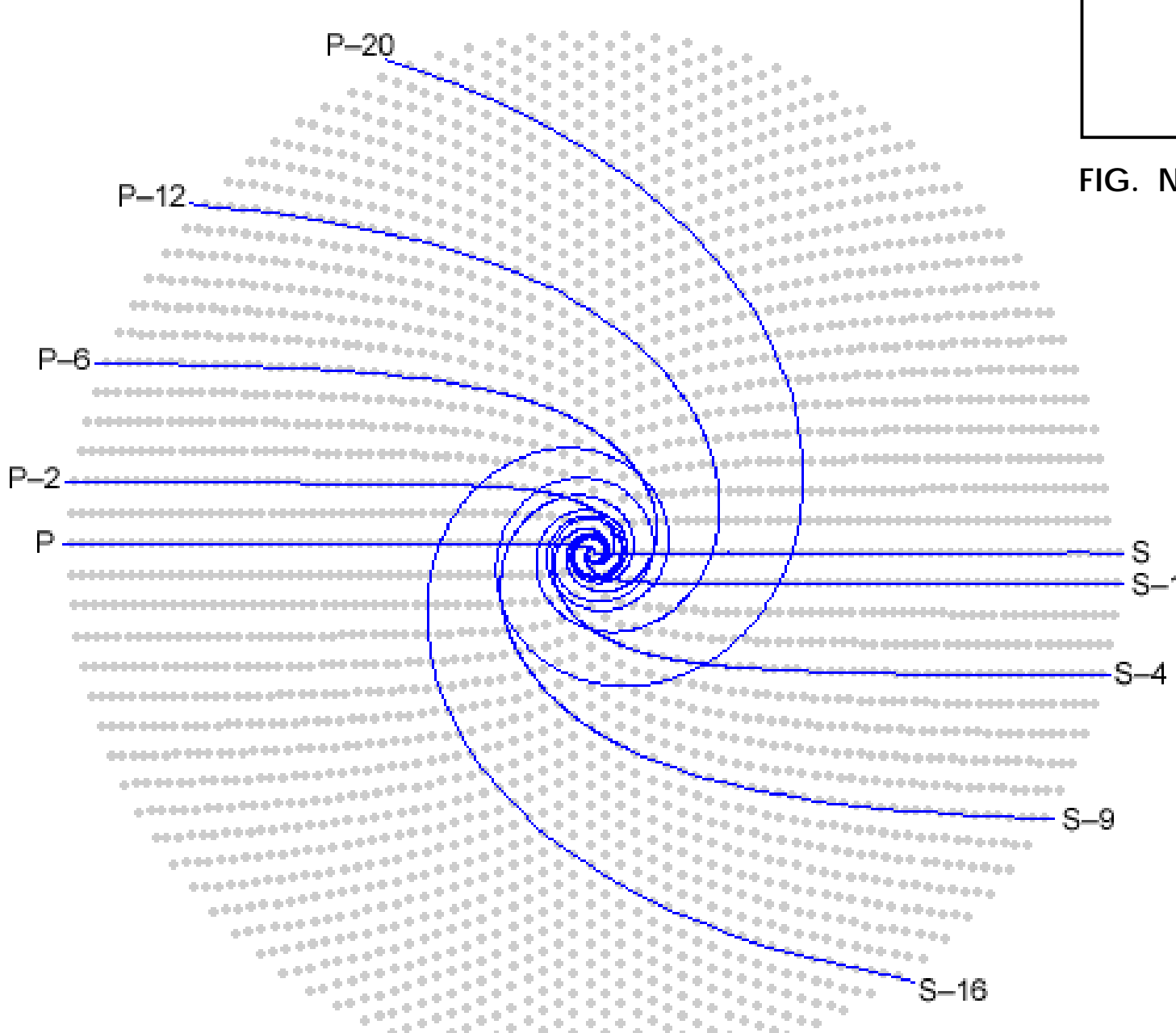

**FIG. NS-8** : Shows the first ten " pronics – curves "

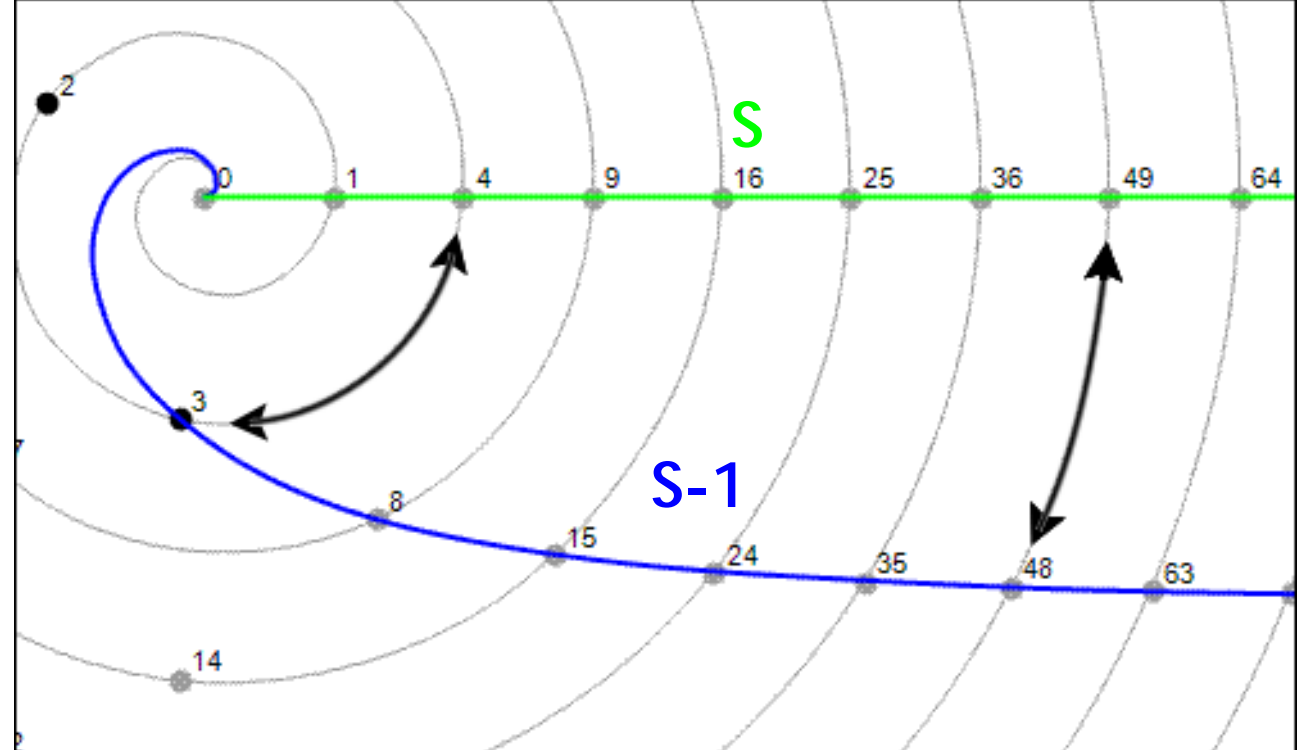

**FIG. NS-9** : Distances between the numbers on curve S and curve S-1 are constant measured along the spiral

## 9.3 Offset Curves

At the bottom of the previous page we saw a picture of the first ten product curves ( pronics – curves ). Product curves are important because every possible way of multiplying one number by another is represented on one of them.

But it turns out that they are only a special case of a more general phenomenon. Their properties are due to the fact that they are located at fixed distances from defined angles. As we will see, other curves which are located at fixed distances from other defined angles have similar properties.
To illustrate these ideas, let's look at product curve S – 1, shown above in blue ( see **Figure NS-9** ). It is located at a fixed distance from angle zero, shown in green.

At first glance the green and blue lines appear to converge at the left. But if we measure the distance between them using the spiral itself as our tape measure, the lines are always one unit apart. The black arrows show how to do this. The left arrow stretches between 3 and 4, a distance of one. The right arrow stretches between 48 and 49, also a distance of one. Even when the blue curve is at zero, it's separated (for measuring purposes) from the green curve by the piece of the spiral that runs from zero to one. No matter where we measure, the distance is always 1.
I call a constant distance of this sort an offset. When a curve is offset from an angle line, I call it an offset curve.
If we zoom out far enough, offset curves look straight. Some of them have so little curvature that we barely have to zoom. For instance:

The blue lines in **Figure NS-10** show the first offset curves ( with offset 0) of the angle lines which have a rotational angle of 0, 1/64, 1/32, 1/16, 1/8, and 1/4 in reference to curve S ( curve which contains the perfect squares). Note : one full rotation = 1.

Offset curves are important because some of them are composite. When I say that a curve is composite I mean that all the integers on them (except for the first, which is always zero, and the second, which may be prime) are non-prime. Moreover, we can predict the factorization of any integer on such curves just from knowing its location.
Every rational angle has composite offset curves.
When I say "rational angle," I mean an angle line whose measurement in rotations ( in reference to curve S ) is a rational number: 1/2 rotation, 1/3 rotation, 1473/2076 rotation, etc.

**Figure NS-11**, for example shows the 1/3 angle line marked in green ( at an angle of 120 degrees in reference to curve S ) :

And **Figure NS-12** on the right shows a few of its offset curves:

The blue offset curves are composite. (To avoid misunderstanding, let me say again that the first integer on a composite curve is zero, the second may or may not be prime, and the rest are composite.)

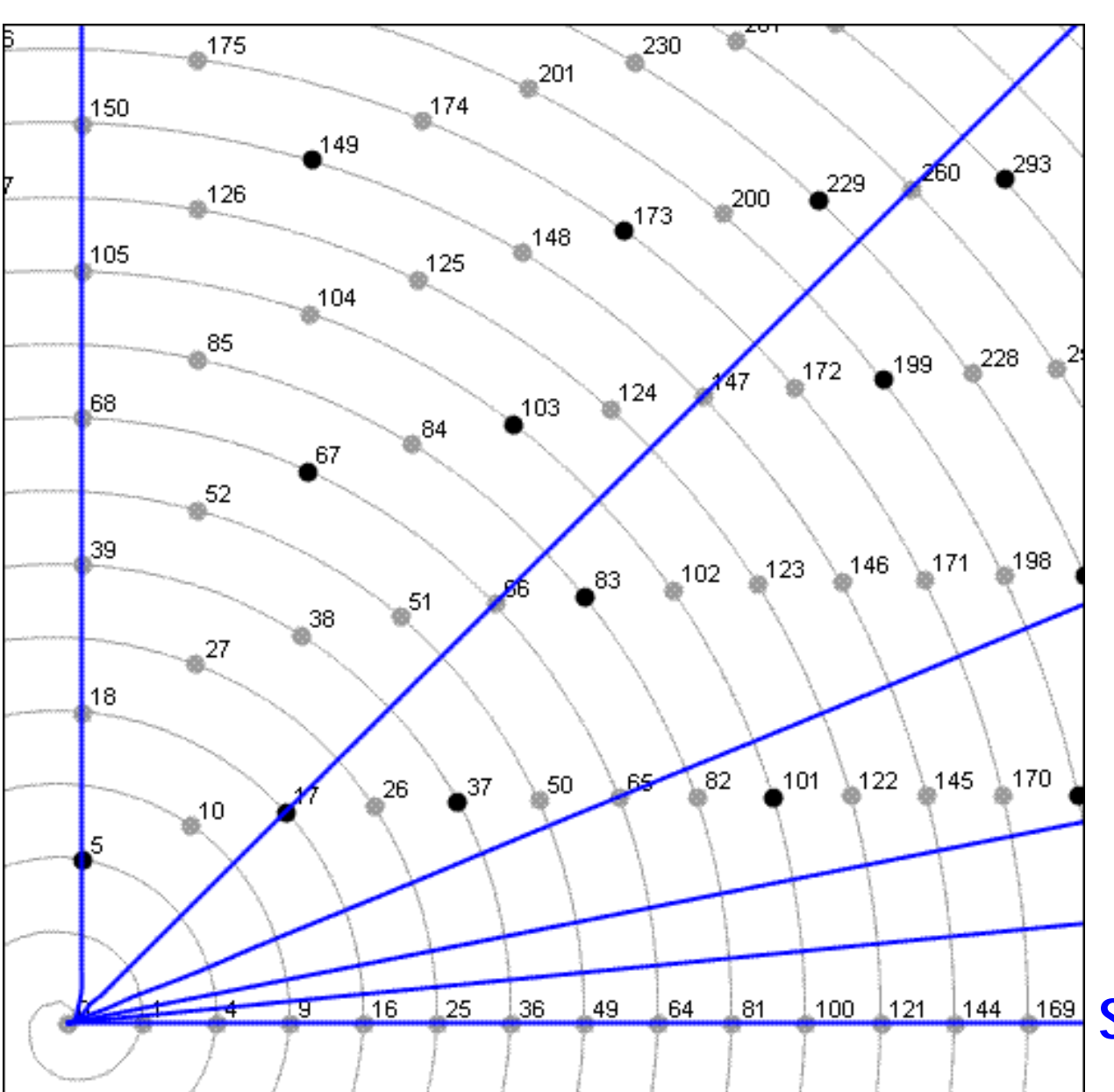

**FIG. NS-10** : First 5 Offset-Curves at rotational angles of 1/64, 1/32, 1/16, 1/8 and 1/4 in reference to Curve S

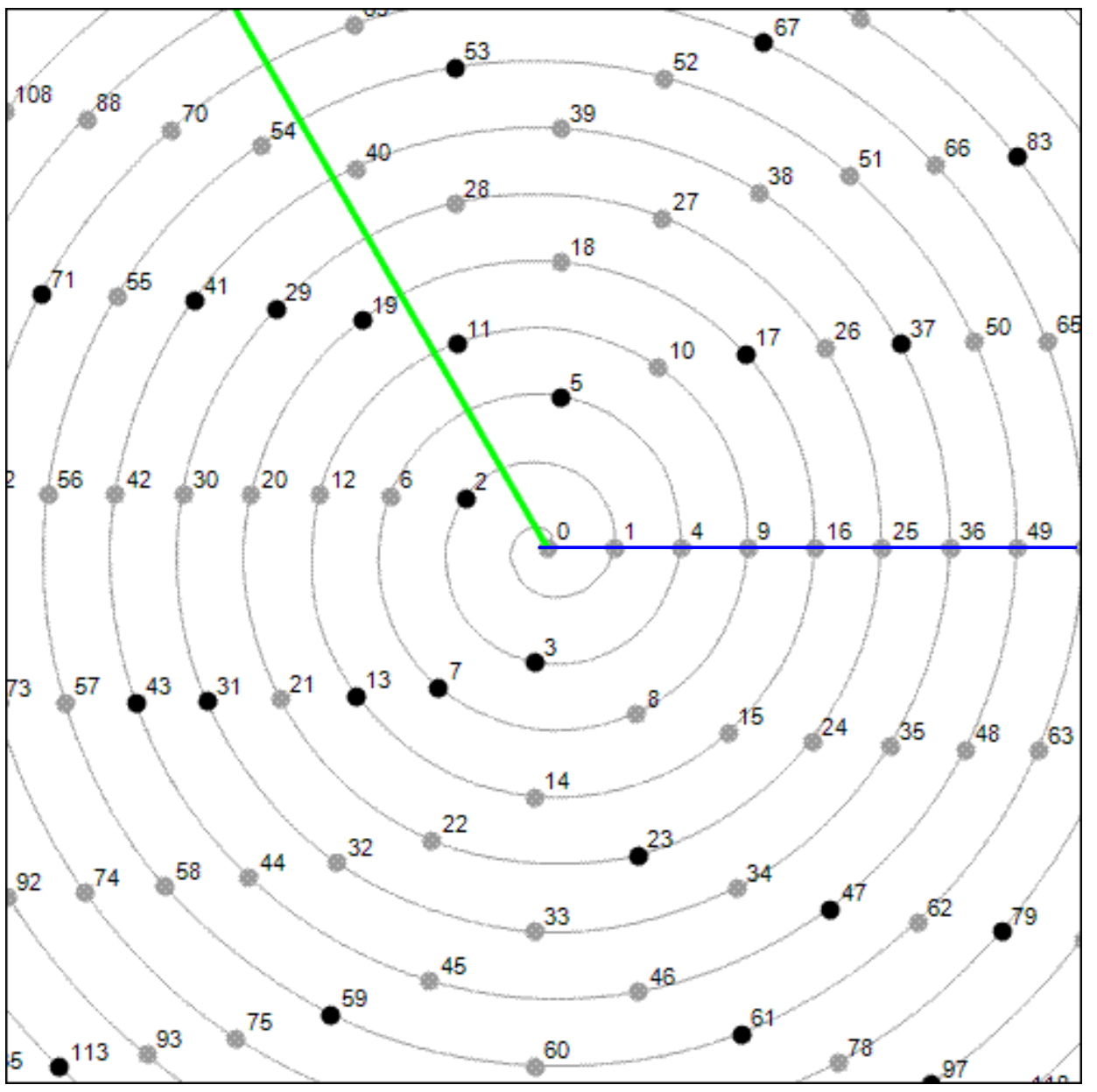

**FIG. NS-11** : Angle Line with rotational angle 1/3 (120 degrees)

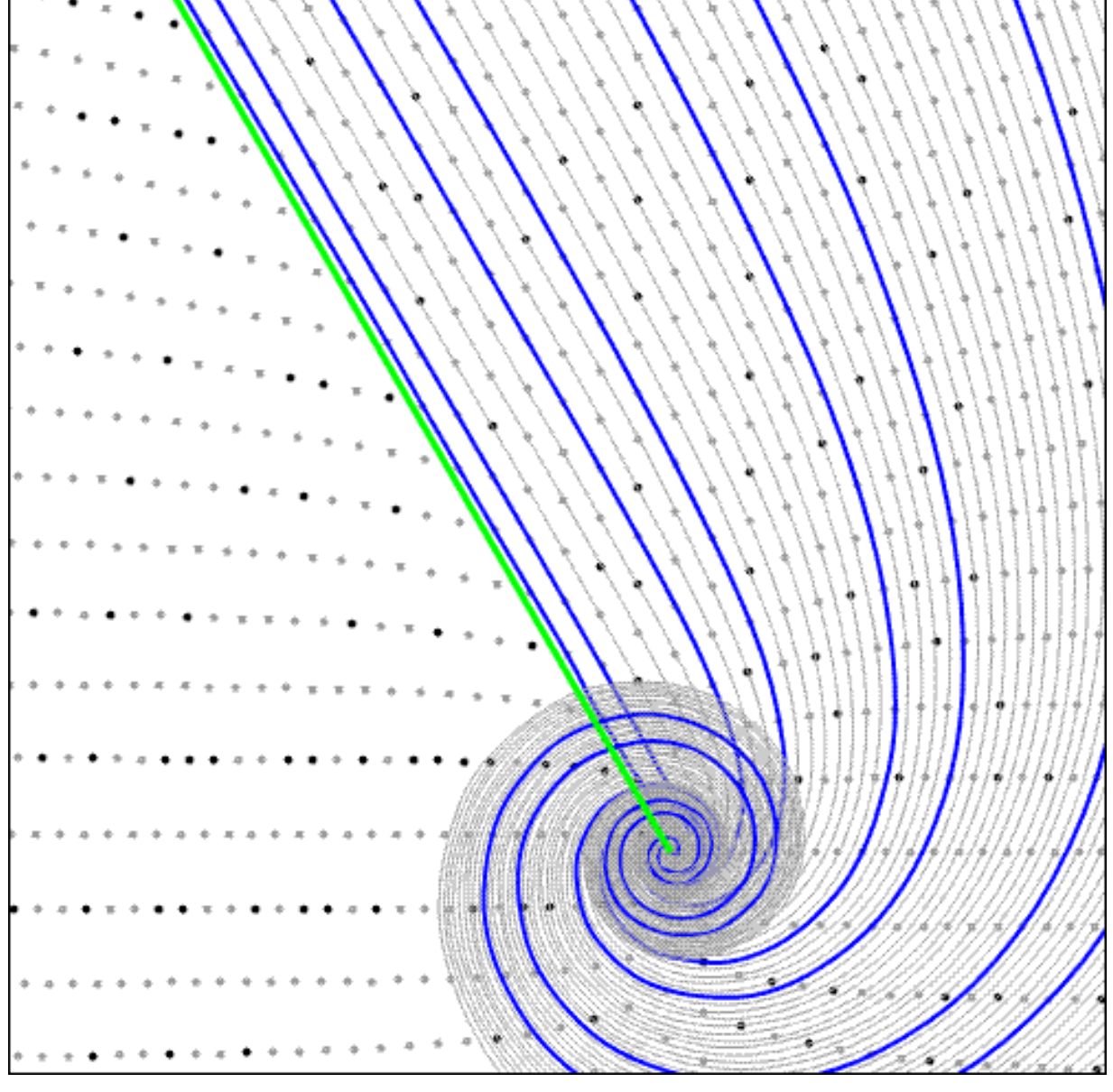

**FIG. NS-12** : Offset-Curves of angle line with angle 1/3 (120 deg.)

As I said a moment ago, every rational angle has blue offset curves of this type. Here in **Figure NS-13** for example the angle 29/33 is shown :

In this Figure you see something new : red curves. Like blue curves, red ones are composite. But red and blue sequences are composite for different reasons, so I distinguish them by color.

One last example: **Figure NS-14** shows the first few offset curves for angle 1/4 (90 deg.) → see below on the right :

I have a special reason for including angle 1/4 in my analysis. When I first wrote my website in 2003, I didn't realize that every rational angle has composite curves. Therefore I was puzzled by the prominent columns of non-prime numbers that march north and south from the center of the spiral, and I published **Figure NS-15** as a mystery to be solved later ( see below ) :

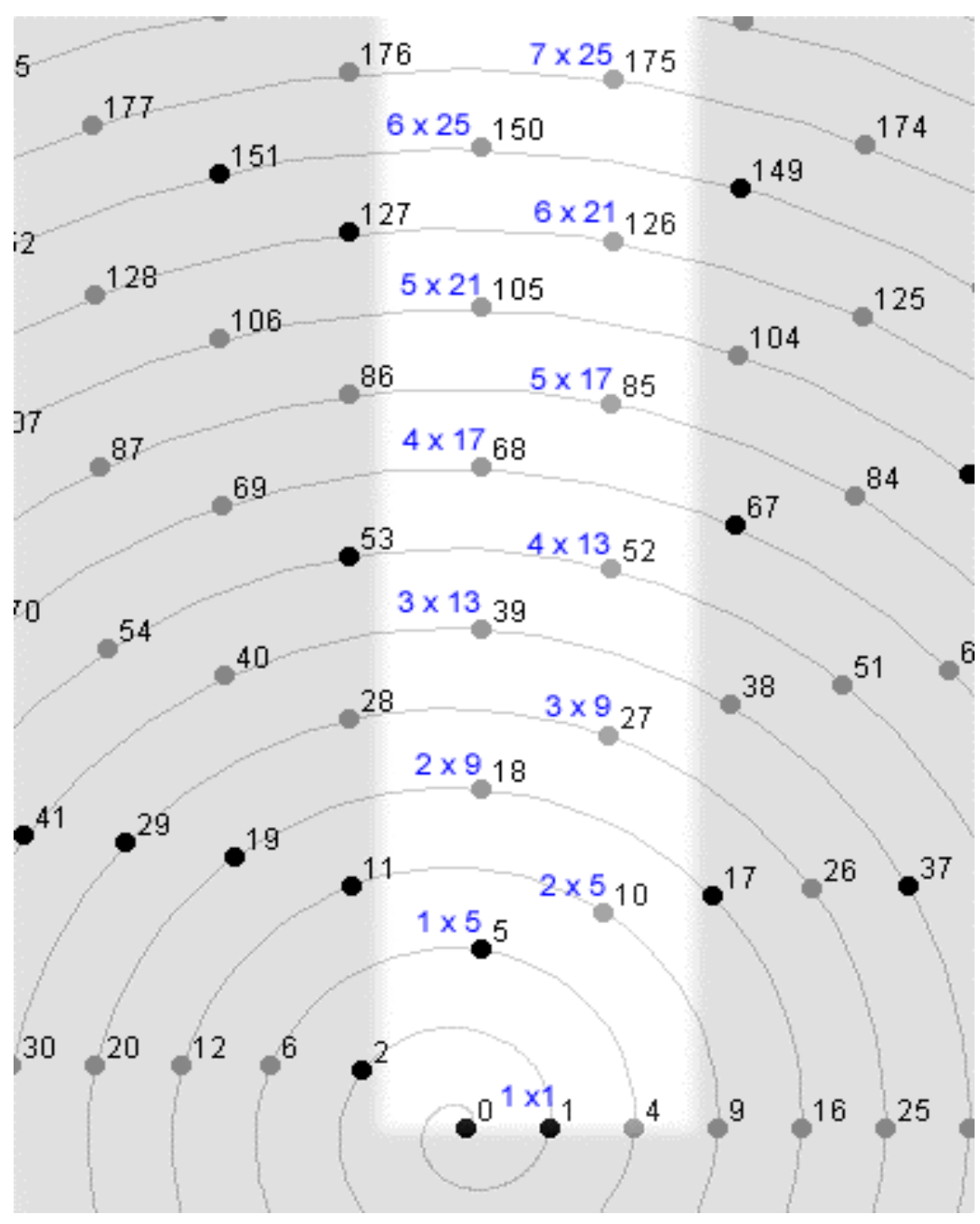

**FIG. NS-15** : A column of non-primes runs north & south

## 9.4 Quadratic Polynomials (Quadratic Functions)

Every offset curve can be generated by a quadratic polynomial of the form:

$$y = ax^2 + bx + c$$

Where $y$ is the number that appears on the graph; $a$, $b$, and $c$ are constants that define the curve; and $x$ is the index of the integer on the curve. In order for a quadratic function to be an offset curve, $a$ must be a perfect square. In order for a function to be a *composite* offset curve, $c$ must be zero. In other words, the function for a composite offset curve looks like this:

$$y = ax^2 + bx$$

Of course $a$ must be a perfect square in the above equation, because composite curves are a subset of offset curves. For composite offset curves of angle $n/d$ where the denominator is even:

$$a = (d/2)^2$$
$$b = n$$
$$\text{offset} = (n/d)^2$$

For angles with odd denominators, multiply the numerator and denominator by two and use the above formulas. The correspondence between angles and coefficients creates an **orderly pattern** on the graph. This can be seen in the illustration on the right **Figure NS-16** , which shows all the composite curves for denominator 6 with the numerator ranging from –9 to 9 inclusive:

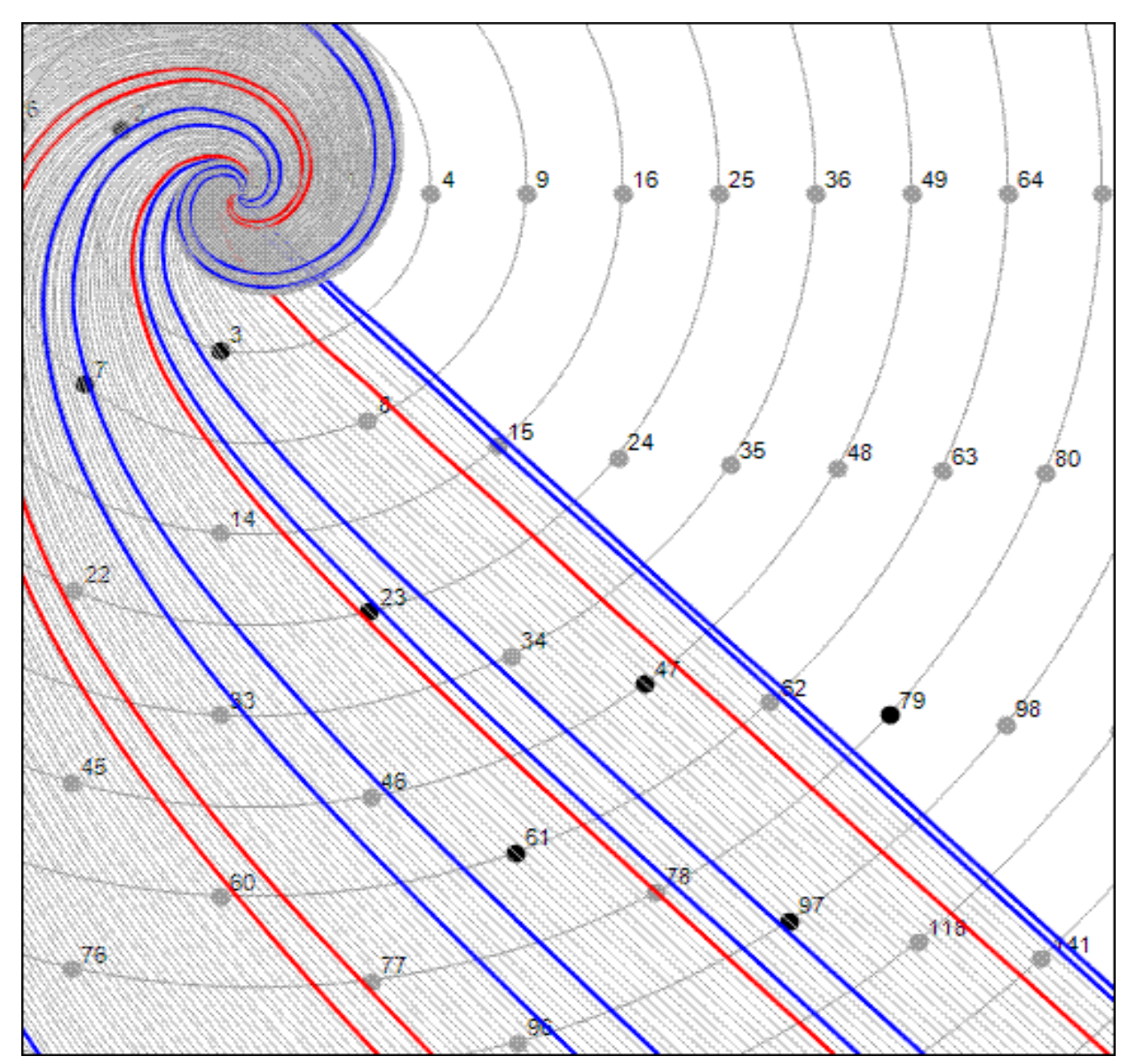

**FIG. NS-13** : Shows the Offset-Curves of angle 29/33

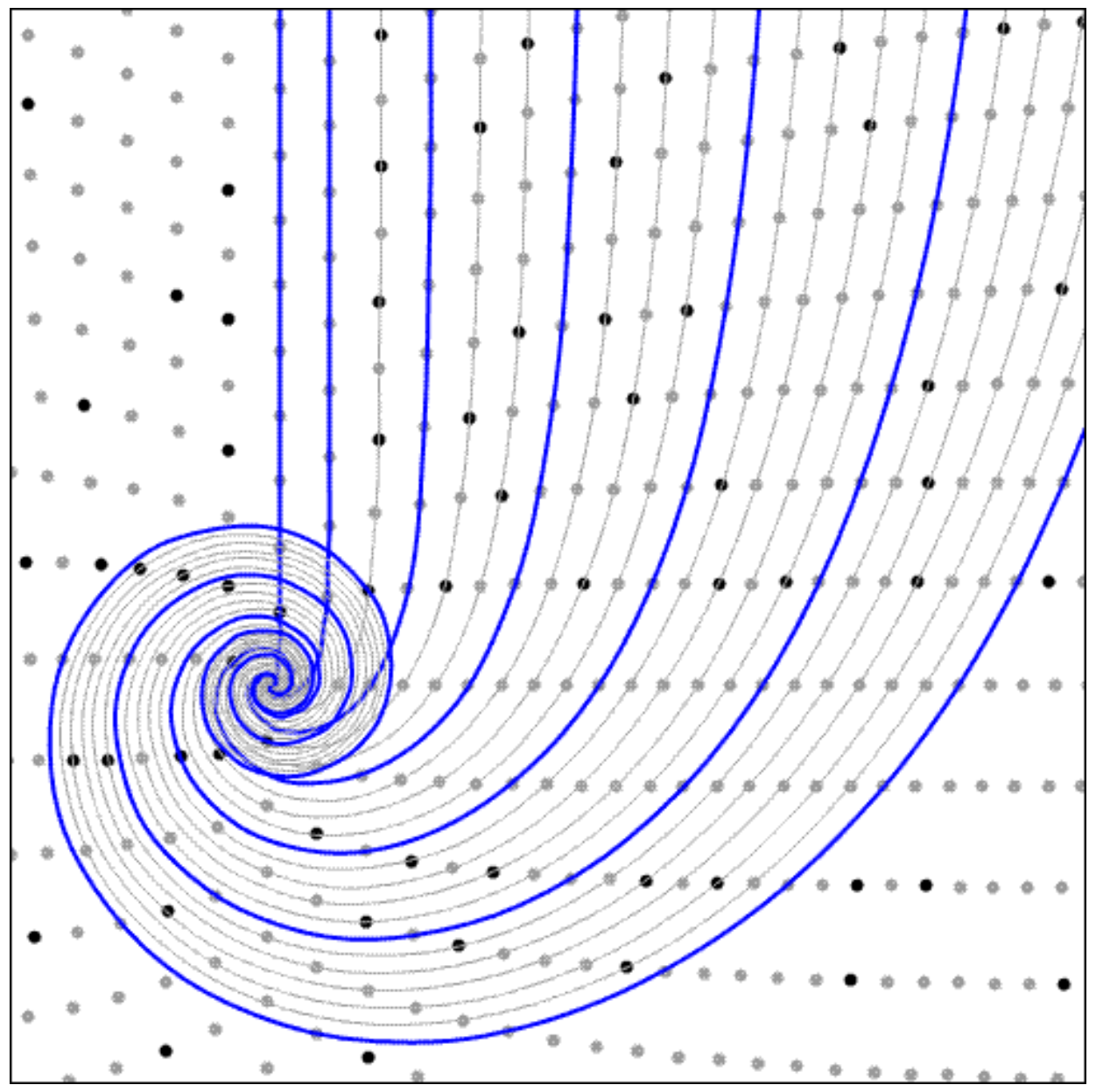

**FIG. NS-14** : Shows the Offset-Curves of angle 1/4 ( 90 deg.)

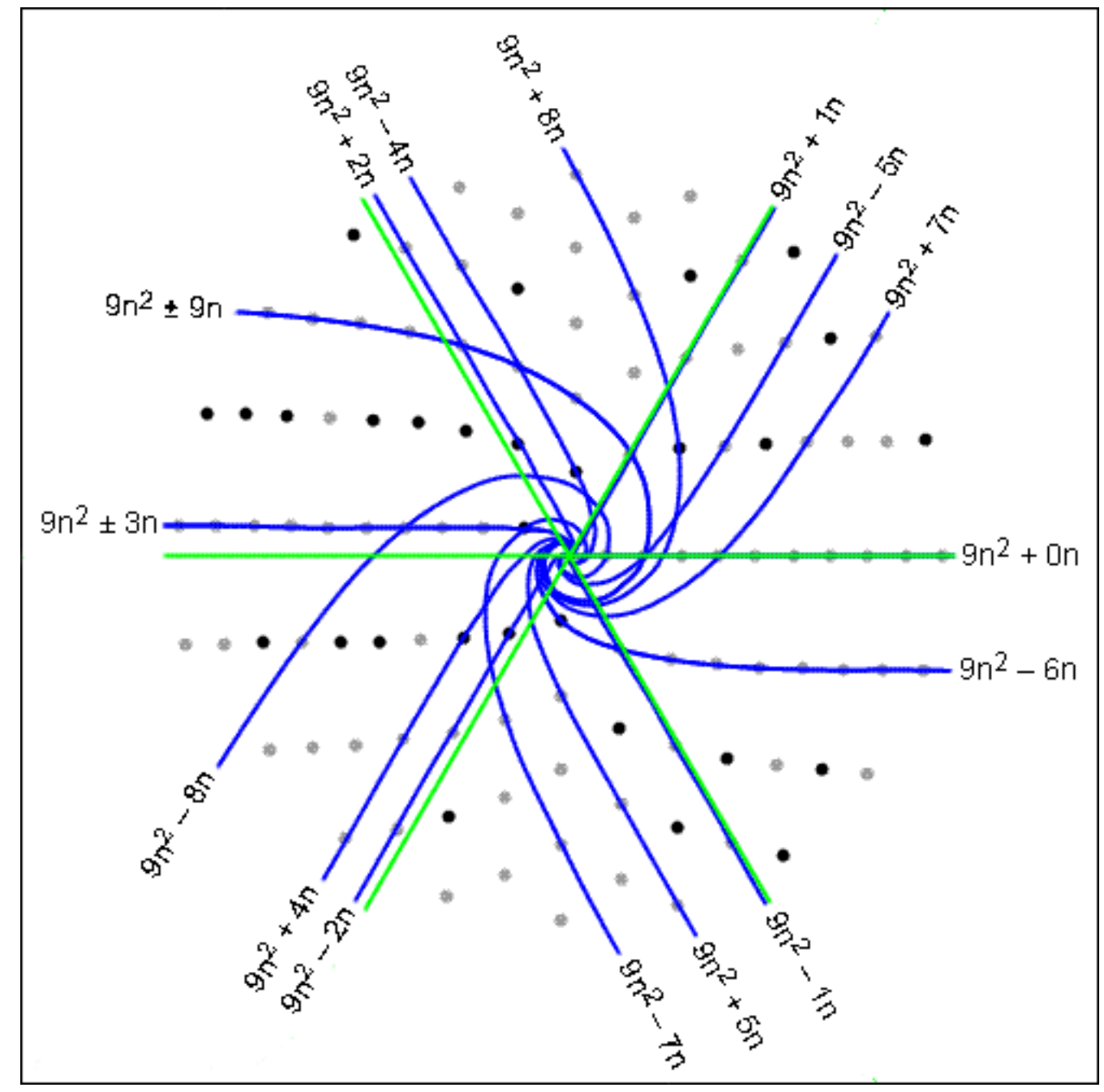

**FIG. NS-16**: offset curves are **orderly** distributed quadratic polynomials

## 9.5 Primes

Now that we have adopted a convention for naming the curves on the number spiral, let's look again at how primes are distributed.

As shown in **Figure NS-17** the densest concentrations of primes seem to occur mainly in curves whose offsets are prime.

When we look at a graph that shows only primes, like **Figure NS-17**, the left-facing curves are much more pronounced than the right-facing curves. However, when we look at a graph that shows all the integers, like **Figure NS-3** the left-facing and right-facing curves are equally salient. The main reason for this seems to be that on the left side, primes can occur only in curves with odd offsets. On the right side, primes can occur in curves with both odd and even offsets. It would be interesting to investigate empirically whether equal quantities of primes occur on both sides.

→ Interresting would also be an analysis of the circular appearances of prime numbers in this diagram !! ( → comment from Harry K. Hahn )

In FIG. NS-17 the first 46656 Prime Numbers are shown.

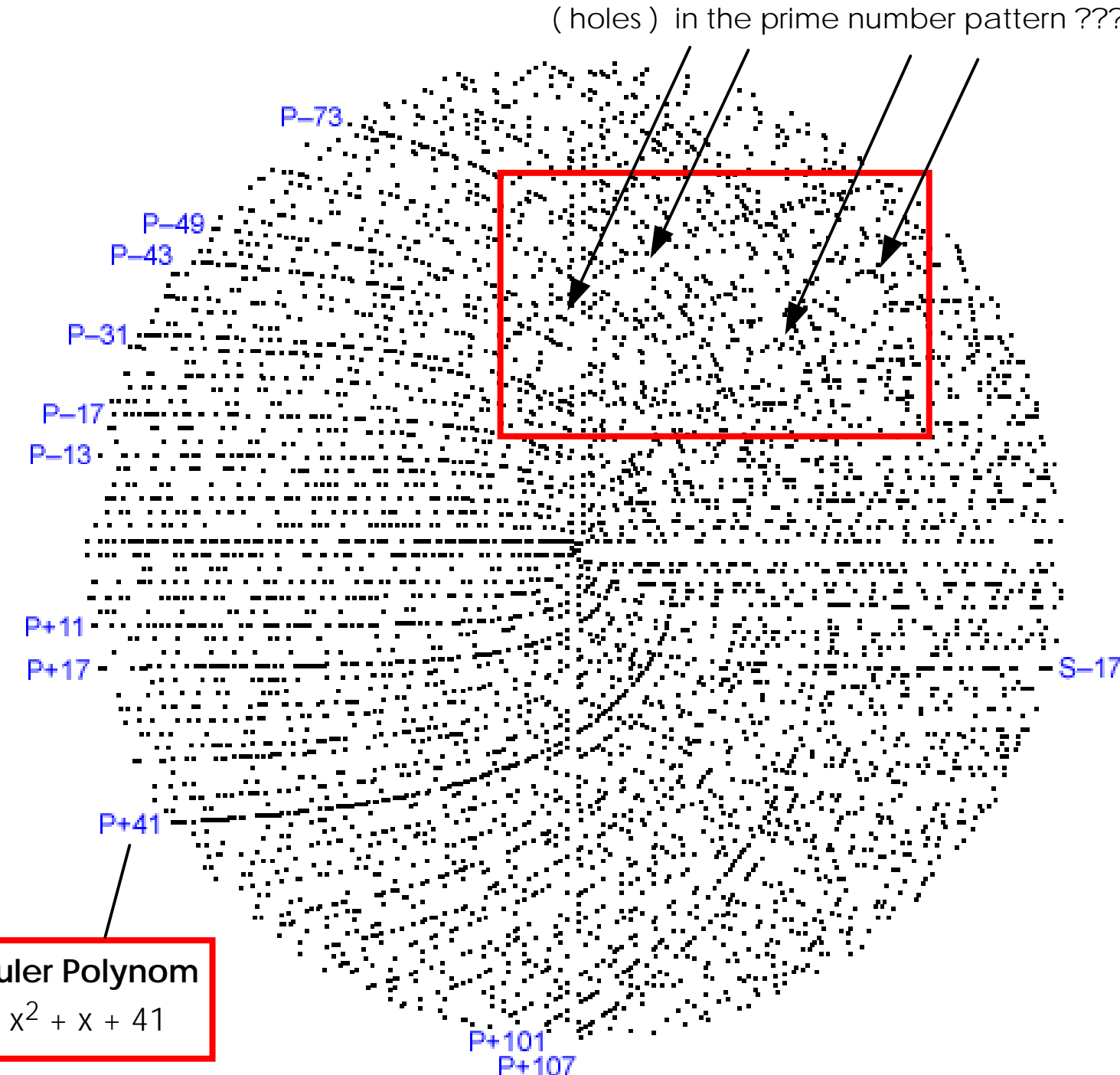

**FIG. NS-17** : Prime Numbers seem to occur mainly in curves whose offsets are prime

Here a comparison of the Number Spiral with the Ulam-Spiral → **see Figure NS-18**

Each diagonal of the Ulam spiral corresponds to a particular curve on the number spiral.

However, the Ulam spiral makes two diagonals out of each curve by allowing both ends of a diagonal to grow in opposite directions and placing alternate members of the sequence at either end. Moreover, the two halves of each diagonal do not usually line up with each other. This sounds terribly confusing in words, but as shown in **Figure NS-18**, it's really pretty simple.

I've labeled only a few diagonals in this Figure to illustrate the pattern, the correspondence extends infinitely:

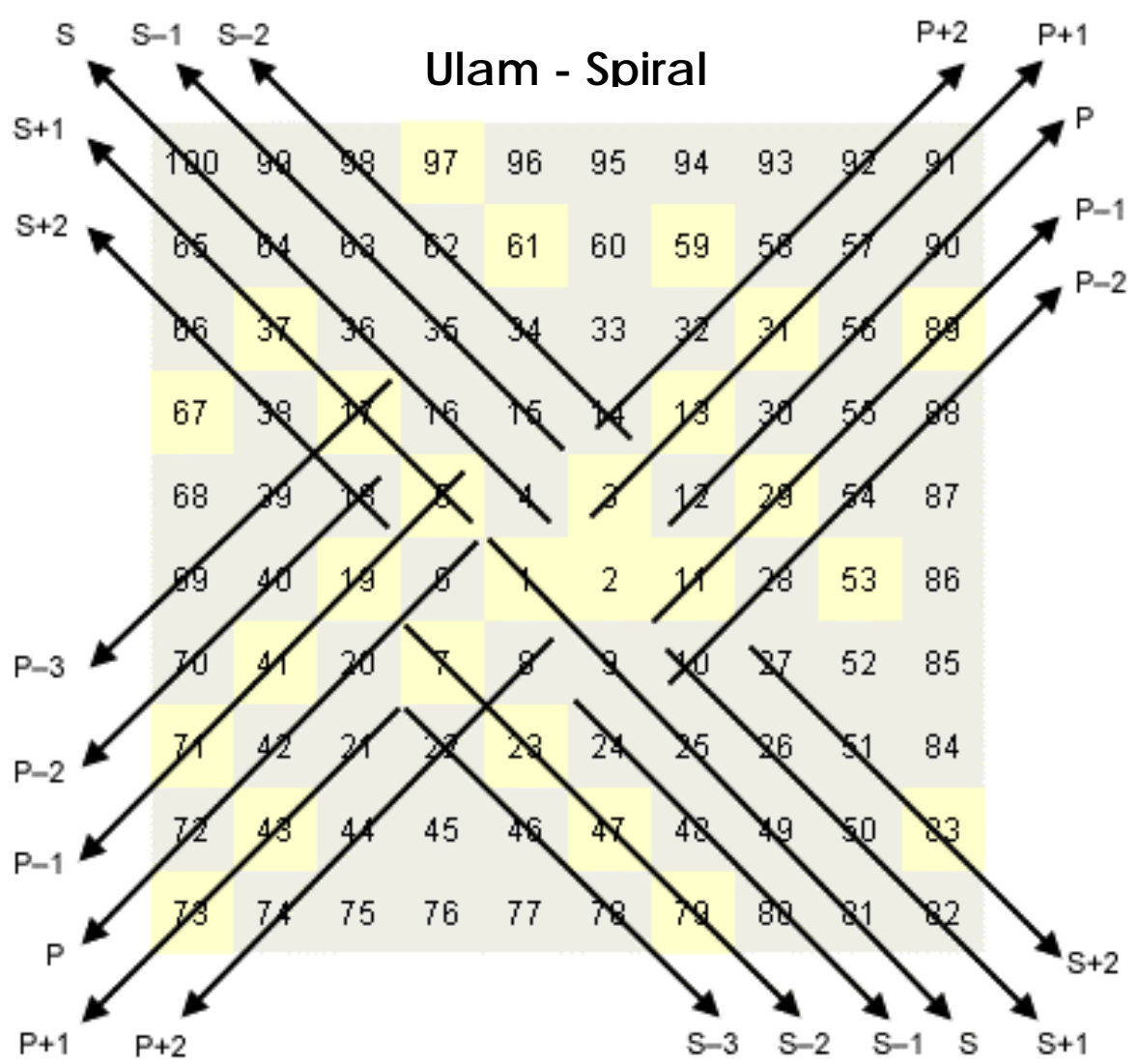

**FIG. NS-18** : Comparison between the Number Spiral and the Ulam-Spiral

---

## 9.6 FORMULAS

To convert from polar coordinates to Windows screen coordinates (in which y increases from top to bottom, unlike conventional graphs):

$$x = r\cos(2\pi\theta)$$

$$y = -r\sin(2\pi\theta)$$

Note: theta is in rotations; *x* and *y* must be scaled and translated before using them for screen display

To plot the spiral, including both the thin gray coiled line and the integers on it: :

$$r = \sqrt{n}$$

$$\theta = \sqrt{n}$$

Note: theta is in rotations

To plot offset curves:

$$r = \sqrt{an^2 + bn + c}$$

$$\theta = r - n\sqrt{a}$$

To derive coefficients *a, b, c* of a quadratic formula from three successive integers *i, j, k* in a quadratic sequence:

$$a = \frac{i - 2j + k}{2}$$

$$b = j - i - 3a$$

$$c = i - a - b$$

For more information about these formulas including their derivation, see Method of Common Differences.

To convert between a composite offset curve of angle *n/d* (measured in rotations) and its related quadratic function $y = ax^2 + bx$:

$$a = \left(\frac{d}{2}\right)^2$$

$$b = n$$

$$offset = \left(\frac{n}{d}\right)^2$$

$$d = 2\sqrt{a}$$

To factor an integer on a composite offset curve:

→ $y = x(ax + b)$

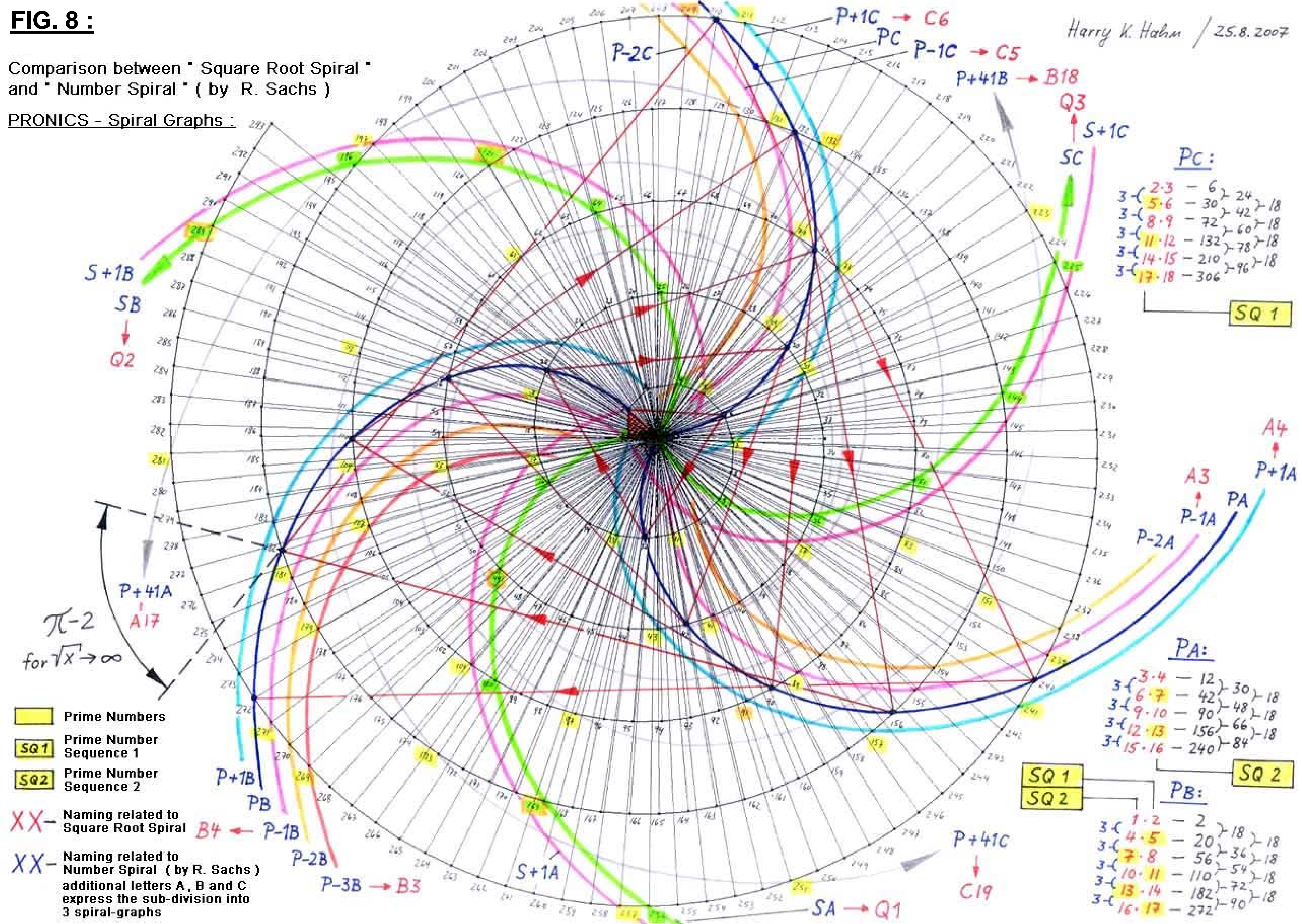

Table 7-A : Comparison of the arrangement of the "PRONICS-Spiral Graphs" ( product curves ) S , P , S+1 and P-1 on the "Number Spiral" & "Square Root Spiral" ( see FIG.: NS-8 at page 29 and FIG.: 8 at page 33 )

**NUMBER SPIRAL** → **SQUARE ROOT SPIRAL**

| factors | SD | S | S' | S" |
|---|---|---|---|---|
| 1x1 | 1 | 1 | | |
| 2x2 | 4 | 4 | 3 | |
| 3x3 | 9 | 9 | 5 | 2 |
| 4x4 | 5 | 16 | 7 | 2 |
| 5x5 | 7 | 25 | 9 | 2 |
| 6x6 | 9 | 36 | 11 | 2 |
| 7x7 | 13 | 49 | 13 | 2 |
| 8x8 | 10 | 64 | 15 | 2 |
| | 9 | 81 | 17 | 2 |
| | 1 | 100 | 19 | 2 |
| | 4 | 121 | 21 | 2 |
| | 9 | 144 | 23 | 2 |
| | 16 | 169 | 25 | 2 |
| | 16 | 196 | 27 | 2 |
| | 9 | 225 | 29 | 2 |
| | 13 | 256 | 31 | 2 |
| | 19 | 289 | 33 | 2 |
| | 9 | 324 | 35 | 2 |
| | 10 | 361 | 37 | 2 |
| | 4 | 400 | 39 | 2 |
| | 9 | 441 | 41 | 2 |
| | 16 | 484 | 43 | 2 |
| | 16 | 529 | 45 | 2 |
| | 18 | 576 | 47 | 2 |
| | 13 | 625 | 49 | 2 |

| factors | SD | SA (Q1) | SA' | SA" |
|---|---|---|---|---|
| 1x1 | 1 | 1 | | |
| 4x4 | 7 | 16 | 15 | |
| 7x7 | 13 | 49 | 33 | 18 |
| 10x10 | 1 | 100 | 51 | 18 |
| 13x13 | 16 | 169 | 69 | 18 |
| 16x16 | 13 | 256 | 87 | 18 |
| 19x19 | 10 | 361 | 105 | 18 |
| 22x22 | 16 | 484 | 123 | 18 |
| | 13 | 625 | 141 | 18 |
| | 19 | 784 | 159 | 18 |
| | 16 | 961 | 177 | 18 |
| | 13 | 1156 | 195 | 18 |
| | 19 | 1369 | 213 | 18 |
| | 7 | 1600 | 231 | 18 |
| | 22 | 1849 | 249 | 18 |
| | 10 | 2116 | 267 | 18 |
| | 7 | 2401 | 285 | 18 |
| | 13 | 2704 | 303 | 18 |
| | 10 | 3025 | 321 | 18 |
| | 16 | 3364 | 339 | 18 |
| | 13 | 3721 | 357 | 18 |
| | 19 | 4096 | 375 | 18 |
| | 25 | 4489 | 393 | 18 |
| | 13 | 4900 | 411 | 18 |
| | 19 | 5329 | 429 | 18 |

| factors | SD | SB (Q2) | SB' | SB" |
|---|---|---|---|---|
| 2x2 | 4 | 4 | | |
| 5x5 | 7 | 25 | 21 | |
| 8x8 | 10 | 64 | 39 | 18 |
| 11x11 | 4 | 121 | 57 | 18 |
| 14x14 | 16 | 196 | 75 | 18 |
| 17x17 | 19 | 289 | 93 | 18 |
| 20x20 | 4 | 400 | 111 | 18 |
| 23x23 | 16 | 529 | 129 | 18 |
| | 19 | 676 | 147 | 18 |
| | 13 | 841 | 165 | 18 |
| | 7 | 1024 | 183 | 18 |
| | 10 | 1225 | 201 | 18 |
| | 13 | 1444 | 219 | 18 |
| | 16 | 1681 | 237 | 18 |
| | 19 | 1936 | 255 | 18 |
| | 13 | 2209 | 273 | 18 |
| | 7 | 2500 | 291 | 18 |
| | 19 | 2809 | 309 | 18 |
| | 13 | 3136 | 327 | 18 |
| | 16 | 3481 | 345 | 18 |
| | 19 | 3844 | 363 | 18 |
| | 13 | 4225 | 381 | 18 |
| | 16 | 4624 | 399 | 18 |
| | 10 | 5041 | 417 | 18 |
| | 22 | 5476 | 435 | 18 |

| factors | SD | SC (Q3) | SC' | SC" |
|---|---|---|---|---|
| 3x3 | 9 | 9 | | |
| 6x6 | 9 | 36 | 27 | |
| 9x9 | 9 | 81 | 45 | 18 |
| 12x12 | 9 | 144 | 63 | 18 |
| 15x15 | 9 | 225 | 81 | 18 |
| 18x18 | 9 | 324 | 99 | 18 |
| 21x21 | 9 | 441 | 117 | 18 |
| 24x24 | 18 | 576 | 135 | 18 |
| | 18 | 729 | 153 | 18 |
| | 9 | 900 | 171 | 18 |
| | 18 | 1089 | 189 | 18 |
| | 18 | 1296 | 207 | 18 |
| | 9 | 1521 | 225 | 18 |
| | 18 | 1764 | 243 | 18 |
| | 9 | 2025 | 261 | 18 |
| | 9 | 2304 | 279 | 18 |
| | 9 | 2601 | 297 | 18 |
| | 18 | 2916 | 315 | 18 |
| | 18 | 3249 | 333 | 18 |
| | 9 | 3600 | 351 | 18 |
| | 27 | 3969 | 369 | 18 |
| | 18 | 4356 | 387 | 18 |
| | 18 | 4761 | 405 | 18 |
| | 18 | 5184 | 423 | 18 |
| | 18 | 5625 | 441 | 18 |

**NUMBER SPIRAL** → **SQUARE ROOT SPIRAL**

| SD | S+1 | S+1' | S+1" |
|---|---|---|---|
| 5 | 5 | | |
| 1 | 10 | 5 | |
| 8 | 17 | 7 | 2 |
| 8 | 26 | 9 | 2 |
| 10 | 37 | 11 | 2 |
| 5 | 50 | 13 | 2 |
| 11 | 65 | 15 | 2 |
| 10 | 82 | 17 | 2 |
| 2 | 101 | 19 | 2 |
| 5 | 122 | 21 | 2 |
| 10 | 145 | 23 | 2 |
| 8 | 170 | 25 | 2 |
| 17 | 197 | 27 | 2 |
| 10 | 226 | 29 | 2 |
| 14 | 257 | 31 | 2 |
| 11 | 290 | 33 | 2 |
| 10 | 325 | 35 | 2 |
| 11 | 362 | 37 | 2 |
| 5 | 401 | 39 | 2 |
| 10 | 442 | 41 | 2 |
| 17 | 485 | 43 | 2 |
| 8 | 530 | 45 | 2 |
| 19 | 577 | 47 | 2 |
| 14 | 626 | 49 | 2 |
| 20 | 677 | 51 | 2 |

| SD | S+1A | S+1A' | S+1A" |
|---|---|---|---|
| 2 | 2 | | |
| 8 | 17 | 15 | |
| 5 | 50 | 33 | 18 |
| 2 | 101 | 51 | 18 |
| 8 | 170 | 69 | 18 |
| 14 | 257 | 87 | 18 |
| 11 | 362 | 105 | 18 |
| 17 | 485 | 123 | 18 |
| 14 | 626 | 141 | 18 |
| 20 | 785 | 159 | 18 |
| 17 | 962 | 177 | 18 |
| 14 | 1157 | 195 | 18 |
| 11 | 1370 | 213 | 18 |
| 8 | 1601 | 231 | 18 |
| 14 | 1850 | 249 | 18 |
| 11 | 2117 | 267 | 18 |
| 8 | 2402 | 285 | 18 |
| 14 | 2705 | 303 | 18 |
| 11 | 3026 | 321 | 18 |
| 17 | 3365 | 339 | 18 |
| 14 | 3722 | 357 | 18 |
| 20 | 4097 | 375 | 18 |
| 17 | 4490 | 393 | 18 |
| 14 | 4901 | 411 | 18 |
| 11 | 5330 | 429 | 18 |
| 26 | 5777 | 447 | 18 |

| SD | S+1B | S+1B' | S+1B" |
|---|---|---|---|
| 5 | 5 | | |
| 8 | 26 | 21 | |
| 11 | 65 | 39 | 18 |
| 5 | 122 | 57 | 18 |
| 17 | 197 | 75 | 18 |
| 11 | 290 | 93 | 18 |
| 5 | 401 | 111 | 18 |
| 8 | 530 | 129 | 18 |
| 20 | 677 | 147 | 18 |
| 14 | 842 | 165 | 18 |
| 8 | 1025 | 183 | 18 |
| 11 | 1226 | 201 | 18 |
| 14 | 1445 | 219 | 18 |
| 17 | 1682 | 237 | 18 |
| 20 | 1937 | 255 | 18 |
| 5 | 2210 | 273 | 18 |
| 8 | 2501 | 291 | 18 |
| 11 | 2810 | 309 | 18 |
| 14 | 3137 | 327 | 18 |
| 17 | 3482 | 345 | 18 |
| 20 | 3845 | 363 | 18 |
| 14 | 4226 | 381 | 18 |
| 17 | 4625 | 399 | 18 |
| 11 | 5042 | 417 | 18 |
| 23 | 5477 | 435 | 18 |

| SD | S+1C | S+1C' | S+1C" |
|---|---|---|---|
| 1 | 1 | | |
| 1 | 10 | 9 | |
| 10 | 37 | 27 | 18 |
| 10 | 82 | 45 | 18 |
| 10 | 145 | 63 | 18 |
| 10 | 226 | 81 | 18 |
| 10 | 325 | 99 | 18 |
| 10 | 442 | 117 | 18 |
| 19 | 577 | 135 | 18 |
| 10 | 730 | 153 | 18 |
| 10 | 901 | 171 | 18 |
| 10 | 1090 | 189 | 18 |
| 19 | 1297 | 207 | 18 |
| 10 | 1522 | 225 | 18 |
| 19 | 1765 | 243 | 18 |
| 10 | 2026 | 261 | 18 |
| 10 | 2305 | 279 | 18 |
| 10 | 2602 | 297 | 18 |
| 19 | 2917 | 315 | 18 |
| 10 | 3250 | 333 | 18 |
| 10 | 3601 | 351 | 18 |
| 19 | 3970 | 369 | 18 |
| 19 | 4357 | 387 | 18 |
| 19 | 4762 | 405 | 18 |
| 19 | 5185 | 423 | 18 |

| factors | SD | P | P' | P" |
|---|---|---|---|---|
| 1x2 | 2 | 2 | | |
| 2x3 | 6 | 6 | 4 | |
| 3x4 | 3 | 12 | 6 | 2 |
| 4x5 | 2 | 20 | 8 | 2 |
| 5x6 | 3 | 30 | 10 | 2 |
| 6x7 | 6 | 42 | 12 | 2 |
| 7x8 | 11 | 56 | 14 | 2 |
| 8x9 | 9 | 72 | 16 | 2 |
| 9x10 | 9 | 90 | 18 | 2 |
| 10x11 | 2 | 110 | 20 | 2 |
| 11x12 | 6 | 132 | 22 | 2 |
| | 12 | 156 | 24 | 2 |
| | 11 | 182 | 26 | 2 |
| | 3 | 210 | 28 | 2 |
| | 6 | 240 | 30 | 2 |
| | 11 | 272 | 32 | 2 |
| | 9 | 306 | 34 | 2 |
| | 9 | 342 | 36 | 2 |
| | 11 | 380 | 38 | 2 |
| | 6 | 420 | 40 | 2 |
| | 12 | 462 | 42 | 2 |
| | 11 | 506 | 44 | 2 |
| | 12 | 552 | 46 | 2 |
| | 6 | 600 | 48 | 2 |
| | 11 | 650 | 50 | 2 |

| factors | SD | PA | PA' | PA" |
|---|---|---|---|---|
| 0x1 | 0 | 0 | | |
| 3x4 | 3 | 12 | 12 | |
| 6x7 | 6 | 42 | 30 | 18 |
| 9x10 | 9 | 90 | 48 | 18 |
| 12x13 | 12 | 156 | 66 | 18 |
| 15x16 | 6 | 240 | 84 | 18 |
| 18x19 | 9 | 342 | 102 | 18 |
| 21x22 | 12 | 462 | 120 | 18 |
| 24x25 | 6 | 600 | 138 | 18 |
| 27x28 | 18 | 756 | 156 | 18 |
| 30x31 | 12 | 930 | 174 | 18 |
| | 6 | 1122 | 192 | 18 |
| | 9 | 1332 | 210 | 18 |
| | 12 | 1560 | 228 | 18 |
| | 15 | 1806 | 246 | 18 |
| | 9 | 2070 | 264 | 18 |
| | 12 | 2352 | 282 | 18 |
| | 15 | 2652 | 300 | 18 |
| | 18 | 2970 | 318 | 18 |
| | 12 | 3306 | 336 | 18 |
| | 15 | 3660 | 354 | 18 |
| | 9 | 4032 | 372 | 18 |
| | 12 | 4422 | 390 | 18 |
| | 15 | 4830 | 408 | 18 |
| | 18 | 5256 | 426 | 18 |

| factors | SD | PB | PB' | PB" |
|---|---|---|---|---|
| 1x2 | 2 | 2 | | |
| 4x5 | 2 | 20 | 18 | |
| 7x8 | 11 | 56 | 36 | 18 |
| 10x11 | 2 | 110 | 54 | 18 |
| 13x14 | 11 | 182 | 72 | 18 |
| 16x17 | 11 | 272 | 90 | 18 |
| 19x20 | 11 | 380 | 108 | 18 |
| 22x23 | 11 | 506 | 126 | 18 |
| 25x26 | 11 | 650 | 144 | 18 |
| 28x29 | 11 | 812 | 162 | 18 |
| 31x32 | 20 | 992 | 180 | 18 |
| | 11 | 1190 | 198 | 18 |
| | 11 | 1406 | 216 | 18 |
| | 11 | 1640 | 234 | 18 |
| | 20 | 1892 | 252 | 18 |
| | 11 | 2162 | 270 | 18 |
| | 11 | 2450 | 288 | 18 |
| | 20 | 2756 | 306 | 18 |
| | 11 | 3080 | 324 | 18 |
| | 11 | 3422 | 342 | 18 |
| | 20 | 3782 | 360 | 18 |
| | 11 | 4160 | 378 | 18 |
| | 20 | 4556 | 396 | 18 |
| | 20 | 4970 | 414 | 18 |
| | 11 | 5402 | 432 | 18 |

| factors | SD | PC | PC' | PC" |
|---|---|---|---|---|
| 2x3 | 6 | 6 | | |
| 5x6 | 3 | 30 | 24 | |
| 8x9 | 9 | 72 | 42 | 18 |
| 11x12 | 6 | 132 | 60 | 18 |
| 14x15 | 3 | 210 | 78 | 18 |
| 17x18 | 9 | 306 | 96 | 18 |
| 20x21 | 6 | 420 | 114 | 18 |
| 23x24 | 12 | 552 | 132 | 18 |
| 26x27 | 9 | 702 | 150 | 18 |
| 29x30 | 15 | 870 | 168 | 18 |
| 32x33 | 12 | 1056 | 186 | 18 |
| | 9 | 1260 | 204 | 18 |
| | 15 | 1482 | 222 | 18 |
| | 12 | 1722 | 240 | 18 |
| | 18 | 1980 | 258 | 18 |
| | 15 | 2256 | 276 | 18 |
| | 12 | 2550 | 294 | 18 |
| | 18 | 2862 | 312 | 18 |
| | 15 | 3192 | 330 | 18 |
| | 12 | 3540 | 348 | 18 |
| | 18 | 3906 | 366 | 18 |
| | 15 | 4290 | 384 | 18 |
| | 21 | 4692 | 402 | 18 |
| | 9 | 5112 | 420 | 18 |
| | 15 | 5550 | 438 | 18 |

| SD | P-1 | P-1' | P-1" |
|---|---|---|---|
| 5 | 5 | | |
| 2 | 11 | 6 | |
| 10 | 19 | 8 | 2 |
| 11 | 29 | 10 | 2 |
| 5 | 41 | 12 | 2 |
| 10 | 55 | 14 | 2 |
| 8 | 71 | 16 | 2 |
| 17 | 89 | 18 | 2 |
| 10 | 109 | 20 | 2 |
| 5 | 131 | 22 | 2 |
| 11 | 155 | 24 | 2 |
| 10 | 181 | 26 | 2 |
| 11 | 209 | 28 | 2 |
| 14 | 239 | 30 | 2 |
| 10 | 271 | 32 | 2 |
| 8 | 305 | 34 | 2 |
| 8 | 341 | 36 | 2 |
| 19 | 379 | 38 | 2 |
| 14 | 419 | 40 | 2 |
| 11 | 461 | 42 | 2 |
| 10 | 505 | 44 | 2 |
| 11 | 551 | 46 | 2 |
| 23 | 599 | 48 | 2 |
| 19 | 649 | 50 | 2 |
| 8 | 701 | 52 | 2 |

| SD | P-1A (A3) | P-1A' | P-1A" |
|---|---|---|---|
| 2 | 11 | | |
| 5 | 41 | 30 | |
| 17 | 89 | 48 | 18 |
| 11 | 155 | 66 | 18 |
| 14 | 239 | 84 | 18 |
| 8 | 341 | 102 | 18 |
| 11 | 461 | 120 | 18 |
| 23 | 599 | 138 | 18 |
| 17 | 755 | 156 | 18 |
| 20 | 929 | 174 | 18 |
| 5 | 1121 | 192 | 18 |
| 8 | 1331 | 210 | 18 |
| 20 | 1559 | 228 | 18 |
| 14 | 1805 | 246 | 18 |
| 17 | 2069 | 264 | 18 |
| 11 | 2351 | 282 | 18 |
| 14 | 2651 | 300 | 18 |
| 26 | 2969 | 318 | 18 |
| 11 | 3305 | 336 | 18 |
| 23 | 3659 | 354 | 18 |
| 8 | 4031 | 372 | 18 |
| 11 | 4421 | 390 | 18 |
| 23 | 4829 | 408 | 18 |
| 17 | 5255 | 426 | 18 |
| 29 | 5699 | 444 | 18 |

| SD | P-1B (B4) | P-1B' | P-1B" |
|---|---|---|---|
| 10 | 19 | | |
| 10 | 55 | 36 | |
| 10 | 109 | 54 | 18 |
| 10 | 181 | 72 | 18 |
| 10 | 271 | 90 | 18 |
| 19 | 379 | 108 | 18 |
| 10 | 505 | 126 | 18 |
| 19 | 649 | 144 | 18 |
| 10 | 811 | 162 | 18 |
| 19 | 991 | 180 | 18 |
| 19 | 1189 | 198 | 18 |
| 10 | 1405 | 216 | 18 |
| 19 | 1639 | 234 | 18 |
| 19 | 1891 | 252 | 18 |
| 10 | 2161 | 270 | 18 |
| 19 | 2449 | 288 | 18 |
| 19 | 2755 | 306 | 18 |
| 19 | 3079 | 324 | 18 |
| 10 | 3421 | 342 | 18 |
| 19 | 3781 | 360 | 18 |
| 19 | 4159 | 378 | 18 |
| 19 | 4555 | 396 | 18 |
| 28 | 4969 | 414 | 18 |
| 10 | 5401 | 432 | 18 |
| 19 | 5851 | 450 | 18 |

| SD | P-1C (C5) | P-1C' | P-1C" |
|---|---|---|---|
| 5 | 5 | | |
| 11 | 29 | 24 | |
| 8 | 71 | 42 | 18 |
| 5 | 131 | 60 | 18 |
| 11 | 209 | 78 | 18 |
| 8 | 305 | 96 | 18 |
| 14 | 419 | 114 | 18 |
| 11 | 551 | 132 | 18 |
| 8 | 701 | 150 | 18 |
| 23 | 869 | 168 | 18 |
| 11 | 1055 | 186 | 18 |
| 17 | 1259 | 204 | 18 |
| 14 | 1481 | 222 | 18 |
| 11 | 1721 | 240 | 18 |
| 26 | 1979 | 258 | 18 |
| 14 | 2255 | 276 | 18 |
| 20 | 2549 | 294 | 18 |
| 17 | 2861 | 312 | 18 |
| 14 | 3191 | 330 | 18 |
| 20 | 3539 | 348 | 18 |
| 17 | 3905 | 366 | 18 |
| 23 | 4289 | 384 | 18 |
| 20 | 4691 | 402 | 18 |
| 8 | 5111 | 420 | 18 |
| 23 | 5549 | 438 | 18 |

XX Prime Number Sequence **SQ1**

XX Prime Number Sequence **SQ2**

Prime Numbers

Numbers not divisible by 2, 3 or 5

Table 7-B : Prime Number Sequences derived from Spiral Graphs **P+41** ( Euler Polynom → see FIG.: NS-17 at page 32 ) and **P+41A** , **P+41B** and **P+41C** ( see FIG.: 8 at page 33 )

## NUMBER SPIRAL

Euler Polynom → P+41

| SD | P+41 | P+41' | P+41'' |
|---|---|---|---|
| 5 | 41 | | |
| 7 | 43 | 2 | |
| 11 | 47 | 4 | 2 |
| 8 | 53 | 6 | 2 |
| 7 | 61 | 8 | 2 |
| 8 | 71 | 10 | 2 |
| 11 | 83 | 12 | 2 |
| 16 | 97 | 14 | 2 |
| 5 | 113 | 16 | 2 |
| 5 | 131 | 18 | 2 |
| 7 | 151 | 20 | 2 |
| 11 | 173 | 22 | 2 |
| 17 | 197 | 24 | 2 |
| 7 | 223 | 26 | 2 |
| 8 | 251 | 28 | 2 |
| 11 | 281 | 30 | 2 |
| 7 | 313 | 32 | 2 |
| 14 | 347 | 34 | 2 |
| 14 | 383 | 36 | 2 |
| 7 | 421 | 38 | 2 |
| 11 | 461 | 40 | 2 |
| 8 | 503 | 42 | 2 |
| 16 | 547 | 44 | 2 |
| 17 | 593 | 46 | 2 |
| 11 | 641 | 48 | 2 |
| 16 | 691 | 50 | 2 |
| 14 | 743 | 52 | 2 |
| 23 | 797 | 54 | 2 |
| 16 | 853 | 56 | 2 |
| 11 | 911 | 58 | 2 |
| 17 | 971 | 60 | 2 |
| 7 | 1033 | 62 | 2 |
| 17 | 1097 | 64 | 2 |
| 11 | 1163 | 66 | 2 |
| 7 | 1231 | 68 | 2 |
| 5 | 1301 | 70 | 2 |
| 14 | 1373 | 72 | 2 |
| 16 | 1447 | 74 | 2 |
| 11 | 1523 | 76 | 2 |
| 8 | 1601 | 78 | 2 |
| 16 | 1681 | 80 | 2 |
| 17 | 1763 | 82 | 2 |
| 20 | 1847 | 84 | 2 |

## SQUARE ROOT SPIRAL

### P+41A

| SD | prime factors | P+41A | P+41A' | P+41A'' |
|---|---|---|---|---|
| 5 | | 41 | | |
| 8 | | 53 | 12 | |
| 11 | | 83 | 30 | 18 |
| 5 | | 131 | 48 | 18 |
| 17 | | 197 | 66 | 18 |
| 11 | | 281 | 84 | 18 |
| 14 | | 383 | 102 | 18 |
| 8 | | 503 | 120 | 18 |
| 11 | | 641 | 138 | 18 |
| 23 | | 797 | 156 | 18 |
| 17 | | 971 | 174 | 18 |
| 11 | | 1163 | 192 | 18 |
| 14 | | 1373 | 210 | 18 |
| 8 | | 1601 | 228 | 18 |
| 20 | | 1847 | 246 | 18 |
| 5 | | 2111 | 264 | 18 |
| 17 | | 2393 | 282 | 18 |
| 20 | | 2693 | 300 | 18 |
| 5 | | 3011 | 318 | 18 |
| 17 | | 3347 | 336 | 18 |
| 11 | | 3701 | 354 | 18 |
| 14 | | 4073 | 372 | 18 |
| 17 | | 4463 | 390 | 18 |
| 20 | | 4871 | 408 | 18 |
| 23 | | 5297 | 426 | 18 |
| 17 | | 5741 | 444 | 18 |
| 11 | | 6203 | 462 | 18 |
| 23 | 41 x 163 | 6683 | 480 | 18 |
| 17 | 43 x 167 | 7181 | 498 | 18 |
| 29 | 43 x 179 | 7697 | 516 | 18 |
| 14 | | 8231 | 534 | 18 |
| 26 | | 8783 | 552 | 18 |
| 20 | 47 x 199 | 9353 | 570 | 18 |
| 23 | | 9941 | 588 | 18 |
| 17 | 53 x 199 | 10547 | 606 | 18 |
| 11 | | 11171 | 624 | 18 |
| 14 | | 11813 | 642 | 18 |
| 17 | | 12473 | 660 | 18 |
| 11 | | 13151 | 678 | 18 |
| 23 | 61 x 227 | 13847 | 696 | 18 |
| 17 | | 14561 | 714 | 18 |
| 20 | 41 x 373 | 15293 | 732 | 18 |

### P+41B

| SD | prime factors | P+41B | P+41B' | P+41B'' |
|---|---|---|---|---|
| 7 | | 43 | | |
| 7 | | 61 | 18 | |
| 16 | | 97 | 36 | 18 |
| 7 | | 151 | 54 | 18 |
| 7 | | 223 | 72 | 18 |
| 7 | | 313 | 90 | 18 |
| 7 | | 421 | 108 | 18 |
| 16 | | 547 | 126 | 18 |
| 16 | | 691 | 144 | 18 |
| 16 | | 853 | 162 | 18 |
| 7 | | 1033 | 180 | 18 |
| 7 | | 1231 | 198 | 18 |
| 16 | | 1447 | 216 | 18 |
| 16 | 41 x 41 | 1681 | 234 | 18 |
| 16 | | 1933 | 252 | 18 |
| 7 | | 2203 | 270 | 18 |
| 16 | 47 x 53 | 2491 | 288 | 18 |
| 25 | | 2797 | 306 | 18 |
| 7 | | 3121 | 324 | 18 |
| 16 | | 3463 | 342 | 18 |
| 16 | | 3823 | 360 | 18 |
| 7 | | 4201 | 378 | 18 |
| 25 | | 4597 | 396 | 18 |
| 7 | | 5011 | 414 | 18 |
| 16 | | 5443 | 432 | 18 |
| 25 | 71 x 83 | 5893 | 450 | 18 |
| 16 | | 6361 | 468 | 18 |
| 25 | 41 x 167 | 6847 | 486 | 18 |
| 16 | | 7351 | 504 | 18 |
| 25 | | 7873 | 522 | 18 |
| 16 | 47 x 179 | 8413 | 540 | 18 |
| 25 | | 8971 | 558 | 18 |
| 25 | | 9547 | 576 | 18 |
| 7 | | 10141 | 594 | 18 |
| 16 | | 10753 | 612 | 18 |
| 16 | | 11383 | 630 | 18 |
| 7 | 53 x 227 | 12031 | 648 | 18 |
| 25 | | 12697 | 666 | 18 |
| 16 | | 13381 | 684 | 18 |
| 16 | | 14083 | 702 | 18 |
| 16 | 113 x 131 | 14803 | 720 | 18 |
| 16 | | 15541 | 738 | 18 |

### P+41C

| SD | prime factors | P+41C | P+41C' | P+41C'' |
|---|---|---|---|---|
| 5 | | 41 | | |
| 11 | | 47 | 6 | |
| 8 | | 71 | 24 | 18 |
| 5 | | 113 | 42 | 18 |
| 11 | | 173 | 60 | 18 |
| 8 | | 251 | 78 | 18 |
| 14 | | 347 | 96 | 18 |
| 11 | | 461 | 114 | 18 |
| 17 | | 593 | 132 | 18 |
| 14 | | 743 | 150 | 18 |
| 11 | | 911 | 168 | 18 |
| 17 | | 1097 | 186 | 18 |
| 5 | | 1301 | 204 | 18 |
| 11 | | 1523 | 222 | 18 |
| 17 | 41 x 43 | 1763 | 240 | 18 |
| 5 | 43 x 47 | 2021 | 258 | 18 |
| 20 | | 2297 | 276 | 18 |
| 17 | | 2591 | 294 | 18 |
| 14 | | 2903 | 312 | 18 |
| 11 | 53 x 61 | 3233 | 330 | 18 |
| 17 | | 3581 | 348 | 18 |
| 23 | | 3947 | 366 | 18 |
| 11 | 61 x 71 | 4331 | 384 | 18 |
| 17 | | 4733 | 402 | 18 |
| 14 | | 5153 | 420 | 18 |
| 20 | | 5591 | 438 | 18 |
| 17 | | 6047 | 456 | 18 |
| 14 | | 6521 | 474 | 18 |
| 11 | | 7013 | 492 | 18 |
| 17 | | 7523 | 510 | 18 |
| 14 | 83 x 97 | 8051 | 528 | 18 |
| 29 | | 8597 | 546 | 18 |
| 17 | | 9161 | 564 | 18 |
| 23 | | 9743 | 582 | 18 |
| 11 | | 10343 | 600 | 18 |
| 17 | 97 x 113 | 10961 | 618 | 18 |
| 23 | | 11597 | 636 | 18 |
| 11 | | 12251 | 654 | 18 |
| 17 | | 12923 | 672 | 18 |
| 14 | | 13613 | 690 | 18 |
| 11 | | 14321 | 708 | 18 |
| 17 | 41 x 367 | 15047 | 726 | 18 |

Prime Numbers ; Numbers not divisible by 2, 3 or 5

Table 7-C : Quadratic Polynomials of the "Prime-Number-Spiral-Graphs" **P+41** ( Number Spiral ) and **P+41A** to **P+41C** {( Square Root Spiral ) --> 2. Differential = 18 }

| | Spiral Graph | Number Sequence of Spiral Graph | Quadratic Polynomial 1 ( calculated with the first 3 numbers of the given sequence ) | Quadratic Polynomial 2 ( calculated with 3 numbers starting with the **2.** Number of the sequence ) | Quadratic Polynomial 3 ( calculated with 3 numbers starting with the **3.** Number of the sequence ) | Quadratic Polynomial 4 ( calculated with 3 numbers starting with the **4.** Number of the sequence ) |
|---|---|---|---|---|---|---|
| Number Spiral | P+41 | 41 , 43 , 47 , 53 , 61 , 71 , ..... | $f_1(x) = 1x^2 - 1x + 41$ | $f_2(x) = 1x^2 + 1x + 41$ | $f_3(x) = 1x^2 + 3x + 43$ | $f_4(x) = 1x^2 + 5x + 47$ |
| Square Root Spiral | P+41A | 41 , 53 , 83 , 131 , 197 , 281 , ..... | $f_1(x) = 9x^2 - 15x + 47$ | $f_2(x) = 9x^2 + 3x + 41$ | $f_3(x) = 9x^2 + 21x + 53$ | $f_4(x) = 9x^2 + 39x + 83$ |
| | P+41B | 43 , 61 , 97 , 151 , 223 , 313 , ..... | $f_1(x) = 9x^2 - 9x + 43$ | $f_2(x) = 9x^2 + 9x + 43$ | $f_3(x) = 9x^2 + 27x + 61$ | $f_4(x) = 9x^2 + 45x + 97$ |
| | P+41C | 41 , 47 , 71 , 113 , 173 , 251 , ..... | $f_1(x) = 9x^2 - 21x + 53$ | $f_2(x) = 9x^2 - 3x + 41$ | $f_3(x) = 9x^2 + 15x + 47$ | $f_4(x) = 9x^2 + 33x + 71$ |

## 10 Comparison of the Ulam Spiral, Number Spiral and Square Root Spiral

Especially interesting should be a direct **comparison of** the **Ulam-Spiral** the **Number Spiral** and the **Square Root Spiral** in regards to the distribution of certain number-groups e.g. Square-Numbers , Pronics , Prime-Numbers etc.

For that purpose I have produced another diagram of the Square Root Spiral which shows some " reference graphs " which will make this comparison easier.

**→ see FIG. 8 - "Comparison between Square Root Spiral & Number Spiral "**

Besides **FIG. 8** ( on page **33** ), the following two images from Mr. Sach's analysis (→ chapter **9** ) should be used for the above mentioned comparison :

Image **NS-8** **: Number Spiral** **-** ( Pronics Curves ) - page **29**

Image **NS-18** **: Ulam Spiral** **-** ( Comparison with Number Spiral ) - page **32**

→ These 3 images : **FIG 8 ,** **NS-8** and **NS-18** of the Square Root Spiral , Number Spiral and Ulam Spiral , show the arrangement of the same number sequences or "Reference Graphs"

The main reference graphs are the ones which contain the square numbers.
These are either named **S** on the Ulam- and Number-Spiral ( FIG. NS-18 / NS-8 ), or **SA, SB, SC** on the Square Root Spiral ( drawn in green color in FIG 8 ).
( Note : the original naming of graphs SA, SB, SC on the Square Root Spiral is actually Q1, Q2, Q3 → see FIG.1 )

The difference in the distribution of the square numbers is as follows :

**Number Spiral** : Square Numbers are located on **1** straight graph
**Ulam Spiral** : Square Numbers are located on **2** straight graphs
**Square Root Spiral** : Square Numbers are located on **3** spiral graphs

The difference is caused through the fact, that on the Number Spiral every wind only contains one square number , whereas the winds of the Ulam Spiral and Square Root Spiral contain either 2 or 3 successive Square Numbers of the Square Number Sequence per wind.

Or put in other words, the Number Spiral has the tightest spiral winding, followed by the Ulam Spiral, which has a wider winding and then followed by the Square Root Spiral which has the widest winding of all three spirals.
There are of course further spiral variants possible with 4, 5, 6 and more successive square numbers per wind of the spiral, which then would have 4, 5, 6 or more ( spiral ) graphs on which the square numbers would be located on.

It would definitely be interesting to analyse all these spiral variants in regards to the distribution of certain number groups like pronics or prime numbers etc. !!

Beside the graphs SA, SB, SC which contain the square numbers in FIG 8 I also drew the graphs **S+1A, S+1B, S+1C** into FIG 8 ( → drawn in pink ).
These graphs are the next parallel graphs to the graphs SA, SB and SC in the negative rotation direction, and the numbers contained in these graphs are the square numbers +1. The same numbers are also contained in the two graphs marked with **S+1** in the Ulam Spiral ( see image NS-18 ) and in the graph **S+1** on the Number Spiral. <u>Note</u> : this graph is not marked in image NS-8 , however it lies on the opposite side of graph S-1 in reference to graph S.

The main graphs of the second group of reference graphs are named with the letter **P** on the Ulam Spiral and Number Spiral ( see FIG. NS-18 and / NS-8 ) and with the letters **PA, PB, PC** on the Square Root Spiral ( graphs drawn in dark blue in FIG 8 ). → <u>Note</u> : letter **P** stays for " pronics"

The distribution of "Pronics-Numbers" ( of the same type ! ) is similar to the distribution of the Square Numbers. It is as follows :

**Number Spiral** : Pronics Numbers are located on **1** spiral graph
**Ulam Spiral** : Pronics Numbers are located on **2** straight graphs
**Square Root Spiral** : Pronics Numbers are located on **3** spiral graphs

The difference in the number of graphs, on which pronics-numbers of the same type are located, has the same cause as already described for the square numbers , which is the difference in the tightness of the winding of the Ulam-, Number- and Square Root-Spiral.

Beside the "main" pronics-graphs PA, PB, PC ( see FIG 8 ) I drew the next parallel graphs to these graphs into FIG 8. These Graphs are named P+1A, P+1B, P+1C and P-1A, P-1B, P-1C and P-2A, P-2B, P-2C.
And their counterparts on the Ulam Spiral and on the Number Spiral are named P+1 as well as P-1 and P-2. → See FIG. NS-18 and / NS-8

<u>Note</u> : From these graphs only graph P-2 is marked on the Number Spiral ( → see NS-8 ). Because only the "pronics-graphs" P, P-2, P-6, P-12, P-20 etc. as well as the "S-graphs" S, S-1, S-4, S-9, S-16 etc. really contain "pronics-numbers".

Other P-graphs like P-4, P-8, P-14, P-16 etc. , which also lie in an even-number distance to graph P , also contain numbers which are composed by factors for which a certain rule applies. However these numbers are no real pronics !

<u>Note</u>: All P-graphs on the Number Spiral with an even-number-distance to graph P contain no Prime Numbers ! (→ see image NS-8 )

It is different with the P-graphs which have an odd-number distance to graph P, like the graphs P-1, P-3, P-5 etc. or P+1, P+3, P+5 etc. ( not marked in NS-8 )

These graphs contain an above-average share of prime numbers !
For example graph P+1 and P-1 ( see image NS-4 in chapter 9 ) or graph P+41 ( see image NS-17 ). These graphs contain high shares in Prime Numbers !

Graph **P+41** is already well known in mathematics as the " **Euler-Polynomial**" .
This graph which contains the number sequence 41, 43, 47, 53, 61, 71,…etc. and which contains a particular high share in Prime Numbers, is defined by the following quadratic polynomial : **$f(x) = x^2 + x + 41$** ( see image NS 17 )

On the Square Root Spiral the Euler Polynom divides into three spiral-graphs which I named **P+41A, P+41B** and **P+41C** ( see FIG 8 ). The number sequences of these graphs are shown in **Table 7-B**. ( → page 35 )

But let's first have a look to **Table 7-A** ( page 34 ). This table shows the number sequences of the "reference-graphs" which I used to draw a comparison between the Ulam-Spiral, the Number-Spiral and the Square Root Spiral.

The lefthand side of Table 7-A shows the number sequences of the two main reference graphs S and P as they appear in image NS-8 of the Number Spiral.

Beside the number sequences S and P, the next three columns show the number sequences of the spiral graphs SA, SB, SC and PA, PB, PC which are the corresponding "S- and P-graphs" on the Square Root Spiral. ( → see FIG 8 )

All above mentioned graphs are defined by "pronics-numbers". The factors of these pronics are shown in the yellow columns. In the number sequence of graph S the factors of the numbers increase by 0 and In the number sequence P the factors of the numbers increase by 1. It is the same with the factors in the numbers of the sequences SA, SB, SC and PA, PB, PC, however the factors of <u>successive</u> numbers in these sequences increase by 3.

By connecting the "pronics" on the spiral graphs PA, PB, PC with a continuous line, in the same order as they appear in the pronics-graph **P** ( with an increase of 1 between the factors of <u>successive</u> numbers ), a rotating triangular line pattern evolves ( → see red line pattern in FIG 8 ).
The angle between two successive lines of this triangular line pattern strives for

**$\pi$ - 2** for PA, PB, PC → $\infty$

Besides, a similar triangular line pattern would evolve if the square numbers in the graphs SA, SB, SC would be connected in the correct order by lines. And the angle in this triangular line pattern would strive for the same constant !

The factors of the pronics which I marked in red or blue in Table 7-A ( see yellow columns ), belong to two special number sequences which contain all existing prime numbers. I called these two special number sequences **SQ1** and **SQ2**
I had a closer look to these two important number sequences in another study which I intend to file with the arXiv – Archiv and which has the titel :

→ " The logic of the prime number distribution "

Number Sequence **SQ1** : 5, 11, 17, 23, 29, ….. ( → numbers marked in blue )

Number Sequence **SQ2** : 1, 7, 13, 19, 25 , 31, ….. ( → numbers marked in red )

Interesting are also the "sums of the digits" which occur in the number sequences of the reference graphs ! If we order the occurring sums of the digits according to their value, then the same sums of the digits sequences appear as already described for the prime number spiral graphs in Table 2.

Noticeable is here the "ordered" sums of the digits sequence which arises from the number sequence belonging to reference graph P ( → image NS-8 ).
The ordered sums of the digits sequence of graph P which is 2, 3, 6, 9, 11, 12, … shows the same periodicity ( …..3, 3, 2, 1,….), of the differences between the numbers in this sequences, as the sums of the digits sequences of the prime-number-sequences with the 2.Differential = 20 shown in Table 2.

And the ordered sums of the digits sequences which arise from the number sequences belonging to reference graphs SA, SB, SC and PA, PB, PC have either a periodicity of 3 or 9.

Worth mentioning is also the periodic occurrence of groups of four numbers which are not divisible by 2, 3 or 5 in the number sequences SA and SB (marked in red ).

On the righthand side of Table 7-A the number sequences S+1 and P-1 are shown. S+1 and P-1 are the two parallel graphs next to the reference graphs S and P in image NS-8 of the Number Spiral. These graphs are not marked in image NS-8 !

The next three columns on the right show the number sequences of the spiral graphs S+1A, S+1B, S+1C and P-1A, P-1B, P-1C which are the corresponding graphs on the Square Root Spiral. ( → see FIG 8 ).
Noticable is here that the number sequences P-1A, P-1B, P-1C are identical to the number sequences A3, B4, C5 shown in Table 6-A1, which are derived from the Prime Number Spiral Graphs A3, B4, C5 shown in FIG. 6-A !!

Worth mentioning are again the "ordered sums of the digits sequences", which arise from the number sequences S+1 and P-1 ( Number Spiral ), as well as from the number sequences S+1A, S+1B, S+1C and P-1A, P-1B, P-1C ( Square Root Spiral ).

The ordered sums of the digits sequences of graphs S+1 and P-1 are :
1, 5, 8, 10, 11, 14, 17, 19, 20, …… and 2, 5, 8, 10, 11, 14, 17, 19,…., which show the same periodicity ( …..3, 3, 2, 1,….), of the differences between the numbers in this sequences, as reference graph P and the sums of the digits sequences of the prime-number-sequences with the 2. Differential = 20 shown in Table 2.

It would definitely be interesting to find the real reason for the often occurrence of this periodicity (...3, 3, 2, 1,…) which not only occurs in the Square Root Spiral but also in the Number Spiral !!

And the ordered sums of the digits sequences which arise from the number sequences belonging to the graphs S+1A, S+1B, S+1C and P-1A, P-1B, P-1C have either a periodicity of 3 or 9.

Noticable is here the periodic occurrence of groups of three numbers which are not divisible by 2, 3 or 5 in the number sequences S+1A, S+1B, S+1C as well as the periodic occurrence of groups of four such numbers in the number sequences P-1A, P-1B, P-1C.

I want to come back now to **Table 7-B** ( → page 35 ), which shows the number sequence of Graph **P+41** ( see image NS 17 on page 32 – "The Number Spiral"), which is known as the " **Euler-Polynomial**" : $f(x) = x^2 + x + 41$

As mentioned before, in the Square Root Spiral the graph **P+41** divides into three spiral-graphs which I named **P+41A, P+41B** and **P+41C** ( see FIG 8 ). The number sequences of these graphs are shown in Table 7-B.

By the way, the spiral graphs and number sequences **P+41A, P+41B** and **P+41C are identical to** the prime number spiral graphs and corresponding number sequences **A17, B18** and **C 19 !!** Unfortunately I haven't carried out my analysis so far in FIG 6-A and Table 6-A1, otherwise these spiral graphs and number sequences would also have appeared here !!

The columns "prime factors" on the left of the number sequences P+41A, P+41B and P+41C show the prime factors of the first non-prime numbers in these sequences. In all three number sequences the smallest prime factor is 41.
Worth mentioning is also the fact, that the smallest number in the number sequences P+41A and P+41C is 41 , whereas the smallest number in the number sequence P+41B is 43.

Noticable is also the periodic occurrence of the prime factors in the non-prime numbers. For example prime factor 41 occurrs in the period 14/27 in number sequence P+41A and P+41B and in the period 27/14 in number sequence P+41C.

And prime factor 43 occurrs in the period 1/42 in number sequence P+41A and P+41C.

Important : As already described in Table 4 and 5A & 5B, the prime factors of the non-prime numbers of all quadratic polynomials, which lie on the Square Root Spiral, occur in clear defined periods !!
In **Table 7-C** ( on the bottom of page 35 ) the quadratic polynomials of the number sequences P+41 (Number Spiral) and P+41A, P+41B, P+41C ( Square Root Spiral ) are shown.

Note: Different quadratic polynomials can be calculated for the above mentioned number-sequences. The resulting quadratic polynomial depends on the selection of the three numbers out of these number-sequences, which are used for the calculation of the quadratic polynomials.

I have calculated the first four quadratic polynomials for each number sequence ( see Table 7-C ).

Worth mentioning is here the following interdependence :

Referring to the general quadratic polynomial : $f(x) = \mathbf{ax^2 + bx + c}$

the following pattern evolves in **Table 7-C** for the coefficients **a, b** and **c** :

**a** → is always equivalent to the **2.Differential** of the Spiral-Graph **divided by two** ( which is 18/2 = 9 for the quadratic polynomials P+41A, P+41B and P+41C)

**b** → this coefficient **increases by the value of the 2.Differential** from column to column shown in Table 7-C

**c** → this coefficient is always **equal to the number which comes before the three numbers, which were chosen to calculate the quadratic polynomial** in the number sequence belonging to this quadratic polynomial.
( Note : this rule doesn's apply to the first quadratic polynomial which was calculated with the first three numbers of the number-sequence, because of course there is no number before the first three numbers ! )

By the way, the same pattern exists in the Tables 6-A2 to 6-C2, which show the first four calculated quadratic polynomials of the Prime Number Spiral Graphs shown in FIG 6-A to 6-C !!

# Appendix :

This image shows exemplary some polar coordinates of the "Prime Number Spiral Graph" **B3** → see **FIG. 6-A**

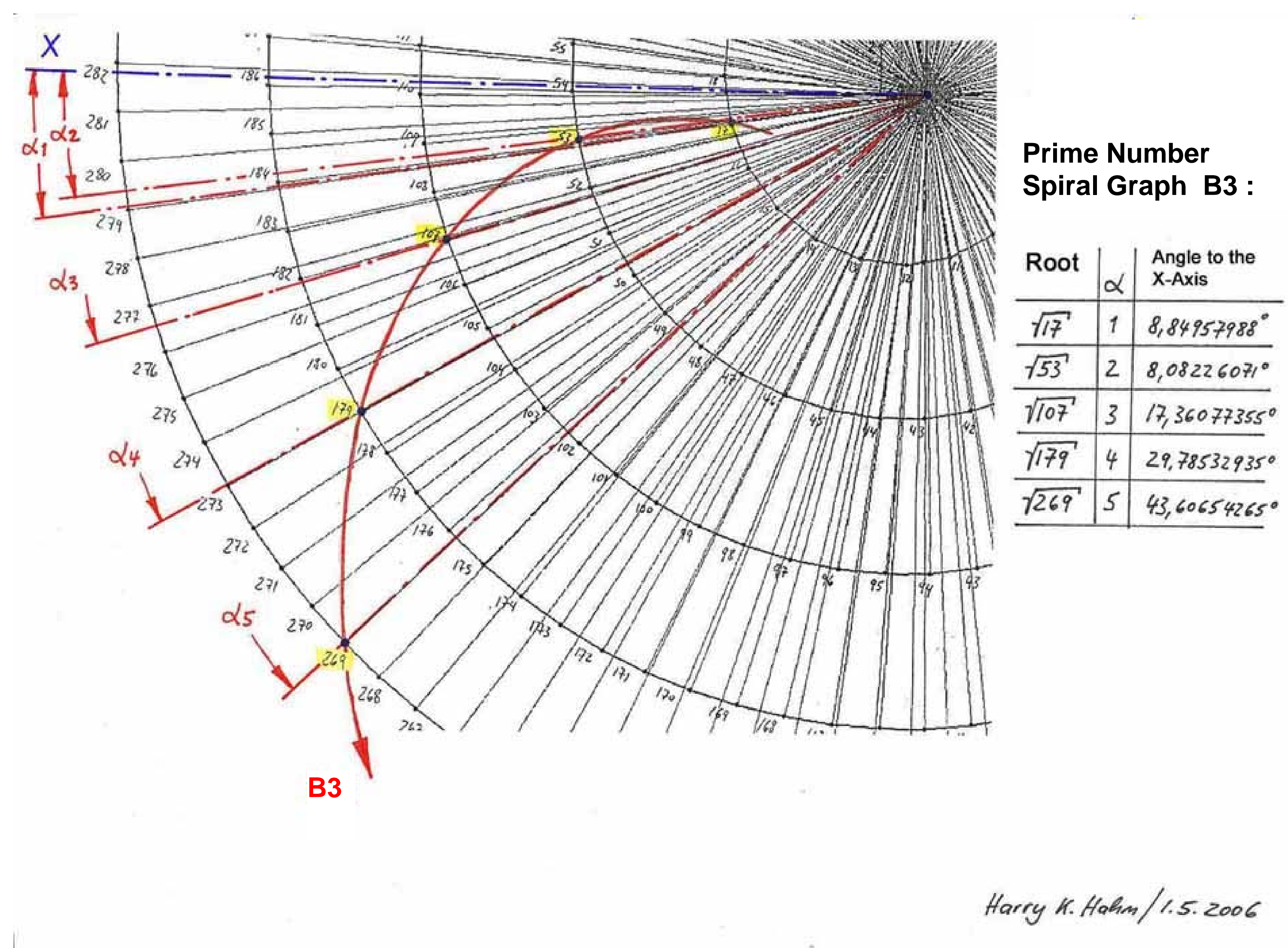

| Root | α | Angle to the X-Axis |
|---|---|---|
| $\sqrt{17}$ | 1 | 8,84957988° |
| $\sqrt{53}$ | 2 | 8,08226071° |
| $\sqrt{107}$ | 3 | 17,36077355° |
| $\sqrt{179}$ | 4 | 29,78532935° |
| $\sqrt{269}$ | 5 | 43,60654265° |

**Table 6-B1 :** "Prime Number Sequences" derived from the graphs schown in the "Prime Number Spiral Systems" **N20-D** to **N20-I** and **P20-D** to **P20-I** → see **FIG. 6-B**

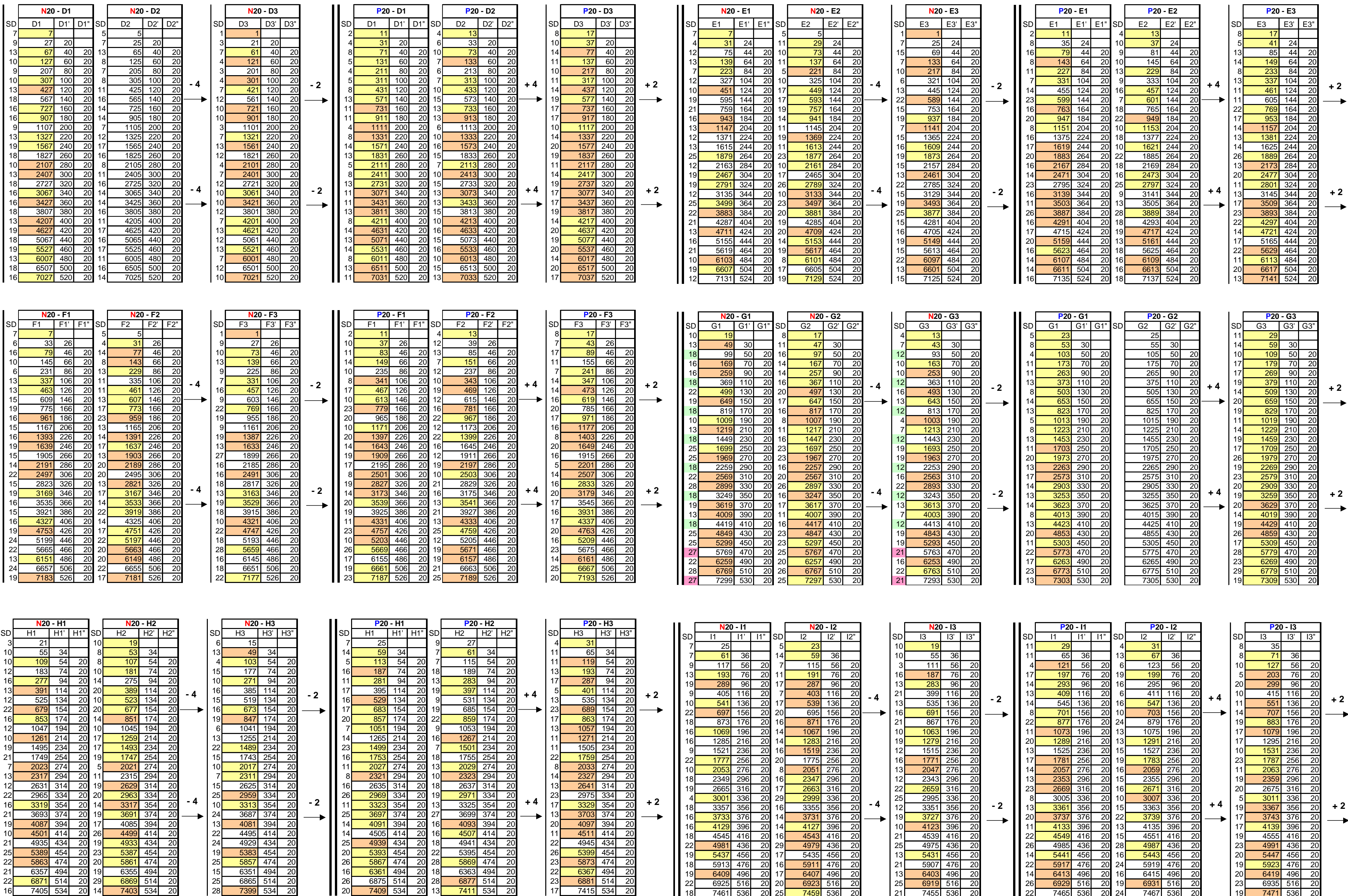

**N20 - D1**

| SD | D1 | D1' | D1" |
|---|---|---|---|
| 7 | 7 | | |
| 9 | 27 | 20 | |
| 13 | 67 | 40 | 20 |
| 10 | 127 | 60 | 20 |
| 9 | 207 | 80 | 20 |
| 10 | 307 | 100 | 20 |
| 13 | 427 | 120 | 20 |
| 18 | 567 | 140 | 20 |
| 16 | 727 | 160 | 20 |
| 16 | 907 | 180 | 20 |
| 9 | 1107 | 200 | 20 |
| 13 | 1327 | 220 | 20 |
| 19 | 1567 | 240 | 20 |
| 18 | 1827 | 260 | 20 |
| 10 | 2107 | 280 | 20 |
| 13 | 2407 | 300 | 20 |
| 18 | 2727 | 320 | 20 |
| 16 | 3067 | 340 | 20 |
| 16 | 3427 | 360 | 20 |
| 18 | 3807 | 380 | 20 |
| 13 | 4207 | 400 | 20 |
| 19 | 4627 | 420 | 20 |
| 18 | 5067 | 440 | 20 |
| 19 | 5527 | 460 | 20 |
| 13 | 6007 | 480 | 20 |
| 18 | 6507 | 500 | 20 |
| 16 | 7027 | 520 | 20 |

**N20 - D2**

| SD | D2 | D2' | D2" |
|---|---|---|---|
| 5 | 5 | | |
| 7 | 25 | 20 | |
| 13 | 65 | 40 | 20 |
| 8 | 125 | 60 | 20 |
| 7 | 205 | 80 | 20 |
| 8 | 305 | 100 | 20 |
| 11 | 425 | 120 | 20 |
| 16 | 565 | 140 | 20 |
| 14 | 725 | 160 | 20 |
| 14 | 905 | 180 | 20 |
| 7 | 1105 | 200 | 20 |
| 12 | 1325 | 220 | 20 |
| 17 | 1565 | 240 | 20 |
| 16 | 1825 | 260 | 20 |
| 8 | 2105 | 280 | 20 |
| 11 | 2405 | 300 | 20 |
| 16 | 2725 | 320 | 20 |
| 14 | 3065 | 340 | 20 |
| 14 | 3425 | 360 | 20 |
| 16 | 3805 | 380 | 20 |
| 11 | 4205 | 400 | 20 |
| 17 | 4625 | 420 | 20 |
| 16 | 5065 | 440 | 20 |
| 17 | 5525 | 460 | 20 |
| 11 | 6005 | 480 | 20 |
| 16 | 6505 | 500 | 20 |
| 14 | 7025 | 520 | 20 |

- 4 →

- 4 →

**N20 - D3**

| SD | D3 | D3' | D3" |
|---|---|---|---|
| 1 | 1 | | |
| 3 | 21 | 20 | |
| 7 | 61 | 40 | 20 |
| 4 | 121 | 60 | 20 |
| 3 | 201 | 80 | 20 |
| 4 | 301 | 100 | 20 |
| 7 | 421 | 120 | 20 |
| 12 | 561 | 140 | 20 |
| 10 | 721 | 160 | 20 |
| 10 | 901 | 180 | 20 |
| 3 | 1101 | 200 | 20 |
| 7 | 1321 | 220 | 20 |
| 13 | 1561 | 240 | 20 |
| 12 | 1821 | 260 | 20 |
| 4 | 2101 | 280 | 20 |
| 7 | 2401 | 300 | 20 |
| 12 | 2721 | 320 | 20 |
| 10 | 3061 | 340 | 20 |
| 10 | 3421 | 360 | 20 |
| 12 | 3801 | 380 | 20 |
| 7 | 4201 | 400 | 20 |
| 13 | 4621 | 420 | 20 |
| 12 | 5061 | 440 | 20 |
| 13 | 5521 | 460 | 20 |
| 7 | 6001 | 480 | 20 |
| 12 | 6501 | 500 | 20 |
| 10 | 7021 | 520 | 20 |

- 2 →

- 2 →

**P20 - D1**

| SD | D1 | D1' | D1" |
|---|---|---|---|
| 2 | 11 | | |
| 4 | 31 | 20 | |
| 8 | 71 | 40 | 20 |
| 5 | 131 | 60 | 20 |
| 4 | 211 | 80 | 20 |
| 5 | 311 | 100 | 20 |
| 8 | 431 | 120 | 20 |
| 13 | 571 | 140 | 20 |
| 11 | 731 | 160 | 20 |
| 11 | 911 | 180 | 20 |
| 4 | 1111 | 200 | 20 |
| 8 | 1331 | 220 | 20 |
| 14 | 1571 | 240 | 20 |
| 13 | 1831 | 260 | 20 |
| 5 | 2111 | 280 | 20 |
| 8 | 2411 | 300 | 20 |
| 13 | 2731 | 320 | 20 |
| 11 | 3071 | 340 | 20 |
| 11 | 3431 | 360 | 20 |
| 13 | 3811 | 380 | 20 |
| 8 | 4211 | 400 | 20 |
| 14 | 4631 | 420 | 20 |
| 13 | 5071 | 440 | 20 |
| 14 | 5531 | 460 | 20 |
| 8 | 6011 | 480 | 20 |
| 13 | 6511 | 500 | 20 |
| 11 | 7031 | 520 | 20 |

**P20 - D2**

| SD | D2 | D2' | D2" |
|---|---|---|---|
| 4 | 13 | | |
| 6 | 33 | 20 | |
| 10 | 73 | 40 | 20 |
| 7 | 133 | 60 | 20 |
| 6 | 213 | 80 | 20 |
| 7 | 313 | 100 | 20 |
| 10 | 433 | 120 | 20 |
| 15 | 573 | 140 | 20 |
| 13 | 733 | 160 | 20 |
| 13 | 913 | 180 | 20 |
| 6 | 1113 | 200 | 20 |
| 10 | 1333 | 220 | 20 |
| 16 | 1573 | 240 | 20 |
| 15 | 1833 | 260 | 20 |
| 7 | 2113 | 280 | 20 |
| 10 | 2413 | 300 | 20 |
| 15 | 2733 | 320 | 20 |
| 13 | 3073 | 340 | 20 |
| 13 | 3433 | 360 | 20 |
| 15 | 3813 | 380 | 20 |
| 10 | 4213 | 400 | 20 |
| 16 | 4633 | 420 | 20 |
| 15 | 5073 | 440 | 20 |
| 16 | 5533 | 460 | 20 |
| 10 | 6013 | 480 | 20 |
| 15 | 6513 | 500 | 20 |
| 13 | 7033 | 520 | 20 |

+ 4 →

+ 4 →

**P20 - D3**

| SD | D3 | D3' | D3" |
|---|---|---|---|
| 8 | 17 | | |
| 10 | 37 | 20 | |
| 14 | 77 | 40 | 20 |
| 11 | 137 | 60 | 20 |
| 10 | 217 | 80 | 20 |
| 11 | 317 | 100 | 20 |
| 14 | 437 | 120 | 20 |
| 19 | 577 | 140 | 20 |
| 17 | 737 | 160 | 20 |
| 17 | 917 | 180 | 20 |
| 10 | 1117 | 200 | 20 |
| 14 | 1337 | 220 | 20 |
| 20 | 1577 | 240 | 20 |
| 19 | 1837 | 260 | 20 |
| 11 | 2117 | 280 | 20 |
| 14 | 2417 | 300 | 20 |
| 19 | 2737 | 320 | 20 |
| 17 | 3077 | 340 | 20 |
| 17 | 3437 | 360 | 20 |
| 19 | 3817 | 380 | 20 |
| 14 | 4217 | 400 | 20 |
| 20 | 4637 | 420 | 20 |
| 19 | 5077 | 440 | 20 |
| 20 | 5537 | 460 | 20 |
| 14 | 6017 | 480 | 20 |
| 20 | 6517 | 500 | 20 |
| 17 | 7037 | 520 | 20 |

+ 2 →

+ 2 →

**N20 - E1**

| SD | E1 | E1' | E1" |
|---|---|---|---|
| 7 | 7 | | |
| 4 | 31 | 24 | |
| 12 | 75 | 44 | 20 |
| 13 | 139 | 64 | 20 |
| 7 | 223 | 84 | 20 |
| 12 | 327 | 104 | 20 |
| 10 | 451 | 124 | 20 |
| 19 | 595 | 144 | 20 |
| 21 | 759 | 164 | 20 |
| 16 | 943 | 184 | 20 |
| 13 | 1147 | 204 | 20 |
| 12 | 1371 | 224 | 20 |
| 13 | 1615 | 244 | 20 |
| 25 | 1879 | 264 | 20 |
| 12 | 2163 | 284 | 20 |
| 19 | 2467 | 304 | 20 |
| 19 | 2791 | 324 | 20 |
| 12 | 3135 | 344 | 20 |
| 25 | 3499 | 364 | 20 |
| 22 | 3883 | 384 | 20 |
| 21 | 4287 | 404 | 20 |
| 13 | 4711 | 424 | 20 |
| 16 | 5155 | 444 | 20 |
| 21 | 5619 | 464 | 20 |
| 10 | 6103 | 484 | 20 |
| 19 | 6607 | 504 | 20 |
| 12 | 7131 | 524 | 20 |

**N20 - E2**

| SD | E2 | E2' | E2" |
|---|---|---|---|
| 5 | 5 | | |
| 11 | 29 | 24 | |
| 10 | 73 | 44 | 20 |
| 11 | 137 | 64 | 20 |
| 5 | 221 | 84 | 20 |
| 10 | 325 | 104 | 20 |
| 17 | 449 | 124 | 20 |
| 17 | 593 | 144 | 20 |
| 19 | 757 | 164 | 20 |
| 14 | 941 | 184 | 20 |
| 11 | 1145 | 204 | 20 |
| 19 | 1369 | 224 | 20 |
| 11 | 1613 | 244 | 20 |
| 23 | 1877 | 264 | 20 |
| 10 | 2161 | 284 | 20 |
| 17 | 2465 | 304 | 20 |
| 26 | 2789 | 324 | 20 |
| 10 | 3133 | 344 | 20 |
| 23 | 3497 | 364 | 20 |
| 20 | 3881 | 384 | 20 |
| 19 | 4285 | 404 | 20 |
| 20 | 4709 | 424 | 20 |
| 14 | 5153 | 444 | 20 |
| 19 | 5617 | 464 | 20 |
| 8 | 6101 | 484 | 20 |
| 17 | 6605 | 504 | 20 |
| 19 | 7129 | 524 | 20 |

- 4 →

- 4 →

**N20 - E3**

| SD | E3 | E3' | E3" |
|---|---|---|---|
| 1 | 1 | | |
| 7 | 25 | 24 | |
| 15 | 69 | 44 | 20 |
| 7 | 133 | 64 | 20 |
| 10 | 217 | 84 | 20 |
| 6 | 321 | 104 | 20 |
| 13 | 445 | 124 | 20 |
| 22 | 589 | 144 | 20 |
| 15 | 753 | 164 | 20 |
| 19 | 937 | 184 | 20 |
| 7 | 1141 | 204 | 20 |
| 15 | 1365 | 224 | 20 |
| 16 | 1609 | 244 | 20 |
| 19 | 1873 | 264 | 20 |
| 15 | 2157 | 284 | 20 |
| 13 | 2461 | 304 | 20 |
| 22 | 2785 | 324 | 20 |
| 15 | 3129 | 344 | 20 |
| 19 | 3493 | 364 | 20 |
| 25 | 3877 | 384 | 20 |
| 15 | 4281 | 404 | 20 |
| 16 | 4705 | 424 | 20 |
| 19 | 5149 | 444 | 20 |
| 15 | 5613 | 464 | 20 |
| 22 | 6097 | 484 | 20 |
| 13 | 6601 | 504 | 20 |
| 15 | 7125 | 524 | 20 |

- 2 →

- 2 →

**P20 - E1**

| SD | E1 | E1' | E1" |
|---|---|---|---|
| 2 | 11 | | |
| 8 | 35 | 24 | |
| 16 | 79 | 44 | 20 |
| 8 | 143 | 64 | 20 |
| 11 | 227 | 84 | 20 |
| 7 | 331 | 104 | 20 |
| 14 | 455 | 124 | 20 |
| 23 | 599 | 144 | 20 |
| 16 | 763 | 164 | 20 |
| 20 | 947 | 184 | 20 |
| 8 | 1151 | 204 | 20 |
| 16 | 1375 | 224 | 20 |
| 17 | 1619 | 244 | 20 |
| 20 | 1883 | 264 | 20 |
| 16 | 2167 | 284 | 20 |
| 14 | 2471 | 304 | 20 |
| 23 | 2795 | 324 | 20 |
| 16 | 3139 | 344 | 20 |
| 11 | 3503 | 364 | 20 |
| 26 | 3887 | 384 | 20 |
| 16 | 4291 | 404 | 20 |
| 17 | 4715 | 424 | 20 |
| 20 | 5159 | 444 | 20 |
| 16 | 5623 | 464 | 20 |
| 14 | 6107 | 484 | 20 |
| 14 | 6611 | 504 | 20 |
| 16 | 7135 | 524 | 20 |

**P20 - E2**

| SD | E2 | E2' | E2" |
|---|---|---|---|
| 4 | 13 | | |
| 10 | 37 | 24 | |
| 9 | 81 | 44 | 20 |
| 10 | 145 | 64 | 20 |
| 13 | 229 | 84 | 20 |
| 9 | 333 | 104 | 20 |
| 16 | 457 | 124 | 20 |
| 7 | 601 | 144 | 20 |
| 18 | 765 | 164 | 20 |
| 22 | 949 | 184 | 20 |
| 10 | 1153 | 204 | 20 |
| 18 | 1377 | 224 | 20 |
| 10 | 1621 | 244 | 20 |
| 22 | 1885 | 264 | 20 |
| 18 | 2169 | 284 | 20 |
| 16 | 2473 | 304 | 20 |
| 25 | 2797 | 324 | 20 |
| 9 | 3141 | 344 | 20 |
| 13 | 3505 | 364 | 20 |
| 28 | 3889 | 384 | 20 |
| 18 | 4293 | 404 | 20 |
| 19 | 4717 | 424 | 20 |
| 13 | 5161 | 444 | 20 |
| 18 | 5625 | 464 | 20 |
| 16 | 6109 | 484 | 20 |
| 16 | 6613 | 504 | 20 |
| 18 | 7137 | 524 | 20 |

+ 4 →

+ 4 →

**P20 - E3**

| SD | E3 | E3' | E3" |
|---|---|---|---|
| 8 | 17 | | |
| 5 | 41 | 24 | |
| 13 | 85 | 44 | 20 |
| 14 | 149 | 64 | 20 |
| 8 | 233 | 84 | 20 |
| 13 | 337 | 104 | 20 |
| 11 | 461 | 124 | 20 |
| 11 | 605 | 144 | 20 |
| 22 | 769 | 164 | 20 |
| 17 | 953 | 184 | 20 |
| 14 | 1157 | 204 | 20 |
| 13 | 1381 | 224 | 20 |
| 14 | 1625 | 244 | 20 |
| 26 | 1889 | 264 | 20 |
| 13 | 2173 | 284 | 20 |
| 20 | 2477 | 304 | 20 |
| 11 | 2801 | 324 | 20 |
| 13 | 3145 | 344 | 20 |
| 17 | 3509 | 364 | 20 |
| 23 | 3893 | 384 | 20 |
| 22 | 4297 | 404 | 20 |
| 14 | 4721 | 424 | 20 |
| 17 | 5165 | 444 | 20 |
| 22 | 5629 | 464 | 20 |
| 11 | 6113 | 484 | 20 |
| 20 | 6617 | 504 | 20 |
| 13 | 7141 | 524 | 20 |

+ 2 →

+ 2 →

**N20 - F1**

| SD | F1 | F1' | F1" |
|---|---|---|---|
| 7 | 7 | | |
| 6 | 33 | 26 | |
| 16 | 79 | 46 | 20 |
| 10 | 145 | 66 | 20 |
| 6 | 231 | 86 | 20 |
| 13 | 337 | 106 | 20 |
| 13 | 463 | 126 | 20 |
| 15 | 609 | 146 | 20 |
| 19 | 775 | 166 | 20 |
| 16 | 961 | 186 | 20 |
| 15 | 1167 | 206 | 20 |
| 16 | 1393 | 226 | 20 |
| 19 | 1639 | 246 | 20 |
| 15 | 1905 | 266 | 20 |
| 14 | 2191 | 286 | 20 |
| 22 | 2497 | 306 | 20 |
| 15 | 2823 | 326 | 20 |
| 19 | 3169 | 346 | 20 |
| 16 | 3535 | 366 | 20 |
| 15 | 3921 | 386 | 20 |
| 16 | 4327 | 406 | 20 |
| 19 | 4753 | 426 | 20 |
| 24 | 5199 | 446 | 20 |
| 22 | 5665 | 466 | 20 |
| 13 | 6151 | 486 | 20 |
| 24 | 6657 | 506 | 20 |
| 19 | 7183 | 526 | 20 |

**N20 - F2**

| SD | F2 | F2' | F2" |
|---|---|---|---|
| 5 | 5 | | |
| 4 | 31 | 26 | |
| 14 | 77 | 46 | 20 |
| 8 | 143 | 66 | 20 |
| 13 | 229 | 86 | 20 |
| 11 | 335 | 106 | 20 |
| 11 | 461 | 126 | 20 |
| 13 | 607 | 146 | 20 |
| 17 | 773 | 166 | 20 |
| 23 | 959 | 186 | 20 |
| 13 | 1165 | 206 | 20 |
| 14 | 1391 | 226 | 20 |
| 17 | 1637 | 246 | 20 |
| 13 | 1903 | 266 | 20 |
| 20 | 2189 | 286 | 20 |
| 20 | 2495 | 306 | 20 |
| 13 | 2821 | 326 | 20 |
| 17 | 3167 | 346 | 20 |
| 14 | 3533 | 366 | 20 |
| 22 | 3919 | 386 | 20 |
| 14 | 4325 | 406 | 20 |
| 17 | 4751 | 426 | 20 |
| 22 | 5197 | 446 | 20 |
| 20 | 5663 | 466 | 20 |
| 20 | 6149 | 486 | 20 |
| 22 | 6655 | 506 | 20 |
| 17 | 7181 | 526 | 20 |

- 4 →

- 4 →

**N20 - F3**

| SD | F3 | F3' | F3" |
|---|---|---|---|
| 1 | 1 | | |
| 9 | 27 | 26 | |
| 10 | 73 | 46 | 20 |
| 13 | 139 | 66 | 20 |
| 9 | 225 | 86 | 20 |
| 7 | 331 | 106 | 20 |
| 16 | 457 | 126 | 20 |
| 9 | 603 | 146 | 20 |
| 22 | 769 | 166 | 20 |
| 19 | 955 | 186 | 20 |
| 9 | 1161 | 206 | 20 |
| 19 | 1387 | 226 | 20 |
| 13 | 1633 | 246 | 20 |
| 27 | 1899 | 266 | 20 |
| 16 | 2185 | 286 | 20 |
| 16 | 2491 | 306 | 20 |
| 18 | 2817 | 326 | 20 |
| 13 | 3163 | 346 | 20 |
| 19 | 3529 | 366 | 20 |
| 18 | 3915 | 386 | 20 |
| 10 | 4321 | 406 | 20 |
| 22 | 4747 | 426 | 20 |
| 18 | 5193 | 446 | 20 |
| 28 | 5659 | 466 | 20 |
| 16 | 6145 | 486 | 20 |
| 18 | 6651 | 506 | 20 |
| 22 | 7177 | 526 | 20 |

- 2 →

- 2 →

**P20 - F1**

| SD | F1 | F1' | F1" |
|---|---|---|---|
| 2 | 11 | | |
| 10 | 37 | 26 | |
| 11 | 83 | 46 | 20 |
| 14 | 149 | 66 | 20 |
| 10 | 235 | 86 | 20 |
| 8 | 341 | 106 | 20 |
| 17 | 467 | 126 | 20 |
| 10 | 613 | 146 | 20 |
| 23 | 779 | 166 | 20 |
| 20 | 965 | 186 | 20 |
| 10 | 1171 | 206 | 20 |
| 20 | 1397 | 226 | 20 |
| 14 | 1643 | 246 | 20 |
| 19 | 1909 | 266 | 20 |
| 17 | 2195 | 286 | 20 |
| 8 | 2501 | 306 | 20 |
| 19 | 2827 | 326 | 20 |
| 14 | 3173 | 346 | 20 |
| 20 | 3539 | 366 | 20 |
| 19 | 3925 | 386 | 20 |
| 11 | 4331 | 406 | 20 |
| 23 | 4757 | 426 | 20 |
| 10 | 5203 | 446 | 20 |
| 26 | 5669 | 466 | 20 |
| 17 | 6155 | 486 | 20 |
| 19 | 6661 | 506 | 20 |
| 23 | 7187 | 526 | 20 |

**P20 - F2**

| SD | F2 | F2' | F2" |
|---|---|---|---|
| 4 | 13 | | |
| 12 | 39 | 26 | |
| 13 | 85 | 46 | 20 |
| 7 | 151 | 66 | 20 |
| 12 | 237 | 86 | 20 |
| 10 | 343 | 106 | 20 |
| 19 | 469 | 126 | 20 |
| 12 | 615 | 146 | 20 |
| 16 | 781 | 166 | 20 |
| 22 | 967 | 186 | 20 |
| 12 | 1173 | 206 | 20 |
| 22 | 1399 | 226 | 20 |
| 16 | 1645 | 246 | 20 |
| 12 | 1911 | 266 | 20 |
| 19 | 2197 | 286 | 20 |
| 10 | 2503 | 306 | 20 |
| 21 | 2829 | 326 | 20 |
| 16 | 3175 | 346 | 20 |
| 13 | 3541 | 366 | 20 |
| 21 | 3927 | 386 | 20 |
| 13 | 4333 | 406 | 20 |
| 25 | 4759 | 426 | 20 |
| 12 | 5205 | 446 | 20 |
| 19 | 5671 | 466 | 20 |
| 19 | 6157 | 486 | 20 |
| 21 | 6663 | 506 | 20 |
| 25 | 7189 | 526 | 20 |

+ 4 →

+ 4 →

**P20 - F3**

| SD | F3 | F3' | F3" |
|---|---|---|---|
| 8 | 17 | | |
| 7 | 43 | 26 | |
| 17 | 89 | 46 | 20 |
| 11 | 155 | 66 | 20 |
| 7 | 241 | 86 | 20 |
| 14 | 347 | 106 | 20 |
| 14 | 473 | 126 | 20 |
| 16 | 619 | 146 | 20 |
| 20 | 785 | 166 | 20 |
| 17 | 971 | 186 | 20 |
| 16 | 1177 | 206 | 20 |
| 8 | 1403 | 226 | 20 |
| 20 | 1649 | 246 | 20 |
| 16 | 1915 | 266 | 20 |
| 5 | 2201 | 286 | 20 |
| 14 | 2507 | 306 | 20 |
| 16 | 2833 | 326 | 20 |
| 20 | 3179 | 346 | 20 |
| 17 | 3545 | 366 | 20 |
| 16 | 3931 | 386 | 20 |
| 17 | 4337 | 406 | 20 |
| 20 | 4763 | 426 | 20 |
| 16 | 5209 | 446 | 20 |
| 23 | 5675 | 466 | 20 |
| 14 | 6161 | 486 | 20 |
| 25 | 6667 | 506 | 20 |
| 20 | 7193 | 526 | 20 |

+ 2 →

+ 2 →

**N20 - G1**

| SD | G1 | G1' | G1" |
|---|---|---|---|
| 10 | 19 | | |
| 13 | 49 | 30 | |
| 18 | 99 | 50 | 20 |
| 16 | 169 | 70 | 20 |
| 16 | 259 | 90 | 20 |
| 18 | 369 | 110 | 20 |
| 22 | 499 | 130 | 20 |
| 19 | 649 | 150 | 20 |
| 18 | 819 | 170 | 20 |
| 10 | 1009 | 190 | 20 |
| 13 | 1219 | 210 | 20 |
| 18 | 1449 | 230 | 20 |
| 25 | 1699 | 250 | 20 |
| 25 | 1969 | 270 | 20 |
| 18 | 2259 | 290 | 20 |
| 22 | 2569 | 310 | 20 |
| 28 | 2899 | 330 | 20 |
| 18 | 3249 | 350 | 20 |
| 19 | 3619 | 370 | 20 |
| 13 | 4009 | 390 | 20 |
| 18 | 4419 | 410 | 20 |
| 25 | 4849 | 430 | 20 |
| 25 | 5299 | 450 | 20 |
| 27 | 5769 | 470 | 20 |
| 22 | 6259 | 490 | 20 |
| 28 | 6769 | 510 | 20 |
| 27 | 7299 | 530 | 20 |

**N20 - G2**

| SD | G2 | G2' | G2" |
|---|---|---|---|
| 8 | 17 | | |
| 11 | 47 | 30 | |
| 16 | 97 | 50 | 20 |
| 14 | 167 | 70 | 20 |
| 14 | 257 | 90 | 20 |
| 16 | 367 | 110 | 20 |
| 20 | 497 | 130 | 20 |
| 17 | 647 | 150 | 20 |
| 16 | 817 | 170 | 20 |
| 8 | 1007 | 190 | 20 |
| 11 | 1217 | 210 | 20 |
| 16 | 1447 | 230 | 20 |
| 23 | 1697 | 250 | 20 |
| 23 | 1967 | 270 | 20 |
| 16 | 2257 | 290 | 20 |
| 20 | 2567 | 310 | 20 |
| 26 | 2897 | 330 | 20 |
| 16 | 3247 | 350 | 20 |
| 17 | 3617 | 370 | 20 |
| 11 | 4007 | 390 | 20 |
| 16 | 4417 | 410 | 20 |
| 23 | 4847 | 430 | 20 |
| 23 | 5297 | 450 | 20 |
| 25 | 5767 | 470 | 20 |
| 20 | 6257 | 490 | 20 |
| 26 | 6767 | 510 | 20 |
| 25 | 7297 | 530 | 20 |

- 4 →

- 4 →

**N20 - G3**

| SD | G3 | G3' | G3" |
|---|---|---|---|
| 4 | 13 | | |
| 7 | 43 | 30 | |
| 12 | 93 | 50 | 20 |
| 10 | 163 | 70 | 20 |
| 10 | 253 | 90 | 20 |
| 12 | 363 | 110 | 20 |
| 16 | 493 | 130 | 20 |
| 13 | 643 | 150 | 20 |
| 12 | 813 | 170 | 20 |
| 4 | 1003 | 190 | 20 |
| 7 | 1213 | 210 | 20 |
| 12 | 1443 | 230 | 20 |
| 19 | 1693 | 250 | 20 |
| 19 | 1963 | 270 | 20 |
| 12 | 2253 | 290 | 20 |
| 16 | 2563 | 310 | 20 |
| 22 | 2893 | 330 | 20 |
| 12 | 3243 | 350 | 20 |
| 13 | 3613 | 370 | 20 |
| 7 | 4003 | 390 | 20 |
| 12 | 4413 | 410 | 20 |
| 19 | 4843 | 430 | 20 |
| 19 | 5293 | 450 | 20 |
| 21 | 5763 | 470 | 20 |
| 16 | 6253 | 490 | 20 |
| 22 | 6763 | 510 | 20 |
| 21 | 7293 | 530 | 20 |

- 2 →

- 2 →

**P20 - G1**

| SD | G1 | G1' | G1" |
|---|---|---|---|
| 5 | 23 | | |
| 8 | 53 | 30 | |
| 4 | 103 | 50 | 20 |
| 11 | 173 | 70 | 20 |
| 11 | 263 | 90 | 20 |
| 13 | 373 | 110 | 20 |
| 8 | 503 | 130 | 20 |
| 14 | 653 | 150 | 20 |
| 13 | 823 | 170 | 20 |
| 5 | 1013 | 190 | 20 |
| 8 | 1223 | 210 | 20 |
| 13 | 1453 | 230 | 20 |
| 11 | 1703 | 250 | 20 |
| 20 | 1973 | 270 | 20 |
| 13 | 2263 | 290 | 20 |
| 17 | 2573 | 310 | 20 |
| 14 | 2903 | 330 | 20 |
| 13 | 3253 | 350 | 20 |
| 14 | 3623 | 370 | 20 |
| 8 | 4013 | 390 | 20 |
| 13 | 4423 | 410 | 20 |
| 20 | 4853 | 430 | 20 |
| 11 | 5303 | 450 | 20 |
| 22 | 5773 | 470 | 20 |
| 17 | 6263 | 490 | 20 |
| 23 | 6773 | 510 | 20 |
| 13 | 7303 | 530 | 20 |

**P20 - G2**

| SD | G2 | G2' | G2" |
|---|---|---|---|
| | 25 | | |
| | 55 | 30 | |
| | 105 | 50 | 20 |
| | 175 | 70 | 20 |
| | 265 | 90 | 20 |
| | 375 | 110 | 20 |
| | 505 | 130 | 20 |
| | 655 | 150 | 20 |
| | 825 | 170 | 20 |
| | 1015 | 190 | 20 |
| | 1225 | 210 | 20 |
| | 1455 | 230 | 20 |
| | 1705 | 250 | 20 |
| | 1975 | 270 | 20 |
| | 2265 | 290 | 20 |
| | 2575 | 310 | 20 |
| | 2905 | 330 | 20 |
| | 3255 | 350 | 20 |
| | 3625 | 370 | 20 |
| | 4015 | 390 | 20 |
| | 4425 | 410 | 20 |
| | 4855 | 430 | 20 |
| | 5305 | 450 | 20 |
| | 5775 | 470 | 20 |
| | 6265 | 490 | 20 |
| | 6775 | 510 | 20 |
| | 7305 | 530 | 20 |

+ 4 →

+ 4 →

**P20 - G3**

| SD | G3 | G3' | G3" |
|---|---|---|---|
| 11 | 29 | | |
| 14 | 59 | 30 | |
| 10 | 109 | 50 | 20 |
| 17 | 179 | 70 | 20 |
| 17 | 269 | 90 | 20 |
| 19 | 379 | 110 | 20 |
| 14 | 509 | 130 | 20 |
| 20 | 659 | 150 | 20 |
| 19 | 829 | 170 | 20 |
| 11 | 1019 | 190 | 20 |
| 14 | 1229 | 210 | 20 |
| 19 | 1459 | 230 | 20 |
| 17 | 1709 | 250 | 20 |
| 26 | 1979 | 270 | 20 |
| 19 | 2269 | 290 | 20 |
| 23 | 2579 | 310 | 20 |
| 20 | 2909 | 330 | 20 |
| 19 | 3259 | 350 | 20 |
| 20 | 3629 | 370 | 20 |
| 14 | 4019 | 390 | 20 |
| 19 | 4429 | 410 | 20 |
| 26 | 4859 | 430 | 20 |
| 17 | 5309 | 450 | 20 |
| 28 | 5779 | 470 | 20 |
| 23 | 6269 | 490 | 20 |
| 29 | 6779 | 510 | 20 |
| 19 | 7309 | 530 | 20 |

+ 2 →

+ 2 →

**N20 - H1**

| SD | H1 | H1' | H1" |
|---|---|---|---|
| 3 | 21 | | |
| 10 | 55 | 34 | |
| 10 | 109 | 54 | 20 |
| 12 | 183 | 74 | 20 |
| 16 | 277 | 94 | 20 |
| 13 | 391 | 114 | 20 |
| 12 | 525 | 134 | 20 |
| 22 | 679 | 154 | 20 |
| 16 | 853 | 174 | 20 |
| 12 | 1047 | 194 | 20 |
| 10 | 1261 | 214 | 20 |
| 19 | 1495 | 234 | 20 |
| 21 | 1749 | 254 | 20 |
| 7 | 2023 | 274 | 20 |
| 13 | 2317 | 294 | 20 |
| 12 | 2631 | 314 | 20 |
| 22 | 2965 | 334 | 20 |
| 16 | 3319 | 354 | 20 |
| 21 | 3693 | 374 | 20 |
| 19 | 4087 | 394 | 20 |
| 10 | 4501 | 414 | 20 |
| 21 | 4935 | 434 | 20 |
| 25 | 5389 | 454 | 20 |
| 22 | 5863 | 474 | 20 |
| 21 | 6357 | 494 | 20 |
| 22 | 6871 | 514 | 20 |
| 16 | 7405 | 534 | 20 |

**N20 - H2**

| SD | H2 | H2' | H2" |
|---|---|---|---|
| 10 | 19 | | |
| 8 | 53 | 34 | |
| 8 | 107 | 54 | 20 |
| 10 | 181 | 74 | 20 |
| 14 | 275 | 94 | 20 |
| 20 | 389 | 114 | 20 |
| 10 | 523 | 134 | 20 |
| 20 | 677 | 154 | 20 |
| 14 | 851 | 174 | 20 |
| 10 | 1045 | 194 | 20 |
| 17 | 1259 | 214 | 20 |
| 17 | 1493 | 234 | 20 |
| 19 | 1747 | 254 | 20 |
| 5 | 2021 | 274 | 20 |
| 11 | 2315 | 294 | 20 |
| 19 | 2629 | 314 | 20 |
| 20 | 2963 | 334 | 20 |
| 14 | 3317 | 354 | 20 |
| 19 | 3691 | 374 | 20 |
| 17 | 4085 | 394 | 20 |
| 26 | 4499 | 414 | 20 |
| 19 | 4933 | 434 | 20 |
| 23 | 5387 | 454 | 20 |
| 20 | 5861 | 474 | 20 |
| 19 | 6355 | 494 | 20 |
| 29 | 6869 | 514 | 20 |
| 14 | 7403 | 534 | 20 |

- 4 →

- 4 →

**N20 - H3**

| SD | H3 | H3' | H3" |
|---|---|---|---|
| 6 | 15 | | |
| 13 | 49 | 34 | |
| 4 | 103 | 54 | 20 |
| 15 | 177 | 74 | 20 |
| 10 | 271 | 94 | 20 |
| 16 | 385 | 114 | 20 |
| 15 | 519 | 134 | 20 |
| 16 | 673 | 154 | 20 |
| 19 | 847 | 174 | 20 |
| 6 | 1041 | 194 | 20 |
| 13 | 1255 | 214 | 20 |
| 22 | 1489 | 234 | 20 |
| 15 | 1743 | 254 | 20 |
| 10 | 2017 | 274 | 20 |
| 7 | 2311 | 294 | 20 |
| 15 | 2625 | 314 | 20 |
| 25 | 2959 | 334 | 20 |
| 10 | 3313 | 354 | 20 |
| 24 | 3687 | 374 | 20 |
| 13 | 4081 | 394 | 20 |
| 22 | 4495 | 414 | 20 |
| 24 | 4929 | 434 | 20 |
| 19 | 5383 | 454 | 20 |
| 25 | 5857 | 474 | 20 |
| 15 | 6351 | 494 | 20 |
| 25 | 6865 | 514 | 20 |
| 28 | 7399 | 534 | 20 |

- 2 →

- 2 →

**P20 - H1**

| SD | H1 | H1' | H1" |
|---|---|---|---|
| 7 | 25 | | |
| 14 | 59 | 34 | |
| 5 | 113 | 54 | 20 |
| 16 | 187 | 74 | 20 |
| 11 | 281 | 94 | 20 |
| 17 | 395 | 114 | 20 |
| 16 | 529 | 134 | 20 |
| 17 | 683 | 154 | 20 |
| 20 | 857 | 174 | 20 |
| 7 | 1051 | 194 | 20 |
| 14 | 1265 | 214 | 20 |
| 23 | 1499 | 234 | 20 |
| 16 | 1753 | 254 | 20 |
| 11 | 2027 | 274 | 20 |
| 8 | 2321 | 294 | 20 |
| 16 | 2635 | 314 | 20 |
| 26 | 2969 | 334 | 20 |
| 11 | 3323 | 354 | 20 |
| 25 | 3697 | 374 | 20 |
| 14 | 4091 | 394 | 20 |
| 14 | 4505 | 414 | 20 |
| 25 | 4939 | 434 | 20 |
| 20 | 5393 | 454 | 20 |
| 26 | 5867 | 474 | 20 |
| 16 | 6361 | 494 | 20 |
| 26 | 6875 | 514 | 20 |
| 20 | 7409 | 534 | 20 |

**P20 - H2**

| SD | H2 | H2' | H2" |
|---|---|---|---|
| 9 | 27 | | |
| 7 | 61 | 34 | |
| 7 | 115 | 54 | 20 |
| 18 | 189 | 74 | 20 |
| 13 | 283 | 94 | 20 |
| 19 | 397 | 114 | 20 |
| 9 | 531 | 134 | 20 |
| 19 | 685 | 154 | 20 |
| 22 | 859 | 174 | 20 |
| 9 | 1053 | 194 | 20 |
| 16 | 1267 | 214 | 20 |
| 7 | 1501 | 234 | 20 |
| 18 | 1755 | 254 | 20 |
| 13 | 2029 | 274 | 20 |
| 10 | 2323 | 294 | 20 |
| 18 | 2637 | 314 | 20 |
| 19 | 2971 | 334 | 20 |
| 13 | 3325 | 354 | 20 |
| 27 | 3699 | 374 | 20 |
| 16 | 4093 | 394 | 20 |
| 16 | 4507 | 414 | 20 |
| 18 | 4941 | 434 | 20 |
| 22 | 5395 | 454 | 20 |
| 28 | 5869 | 474 | 20 |
| 18 | 6363 | 494 | 20 |
| 28 | 6877 | 514 | 20 |
| 13 | 7411 | 534 | 20 |

+ 4 →

+ 4 →

**P20 - H3**

| SD | H3 | H3' | H3" |
|---|---|---|---|
| 4 | 31 | | |
| 11 | 65 | 34 | |
| 11 | 119 | 54 | 20 |
| 13 | 193 | 74 | 20 |
| 17 | 287 | 94 | 20 |
| 5 | 401 | 114 | 20 |
| 13 | 535 | 134 | 20 |
| 23 | 689 | 154 | 20 |
| 17 | 863 | 174 | 20 |
| 13 | 1057 | 194 | 20 |
| 11 | 1271 | 214 | 20 |
| 11 | 1505 | 234 | 20 |
| 22 | 1759 | 254 | 20 |
| 8 | 2033 | 274 | 20 |
| 14 | 2327 | 294 | 20 |
| 13 | 2641 | 314 | 20 |
| 23 | 2975 | 334 | 20 |
| 17 | 3329 | 354 | 20 |
| 13 | 3703 | 374 | 20 |
| 20 | 4097 | 394 | 20 |
| 11 | 4511 | 414 | 20 |
| 22 | 4945 | 434 | 20 |
| 26 | 5399 | 454 | 20 |
| 23 | 5873 | 474 | 20 |
| 22 | 6367 | 494 | 20 |
| 23 | 6881 | 514 | 20 |
| 17 | 7415 | 534 | 20 |

+ 2 →

+ 2 →

**N20 - I1**

| SD | I1 | I1' | I1" |
|---|---|---|---|
| 7 | 25 | | |
| 7 | 61 | 36 | |
| 9 | 117 | 56 | 20 |
| 13 | 193 | 76 | 20 |
| 19 | 289 | 96 | 20 |
| 9 | 405 | 116 | 20 |
| 10 | 541 | 136 | 20 |
| 22 | 697 | 156 | 20 |
| 18 | 873 | 176 | 20 |
| 16 | 1069 | 196 | 20 |
| 16 | 1285 | 216 | 20 |
| 9 | 1521 | 236 | 20 |
| 22 | 1777 | 256 | 20 |
| 10 | 2053 | 276 | 20 |
| 18 | 2349 | 296 | 20 |
| 19 | 2665 | 316 | 20 |
| 4 | 3001 | 336 | 20 |
| 18 | 3357 | 356 | 20 |
| 16 | 3733 | 376 | 20 |
| 16 | 4129 | 396 | 20 |
| 18 | 4545 | 416 | 20 |
| 22 | 4981 | 436 | 20 |
| 19 | 5437 | 456 | 20 |
| 18 | 5913 | 476 | 20 |
| 19 | 6409 | 496 | 20 |
| 22 | 6925 | 516 | 20 |
| 18 | 7461 | 536 | 20 |

**N20 - I2**

| SD | I2 | I2' | I2" |
|---|---|---|---|
| 5 | 23 | | |
| 14 | 59 | 36 | |
| 7 | 115 | 56 | 20 |
| 11 | 191 | 76 | 20 |
| 17 | 287 | 96 | 20 |
| 7 | 403 | 116 | 20 |
| 17 | 539 | 136 | 20 |
| 20 | 695 | 156 | 20 |
| 16 | 871 | 176 | 20 |
| 14 | 1067 | 196 | 20 |
| 14 | 1283 | 216 | 20 |
| 16 | 1519 | 236 | 20 |
| 20 | 1775 | 256 | 20 |
| 8 | 2051 | 276 | 20 |
| 16 | 2347 | 296 | 20 |
| 17 | 2663 | 316 | 20 |
| 29 | 2999 | 336 | 20 |
| 16 | 3355 | 356 | 20 |
| 14 | 3731 | 376 | 20 |
| 14 | 4127 | 396 | 20 |
| 16 | 4543 | 416 | 20 |
| 29 | 4979 | 436 | 20 |
| 17 | 5435 | 456 | 20 |
| 16 | 5911 | 476 | 20 |
| 17 | 6407 | 496 | 20 |
| 20 | 6923 | 516 | 20 |
| 25 | 7459 | 536 | 20 |

- 4 →

- 4 →

**N20 - I3**

| SD | I3 | I3' | I3" |
|---|---|---|---|
| 10 | 19 | | |
| 10 | 55 | 36 | |
| 3 | 111 | 56 | 20 |
| 16 | 187 | 76 | 20 |
| 13 | 283 | 96 | 20 |
| 21 | 399 | 116 | 20 |
| 13 | 535 | 136 | 20 |
| 16 | 691 | 156 | 20 |
| 21 | 867 | 176 | 20 |
| 10 | 1063 | 196 | 20 |
| 19 | 1279 | 216 | 20 |
| 12 | 1515 | 236 | 20 |
| 16 | 1771 | 256 | 20 |
| 13 | 2047 | 276 | 20 |
| 12 | 2343 | 296 | 20 |
| 22 | 2659 | 316 | 20 |
| 25 | 2995 | 336 | 20 |
| 12 | 3351 | 356 | 20 |
| 19 | 3727 | 376 | 20 |
| 10 | 4123 | 396 | 20 |
| 21 | 4539 | 416 | 20 |
| 25 | 4975 | 436 | 20 |
| 13 | 5431 | 456 | 20 |
| 21 | 5907 | 476 | 20 |
| 13 | 6403 | 496 | 20 |
| 25 | 6919 | 516 | 20 |
| 21 | 7455 | 536 | 20 |

- 2 →

- 2 →

**P20 - I1**

| SD | I1 | I1' | I1" |
|---|---|---|---|
| 11 | 29 | | |
| 11 | 65 | 36 | |
| 4 | 121 | 56 | 20 |
| 17 | 197 | 76 | 20 |
| 14 | 293 | 96 | 20 |
| 13 | 409 | 116 | 20 |
| 14 | 545 | 136 | 20 |
| 8 | 701 | 156 | 20 |
| 22 | 877 | 176 | 20 |
| 11 | 1073 | 196 | 20 |
| 20 | 1289 | 216 | 20 |
| 13 | 1525 | 236 | 20 |
| 17 | 1781 | 256 | 20 |
| 14 | 2057 | 276 | 20 |
| 13 | 2353 | 296 | 20 |
| 23 | 2669 | 316 | 20 |
| 8 | 3005 | 336 | 20 |
| 13 | 3361 | 356 | 20 |
| 20 | 3737 | 376 | 20 |
| 11 | 4133 | 396 | 20 |
| 22 | 4549 | 416 | 20 |
| 26 | 4985 | 436 | 20 |
| 14 | 5441 | 456 | 20 |
| 22 | 5917 | 476 | 20 |
| 14 | 6413 | 496 | 20 |
| 26 | 6929 | 516 | 20 |
| 22 | 7465 | 536 | 20 |

**P20 - I2**

| SD | I2 | I2' | I2" |
|---|---|---|---|
| 4 | 31 | | |
| 13 | 67 | 36 | |
| 6 | 123 | 56 | 20 |
| 19 | 199 | 76 | 20 |
| 16 | 295 | 96 | 20 |
| 6 | 411 | 116 | 20 |
| 16 | 547 | 136 | 20 |
| 10 | 703 | 156 | 20 |
| 24 | 879 | 176 | 20 |
| 13 | 1075 | 196 | 20 |
| 13 | 1291 | 216 | 20 |
| 15 | 1527 | 236 | 20 |
| 19 | 1783 | 256 | 20 |
| 16 | 2059 | 276 | 20 |
| 15 | 2355 | 296 | 20 |
| 16 | 2671 | 316 | 20 |
| 10 | 3007 | 336 | 20 |
| 15 | 3363 | 356 | 20 |
| 22 | 3739 | 376 | 20 |
| 13 | 4135 | 396 | 20 |
| 15 | 4551 | 416 | 20 |
| 28 | 4987 | 436 | 20 |
| 16 | 5443 | 456 | 20 |
| 24 | 5919 | 476 | 20 |
| 16 | 6415 | 496 | 20 |
| 19 | 6931 | 516 | 20 |
| 24 | 7467 | 536 | 20 |

+ 4 →

+ 4 →

**P20 - I3**

| SD | I3 | I3' | I3" |
|---|---|---|---|
| 8 | 35 | | |
| 8 | 71 | 36 | |
| 10 | 127 | 56 | 20 |
| 5 | 203 | 76 | 20 |
| 20 | 299 | 96 | 20 |
| 10 | 415 | 116 | 20 |
| 11 | 551 | 136 | 20 |
| 14 | 707 | 156 | 20 |
| 19 | 883 | 176 | 20 |
| 17 | 1079 | 196 | 20 |
| 17 | 1295 | 216 | 20 |
| 10 | 1531 | 236 | 20 |
| 23 | 1787 | 256 | 20 |
| 11 | 2063 | 276 | 20 |
| 19 | 2359 | 296 | 20 |
| 20 | 2675 | 316 | 20 |
| 5 | 3011 | 336 | 20 |
| 19 | 3367 | 356 | 20 |
| 17 | 3743 | 376 | 20 |
| 17 | 4139 | 396 | 20 |
| 19 | 4555 | 416 | 20 |
| 23 | 4991 | 436 | 20 |
| 20 | 5447 | 456 | 20 |
| 19 | 5923 | 476 | 20 |
| 20 | 6419 | 496 | 20 |
| 23 | 6935 | 516 | 20 |
| 19 | 7471 | 536 | 20 |

+ 2 →

+ 2 →

**Table 6-C1 :** “Prime Number Sequences” derived from the graphs schown in the “Prime Number Spiral Systems” **N22-J** to **N22-K** → see **FIG. 6-C**

| | N22 - J1 | | | | N22 - J2 | | | | N22 - J3 | | |
|---|---|---|---|---|---|---|---|---|---|---|---|
| SD | J1 | J1' | J1" | SD | J2 | J2' | J2" | SD | J3 | J3' | J3" |
| 4 | 31 | | | 11 | 29 | | | 7 | 25 | | |
| 15 | 69 | 38 | | 13 | 67 | 38 | | 9 | 63 | 38 | |
| 12 | 129 | 60 | 22 | 10 | 127 | 60 | 22 | 6 | 123 | 60 | 22 |
| 4 | 211 | 82 | 22 | 11 | 209 | 82 | 22 | 7 | 205 | 82 | 22 |
| 9 | 315 | 104 | 22 | 7 | 313 | 104 | 22 | 12 | 309 | 104 | 22 |
| 9 | 441 | 126 | 22 | 16 | 439 | 126 | 22 | 12 | 435 | 126 | 22 |
| 22 | 589 | 148 | 22 | 20 | 587 | 148 | 22 | 16 | 583 | 148 | 22 |
| 21 | 759 | 170 | 22 | 19 | 757 | 170 | 22 | 15 | 753 | 170 | 22 |
| 15 | 951 | 192 | 22 | 22 | 949 | 192 | 22 | 18 | 945 | 192 | 22 |
| 13 | 1165 | 214 | 22 | 11 | 1163 | 214 | 22 | 16 | 1159 | 214 | 22 |
| 6 | 1401 | 236 | 22 | 22 | 1399 | 236 | 22 | 18 | 1395 | 236 | 22 |
| 21 | 1659 | 258 | 22 | 19 | 1657 | 258 | 22 | 15 | 1653 | 258 | 22 |
| 22 | 1939 | 280 | 22 | 20 | 1937 | 280 | 22 | 16 | 1933 | 280 | 22 |
| 9 | 2241 | 302 | 22 | 16 | 2239 | 302 | 22 | 12 | 2235 | 302 | 22 |
| 18 | 2565 | 324 | 22 | 16 | 2563 | 324 | 22 | 21 | 2559 | 324 | 22 |
| 13 | 2911 | 346 | 22 | 20 | 2909 | 346 | 22 | 16 | 2905 | 346 | 22 |
| 21 | 3279 | 368 | 22 | 19 | 3277 | 368 | 22 | 15 | 3273 | 368 | 22 |
| 24 | 3669 | 390 | 22 | 22 | 3667 | 390 | 22 | 18 | 3663 | 390 | 22 |
| 13 | 4081 | 412 | 22 | 20 | 4079 | 412 | 22 | 16 | 4075 | 412 | 22 |
| 15 | 4515 | 434 | 22 | 13 | 4513 | 434 | 22 | 18 | 4509 | 434 | 22 |
| 21 | 4971 | 456 | 22 | 28 | 4969 | 456 | 22 | 24 | 4965 | 456 | 22 |
| 22 | 5449 | 478 | 22 | 20 | 5447 | 478 | 22 | 16 | 5443 | 478 | 22 |
| 27 | 5949 | 500 | 22 | 25 | 5947 | 500 | 22 | 21 | 5943 | 500 | 22 |
| 18 | 6471 | 522 | 22 | 25 | 6469 | 522 | 22 | 21 | 6465 | 522 | 22 |
| 13 | 7015 | 544 | 22 | 11 | 7013 | 544 | 22 | 16 | 7009 | 544 | 22 |
| 21 | 7581 | 566 | 22 | 28 | 7579 | 566 | 22 | 24 | 7575 | 566 | 22 |
| 24 | 8169 | 588 | 22 | 22 | 8167 | 588 | 22 | 18 | 8163 | 588 | 22 |
| 31 | 8779 | 610 | 22 | 29 | 8777 | 610 | 22 | 29 | 8773 | 610 | 22 |

- 4 → - 4 → - 2 → - 2 →

| | N22 - K1 | | | | N22 - K2 | | | | N22 - K3 | | |
|---|---|---|---|---|---|---|---|---|---|---|---|
| SD | K1 | K1' | K1" | SD | K2 | K2' | K2" | SD | K3 | K3' | K3" |
| 7 | 25 | | | 5 | 23 | | | 10 | 19 | | |
| 7 | 61 | 36 | | 14 | 59 | 36 | | 10 | 55 | 36 | |
| 11 | 119 | 58 | 22 | 9 | 117 | 58 | 22 | 5 | 113 | 58 | 22 |
| 19 | 199 | 80 | 22 | 17 | 197 | 80 | 22 | 13 | 193 | 80 | 22 |
| 4 | 301 | 102 | 22 | 20 | 299 | 102 | 22 | 16 | 295 | 102 | 22 |
| 11 | 425 | 124 | 22 | 9 | 423 | 124 | 22 | 14 | 419 | 124 | 22 |
| 13 | 571 | 146 | 22 | 20 | 569 | 146 | 22 | 16 | 565 | 146 | 22 |
| 19 | 739 | 168 | 22 | 17 | 737 | 168 | 22 | 13 | 733 | 168 | 22 |
| 20 | 929 | 190 | 22 | 18 | 927 | 190 | 22 | 14 | 923 | 190 | 22 |
| 7 | 1141 | 212 | 22 | 14 | 1139 | 212 | 22 | 10 | 1135 | 212 | 22 |
| 16 | 1375 | 234 | 22 | 14 | 1373 | 234 | 22 | 19 | 1369 | 234 | 22 |
| 11 | 1631 | 256 | 22 | 18 | 1629 | 256 | 22 | 14 | 1625 | 256 | 22 |
| 19 | 1909 | 278 | 22 | 17 | 1907 | 278 | 22 | 13 | 1903 | 278 | 22 |
| 13 | 2209 | 300 | 22 | 11 | 2207 | 300 | 22 | 7 | 2203 | 300 | 22 |
| 11 | 2531 | 322 | 22 | 18 | 2529 | 322 | 22 | 14 | 2525 | 322 | 22 |
| 22 | 2875 | 344 | 22 | 20 | 2873 | 344 | 22 | 25 | 2869 | 344 | 22 |
| 10 | 3241 | 366 | 22 | 17 | 3239 | 366 | 22 | 13 | 3235 | 366 | 22 |
| 20 | 3629 | 388 | 22 | 18 | 3627 | 388 | 22 | 14 | 3623 | 388 | 22 |
| 16 | 4039 | 410 | 22 | 14 | 4037 | 410 | 22 | 10 | 4033 | 410 | 22 |
| 16 | 4471 | 432 | 22 | 23 | 4469 | 432 | 22 | 19 | 4465 | 432 | 22 |
| 20 | 4925 | 454 | 22 | 18 | 4923 | 454 | 22 | 23 | 4919 | 454 | 22 |
| 10 | 5401 | 476 | 22 | 26 | 5399 | 476 | 22 | 22 | 5395 | 476 | 22 |
| 31 | 5899 | 498 | 22 | 29 | 5897 | 498 | 22 | 25 | 5893 | 498 | 22 |
| 20 | 6419 | 520 | 22 | 18 | 6417 | 520 | 22 | 14 | 6413 | 520 | 22 |
| 22 | 6961 | 542 | 22 | 29 | 6959 | 542 | 22 | 25 | 6955 | 542 | 22 |
| 19 | 7525 | 564 | 22 | 17 | 7523 | 564 | 22 | 22 | 7519 | 564 | 22 |
| 11 | 8111 | 586 | 22 | 18 | 8109 | 586 | 22 | 14 | 8105 | 586 | 22 |
| 25 | 8719 | 608 | 22 | 23 | 8717 | 608 | 22 | 19 | 8713 | 608 | 22 |

- 4 → - 4 → - 2 → - 2 →

| | N22 - L1 | | | | N22 - L2 | | | | N22 - L3 | | |
|---|---|---|---|---|---|---|---|---|---|---|---|
| SD | L1 | L1' | L1" | SD | L2 | L2' | L2" | SD | L3 | L3' | L3" |
| 7 | 25 | | | 5 | 23 | | | 10 | 19 | | |
| 14 | 59 | 34 | | 12 | 57 | 34 | | 8 | 53 | 34 | |
| 7 | 115 | 56 | 22 | 5 | 113 | 56 | 22 | 10 | 109 | 56 | 22 |
| 13 | 193 | 78 | 22 | 11 | 191 | 78 | 22 | 16 | 187 | 78 | 22 |
| 14 | 293 | 100 | 22 | 12 | 291 | 100 | 22 | 17 | 287 | 100 | 22 |
| 10 | 415 | 122 | 22 | 8 | 413 | 122 | 22 | 13 | 409 | 122 | 22 |
| 19 | 559 | 144 | 22 | 17 | 557 | 144 | 22 | 13 | 553 | 144 | 22 |
| 14 | 725 | 166 | 22 | 12 | 723 | 166 | 22 | 17 | 719 | 166 | 22 |
| 13 | 913 | 188 | 22 | 11 | 911 | 188 | 22 | 16 | 907 | 188 | 22 |
| 7 | 1123 | 210 | 22 | 5 | 1121 | 210 | 22 | 10 | 1117 | 210 | 22 |
| 14 | 1355 | 232 | 22 | 12 | 1353 | 232 | 22 | 17 | 1349 | 232 | 22 |
| 16 | 1609 | 254 | 22 | 14 | 1607 | 254 | 22 | 10 | 1603 | 254 | 22 |
| 22 | 1885 | 276 | 22 | 20 | 1883 | 276 | 22 | 25 | 1879 | 276 | 22 |
| 14 | 2183 | 298 | 22 | 12 | 2181 | 298 | 22 | 17 | 2177 | 298 | 22 |
| 10 | 2503 | 320 | 22 | 8 | 2501 | 320 | 22 | 22 | 2497 | 320 | 22 |
| 19 | 2845 | 342 | 22 | 17 | 2843 | 342 | 22 | 22 | 2839 | 342 | 22 |
| 14 | 3209 | 364 | 22 | 12 | 3207 | 364 | 22 | 8 | 3203 | 364 | 22 |
| 22 | 3595 | 386 | 22 | 20 | 3593 | 386 | 22 | 25 | 3589 | 386 | 22 |
| 7 | 4003 | 408 | 22 | 5 | 4001 | 408 | 22 | 28 | 3997 | 408 | 22 |
| 14 | 4433 | 430 | 22 | 12 | 4431 | 430 | 22 | 17 | 4427 | 430 | 22 |
| 25 | 4885 | 452 | 22 | 23 | 4883 | 452 | 22 | 28 | 4879 | 452 | 22 |
| 22 | 5359 | 474 | 22 | 20 | 5357 | 474 | 22 | 16 | 5353 | 474 | 22 |
| 23 | 5855 | 496 | 22 | 21 | 5853 | 496 | 22 | 26 | 5849 | 496 | 22 |
| 19 | 6373 | 518 | 22 | 17 | 6371 | 518 | 22 | 22 | 6367 | 518 | 22 |
| 19 | 6913 | 540 | 22 | 17 | 6911 | 540 | 22 | 22 | 6907 | 540 | 22 |
| 23 | 7475 | 562 | 22 | 21 | 7473 | 562 | 22 | 26 | 7469 | 562 | 22 |
| 22 | 8059 | 584 | 22 | 20 | 8057 | 584 | 22 | 16 | 8053 | 584 | 22 |
| 25 | 8665 | 606 | 22 | 23 | 8663 | 606 | 22 | 28 | 8659 | 606 | 22 |

- 4 → - 4 → - 2 → - 2 →

| | N22 - M1 | | | | N22 - M2 | | | | N22 - M3 | | |
|---|---|---|---|---|---|---|---|---|---|---|---|
| SD | M1 | M1' | M1" | SD | M2 | M2' | M2" | SD | M3 | M3' | M3" |
| 5 | 23 | | | 3 | 21 | | | 8 | 17 | | |
| 10 | 55 | 32 | | 8 | 53 | 32 | | 13 | 49 | 32 | |
| 10 | 109 | 54 | 22 | 8 | 107 | 54 | 22 | 4 | 103 | 54 | 22 |
| 14 | 185 | 76 | 22 | 12 | 183 | 76 | 22 | 17 | 179 | 76 | 22 |
| 13 | 283 | 98 | 22 | 11 | 281 | 98 | 22 | 16 | 277 | 98 | 22 |
| 7 | 403 | 120 | 22 | 5 | 401 | 120 | 22 | 19 | 397 | 120 | 22 |
| 14 | 545 | 142 | 22 | 12 | 543 | 142 | 22 | 17 | 539 | 142 | 22 |
| 16 | 709 | 164 | 22 | 14 | 707 | 164 | 22 | 10 | 703 | 164 | 22 |
| 22 | 895 | 186 | 22 | 20 | 893 | 186 | 22 | 25 | 889 | 186 | 22 |
| 5 | 1103 | 208 | 22 | 3 | 1101 | 208 | 22 | 17 | 1097 | 208 | 22 |
| 10 | 1333 | 230 | 22 | 8 | 1331 | 230 | 22 | 13 | 1327 | 230 | 22 |
| 19 | 1585 | 252 | 22 | 17 | 1583 | 252 | 22 | 22 | 1579 | 252 | 22 |
| 23 | 1859 | 274 | 22 | 21 | 1857 | 274 | 22 | 17 | 1853 | 274 | 22 |
| 13 | 2155 | 296 | 22 | 11 | 2153 | 296 | 22 | 16 | 2149 | 296 | 22 |
| 16 | 2473 | 318 | 22 | 14 | 2471 | 318 | 22 | 19 | 2467 | 318 | 22 |
| 14 | 2813 | 340 | 22 | 12 | 2811 | 340 | 22 | 17 | 2807 | 340 | 22 |
| 16 | 3175 | 362 | 22 | 14 | 3173 | 362 | 22 | 19 | 3169 | 362 | 22 |
| 22 | 3559 | 384 | 22 | 20 | 3557 | 384 | 22 | 16 | 3553 | 384 | 22 |
| 23 | 3965 | 406 | 22 | 21 | 3963 | 406 | 22 | 26 | 3959 | 406 | 22 |
| 19 | 4393 | 428 | 22 | 17 | 4391 | 428 | 22 | 22 | 4387 | 428 | 22 |
| 19 | 4843 | 450 | 22 | 17 | 4841 | 450 | 22 | 22 | 4837 | 450 | 22 |
| 14 | 5315 | 472 | 22 | 12 | 5313 | 472 | 22 | 17 | 5309 | 472 | 22 |
| 22 | 5809 | 494 | 22 | 20 | 5807 | 494 | 22 | 16 | 5803 | 494 | 22 |
| 16 | 6325 | 516 | 22 | 14 | 6323 | 516 | 22 | 19 | 6319 | 516 | 22 |
| 23 | 6863 | 538 | 22 | 21 | 6861 | 538 | 22 | 26 | 6857 | 538 | 22 |
| 16 | 7423 | 560 | 22 | 14 | 7421 | 560 | 22 | 19 | 7417 | 560 | 22 |
| 13 | 8005 | 582 | 22 | 11 | 8003 | 582 | 22 | 34 | 7999 | 582 | 22 |
| 23 | 8609 | 604 | 22 | 21 | 8607 | 604 | 22 | 17 | 8603 | 604 | 22 |

- 4 → - 4 → - 2 → - 2 →

| | N22 - N1 | | | | N22 - N2 | | | | N22 - N3 | | |
|---|---|---|---|---|---|---|---|---|---|---|---|
| SD | N1 | N1' | N1" | SD | N2 | N2' | N2" | SD | N3 | N3' | N3" |
| 10 | 19 | | | 8 | 17 | | | 4 | 13 | | |
| 13 | 49 | 30 | | 11 | 47 | 30 | | 7 | 43 | 30 | |
| 2 | 101 | 52 | 22 | 18 | 99 | 52 | 22 | 14 | 95 | 52 | 22 |
| 13 | 175 | 74 | 22 | 11 | 173 | 74 | 22 | 16 | 169 | 74 | 22 |
| 10 | 271 | 96 | 22 | 17 | 269 | 96 | 22 | 13 | 265 | 96 | 22 |
| 20 | 389 | 118 | 22 | 18 | 387 | 118 | 22 | 14 | 383 | 118 | 22 |
| 16 | 529 | 140 | 22 | 14 | 527 | 140 | 22 | 10 | 523 | 140 | 22 |
| 16 | 691 | 162 | 22 | 23 | 689 | 162 | 22 | 19 | 685 | 162 | 22 |
| 20 | 875 | 184 | 22 | 18 | 873 | 184 | 22 | 23 | 869 | 184 | 22 |
| 10 | 1081 | 206 | 22 | 17 | 1079 | 206 | 22 | 13 | 1075 | 206 | 22 |
| 13 | 1309 | 228 | 22 | 11 | 1307 | 228 | 22 | 7 | 1303 | 228 | 22 |
| 20 | 1559 | 250 | 22 | 18 | 1557 | 250 | 22 | 14 | 1553 | 250 | 22 |
| 13 | 1831 | 272 | 22 | 20 | 1829 | 272 | 22 | 16 | 1825 | 272 | 22 |
| 10 | 2125 | 294 | 22 | 8 | 2123 | 294 | 22 | 13 | 2119 | 294 | 22 |
| 11 | 2441 | 316 | 22 | 18 | 2439 | 316 | 22 | 14 | 2435 | 316 | 22 |
| 25 | 2779 | 338 | 22 | 23 | 2777 | 338 | 22 | 19 | 2773 | 338 | 22 |
| 16 | 3139 | 360 | 22 | 14 | 3137 | 360 | 22 | 10 | 3133 | 360 | 22 |
| 11 | 3521 | 382 | 22 | 18 | 3519 | 382 | 22 | 14 | 3515 | 382 | 22 |
| 19 | 3925 | 404 | 22 | 17 | 3923 | 404 | 22 | 22 | 3919 | 404 | 22 |
| 13 | 4351 | 426 | 22 | 20 | 4349 | 426 | 22 | 16 | 4345 | 426 | 22 |
| 29 | 4799 | 448 | 22 | 27 | 4797 | 448 | 22 | 23 | 4793 | 448 | 22 |
| 22 | 5269 | 470 | 22 | 20 | 5267 | 470 | 22 | 16 | 5263 | 470 | 22 |
| 19 | 5761 | 492 | 22 | 26 | 5759 | 492 | 22 | 22 | 5755 | 492 | 22 |
| 20 | 6275 | 514 | 22 | 18 | 6273 | 514 | 22 | 23 | 6269 | 514 | 22 |
| 16 | 6811 | 536 | 22 | 23 | 6809 | 536 | 22 | 19 | 6805 | 536 | 22 |
| 25 | 7369 | 558 | 22 | 23 | 7367 | 558 | 22 | 19 | 7363 | 558 | 22 |
| 29 | 7949 | 580 | 22 | 27 | 7947 | 580 | 22 | 23 | 7943 | 580 | 22 |

- 4 → - 4 → - 2 → - 2 →

| | N22 - O1 | | | | N22 - O2 | | | | N22 - O3 | | |
|---|---|---|---|---|---|---|---|---|---|---|---|
| SD | O1 | O1' | O1" | SD | O2 | O2' | O2" | SD | O3 | O3' | O3" |
| 7 | 25 | | | 5 | 23 | | | 10 | 19 | | |
| 8 | 53 | 28 | | 6 | 51 | 28 | | 11 | 47 | 28 | |
| 4 | 103 | 50 | 22 | 2 | 101 | 50 | 22 | 16 | 97 | 50 | 22 |
| 13 | 175 | 72 | 22 | 11 | 173 | 72 | 22 | 16 | 169 | 72 | 22 |
| 17 | 269 | 94 | 22 | 15 | 267 | 94 | 22 | 11 | 263 | 94 | 22 |
| 16 | 385 | 116 | 22 | 14 | 383 | 116 | 22 | 19 | 379 | 116 | 22 |
| 10 | 523 | 138 | 22 | 8 | 521 | 138 | 22 | 13 | 517 | 138 | 22 |
| 17 | 683 | 160 | 22 | 15 | 681 | 160 | 22 | 20 | 677 | 160 | 22 |
| 19 | 865 | 182 | 22 | 17 | 863 | 182 | 22 | 22 | 859 | 182 | 22 |
| 16 | 1069 | 204 | 22 | 14 | 1067 | 204 | 22 | 10 | 1063 | 204 | 22 |
| 17 | 1295 | 226 | 22 | 15 | 1293 | 226 | 22 | 20 | 1289 | 226 | 22 |
| 13 | 1543 | 248 | 22 | 11 | 1541 | 248 | 22 | 16 | 1537 | 248 | 22 |
| 13 | 1813 | 270 | 22 | 11 | 1811 | 270 | 22 | 16 | 1807 | 270 | 22 |
| 8 | 2105 | 292 | 22 | 6 | 2103 | 292 | 22 | 20 | 2099 | 292 | 22 |
| 16 | 2419 | 314 | 22 | 14 | 2417 | 314 | 22 | 10 | 2413 | 314 | 22 |
| 19 | 2755 | 336 | 22 | 17 | 2753 | 336 | 22 | 22 | 2749 | 336 | 22 |
| 8 | 3113 | 358 | 22 | 6 | 3111 | 358 | 22 | 11 | 3107 | 358 | 22 |
| 19 | 3493 | 380 | 22 | 17 | 3491 | 380 | 22 | 22 | 3487 | 380 | 22 |
| 25 | 3895 | 402 | 22 | 23 | 3893 | 402 | 22 | 28 | 3889 | 402 | 22 |
| 17 | 4319 | 424 | 22 | 15 | 4317 | 424 | 22 | 11 | 4313 | 424 | 22 |
| 22 | 4765 | 446 | 22 | 20 | 4763 | 446 | 22 | 25 | 4759 | 446 | 22 |
| 13 | 5233 | 468 | 22 | 11 | 5231 | 468 | 22 | 16 | 5227 | 468 | 22 |
| 17 | 5723 | 490 | 22 | 15 | 5721 | 490 | 22 | 20 | 5717 | 490 | 22 |
| 16 | 6235 | 512 | 22 | 14 | 6233 | 512 | 22 | 22 | 6229 | 512 | 22 |
| 28 | 6769 | 534 | 22 | 26 | 6767 | 534 | 22 | 22 | 6763 | 534 | 22 |
| 17 | 7325 | 556 | 22 | 15 | 7323 | 556 | 22 | 20 | 7319 | 556 | 22 |
| 19 | 7903 | 578 | 22 | 17 | 7901 | 578 | 22 | 31 | 7897 | 578 | 22 |

- 4 → - 4 → - 2 → - 2 →

| | N22 - P1 | | | | N22 - P2 | | | | N22 - P3 | | |
|---|---|---|---|---|---|---|---|---|---|---|---|
| SD | P1 | P1' | P1" | SD | P2 | P2' | P2" | SD | P3 | P3' | P3" |
| 8 | 17 | | | 6 | 15 | | | 2 | 11 | | |
| 7 | 43 | 26 | | 5 | 41 | 26 | | 10 | 37 | 26 | |
| 10 | 91 | 48 | 22 | 17 | 89 | 48 | 22 | 13 | 85 | 48 | 22 |
| 8 | 161 | 70 | 22 | 15 | 159 | 70 | 22 | 11 | 155 | 70 | 22 |
| 10 | 253 | 92 | 22 | 8 | 251 | 92 | 22 | 13 | 247 | 92 | 22 |
| 16 | 367 | 114 | 22 | 14 | 365 | 114 | 22 | 10 | 361 | 114 | 22 |
| 8 | 503 | 136 | 22 | 6 | 501 | 136 | 22 | 20 | 497 | 136 | 22 |
| 13 | 661 | 158 | 22 | 20 | 659 | 158 | 22 | 16 | 655 | 158 | 22 |
| 13 | 841 | 180 | 22 | 21 | 839 | 180 | 22 | 16 | 835 | 180 | 22 |
| 8 | 1043 | 202 | 22 | 6 | 1041 | 202 | 22 | 11 | 1037 | 202 | 22 |
| 16 | 1267 | 224 | 22 | 14 | 1265 | 224 | 22 | 10 | 1261 | 224 | 22 |
| 10 | 1513 | 246 | 22 | 8 | 1511 | 246 | 22 | 13 | 1507 | 246 | 22 |
| 17 | 1781 | 268 | 22 | 24 | 1779 | 268 | 22 | 20 | 1775 | 268 | 22 |
| 10 | 2071 | 290 | 22 | 17 | 2069 | 290 | 22 | 13 | 2065 | 290 | 22 |
| 16 | 2383 | 312 | 22 | 14 | 2381 | 312 | 22 | 19 | 2377 | 312 | 22 |
| 17 | 2717 | 334 | 22 | 15 | 2715 | 334 | 22 | 11 | 2711 | 334 | 22 |
| 13 | 3073 | 356 | 22 | 11 | 3071 | 356 | 22 | 16 | 3067 | 356 | 22 |
| 13 | 3451 | 378 | 22 | 20 | 3449 | 378 | 22 | 16 | 3445 | 378 | 22 |
| 17 | 3851 | 400 | 22 | 24 | 3849 | 400 | 22 | 20 | 3845 | 400 | 22 |
| 16 | 4273 | 422 | 22 | 14 | 4271 | 422 | 22 | 19 | 4267 | 422 | 22 |
| 19 | 4717 | 444 | 22 | 17 | 4715 | 444 | 22 | 13 | 4711 | 444 | 22 |
| 17 | 5183 | 466 | 22 | 15 | 5181 | 466 | 22 | 20 | 5177 | 466 | 22 |
| 19 | 5671 | 488 | 22 | 26 | 5669 | 488 | 22 | 22 | 5665 | 488 | 22 |
| 16 | 6181 | 510 | 22 | 23 | 6179 | 510 | 22 | 19 | 6175 | 510 | 22 |
| 17 | 6713 | 532 | 22 | 15 | 6711 | 532 | 22 | 20 | 6707 | 532 | 22 |
| 22 | 7267 | 554 | 22 | 20 | 7265 | 554 | 22 | 16 | 7261 | 554 | 22 |
| 22 | 7843 | 576 | 22 | 20 | 7841 | 576 | 22 | 25 | 7837 | 576 | 22 |

- 4 → - 4 → - 2 → - 2 →

| | N22 - Q1 | | | | N22 - Q2 | | | | N22 - Q3 | | |
|---|---|---|---|---|---|---|---|---|---|---|---|
| SD | Q1 | Q1' | Q1" | SD | Q2 | Q2' | Q2" | SD | Q3 | Q3' | Q3" |
| 10 | 19 | | | 8 | 17 | | | 4 | 13 | | |
| 7 | 43 | 24 | | 5 | 41 | 24 | | 10 | 37 | 24 | |
| 17 | 89 | 46 | 22 | 15 | 87 | 46 | 22 | 11 | 83 | 46 | 22 |
| 13 | 157 | 68 | 22 | 11 | 155 | 68 | 22 | 7 | 151 | 68 | 22 |
| 13 | 247 | 90 | 22 | 11 | 245 | 90 | 22 | 7 | 241 | 90 | 22 |
| 17 | 359 | 112 | 22 | 15 | 357 | 112 | 22 | 11 | 353 | 112 | 22 |
| 16 | 493 | 134 | 22 | 14 | 491 | 134 | 22 | 19 | 487 | 134 | 22 |
| 19 | 649 | 156 | 22 | 17 | 647 | 156 | 22 | 13 | 643 | 156 | 22 |
| 17 | 827 | 178 | 22 | 15 | 825 | 178 | 22 | 11 | 821 | 178 | 22 |
| 10 | 1027 | 200 | 22 | 8 | 1025 | 200 | 22 | 4 | 1021 | 200 | 22 |
| 16 | 1249 | 222 | 22 | 14 | 1247 | 222 | 22 | 10 | 1243 | 222 | 22 |
| 17 | 1493 | 244 | 22 | 15 | 1491 | 244 | 22 | 20 | 1487 | 244 | 22 |
| 22 | 1759 | 266 | 22 | 20 | 1757 | 266 | 22 | 16 | 1753 | 266 | 22 |
| 13 | 2047 | 288 | 22 | 11 | 2045 | 288 | 22 | 7 | 2041 | 288 | 22 |
| 17 | 2357 | 310 | 22 | 15 | 2355 | 310 | 22 | 11 | 2351 | 310 | 22 |
| 25 | 2689 | 332 | 22 | 23 | 2687 | 332 | 22 | 19 | 2683 | 332 | 22 |
| 10 | 3043 | 354 | 22 | 8 | 3041 | 354 | 22 | 13 | 3037 | 354 | 22 |
| 17 | 3419 | 376 | 22 | 15 | 3417 | 376 | 22 | 11 | 3413 | 376 | 22 |
| 19 | 3817 | 398 | 22 | 17 | 3815 | 398 | 22 | 13 | 3811 | 398 | 22 |
| 16 | 4237 | 420 | 22 | 14 | 4235 | 420 | 22 | 10 | 4231 | 420 | 22 |
| 26 | 4679 | 442 | 22 | 24 | 4677 | 442 | 22 | 20 | 4673 | 442 | 22 |
| 13 | 5143 | 464 | 22 | 11 | 5141 | 464 | 22 | 16 | 5137 | 464 | 22 |
| 22 | 5629 | 486 | 22 | 20 | 5627 | 486 | 22 | 16 | 5623 | 486 | 22 |
| 17 | 6137 | 508 | 22 | 15 | 6135 | 508 | 22 | 11 | 6131 | 508 | 22 |
| 25 | 6667 | 530 | 22 | 23 | 6665 | 530 | 22 | 19 | 6661 | 530 | 22 |
| 19 | 7219 | 552 | 22 | 17 | 7217 | 552 | 22 | 13 | 7213 | 552 | 22 |
| 26 | 7793 | 574 | 22 | 24 | 7791 | 574 | 22 | 29 | 7787 | 574 | 22 |

- 4 → - 4 → - 2 → - 2 →

| | N22 - R1 | | | | N22 - R2 | | | | N22 - R3 | | |
|---|---|---|---|---|---|---|---|---|---|---|---|
| SD | R1 | R1' | R1" | SD | R2 | R2' | R2" | SD | R3 | R3' | R3" |
| 4 | 13 | | | 2 | 11 | | | 7 | 7 | | |
| 8 | 35 | 22 | | 6 | 33 | 22 | | 11 | 29 | 22 | |
| 16 | 79 | 44 | 22 | 14 | 77 | 44 | 22 | 10 | 73 | 44 | 22 |
| 10 | 145 | 66 | 22 | 8 | 143 | 66 | 22 | 13 | 139 | 66 | 22 |
| 8 | 233 | 88 | 22 | 6 | 231 | 88 | 22 | 11 | 227 | 88 | 22 |
| 10 | 343 | 110 | 22 | 8 | 341 | 110 | 22 | 13 | 337 | 110 | 22 |
| 16 | 475 | 132 | 22 | 14 | 473 | 132 | 22 | 19 | 469 | 132 | 22 |
| 17 | 629 | 154 | 22 | 15 | 627 | 154 | 22 | 11 | 623 | 154 | 22 |
| 13 | 805 | 176 | 22 | 11 | 803 | 176 | 22 | 25 | 799 | 176 | 22 |
| 4 | 1003 | 198 | 22 | 2 | 1001 | 198 | 22 | 25 | 997 | 198 | 22 |
| 8 | 1223 | 220 | 22 | 6 | 1221 | 220 | 22 | 11 | 1217 | 220 | 22 |
| 16 | 1465 | 242 | 22 | 14 | 1463 | 242 | 22 | 19 | 1459 | 242 | 22 |
| 19 | 1729 | 264 | 22 | 17 | 1727 | 264 | 22 | 13 | 1723 | 264 | 22 |
| 8 | 2015 | 286 | 22 | 6 | 2013 | 286 | 22 | 11 | 2009 | 286 | 22 |
| 10 | 2323 | 308 | 22 | 8 | 2321 | 308 | 22 | 13 | 2317 | 308 | 22 |
| 16 | 2653 | 330 | 22 | 14 | 2651 | 330 | 22 | 19 | 2647 | 330 | 22 |
| 8 | 3005 | 352 | 22 | 6 | 3003 | 352 | 22 | 29 | 2999 | 352 | 22 |
| 22 | 3379 | 374 | 22 | 20 | 3377 | 374 | 22 | 16 | 3373 | 374 | 22 |
| 22 | 3775 | 396 | 22 | 20 | 3773 | 396 | 22 | 25 | 3769 | 396 | 22 |
| 17 | 4193 | 418 | 22 | 15 | 4191 | 418 | 22 | 20 | 4187 | 418 | 22 |
| 16 | 4633 | 440 | 22 | 14 | 4631 | 440 | 22 | 19 | 4627 | 440 | 22 |
| 19 | 5095 | 462 | 22 | 17 | 5093 | 462 | 22 | 22 | 5089 | 462 | 22 |
| 26 | 5579 | 484 | 22 | 24 | 5577 | 484 | 22 | 20 | 5573 | 484 | 22 |
| 19 | 6085 | 506 | 22 | 17 | 6083 | 506 | 22 | 22 | 6079 | 506 | 22 |
| 16 | 6613 | 528 | 22 | 14 | 6611 | 528 | 22 | 19 | 6607 | 528 | 22 |
| 17 | 7163 | 550 | 22 | 15 | 7161 | 550 | 22 | 20 | 7157 | 550 | 22 |
| 22 | 7735 | 572 | 22 | 20 | 7733 | 572 | 22 | 25 | 7729 | 572 | 22 |

- 4 → - 4 → - 2 → - 2 →

| | N22 - S1 | | | | N22 - S2 | | | | N22 - S3 | | |
|---|---|---|---|---|---|---|---|---|---|---|---|
| SD | S1 | S1' | S1" | SD | S2 | S2' | S2" | SD | S3 | S3' | S3" |
| 2 | 11 | | | 9 | 9 | | | 5 | 5 | | |
| 4 | 31 | 20 | | 11 | 29 | 20 | | 7 | 25 | 20 | |
| 10 | 73 | 42 | 22 | 8 | 71 | 42 | 22 | 13 | 67 | 42 | 22 |
| 11 | 137 | 64 | 22 | 9 | 135 | 64 | 22 | 5 | 131 | 64 | 22 |
| 7 | 223 | 86 | 22 | 5 | 221 | 86 | 22 | 10 | 217 | 86 | 22 |
| 7 | 331 | 108 | 22 | 14 | 329 | 108 | 22 | 10 | 325 | 108 | 22 |
| 11 | 461 | 130 | 22 | 18 | 459 | 130 | 22 | 14 | 455 | 130 | 22 |
| 10 | 613 | 152 | 22 | 8 | 611 | 152 | 22 | 13 | 607 | 152 | 22 |
| 22 | 787 | 174 | 22 | 20 | 785 | 174 | 22 | 16 | 781 | 174 | 22 |
| 20 | 983 | 196 | 22 | 18 | 981 | 196 | 22 | 23 | 977 | 196 | 22 |
| 4 | 1201 | 218 | 22 | 20 | 1199 | 218 | 22 | 16 | 1195 | 218 | 22 |
| 10 | 1441 | 240 | 22 | 17 | 1439 | 240 | 22 | 13 | 1435 | 240 | 22 |
| 11 | 1703 | 262 | 22 | 9 | 1701 | 262 | 22 | 23 | 1697 | 262 | 22 |
| 25 | 1987 | 284 | 22 | 23 | 1985 | 284 | 22 | 19 | 1981 | 284 | 22 |
| 16 | 2293 | 306 | 22 | 14 | 2291 | 306 | 22 | 19 | 2287 | 306 | 22 |
| 11 | 2621 | 328 | 22 | 18 | 2619 | 328 | 22 | 14 | 2615 | 328 | 22 |
| 19 | 2971 | 350 | 22 | 26 | 2969 | 350 | 22 | 22 | 2965 | 350 | 22 |
| 13 | 3343 | 372 | 22 | 11 | 3341 | 372 | 22 | 16 | 3337 | 372 | 22 |
| 20 | 3737 | 394 | 22 | 18 | 3735 | 394 | 22 | 14 | 3731 | 394 | 22 |
| 13 | 4153 | 416 | 22 | 11 | 4151 | 416 | 22 | 16 | 4147 | 416 | 22 |
| 19 | 4591 | 438 | 22 | 26 | 4589 | 438 | 22 | 22 | 4585 | 438 | 22 |
| 11 | 5051 | 460 | 22 | 18 | 5049 | 460 | 22 | 14 | 5045 | 460 | 22 |
| 16 | 5533 | 482 | 22 | 14 | 5531 | 482 | 22 | 20 | 5527 | 482 | 22 |
| 16 | 6037 | 504 | 22 | 14 | 6035 | 504 | 22 | 10 | 6031 | 504 | 22 |
| 20 | 6563 | 526 | 22 | 18 | 6561 | 526 | 22 | 23 | 6557 | 526 | 22 |
| 10 | 7111 | 548 | 22 | 17 | 7109 | 548 | 22 | 13 | 7105 | 548 | 22 |
| 22 | 7681 | 570 | 22 | 29 | 7679 | 570 | 22 | 25 | 7675 | 570 | 22 |

- 4 → - 4 → - 2 → - 2 →

| | N22 - T1 | | | | N22 - T2 | | | | N22 - T3 | | |
|---|---|---|---|---|---|---|---|---|---|---|---|
| SD | T1 | T1' | T1" | SD | T2 | T2' | T2" | SD | T3 | T3' | T3" |
| 4 | 13 | | | 2 | 11 | | | 7 | 7 | | |
| 4 | 31 | 18 | | 11 | 29 | 18 | | 7 | 25 | 18 | |
| 8 | 71 | 40 | 22 | 15 | 69 | 40 | 22 | 11 | 65 | 40 | 22 |
| 7 | 133 | 62 | 22 | 5 | 131 | 62 | 22 | 10 | 127 | 62 | 22 |
| 10 | 217 | 84 | 22 | 8 | 215 | 84 | 22 | 4 | 211 | 84 | 22 |
| 8 | 323 | 106 | 22 | 6 | 321 | 106 | 22 | 11 | 317 | 106 | 22 |
| 10 | 451 | 128 | 22 | 17 | 449 | 128 | 22 | 13 | 445 | 128 | 22 |
| 7 | 601 | 150 | 22 | 23 | 599 | 150 | 22 | 19 | 595 | 150 | 22 |
| 17 | 773 | 172 | 22 | 15 | 771 | 172 | 22 | 20 | 767 | 172 | 22 |
| 22 | 967 | 194 | 22 | 20 | 965 | 194 | 22 | 16 | 961 | 194 | 22 |
| 13 | 1183 | 216 | 22 | 11 | 1181 | 216 | 22 | 16 | 1177 | 216 | 22 |
| 8 | 1421 | 238 | 22 | 15 | 1419 | 238 | 22 | 11 | 1415 | 238 | 22 |
| 16 | 1681 | 260 | 22 | 23 | 1679 | 260 | 22 | 19 | 1675 | 260 | 22 |
| 19 | 1963 | 282 | 22 | 17 | 1961 | 282 | 22 | 22 | 1957 | 282 | 22 |
| 17 | 2267 | 304 | 22 | 15 | 2265 | 304 | 22 | 11 | 2261 | 304 | 22 |
| 19 | 2593 | 326 | 22 | 17 | 2591 | 326 | 22 | 22 | 2587 | 326 | 22 |
| 16 | 2941 | 348 | 22 | 23 | 2939 | 348 | 22 | 19 | 2935 | 348 | 22 |
| 8 | 3311 | 370 | 22 | 15 | 3309 | 370 | 22 | 11 | 3305 | 370 | 22 |
| 13 | 3703 | 392 | 22 | 11 | 3701 | 392 | 22 | 25 | 3697 | 392 | 22 |
| 13 | 4117 | 414 | 22 | 11 | 4115 | 414 | 22 | 7 | 4111 | 414 | 22 |
| 17 | 4553 | 436 | 22 | 15 | 4551 | 436 | 22 | 20 | 4547 | 436 | 22 |
| 7 | 5011 | 458 | 22 | 14 | 5009 | 458 | 22 | 10 | 5005 | 458 | 22 |
| 19 | 5491 | 480 | 22 | 26 | 5489 | 480 | 22 | 22 | 5485 | 480 | 22 |
| 26 | 5993 | 502 | 22 | 24 | 5991 | 502 | 22 | 29 | 5987 | 502 | 22 |
| 19 | 6517 | 524 | 22 | 17 | 6515 | 524 | 22 | 13 | 6511 | 524 | 22 |
| 16 | 7063 | 546 | 22 | 14 | 7061 | 546 | 22 | 19 | 7057 | 546 | 22 |
| 17 | 7631 | 568 | 22 | 24 | 7629 | 568 | 22 | 20 | 7625 | 568 | 22 |

- 4 → - 4 → - 2 → - 2 →

**Table 6-B2 :** Quadratic Polynomials of the Spiral-Graphs belonging to the "Prime-Number-Spiral-Systems" **P20-D** to **P20-I** and **N20-D** to **N20-I** ( with the **2. Differential = 20** )

| Spiral Graph System | Spiral Graph | Number Sequence of Spiral Graph | Quadratic Polynomial 1 ( calculated with the first 3 numbers of the given sequence ) | Quadratic Polynomial 2 ( calculated with 3 numbers starting with the **2.** Number of the sequence ) | Quadratic Polynomial 3 ( calculated with 3 numbers starting with the **3.** Number of the sequence ) | Quadratic Polynomial 4 ( calculated with 3 numbers starting with the **4.** Number of the sequence ) |
|---|---|---|---|---|---|---|
| | D3 | 1 , 21 , 61 , 121 , 201 , 301 ,...... | $f_1(x) = 10x^2 - 10x + 1$ | $f_2(x) = 10x^2 + 10x + 1$ | $f_3(x) = 10x^2 + 30x + 21$ | $f_4(x) = 10x^2 + 50x + 61$ |
| N20-D | D2 | 5 , 25 , 65 , 125 , 205 , 305 ,...... | $f_1(x) = 10x^2 - 10x + 5$ | $f_2(x) = 10x^2 + 10x + 5$ | $f_3(x) = 10x^2 + 30x + 25$ | $f_4(x) = 10x^2 + 50x + 65$ |
| | D1 | 7 , 27 , 67 , 127 , 207 , 307 ,...... | $f_1(x) = 10x^2 - 10x + 7$ | $f_2(x) = 10x^2 + 10x + 7$ | $f_3(x) = 10x^2 + 30x + 27$ | $f_4(x) = 10x^2 + 50x + 67$ |
| | D1 | 11 , 31 , 71 , 131 , 211 , 311 ,...... | $f_1(x) = 10x^2 - 10x + 11$ | $f_2(x) = 10x^2 + 10x + 11$ | $f_3(x) = 10x^2 + 30x + 31$ | $f_4(x) = 10x^2 + 50x + 71$ |
| P20-D | D2 | 13 , 33 , 73 , 133 , 213 , 313 ,...... | $f_1(x) = 10x^2 - 10x + 13$ | $f_2(x) = 10x^2 + 10x + 13$ | $f_3(x) = 10x^2 + 30x + 33$ | $f_4(x) = 10x^2 + 50x + 73$ |
| | D3 | 17 , 37 , 77 , 137 , 217 , 317 ,...... | $f_1(x) = 10x^2 - 10x + 17$ | $f_2(x) = 10x^2 + 10x + 17$ | $f_3(x) = 10x^2 + 30x + 37$ | $f_4(x) = 10x^2 + 50x + 77$ |
| | E3 | 1 , 25 , 69 , 133 , 217 , 321 ,...... | $f_1(x) = 10x^2 - 6x - 3$ | $f_2(x) = 10x^2 + 14x + 1$ | $f_3(x) = 10x^2 + 34x + 25$ | $f_4(x) = 10x^2 + 54x + 69$ |
| N20-E | E2 | 1 , 5 , 29 , 73 , 137 , 221 ,...... | $f_1(x) = 10x^2 - 26x + 17$ | $f_2(x) = 10x^2 - 6x + 1$ | $f_3(x) = 10x^2 + 14x + 5$ | $f_4(x) = 10x^2 + 34x + 29$ |
| | E1 | 3 , 7 , 31 , 75 , 139 , 223 ,...... | $f_1(x) = 10x^2 - 26x + 19$ | $f_2(x) = 10x^2 - 6x + 3$ | $f_3(x) = 10x^2 + 14x + 7$ | $f_4(x) = 10x^2 + 34x + 31$ |
| | E1 | 7 , 11 , 35 , 79 , 143 , 227 ,...... | $f_1(x) = 10x^2 - 26x + 23$ | $f_2(x) = 10x^2 - 6x + 7$ | $f_3(x) = 10x^2 + 14x + 11$ | $f_4(x) = 10x^2 + 34x + 35$ |
| P20-E | E2 | 9 , 13 , 37 , 81 , 145 , 229 ,...... | $f_1(x) = 10x^2 - 26x + 25$ | $f_2(x) = 10x^2 - 6x + 9$ | $f_3(x) = 10x^2 + 14x + 13$ | $f_4(x) = 10x^2 + 34x + 37$ |
| | E3 | 13 , 17 , 41 , 85 , 149 , 233 ,...... | $f_1(x) = 10x^2 - 26x + 29$ | $f_2(x) = 10x^2 - 6x + 13$ | $f_3(x) = 10x^2 + 14x + 17$ | $f_4(x) = 10x^2 + 34x + 41$ |
| | F3 | 1 , 27 , 73 , 139 , 225 , 331 ,...... | $f_1(x) = 10x^2 - 4x - 5$ | $f_2(x) = 10x^2 + 16x + 1$ | $f_3(x) = 10x^2 + 36x + 27$ | $f_4(x) = 10x^2 + 56x + 73$ |
| N20-F | F2 | 5 , 31 , 77 , 143 , 229 , 335 ,...... | $f_1(x) = 10x^2 - 4x - 1$ | $f_2(x) = 10x^2 + 16x + 5$ | $f_3(x) = 10x^2 + 36x + 31$ | $f_4(x) = 10x^2 + 56x + 77$ |
| | F1 | 7 , 33 , 79 , 145 , 231 , 337 ,...... | $f_1(x) = 10x^2 - 4x + 1$ | $f_2(x) = 10x^2 + 16x + 7$ | $f_3(x) = 10x^2 + 36x + 33$ | $f_4(x) = 10x^2 + 56x + 79$ |
| | F1 | 11 , 37 , 83 , 149 , 235 , 341 ,...... | $f_1(x) = 10x^2 - 4x + 5$ | $f_2(x) = 10x^2 + 16x + 11$ | $f_3(x) = 10x^2 + 36x + 37$ | $f_4(x) = 10x^2 + 56x + 83$ |
| P20-F | F2 | 13 , 39 , 85 , 151 , 237 , 343 ,...... | $f_1(x) = 10x^2 - 4x + 7$ | $f_2(x) = 10x^2 + 16x + 13$ | $f_3(x) = 10x^2 + 36x + 39$ | $f_4(x) = 10x^2 + 56x + 85$ |
| | F3 | 17 , 43 , 89 , 155 , 241 , 347 ,...... | $f_1(x) = 10x^2 - 4x + 11$ | $f_2(x) = 10x^2 + 16x + 17$ | $f_3(x) = 10x^2 + 36x + 43$ | $f_4(x) = 10x^2 + 56x + 89$ |
| | G3 | 3 , 13 , 43 , 93 , 163 , 253 ,...... | $f_1(x) = 10x^2 - 20x + 13$ | $f_2(x) = 10x^2 + 0x + 3$ | $f_3(x) = 10x^2 + 20x + 13$ | $f_4(x) = 10x^2 + 40x + 43$ |
| N20-G | G2 | 7 , 17 , 47 , 97 , 167 , 257 ,...... | $f_1(x) = 10x^2 - 20x + 17$ | $f_2(x) = 10x^2 + 0x + 7$ | $f_3(x) = 10x^2 + 20x + 17$ | $f_4(x) = 10x^2 + 40x + 47$ |
| | G1 | 9 , 19 , 49 , 99 , 169 , 259 ,...... | $f_1(x) = 10x^2 - 20x + 19$ | $f_2(x) = 10x^2 + 0x + 9$ | $f_3(x) = 10x^2 + 20x + 19$ | $f_4(x) = 10x^2 + 40x + 49$ |
| | G1 | 13 , 23 , 53 , 103 , 173 , 263 ,...... | $f_1(x) = 10x^2 - 20x + 23$ | $f_2(x) = 10x^2 + 0x + 13$ | $f_3(x) = 10x^2 + 20x + 23$ | $f_4(x) = 10x^2 + 40x + 53$ |
| P20-G | G2 | 15 , 25 , 55 , 105 , 175 , 265 ,...... | $f_1(x) = 10x^2 - 20x + 25$ | $f_2(x) = 10x^2 + 0x + 15$ | $f_3(x) = 10x^2 + 20x + 25$ | $f_4(x) = 10x^2 + 40x + 55$ |
| | G3 | 19 , 29 , 59 , 109 , 179 , 269 ,...... | $f_1(x) = 10x^2 - 20x + 29$ | $f_2(x) = 10x^2 + 0x + 19$ | $f_3(x) = 10x^2 + 20x + 29$ | $f_4(x) = 10x^2 + 40x + 59$ |
| | H3 | 1 , 15 , 49 , 103 , 177 , 271 ,...... | $f_1(x) = 10x^2 - 16x + 7$ | $f_2(x) = 10x^2 + 4x + 1$ | $f_3(x) = 10x^2 + 24x + 15$ | $f_4(x) = 10x^2 + 44x + 49$ |
| N20-H | H2 | 5 , 19 , 53 , 107 , 181 , 275 ,...... | $f_1(x) = 10x^2 - 16x + 11$ | $f_2(x) = 10x^2 + 4x + 5$ | $f_3(x) = 10x^2 + 24x + 19$ | $f_4(x) = 10x^2 + 44x + 53$ |
| | H1 | 7 , 21 , 55 , 109 , 183 , 277 ,...... | $f_1(x) = 10x^2 - 16x + 13$ | $f_2(x) = 10x^2 + 4x + 7$ | $f_3(x) = 10x^2 + 24x + 21$ | $f_4(x) = 10x^2 + 44x + 55$ |
| | H1 | 11 , 25 , 59 , 113 , 187 , 281 ,...... | $f_1(x) = 10x^2 - 16x + 17$ | $f_2(x) = 10x^2 + 4x + 11$ | $f_3(x) = 10x^2 + 24x + 25$ | $f_4(x) = 10x^2 + 44x + 59$ |
| P20-H | H2 | 13 , 27 , 61 , 115 , 189 , 283 ,...... | $f_1(x) = 10x^2 - 16x + 19$ | $f_2(x) = 10x^2 + 4x + 13$ | $f_3(x) = 10x^2 + 24x + 27$ | $f_4(x) = 10x^2 + 44x + 61$ |
| | H3 | 17 , 31 , 65 , 119 , 193 , 287 ,...... | $f_1(x) = 10x^2 - 16x + 23$ | $f_2(x) = 10x^2 + 4x + 17$ | $f_3(x) = 10x^2 + 24x + 31$ | $f_4(x) = 10x^2 + 44x + 65$ |
| | I3 | 3 , 19 , 55 , 111 , 187 , 283 ,...... | $f_1(x) = 10x^2 - 14x + 7$ | $f_2(x) = 10x^2 + 6x + 3$ | $f_3(x) = 10x^2 + 26x + 19$ | $f_4(x) = 10x^2 + 46x + 55$ |
| N20-I | I2 | 7 , 23 , 59 , 115 , 191 , 287 ,...... | $f_1(x) = 10x^2 - 14x + 11$ | $f_2(x) = 10x^2 + 6x + 7$ | $f_3(x) = 10x^2 + 26x + 23$ | $f_4(x) = 10x^2 + 46x + 59$ |
| | I1 | 9 , 25 , 61 , 117 , 193 , 289 ,...... | $f_1(x) = 10x^2 - 14x + 13$ | $f_2(x) = 10x^2 + 6x + 9$ | $f_3(x) = 10x^2 + 26x + 25$ | $f_4(x) = 10x^2 + 46x + 61$ |
| | I1 | 13 , 29 , 65 , 121 , 197 , 293 ,...... | $f_1(x) = 10x^2 - 14x + 17$ | $f_2(x) = 10x^2 + 6x + 13$ | $f_3(x) = 10x^2 + 26x + 29$ | $f_4(x) = 10x^2 + 46x + 65$ |
| P20-I | I2 | 15 , 31 , 67 , 123 , 199 , 295 ,...... | $f_1(x) = 10x^2 - 14x + 19$ | $f_2(x) = 10x^2 + 6x + 15$ | $f_3(x) = 10x^2 + 26x + 31$ | $f_4(x) = 10x^2 + 46x + 67$ |
| | I3 | 19 , 35 , 71 , 127 , 203 , 299 ,...... | $f_1(x) = 10x^2 - 14x + 23$ | $f_2(x) = 10x^2 + 6x + 19$ | $f_3(x) = 10x^2 + 26x + 35$ | $f_4(x) = 10x^2 + 46x + 71$ |

**Table 6-C2 :** Quadratic Polynomials of the Spiral-Graphs belonging to the "Prime-Number-Spiral-Systems" **N22-J** to **N22-T** ( with the **2. Differential = 22** )

| Spiral Graph System | Spiral Graph | Number Sequence of Spiral Graph | Quadratic Polynomial 1 ( calculated with the first 3 numbers of the given sequence ) | Quadratic Polynomial 2 ( calculated with 3 numbers starting with the **2.** Number of the sequence ) | Quadratic Polynomial 3 ( calculated with 3 numbers starting with the **3.** Number of the sequence ) | Quadratic Polynomial 4 ( calculated with 3 numbers starting with the **4.** Number of the sequence ) |
|---|---|---|---|---|---|---|
| **N22-J** | J1 | 15 , 31 , 69 , 129 , 211 , 315 ,...... | $f_1(x) = 11x^2 - 17x + 21$ | $f_2(x) = 11x^2 + 5x + 15$ | $f_3(x) = 11x^2 + 27x + 31$ | $f_4(x) = 11x^2 + 49x + 69$ |
| | J2 | 13 , 29 , 67 , 127 , 209 , 313 ,...... | $f_1(x) = 11x^2 - 17x + 19$ | $f_2(x) = 11x^2 + 5x + 13$ | $f_3(x) = 11x^2 + 27x + 29$ | $f_4(x) = 11x^2 + 49x + 67$ |
| | J3 | 9 , 25 , 63 , 123 , 205 , 309 ,...... | $f_1(x) = 11x^2 - 17x + 15$ | $f_2(x) = 11x^2 + 5x + 9$ | $f_3(x) = 11x^2 + 27x + 25$ | $f_4(x) = 11x^2 + 49x + 63$ |
| **N22-K** | K1 | 11 , 25 , 61 , 119 , 199 , 301 ,...... | $f_1(x) = 11x^2 - 19x + 19$ | $f_2(x) = 11x^2 + 3x + 11$ | $f_3(x) = 11x^2 + 25x + 25$ | $f_4(x) = 11x^2 + 47x + 61$ |
| | K2 | 9 , 23 , 59 , 117 , 197 , 299 ,...... | $f_1(x) = 11x^2 - 19x + 17$ | $f_2(x) = 11x^2 + 3x + 9$ | $f_3(x) = 11x^2 + 25x + 23$ | $f_4(x) = 11x^2 + 47x + 59$ |
| | K3 | 5 , 19 , 55 , 113 , 193 , 295 ,...... | $f_1(x) = 11x^2 - 19x + 13$ | $f_2(x) = 11x^2 + 3x + 5$ | $f_3(x) = 11x^2 + 25x + 19$ | $f_4(x) = 11x^2 + 47x + 55$ |
| **N22-L** | L1 | 13 , 25 , 59 , 115 , 193 , 293 ,...... | $f_1(x) = 11x^2 - 21x + 23$ | $f_2(x) = 11x^2 + 1x + 13$ | $f_3(x) = 11x^2 + 23x + 25$ | $f_4(x) = 11x^2 + 45x + 59$ |
| | L2 | 11 , 23 , 57 , 113 , 191 , 291 ,...... | $f_1(x) = 11x^2 - 21x + 21$ | $f_2(x) = 11x^2 + 1x + 11$ | $f_3(x) = 11x^2 + 23x + 23$ | $f_4(x) = 11x^2 + 45x + 57$ |
| | L3 | 7 , 19 , 53 , 109 , 187 , 287 ,...... | $f_1(x) = 11x^2 - 21x + 17$ | $f_2(x) = 11x^2 + 1x + 7$ | $f_3(x) = 11x^2 + 23x + 19$ | $f_4(x) = 11x^2 + 45x + 53$ |
| **N22-M** | M1 | 13 , 23 , 55 , 109 , 185 , 283 ,...... | $f_1(x) = 11x^2 - 23x + 25$ | $f_2(x) = 11x^2 - 1x + 13$ | $f_3(x) = 11x^2 + 21x + 23$ | $f_4(x) = 11x^2 + 43x + 55$ |
| | M2 | 11 , 21 , 53 , 107 , 183 , 281 ,...... | $f_1(x) = 11x^2 - 23x + 23$ | $f_2(x) = 11x^2 - 1x + 11$ | $f_3(x) = 11x^2 + 21x + 21$ | $f_4(x) = 11x^2 + 43x + 53$ |
| | M3 | 7 , 17 , 49 , 103 , 179 , 277 ,...... | $f_1(x) = 11x^2 - 23x + 19$ | $f_2(x) = 11x^2 - 1x + 7$ | $f_3(x) = 11x^2 + 21x + 17$ | $f_4(x) = 11x^2 + 43x + 49$ |
| **N22-N** | N1 | 11 , 19 , 49 , 101 , 175 , 271 ,...... | $f_1(x) = 11x^2 - 25x + 25$ | $f_2(x) = 11x^2 - 3x + 11$ | $f_3(x) = 11x^2 + 19x + 19$ | $f_4(x) = 11x^2 + 41x + 49$ |
| | N2 | 9 , 17 , 47 , 99 , 173 , 269 ,...... | $f_1(x) = 11x^2 - 25x + 23$ | $f_2(x) = 11x^2 - 3x + 9$ | $f_3(x) = 11x^2 + 19x + 17$ | $f_4(x) = 11x^2 + 41x + 47$ |
| | N3 | 5 , 13 , 43 , 95 , 169 , 265 ,...... | $f_1(x) = 11x^2 - 25x + 19$ | $f_2(x) = 11x^2 - 3x + 5$ | $f_3(x) = 11x^2 + 19x + 13$ | $f_4(x) = 11x^2 + 41x + 43$ |
| **N22-O** | O1 | 19 , 25 , 53 , 103 , 175 , 269 ,...... | $f_1(x) = 11x^2 - 27x + 35$ | $f_2(x) = 11x^2 - 5x + 19$ | $f_3(x) = 11x^2 + 17x + 25$ | $f_4(x) = 11x^2 + 39x + 53$ |
| | O2 | 17 , 23 , 51 , 101 , 173 , 267 ,...... | $f_1(x) = 11x^2 - 27x + 33$ | $f_2(x) = 11x^2 - 5x + 17$ | $f_3(x) = 11x^2 + 17x + 23$ | $f_4(x) = 11x^2 + 39x + 51$ |
| | O3 | 13 , 19 , 47 , 97 , 169 , 263 ,...... | $f_1(x) = 11x^2 - 27x + 29$ | $f_2(x) = 11x^2 - 5x + 13$ | $f_3(x) = 11x^2 + 17x + 19$ | $f_4(x) = 11x^2 + 39x + 47$ |
| **N22-P** | P1 | 13 , 17 , 43 , 91 , 161 , 253 ,...... | $f_1(x) = 11x^2 - 29x + 31$ | $f_2(x) = 11x^2 - 7x + 13$ | $f_3(x) = 11x^2 + 15x + 17$ | $f_4(x) = 11x^2 + 37x + 43$ |
| | P2 | 11 , 15 , 41 , 89 , 159 , 251 ,...... | $f_1(x) = 11x^2 - 29x + 29$ | $f_2(x) = 11x^2 - 7x + 11$ | $f_3(x) = 11x^2 + 15x + 15$ | $f_4(x) = 11x^2 + 37x + 41$ |
| | P3 | 7 , 11 , 37 , 85 , 155 , 247 ,...... | $f_1(x) = 11x^2 - 29x + 25$ | $f_2(x) = 11x^2 - 7x + 7$ | $f_3(x) = 11x^2 + 15x + 11$ | $f_4(x) = 11x^2 + 37x + 37$ |
| **N22-Q** | Q1 | 17 , 19 , 43 , 89 , 157 , 247 ,...... | $f_1(x) = 11x^2 - 31x + 37$ | $f_2(x) = 11x^2 - 9x + 17$ | $f_3(x) = 11x^2 + 13x + 19$ | $f_4(x) = 11x^2 + 35x + 43$ |
| | Q2 | 15 , 17 , 41 , 87 , 155 , 245 ,...... | $f_1(x) = 11x^2 - 31x + 35$ | $f_2(x) = 11x^2 - 9x + 15$ | $f_3(x) = 11x^2 + 13x + 17$ | $f_4(x) = 11x^2 + 35x + 41$ |
| | Q3 | 11 , 13 , 37 , 83 , 151 , 241 ,...... | $f_1(x) = 11x^2 - 31x + 31$ | $f_2(x) = 11x^2 - 9x + 11$ | $f_3(x) = 11x^2 + 13x + 13$ | $f_4(x) = 11x^2 + 35x + 37$ |
| **N22-R** | R1 | 13 , 35 , 79 , 145 , 233 , 343 ,...... | $f_1(x) = 11x^2 - 11x + 13$ | $f_2(x) = 11x^2 + 11x + 13$ | $f_3(x) = 11x^2 + 33x + 35$ | $f_4(x) = 11x^2 + 55x + 79$ |
| | R2 | 11 , 33 , 77 , 143 , 231 , 341 ,...... | $f_1(x) = 11x^2 - 11x + 11$ | $f_2(x) = 11x^2 + 11x + 11$ | $f_3(x) = 11x^2 + 33x + 33$ | $f_4(x) = 11x^2 + 55x + 77$ |
| | R3 | 7 , 29 , 73 , 139 , 227 , 337 ,...... | $f_1(x) = 11x^2 - 11x + 7$ | $f_2(x) = 11x^2 + 11x + 7$ | $f_3(x) = 11x^2 + 33x + 29$ | $f_4(x) = 11x^2 + 55x + 73$ |
| **N22-S** | S1 | 11 , 31 , 73 , 137 , 223 , 331 ,...... | $f_1(x) = 11x^2 - 13x + 13$ | $f_2(x) = 11x^2 + 9x + 11$ | $f_3(x) = 11x^2 + 31x + 31$ | $f_4(x) = 11x^2 + 53x + 73$ |
| | S2 | 9 , 29 , 71 , 135 , 221 , 329 ,...... | $f_1(x) = 11x^2 - 13x + 11$ | $f_2(x) = 11x^2 + 9x + 9$ | $f_3(x) = 11x^2 + 31x + 29$ | $f_4(x) = 11x^2 + 53x + 71$ |
| | S3 | 5 , 25 , 67 , 131 , 217 , 325 ,...... | $f_1(x) = 11x^2 - 13x + 7$ | $f_2(x) = 11x^2 + 9x + 5$ | $f_3(x) = 11x^2 + 31x + 25$ | $f_4(x) = 11x^2 + 53x + 67$ |
| **N22-T** | T1 | 13 , 31 , 71 , 133 , 217 , 323 ,...... | $f_1(x) = 11x^2 - 15x + 17$ | $f_2(x) = 11x^2 + 7x + 13$ | $f_3(x) = 11x^2 + 29x + 31$ | $f_4(x) = 11x^2 + 51x + 71$ |
| | T2 | 11 , 29 , 69 , 131 , 215 , 321 ,...... | $f_1(x) = 11x^2 - 15x + 15$ | $f_2(x) = 11x^2 + 7x + 11$ | $f_3(x) = 11x^2 + 29x + 29$ | $f_4(x) = 11x^2 + 51x + 69$ |
| | T3 | 7 , 25 , 65 , 127 , 211 , 317 ,...... | $f_1(x) = 11x^2 - 15x + 11$ | $f_2(x) = 11x^2 + 7x + 7$ | $f_3(x) = 11x^2 + 29x + 25$ | $f_4(x) = 11x^2 + 51x + 65$ |